\documentclass[12pt, makeidx]{amsbook}  
 
\usepackage{amssymb} 
\usepackage{enumerate, amsfonts, latexsym, color,graphicx} 

\usepackage{epsfig, wasysym, fancyhdr}

\fancypagestyle{plain}
\def\widecheck{\check}

\newtheorem{thmA}{Theorem} 
 
\newtheorem{corA}[thmA]{Corollary}

\newtheorem{theorem}{Theorem}[section]

\newtheorem{prop}[theorem]{Proposition} 
\newtheorem{proposition}[theorem]{Proposition} 
\newtheorem{addendum}[theorem]{Addendum}

\newtheorem*{unnumb-def}{Definition}

\newtheorem{convention}[theorem]{Convention}
\newtheorem{workingassumption}[theorem]{Working Assumption}
 \newtheorem {coroll} [theorem] {Corollary}

\newtheorem{lemma} [theorem]{Lemma} 
\newtheorem{corollary} [theorem] {Corollary}

\newtheorem{observation} {Observation} 
\newtheorem{notation}[theorem]{Notation} 
 
\theoremstyle{definition} 
\newtheorem{definition}[theorem]{Definition} 
\newtheorem{example}[theorem]{Example} 
 
\theoremstyle{remark} 
\newtheorem{remark}[theorem]{Remark} 
\newtheorem{RIproof}[theorem]{}
 
\newtheorem{remarks}[theorem]{Remarks} 
 
\numberwithin{equation}{section} 
 
\newfont{\msb}{msbm10 scaled 1200} 
\newfont{\euf}{eufm10 scaled 1200} 
\def\Prf:       {\it Proof. \rm} 
 
\def\Gthree{\mathcal G_3}
\def\bt {\begin{theorem}} 
\def\et         {\hfill $\square$\end{theorem}} 
\def\bl {\begin{lemma}} 
\def\el         {\hfill $\square$\end{lemma}} 
\def\bd         {\begin{definition}} 
\def\ed{\end{definition}} 
\def\bc         {\begin{corollary}} 
\def\ec         {\hfill $\square$ \end{corollary}} 
\def\bea {\begin{eqnarray}} 
\def\eea        {\end{eqnarray}} 
\def\be {\begin{eqnarray*}} 
\def\ee {\end{eqnarray*}} 
\def\bp {\begin{proposition}} 
\def\ep {\end{proposition}}  
\def\bex {\begin{example}} 
\def\eex {\end{example}} 
\def\height{\text{\rm{time}}} 
\def\Pin{\Pi} 
  
\def \Cal{\mathcal} 
\def \Bbb{\mathbb} 
\def \S{\Sigma} 
\def\D{\Delta} 
\def\E{\cal E} 
\def\ssm{\smallsetminus} 
 
\def\area{\text{\rm{Area}}}

\def\<{\langle} 
\def\>{\rangle} 
\def\|{{\,|\! |\, }}

\def\P{\Cal P} 
\def\cal{\Cal} 
\def\e{\varepsilon}

\def\-{\underline} 
\def\wt{\widetilde} 
\def\N{\Bbb N} 
\def\Z{\Bbb Z} 
 
\def\G{\Gamma} 
\def\A{\Cal A} 
 
\def\time{\text{\rm{time}}}

\def\On{\text{\rm{Out}} (F_n)}

\def\iso{\cong} 
\def\supp{\text{\rm{Supp}}} 
\def\pre{\text{\rm{Pre}}} 
\def\lpl{left para-linear } 
\def\rpl{right para-linear } 
\def\n  {\,|\partial\Delta|} 
\def\vecZ{\mathcal Z} 
\def\tz{{T_0}}  
\def\life{\text{\rm Life}} 
\def\ttt{T_1} 
\def\mess  {\text{\rm Mess}} 
\def\R{\Cal R} 
\def\AFourC     {2C_1 + 6\ll + 2B(5T_0 + 6\ttt + 2) + 2LC_4(6\ttt +
8T_0 + 3) + (B+3)(3\ttt + 2T_0)L +5L+2} 

\def\T{\mathcal T} 
\def\cmm{C_{(\mu,\mu')}(2)} 
\def\vin        {\in_v}  
\def\subT{\chi(\Pin_{\T})} 
\def\CT{\chi_c(\T)} 
\def\down{\text{\rm{down}}}

\def\bonus{\text{\rm{bonus}}} 
\def\dlong{\text{\rm{D}}\Lambda} 
 
\def\plT{p_l({\mathcal T})} 
\def\tplT{\tilde p_l(\T)}
\def\Bb         {(B+3)(3\ttt + 2T_0)L + 6B\ttt + 4BT_0 + 2\ll + 2B + 5L + 1} 
\def\F{\mathcal F} 
\def\ptmm{\hat t_1(\mu,\mu')} 
\def\prmm{\hat\rho(\mu,\mu')} 
\def\pEmm{\hat{\text{\euf T}}  (\mu,\mu')} 

\def\pTmm{\hat\T(\mu,\mu')} 

\def\pT{\hat\T}

\def\xT{x({\mathcal T})}  
\def\prT{p_r({\mathcal T})}  
\def\tprT{\tilde p_r(\T)}
\def\ttwo{t_2(\T)} 
\def\ET{{\text{\euf T}}}
\def\Emm{{\text{\euf T}}(\mu,\mu')}
\def\tone{t_1(\T)} 
\def\I{\mathcal I} 
\def\r{\rho} 
\def\rT{\rho_\T} 
 
\def\bonusT{\text{\rm{bonus}}(\T)} 
\def\K  { 2C_0 + 2K_1 + 2B + 1}  

\def\ll {\lambda_0} 
 
\def\L{\Cal L} 
\def\R{\Cal R}
\def\QT{Q(\T)} 
 
 \def\On{\text{\rm{Out}} (F_n)}

\def\Aut{${\rm{Aut}} (F_n)$}

\def\gep{{\rm{GEP}}}
\def\pep{{\rm{$\Psi$EP}}}
\def\atom{atom}
\def\supp{\text{\rm{Supp}}}

 \def\NF {Nibbled Future}
 \def\nf {nibbled future}
 
 \def\nib {nibbling}

\def\gep{{\rm GEP}}
\def\pep{{\rm $\Psi$EP}}
\def\rpep{{\rm PEP}}
\def\atom{atom}
\def\lbite{\LEFTcircle}
\def\rbite{\RIGHTcircle}
\def\lhnp{\Leftcircle}
\def\rhnp{\Rightcircle}

\def\G{\Gamma}
\def\AR{\langle\Cal A\mid\Cal R\rangle}

\def\At{(\mathcal A\cup\{t\})^{\pm 1}}
\def\A{\mathcal A}
\def\Apm{\mathcal A^{\pm}}
\def\FAt{(\mathcal A^{\pm}\cup\{t^{\pm 1}\})^*}
\def\FAtp{(\mathcal A^{\pm}\cup\{t^{\pm p}\})^*}
\def\FA{F(A)}

\def\ssm{\smallsetminus}

\catcode`\@=11 
\def\serieslogo@{\relax} 
\def\@setcopyright{\relax} 
\catcode`\@=12

\makeindex
 
\begin{document}

\frontmatter

\title[Mapping tori of free group automorphisms]
{The quadratic isoperimetric inequality
for mapping tori of free group automorphisms.}
 
\author[Martin R. Bridson]{Martin R.~Bridson} 
\address{Martin R.~Bridson\\
Mathematical Institute\\
24--29 St Giles'\\
Oxford, OX1 3LB, UK}
\email{bridson@maths.ox.ac.uk}

\author[Daniel Groves]{Daniel Groves}
\address{Daniel Groves\\
MSCS 322 SEO, \textsc{M/C} 249\\
University of Illinois at Chicago\\
851 S. Morgan St.\\
Chicago, IL 60607-7045, USA}
\email{groves@math.uic.edu}

\date{8 February, 2008}

 

\begin{abstract} {If $F$ is a finitely generated free group
and $\phi$ is an automorphism of $F$
then  $F\rtimes_\phi\mathbb Z$
satisfies a quadratic isoperimetric inequality.} 
\end{abstract}

\maketitle

\setcounter{tocdepth}{1}
\tableofcontents 

\chapter*{Introduction}

Associated to an  automorphism $\phi$
of any group $G$ one has the algebraic {\em mapping torus}
$G\rtimes_\phi\mathbb Z$. In this paper we shall be concerned
with the case where $G$ is a finitely generated
free group, denoted $F$. We seek to understand the complexity
of the  word problem in the
 groups $F\rtimes_\phi\mathbb Z$ as measured
by their Dehn functions.

The class of groups of the form $F\rtimes_\phi\mathbb Z$
has been the subject of intensive investigation 
in recent years and a rich structure has begun to
emerge  in keeping with the subtlety of the
classification of free group automorphisms \cite{BFH}, \cite{BFH2}
\cite{BH2}, \cite{FH}, \cite{Lu}, \cite{Sela}.  (See \cite{BestICM} and
the references therein.)
Bestvina--Feighn and Brinkmann proved that if $F\rtimes_\phi\mathbb Z$
doesn't contain
a free abelian subgroup of rank two then
it is hyperbolic \cite{BF}, \cite{Brink},
i.e. its Dehn function
is linear. Epstein and Thurston \cite{E+}
proved that if $\phi$ is induced by a surface automorphism
(in the sense discussed below) then $F\rtimes_\phi\mathbb Z$ is automatic
and hence
has
a quadratic Dehn function. The question
of whether or not all non-hyperbolic groups of the form
$F\rtimes_\phi\mathbb Z$  have quadratic Dehn functions has
attracted a good deal of attention.

\renewcommand{\thethmspec}{{\bf{Main Theorem. \kern-.3em}}}

\smallskip

\noindent
\begin{thethmspec} \label{MainThm}
{\em If $F$ is a finitely generated
free group and $\phi$ is
an automorphism of $F$ then $F \rtimes_{\phi} \mathbb Z$ satisfies a
quadratic isoperimetric inequality.}
\end{thethmspec}

\smallskip
Papasoglu \cite{papa} proved that if a finitely presented group satisfies a quadratic isoperimetric
inequality, then all of its asymptotic cones are simply connected.

\begin{corA}\label{asCones} If $F$ is a finitely generated
free group and $\phi$ is
an automorphism of $F$ then,  then every asymptotic cone of
$F \rtimes_{\phi} \mathbb Z$ is simply connected.
\end{corA}

Ol'shanskii and Sapir \cite[Theorem 2.5]{OS} proved that if a
multiple HNN extension of a free group has Dehn function less
than $n^2 \log n$ (with a somewhat
technical definition of `less than') then
it has a solvable conjugacy problem. Our Main Theorem shows
that free-by-cyclic groups fall into this class.

\begin{corA} \label{Conjug}
If $F$ is a finitely generated
free group and $\phi$ is
an automorphism of $F$, then the
conjugacy problem for $F \rtimes_{\phi} \mathbb Z$ is solvable.
\end{corA}

Corollary \ref{Conjug} was first proved in \cite{BMMV} using
different methods. 

Gromov \cite{Gromov} proved that a finitely presented group is hyperbolic if and only if
its Dehn function is linear. He also proved that if a Dehn function is subquadratic then it
must be linear. Thus if one ranks groups according to the complexity of their Dehn functions,
the groups that have a quadratic Dehn function demand particular attention. The
nature of these groups is far from clear for the moment; in particular it is unclear
what they have in common. It is not known, for example,
whether they all have a solvable conjugacy problem. Nor is it
known  whether the isomorphism problem is solvable amongst them.  Our Main Theorem provides
a rich source of new examples on which to test such questions. 

\medskip

Much of our modern understanding of the automorphisms of
free groups has been guided by the  analogies with
automorphisms of free-abelian groups and surface groups \cite{mb-kv}.
The former analogy will prove useful  is our analysis of
how elements of a free group grow when one repeatedly applies an
automorphism,  but it offers
offers us poor guidance at the level of Dehn functions: the
Dehn function of $\Z^d\rtimes_\phi\Z$ can be polynomial
of degree $2,3,\dots,d+1$ or it can be  exponential; it depends
on the growth rate of $\phi$ and is quadratic only if $\phi\in{\rm{GL}}(n,\Z)$
has finite order \cite{BG}.

The analogy with surface automorphisms is more apt.
A self-homeomorphism of a compact surface $S$ defines an outer
automorphism of $\pi_1S$ and hence a semidirect product
$\pi_1S\rtimes_\phi\mathbb Z$. This group  is the fundamental
group of a compact   3-manifold, namely the
mapping torus $M_\phi$ of the homeomorphism. By using
Thurston's Geometrization
Theorem for Haken manifolds, Epstein and Thurston \cite{E+} were
able to prove that
$\pi_1S\rtimes_\phi\mathbb Z$
 is an automatic group; hence its Dehn function is
either linear or quadratic. If $S$ has boundary then
only the quadratic case arises. A more geometric explanation for the existence of
a quadratic isoperimetric inequality in the bounded case
comes from the fact that
 $M_\phi$ supports a
metric of non-positive curvature, as does any irreducible
3-manifold with non-empty boundary \cite{mb-shs}, \cite{leeb}.

If $S$ has boundary, then $\pi_1S$ is free. Thus the
foregoing considerations give many examples of free-by-cyclic
groups that have quadratic Dehn functions. But
there are many types of  free group automorphisms that
do not arise from surface automorphisms, for example
those $\phi$ that do not have a power leaving any
non-trivial conjugacy class invariant, and those $\phi$
for which there is a word $w\in F$ such that the function $n\mapsto |\phi^n(w)|$
grows like a super-linear polynomial.

The non-automaticity of certain  $F\rtimes_\phi\mathbb Z$ provides
a more subtle obstruction to realising $\phi$ as a surface automorphism:
in contrast to the
Epstein-Thurston Theorem, Brady, Bridson and Reeves  \cite{BB},
\cite{BR} showed that certain mapping tori $F_3\rtimes \mathbb Z$ are
not automatic, for example that associated to the  automorphism
$[a\mapsto a,\, b\mapsto ab,\, c\mapsto a^2c]$.
Such examples show that one cannot proceed  via automaticity in order
to prove the Main Theorem. Nor can one rely on non-positive
curvature,  because Gersten \cite{Ge} showed that the above example
$F_3\rtimes\mathbb Z$ is not the fundamental group of any compact
non-positively curved space. Thus one needs a new approach to the
quadratic isoperimetric inequality.

A technique for dealing with classes of linearly growing automorphisms is described by
Brady and Bridson in \cite{BB}, while Macura  \cite{Mac} developed techniques for
dealing with polynomially growing automorphisms.  But these techniques apply only
to  restricted classes of automorphisms and do not speak to the core
problem of establishing the quadratic isoperimetric inequality
for mapping tori of general free group automorphisms.
In the present work
we attack  this core problem
directly, undertaking a detailed analysis of the geometry of van Kampen diagrams over the
natural presentations of free-by-cyclic groups.  

The focus of this analysis is on the dynamics of the {\em{time flow of 
$t$-corridors}}, which is closely related to the dynamics of
the given free group automorphism.
Here, $t$ is the generator of the $\Z$ factor in $F\rtimes_\phi\Z$
and a $t$-corridor is a chain of 2-cells extending across
a van Kampen diagram with adjacent 2-cells abutting
along an edge labelled $t$ (see Subsection \ref{sec:Time,corridors}).

The key estimate -- a linear
bound on the length
of $t$-corridors (Theorem \ref{BoundS:BG3}) -- admits the following algebraic formulation. This
clarifies the manner in which our results concerning the
geometry of van Kampen diagrams give
rise to a non-deterministic quadratic time algorithm for the word
problem in free-by-cyclic groups (for an alternative approach see
\cite{schleimer}).

Fix a set of generators $\mathcal A$ for $F$ and let
$d_{F}$ be the corresponding word metric.
We
consider words over the  $e_i\in(\mathcal A\cup\{t\})^{\pm 1}$,
where $t$ is a generator of the righthand factor of
$F\rtimes_\phi\Z$.
A {\em{bracket}} $\beta$ in a word $w$ is a decomposition
$w
\equiv
w_1(w_2)w_3$; the subword $w_2$ is the {\em{content}} of
$\beta$, and the
initial and terminal letters of $w_2$ are its {\em{sentinels}}.
A  second bracket $\beta'$, giving $w \equiv w_1'
(w_2') w_3'$
is {\em{compatible}}
with $\beta$ if $w_2' \subset w_i$ for some
$i \in\{1, 2, 3\}$ or $w_2\subset w_i'$.
A
{\em{$t$-complete}} bracketing is a set of pairwise compatible
brackets $\beta_1,\dots,\beta_m$ such that the sentinels of
 each $\beta_i$
are $\{t, t^{-1}\}$ and every $t^{\pm 1}$ in $w$ is a sentinel of a
unique bracket.  In such
a bracketing, the content of each bracket is equal in
$F\rtimes_\phi \Z$ to an element of $F$.

\renewcommand{\thethmspec}{{\bf{Bracketing Theorem. \kern-.3em}}}

\smallskip
\noindent
\begin{thethmspec}\label{Bracket}\index{Bracketing Theorem}
{\em{There exists a constant $K = K(\phi,\mathcal B)$
such that any word $w \equiv e_1 \dots e_n$ that represents
the identity in $F\rtimes_\phi\Z$
admits a $t$-complete bracketing $\beta_1,\dots , \beta_m$
such that the content $c_i$ of each $\beta_i$
satisfies $d_F(1,c_i)\le Kn$.}}
\end{thethmspec}
\smallskip

In order to prove the above theorems one  has to delve deeply into the nature
of free-group automorphisms. In particular, one needs a precise understanding
of how  the iterated images $\phi^n(w)$ of an arbitrary element  $w\in F$ can
evolve. This delicate task is made possible by the existence of informative
geometric representatives for $\phi$.

We already alluded to the fact that
the study of automorphisms of free groups is  informed greatly by the 
analogies with automorphisms of free-abelian groups and surface groups.
However, one
often has to work considerably harder in the free group case in order to
obtain the appropriate analogues of familiar results from these other
contexts. Nowhere is this more true than in the quest for suitable normal
forms and geometric representatives. One can gain insight into the nature of
individual elements of $\hbox{GL}(n,\mathbb Z)$ by realizing them
as diffeomorphisms  of the $n$-torus. Likewise, one analyzes individual
elements   of the mapping class group by realizing them
as diffeomorphisms  of a surface.  The situation for ${\rm{Aut}}(F)$
 and ${\rm{Out}}(F)$
is more complicated: the natural choices of classifying space $K(F_n,1)$
are finite graphs of genus $n$, and no element of infinite order in
${\rm{Out}}(F)$
is induced by the action on $\pi_1(Y)$ of a
homeomorphism of $Y$. Thus the best that
one can hope for in this situation is to identify a graph $Y_\phi$ that admits a
{\em{homotopy equivalence}} inducing $\phi$ and has additional structure
well-adapted to $\phi$.
This is the context of the {\em train track technology} of
 Bestvina, Feighn and Handel \cite{BH2, BFH, BFH2}.

Their work results in a decomposition
theory for elements of ${\rm{Out}}(F)$ that is closely analogous to (but
more complicated than) the Nielsen-Thurston theory for surface
automorphisms \cite{casson}.  The finer features of
the topological normal forms that they obtain are adapted to the problems that
they wished to solve in each of their papers: the Scott conjecture in \cite{BH2}
and the Tits alternative in the series of papers \cite{BFH, BFH2, BFHGeomDed}.
The problem that we solve in this book, that of determining
the Dehn functions of  all free-by-cyclic groups, requires
a further refinement of the train-track technology. Specifically,
we must adapt our topological representatives
so as to make tractable the problem of determining the isoperimetric properties
of the mapping torus of the homotopy equivalence $f: Y_\phi \to Y_\phi$ realizing
an iterate of $\phi$.

Recall that an automorphism $\phi$ of a finitely generated free group $F$ is
called {\em{positive}} if there is a basis $a_1,\dots,a_n$ for $F$ such that the
reduced word representing each $\phi(a_i)\in F$ contains no inverses $a_j^{-1}$.
On the rose (1-vertex graph) with directed edges labelled $a_i$, one has a
natural representative for any automorphism of $F$.  The key feature of positive
automorphisms is the fact that the positive iterates of
this representative  restrict to  injections on each edge
of the graph. Such maps are the prototypes for  train-track representatives.

This discussion suggests a strategy that one might follow in order
to prove one Main Theorem:
first, one should  prove it
in the case of positive automorphisms, relying on the simplifications
afforded by the positivity hypothesis to confront
the web of large-scale cancellation phenomena that must be understood
if one is to have any chance of proving the theorem in general.
Then, in the general case, one should attempt to
follow the architecture of the proof in the positive case,
using a suitably refined
train-track description of the automorphism in
place of the positivity assumption. We shall implement the two stages of
this plan in Parts 1 and 3 of this monograph, respectively.
Ultimately, this strategy works. However, in Part 3,
in order to bring our plan to fruition we have to deal with myriad additional complexities
arising from  intricate cancellations that do not arise in the positive case. 

Roughly speaking, these additional complexities correspond to the fact that
most free group automorphisms do not have train track representatives, 
only {\em{relative}} train track representatives.
In Part 2 of  this monograph, we refine the theory of \index{train track map!improved relative}
{\em improved relative train track maps} due
Bestvina, Feighn and Handel
\cite{BFH}, so as to tease-out features that 
allow us to adapt the crucial arguments from Part 1. A
vital ingredient in this approach is the identification
of basic units that will play the
role in the general case that single edges (letters)
played in the positive case. To this end, we develop  a theory
of {\em{beads}}, whose claim to the
role is clinched by the
{\em{Beaded Decomposition Theorem \ref{BDT}}}. This theorem is the main objective of
Part 2. Indeed we have gone to considerable lengths to distill the entire contribution of
Part 2 to Part 3 into this single statement and the important technical refinement
of it described in Addendum \ref {add:SameJ}. We have done so
in order that the reader who is willing to accept
it as an article of faith may proceed directly from Part 1 to Part 3.

The introduction to each part of the book contains a more detailed explanation of its
contents.

\medskip

\noindent{\bf{Acknowledgements.}}  The first author's work was supported in part by
 Research Fellowships from the EPSRC of Great Britain and by a Royal Society
Wolfson Research Merit Award.  Much of this work was undertaken
whilst he was a Professor at Imperial College London, from which
he was granted two terms of sabbatical leave. The second author was supported
in part by a Junior Research Fellowship at Merton College, Oxford, by a Taussky-Todd
Instructorship and a Senior Research Fellowship at the California Institute of
Technology and by NSF Grant DMS-0504251.  We thank these organisations for their support.
We also thank the anonymous referee for his careful reading and helpful comments.
  
\mainmatter

\part{Positive Automorphisms}\label{Part:BG1}

\renewcommand{\thesection}{\ref{Part:BG1}.\arabic{section}}

\setcounter{section}{0}

 An automorphism $\phi$ of a finitely generated free
group $F$ is called \index{positive automorphism}{\em positive} if there is a 
basis $a_1,\dots,a_n$ for $F$ such that the reduced
word representing each $\phi(a_i)\in F$ contains no inverses $a_j^{-1}$. 
Part \ref{Part:BG1} of this work
is dedicated entirely to proving the following special case of the Main Theorem.

 \smallskip

\begin{thmA}\label{MainThmPos}
Let $F$ be a finitely generated free group.
If $\phi$ is a positive  automorphism of $F$, then
$F\rtimes_\phi\mathbb Z$ satisfies a  quadratic isoperimetric inequality.
\end{thmA}

This part of the book is organised as follows.
In Section \ref{vanKampSection} we
recall some basic definitions associated to Dehn functions.
In Sections \ref{BCSection} and \ref{time}
 we record some simple but important observations
concerning the large-scale behaviour of  
the van Kampen
diagrams associated to free-by-cyclic groups and in
particular the geometry of {\em corridor} subdiagrams. (The automorphisms
considered up to this point are not assumed to be positive.)
These observations lead us to a strategy for proving Theorem  \ref{MainThmPos} based
on the geometry of the {\em time flow of corridors}. In Section  \ref{StrategySection}
we state a sharper version of Theorem  \ref{MainThmPos} adapted to this strategy
and reduce to the study of automorphisms with  stability
properties  that regulate the evolution of corridors. 
  In Section \ref{PrefFutSec} we develop the notion of
 {\em preferred future} which allows
 us to trace the trajectory of
  $1$-cells in the corridor flow.
  
The estimates that we establish in Sections \ref{PrefFutSec} and 
\ref{NonConstantSection} reduce us to the
nub of the difficulties that one faces in trying to prove Theorem \ref{MainThmPos},
 namely
the possible existence of  large blocks of ``constant letters". A
sketch of the strategy that we shall use to overcome this problem is presented
in Section \ref{A4sec}. The three main ingredients in this strategy are
the elaborate global cancellation arguments in Section \ref{ConstantSection}, the machinery
of {\em teams} developed in Section \ref{teamSec}, and the {\em bonus scheme}
developed in Section \ref{BonusScheme} to accommodate a final tranche of cancellation
phenomena whose quirkiness eludes the grasp of teams. In a brief final
section we gather our many estimates to establish the bound required for
Theorem \ref{MainThmPos}. A glossary of constants is included for the reader's convenience. 

\section{Van Kampen Diagrams} \label{vanKampSection}

We recall some basic definitions and facts concerning 
Dehn functions and van Kampen diagrams.

\subsection{Dehn Functions and Isoperimetric Inequalities} 
 
Given a finitely presented group $G=\langle \mathcal A \mid 
\mathcal R \rangle$
and a word $w $ in the generators $\mathcal A^{\pm 1}$ that
represents $1\in G$,
one defines
$$
\area (w)  = \\
 \min\big\{ N \in {\mathbb N}^+ \; | \;
\exists\text{ equality }
w = \prod^N_{j=1}u_j^{-1}r_j u_j \text{ in $F(\mathcal A)$
with } r_j \in \mathcal R^{\pm 1} \big\}\, .
$$

The \index{Dehn function}{\it Dehn function} $\delta(n)$ of the finite 
presentation $\langle \mathcal A \mid \mathcal R\rangle$ is defined by 
$$ 
\delta(n) \; = \; \max\{\text{\rm{Area}}(w) \; 
|\; w \in \text{\rm{ker}}(F(\mathcal A) \twoheadrightarrow G), 
\; |w| \leq n \, \} \, , 
$$ 
where $|w|$ denotes the length of the   word $w$. 
Whenever two presentations 
define isomorphic (or indeed  quasi-isometric) 
groups, the Dehn functions of 
the finite presentations   
are equivalent under the relation 
$\simeq$ that identifies functions  
$[0,\infty)\to [0,\infty)$ that only differ by a quasi-Lipschitz 
distortion of their domain and their range.

For any constants $p,q\ge 1$, one sees that 
$n\mapsto n^p$ is $\simeq$ equivalent to $n\mapsto n^q$ 
only if $p=q$. Thus it makes sense to say that the 
``Dehn function of a group" is $\simeq n^p$. 
 
A group $\G$ is said to \index{isoperimetric 
inequality} {\em satisfy a quadratic isoperimetric 
inequality} if its Dehn function is $\simeq n$ or 
$\simeq n^2$. A result of  Gromov \cite{Gromov}, detailed 
proofs of which were given by several authors, states that if  
a  Dehn function is subquadratic, then it is linear --- 
see \cite[III.H]{BH} for a discussion, proof and references.

See \cite{steer} for a thorough and  elementary account of 
what is known about  Dehn functions and an  
explanation of their connection  with filling 
problems in Riemannian geometry.  
\smallskip  
 
\subsection{Van Kampen diagrams}\label{vkD} 
 
According to van Kampen's lemma (see \cite{vK}, 
\cite{LS} or \cite{steer})   
an equality $w = \prod^N_{j=1}u_jr_ju_j^{-1}$ in the 
free group $\mathcal A$, with $N=\area(w)$, 
can be portrayed by  
a finite, 1-connected, 
combinatorial 2-complex with basepoint, embedded in $\mathbb R^2$. Such a complex is 
called a \index{van Kampen diagram}{\em van Kampen diagram} for $w$; its oriented 1-cells   
are labelled by elements of $\mathcal A^{\pm 1}$; 
the boundary label on each 2-cell (read with clockwise 
orientation from one of its vertices) is an element  
of $\mathcal R^{\pm 1}$;  and the boundary cycle of the 
complex (read with positive orientation from the basepoint) 
is the word $w$; 
the number 
of 2-cells in the  diagram   is $N$.  Conversely, any van Kampen diagram with $M$ 
2-cells gives rise 
to an equality in $F(\mathcal A)$ expressing the word 
labelling the boundary cycle 
of the diagram as a product of $M$  
conjugates of the defining relations.  
Thus 
$\text{\rm{Area}}(w)$ is the minimum number of 2-cells among all 
van Kampen diagrams  
for $w$. If a van Kampen diagram $\Delta$ for $w$ has $\text{\rm{Area}}(w)$ 
2-cells, then $\Delta$ is a called a {\em least-area} diagram. If 
the underlying 2-complex is homeomorphic to a 2-dimensional 
disc, then the van Kampen diagram is called a {\em disc diagram}. 
 
We use the term  \index{van Kampen diagram!area of}{\em area} to describe the number of 2-cells in a 
van Kampen diagram, and write $\text{\rm{Area }} \Delta$. We write 
$\partial \Delta$ to denote the boundary cycle of the diagram; we write 
$|\partial\Delta|$ to denote the length of this cycle.

Note that associated to a van Kampen diagram $\Delta$ with basepoint $p$ 
one has a morphism of 
labelled, oriented graphs $h_\Delta: (\Delta^{(1)},p)\to (\mathcal C_\A, 1)$,  where 
$\mathcal C_\A$ is the Cayley graph associated to the choice of 
generators $\A$ for $G$. The map $h_\Delta$ takes $p$ to the 
identity vertex  $1\in \mathcal C_\A$ and preserves the labels on oriented edges. 
 
We shall need the following simple observations. 
 
\begin{lemma} If a van Kampen diagram 
$\Delta$ is least-area, then every simply-connected 
subdiagram of $\Delta$ is also least-area. 
\end{lemma} 
 
Recall that a function $f:\mathbb N\to [0,\infty)$ 
is {\em sub-additive} if $f(n+m)\le f(n) + f(m)$ 
for all $n,m\in\mathbb N$. For example, given $r\ge 1,\, k>0$, 
the function $n\mapsto kn^r$ is sub-additive. 
 
\begin{lemma}  \label{disc} 
Let $f:\mathbb N\to [0,\infty)$ be a sub-additive function and let $\P$ 
be a finite presentation of a group. 
If  $\text{\rm{Area }} \Delta\le f(|\partial\Delta|)$ for every 
least-area disc diagram $\Delta$ over $\P$, then the Dehn 
function of $\P$ is $\le f(n)$. 
\end{lemma}

\subsection{Presenting $F\rtimes\mathbb Z$} 
 
We shall establish the quadratic bound required for 
the Theorem \ref{MainThmPos} by examining the nature of van Kampen 
diagrams over the following natural (aspherical) presentations of
free-by-cyclic groups. 
 
Given a finitely generated free group $F$ and 
an automorphism $\phi$ of $F$, we fix a basis 
$a_1,\dots,a_m$ for $F$, write $u_i$ to denote 
the reduced word equal to $\phi(a_i)$ in $F$, and 
present $ 
F\rtimes_\phi\mathbb Z$ by 
\begin{equation}\label{presentation} 
\P\iso \langle a_1,\dots,a_m,t\mid 
t^{-1}a_1tu_1^{-1},\dots, t^{-1}a_mtu_m^{-1}\rangle. 
\end{equation} 
Throughout Part \ref{Part:BG1}, we shall work exclusively with this presentation. 
\medskip

\begin{figure}[htbp] 
\begin{center} 
  
  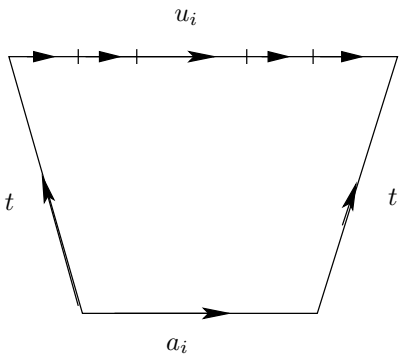  
\caption{A $2$-cell in a van Kampen diagram for $F \rtimes_{\phi} \mathbb Z$.} 
\label{figure:2cell} 
\end{center} 
\end{figure}

\subsection{Time and $t$-Corridors with naive tops} \label{sec:Time,corridors}
 
The use of $t$-corridors as a tool for investigating 
van Kampen diagrams has 
become well-established in recent years. In the 
setting of van Kampen diagrams over the above presentation,  
$t$-corridors are easily described. 
 
Consider a van Kampen diagram $\Delta$ over 
the above presentation $\P$ and focus on an edge in the boundary 
$\partial \Delta$ that is labelled 
$t^{\pm 1}$ (read with positive orientation from the basepoint).  
If this edge lies in the boundary 
of a 2-cell, then the boundary  cycle of this 2-cell has 
the form $t^{-1}a_itu_i^{-1}$ (read with suitable orientation from 
a suitable point, see Figure \ref{figure:2cell}). In particular, there 
is  exactly one other edge  
in the boundary of the 2-cell that is labelled $t$; crossing 
this edge we enter another 2-cell with a similar boundary 
label, and iterating the argument we get a chain of 2-cells 
running across the diagram; this chain terminates at an edge of 
$\partial \Delta$ which (following the orientation of $\partial \Delta$ 
in the direction of our original edge labelled $t^{\pm 1}$) is labelled 
$t^{\mp 1}$. This chain of 2-cells is called a \index{corridor}{\it{$t$-corridor}}. 
The edges labelled $t$ that we crossed in the above description 
are called the {\em vertical} edges of the corridor.  
The vertical edge on $\partial \Delta$ labelled $t^{-1}$ is 
called the {\em initial} end of the corridor, and at the other end one 
has the {\em terminal} edge. 
 
Formally, one should define a $t$-corridor to be a combinatorial map 
to $\Delta$ from a suitable subdivision of $[0,1]\times [0,1]$: the 
initial edge is the restriction of this map to $\{0\}\times [0,1]$; the 
vertical edges are the images of the 1-cells of the form $\{s\}\times [0,1]$, 
oriented so that the edge joining $(s,0)$ to $(s,1)$ is labelled $t$.
The {\em naive top} of the corridor is the edge-path obtained by restricting 
the above map to $[0,1]\times\{1\}$, and the {\em bottom} is the restriction 
to  $[0,1]\times\{0\}$. 
\smallskip 
 
\noindent{\bf Left/Right Terminology:} The orientation of a disc 
diagram induces an orientation on its corridors. Whenever 
we focus on an individual corridor, we shall regard its 
 initial edge as being {\em left}most  and its terminal 
edge as being {\em right}most. (This is just a suggestive way of saying 
 that the 
corridor map from $[0,1]\times (0,1)\subset \mathbb R^2$ to $\Delta\subset 
\mathbb R^2$ is 
orientation-preserving.) 
 
\smallskip

\medskip 
 
\begin{figure}[htbp] 
\begin{center} 
  
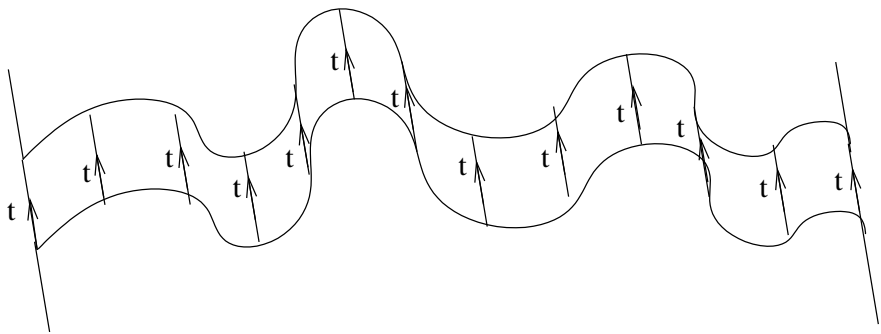 
  
\caption{A $t$-corridor} 
\label{figure:corridor} 
\end{center} 
\end{figure}

See \cite{BG} for a  detailed account of $t$-corridors. 
Here we shall need only the following easy facts: 
\begin{enumerate} 
\item 
distinct $t$-corridors 
have disjoint interiors;  
\item 
if $\sigma$ is the edge-path in $\Delta$ running along 
the (naive) top or bottom  of a $t$-corridor, then $\sigma$ is 
 labelled 
by a word in the letters $\mathcal A^{\pm 1}$ 
that is equal in $F\rtimes\mathbb Z$   to the 
words labelling the subarcs of $\partial \Delta$ 
which share the  endpoints of $\sigma$ (given appropriate 
orientations); 
\item if we are in a least-area diagram  
then the word on the bottom of the corridor is freely reduced; 
\item the number of 2-cells in the 
$t$-corridor is the length of the word labelling the 
bottom side. 
\item In subsection 1.2 we described the map $h_\Delta$ 
associated to a van Kampen diagram. This map 
sends vertices of $\Delta$ to vertices of the Cayley graph $\mathcal C_\A$, 
i.e. elements of $F\rtimes \langle t \rangle$.  If the initial vertex of a  
directed edge in $\D$ is sent to an element of the form $wt^j$, with 
$w\in F$, then the edge is defined to occur at \index{time}{\bf time} $j$. Note that the 
vertical edges of a fixed corridor all occur at the same time. 
\end{enumerate}     
 
We will
consider the {\em dynamics} 
of the automorphism $\phi$ with respect to this notion of time.

\begin{definition} 
[Time and Length]  Item (5) above 
implies that the time of each 
$t$-corridor $S$ is well-defined;
we denote it $\height(S)$.  
 
 We define the \index{corridor!length of}{\em length} of a  
corridor $S$ to be the number of 2-cells that it 
contains, which is equal to the number of 1-cells along its bottom. 
We write $|S|$ to denote the length of $S$. 
\end{definition}

\subsection{Conditioning the Diagram} \label{ss:Cond}
 
We are working with the following presentation of $F\rtimes_\phi\mathbb Z$ 
$$ 
\mathcal P = \langle a_1,\dots,a_m,t\mid 
t^{-1}a_1tu_1^{-1},\dots, t^{-1}a_mtu_m^{-1}\rangle. 
$$ 
 
In the light of Lemma \ref{disc}, in order to prove the  
Theorem \ref{MainThmPos} it suffices to consider only {\em disc diagrams}. Therefore, 
henceforth we shall assume that all diagrams are topological discs. 
We shall also assume that all of the discs considered are 
{\em least-area} diagrams for freely reduced words.  
 
\begin{lemma} 
Every least-area disc diagram  over $\P$ is  
the union of its $t$-corridors. 
\end{lemma} 
 
\begin{proof} Since the diagram is a disc, every 1-cell lies in 
the boundary of some 2-cell. The boundary of each 2-cell 
contains two edges labelled $t$. Consider the equivalence relation 
on 2-cells generated by $e\sim e'$ if the boundaries of $e$ and $e'$ share an 
edge labelled $t$. Each equivalence class forms either a $t$-corridor 
or else a $t$-ring, i.e. the closure of an annular sub-diagram 
whose internal and external cycles are labelled by a word in the 
generators of $F$. If the latter case arose, then since 
$F$ is a free group, the word $u$ on the external 
cycle would be freely equal to the empty word (since it contains no edges 
 labelled $t$). This would contradict the hypothesis that the diagram 
  is least-area, because one could reduce its area 
by excising the simply-connected sub-diagram bounded by this cycle, 
replacing it  with the zero-area diagram for $u$ over the free 
presentation of $F$. 
\end{proof}

\subsection{Folded Corridors} \label{ss:Fold}

In the light of the above lemma, we see that the diagrams $\Delta$ that we 
need to consider are essentially determined once one knows which 
pairs of boundary edges are connected by $t$-corridors. However, there 
remains a slight ambiguity arising from the fact that free-reduction in 
the free group is not a canonical process (e.g.  $x = (xx^{-1})x = x (x^{-1}x)$).

\begin{figure}[htbp] 
\begin{center} 
  
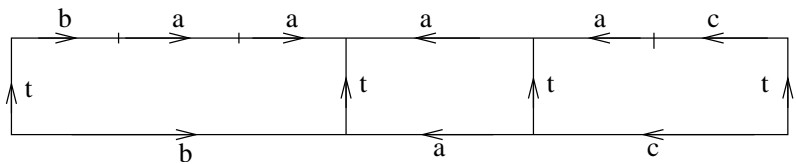 
  
\caption{An unfolded corridor} 
\label{figure:Fold1} 
\end{center} 
\end{figure}

\begin{figure}[htbp] 
\begin{center} 
  
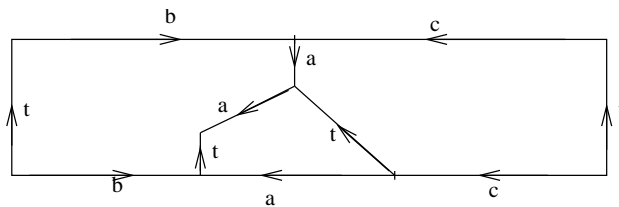 
 
\caption{The  corresponding folded corridor.} 
\label{figure:Fold2} 
\end{center} 
\end{figure}

To avoid this ambiguity, we  fix a least area disc diagram $\Delta$  
and assume that its corridors are {\em folded} in the sense of \cite{B-plms}. 
The topological closure  $T\subset\Delta$ of each corridor is a combinatorial 
disc. 
The hypothesis ``least area" alone 
forces the label on the {\em bottom} of the corridor 
to be a {\em freely reduced} word in the letters $a_i^{\pm 1}$. 
We define the 
 {\em top} of the (folded) corridor to be the  
 injective edge-path that remains when one deletes from the 
frontier of $T$ 
 the bottom and ends of the corridor. The word labelling 
 this path is the freely reduced word in $F$ that equals the 
 label on the naive top of the corridor. Note that, unlike the 
bottom of the corridor, the top may fail to intersect  the closure of some 
2-cells --- see Figures \ref{figure:Fold1} and \ref{figure:Fold2} (where the automorphism is $a \mapsto a, 
b \mapsto ba^2, c \mapsto ca$). 
 
\begin{notation} 
We write $\top (S)$ and $\bot (S)$, respectively, to denote the top and bottom 
of a folded corridor $S$. 
\end{notation} 
 
\smallskip 
 \index{corridor}
{\centerline{ 
{\em Henceforth we shall refer to folded $t$-corridors simply as ``corridors".}}}

\bigskip

\medskip 
 
\subsection{Naive Expansion and Death} 
 
For each generator $a_i\in F$ we have the reduced word 
$u_i=\phi(a_i)$. Given a reduced  word $v=a_{i(1)}\dots a_{i(m)}$ 
we define the \index{naive expansion}{\em naive expansion} of $\phi(v)$ to be  
the (unreduced) concatenation $u_{i(1)}\dots u_{i(m)}$. 
 
Note that if $v$ is the label on an interval of the bottom of a corridor, 
then the naive expansion of $\phi(v)$ is the label on the 
corresponding arc of the naive top of the corridor.

An edge $\e$ on the bottom of a  corridor $S$ is said to \index{death (of an edge)}{\em die} in $S$ 
if the 2-cell containing that edge  does not contain any edge of  
$\top(S)$.  (Equivalently, if $w$ is the label on  $\bot(S)$ and $a_i$ is 
the label on $\e$, then the subword  $u_i=\phi(a_i)$ in 
the naive expansion of $\phi(w)$ is cancelled 
completely during the free reduction encoded in $\Delta$.) In Figure \ref{figure:Fold2}
the edge labelled $a$ on the bottom of the corridor dies.

\section{Singularities and Bounded Cancellation} \label{BCSection} 
 
We have noted that the structure  of a (folded, least-area disc) diagram 
over the natural presentation of a free-by-cyclic group 
is  the union of its 
corridors. 
In this section we pursue an 
understanding of how these corridors meet. 
\smallskip

\medskip

\begin{figure}[htbp] 
\begin{center} 
  
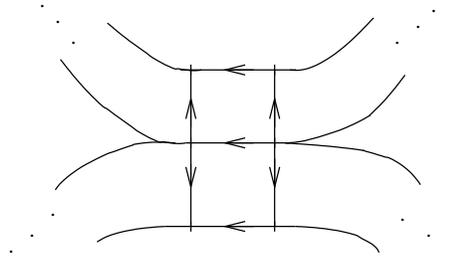 
  
\caption{Corridors cannot meet this way in a least-area diagram} 
\label{figure:wrongway} 
\end{center} 
\end{figure}

The first observation to make is that corridors cannot meet as in 
Figure \ref{figure:wrongway}. 
 
\begin{lemma} \label{l:onepoint}
 If $S\neq S'$, then $\bot(S)\cap\bot(S')$ 
consists of at most one point. 
\end{lemma} 
 
\begin{proof} For each letter $a$, there is only one type 
of 2-cell which has the label $a$ on its bottom side.  Thus, if two corridors 
were to meet in the manner of  Figure \ref{figure:wrongway}, then we would have a pair 
of 2-cells whose union was bounded by a loop labelled  
$u_it^{-1}tu_i^{-1}t^{-1}t$, which is 
freely equal to the identity. By excising this pair of 2-cells and 
filling the loop with a diagram of zero area, we would 
reduce the area of $\Delta$ without altering its boundary label --- 
but  $\Delta$  is assumed to be a least-area diagram. 
 
Thus $\bot(S)\cap\bot(S')$ contains no edges. To see that it cannot 
contain more than one vertex, follow the proof of Proposition 
\ref{SingularityProp}(1). 
\end{proof}

\begin{definition}  
A \index{singularity}{\em singularity} in $\Delta$ is a non-empty connected  component of the intersection 
of the tops of two 
distinct folded corridors. A 2-cell  is said to {\em hit} the 
singularity if 
it contains an edge of the singularity.  

The singularity   is said to be {\em degenerate} if it consists of a single point, and
otherwise it is {\em non-degenerate}.
\end{definition}  
 
\medskip 
 
\begin{figure}[htbp] 
\begin{center} 
  
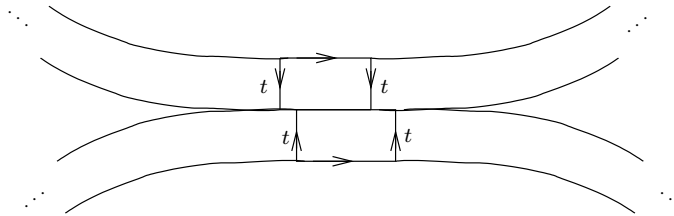 
  
\caption{A `singularity'} 
\label{figure:Singularity} 
\end{center} 
\end{figure}

Let $L$ be the maximum of the lengths of the words $u_i$ in 
our fixed presentation $\mathcal P$ of $F\rtimes_\phi \mathbb Z$. 
 
\begin{proposition}[Bounded singularities] \label{SingularityProp}$\ $ 
\index{Bounded Singularities Lemma}
 
\begin{enumerate} 
\item[1.] If the tops of two corridors in a  least-area 
diagram meet, then their intersection is a singularity. 
\item[2.] 
There exists a constant $B$ depending only 
on $\phi$  such that less than $B$ 2-cells 
hit each singularity in a  least-area diagram over $\P$. 
\item[3.]  
If $\Delta$ is a least-area diagram over $\P$, 
then there are less than $2|\partial \Delta|$ non-degenerate singularities 
in $\Delta$, and each has length at most $LB$.
\end{enumerate} 
\end{proposition}

\begin{proof} Suppose that the intersection of the tops of two corridors $S$ and $S'$ 
contains two distinct vertices, $p$ and $q$ say. Consider the unique subarcs 
of $\top(S)$ and $\top(S')$ connecting $p$ to $q$. 
 Each of these arcs is labelled by a reduced word in 
the generators of $F$; since the arcs have the same endpoints in $\Delta$, 
these words must be identical.  If the arcs did not coincide, then 
we could excise the subdiagram that they bounded and replace it with 
a zero-area diagram, contradicting our least-area hypothesis. This proves 
(1).

\begin{figure}[htbp] 
\begin{center} 
  
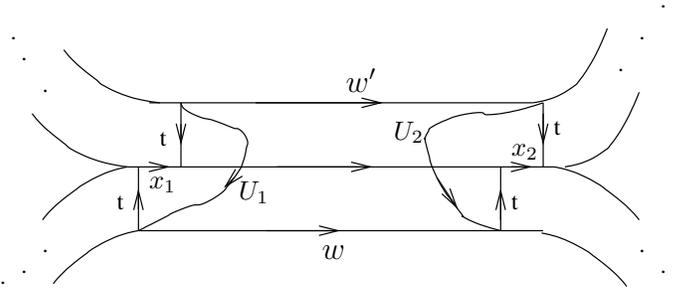 
  
\caption{The proof of Proposition 1.2.3}
\label{figure:BoundedSing} 
\end{center} 
\end{figure}

Figure \ref{figure:BoundedSing} portrays the argument we use to prove (2). In $S$ 
(respectively $S'$), we choose  
an outermost pair of oriented edges $\e_1, \e_2$ (resp. 
$\e_1',\e_2'$) labelled $t$ whose termini lie on the  
singularity. We then connect their endpoints by shortest 
arcs in the singularity as shown. Note that 
each of the arcs labelled $x_1$ and $x_2$ is contained in the top 
of a single 2-cell, and hence has length at most $L$. 
We write $\alpha_i$ to denote the concatenation of $\e_i$, the arc labelled $x_i$ 
and the inverse of $\e_i'$.

Let $U^{-1}_i\in F$ be the reduced word representing 
$\phi^{-1}(x_i)$. In $F\rtimes_\phi\mathbb Z$ we have $tx_it^{-1}U_i=1$; 
let  $\Delta_i$ be a least-area van Kampen diagram portraying  
this equality.   
 
Let $w$ (resp. $w'$) be the label on the edge-path 
in $\bot (S)$ (resp. $\bot(S')$) that connects 
the initial point of $\e_1$ (resp. $\e_1'$) to 
the initial point of $\e_2$ (resp. $\e_2'$).

If we excise from $\Delta$ the subdiagram bounded by the loop whose label 
is 
$t^{-1}wtx_2t^{-1}{w'}^{-1}tx_1^{-1}$, then we reduce the area of  
$\Delta$ by $|w| + |w'|$. (Recall that the edges on the bottom 
of a corridor are in 1-1 correspondence with the 2-cells of the 
corridor.) We may then attach a copy of $\Delta_i$ along $\alpha_i$ 
and fill the resulting loop labelled $U_1wU_2^{-1}{w'}^{-1}$ with 
a diagram of zero area, because this word is equal to $1$ 
in the free group $F$.  
Thus we obtain a new van Kampen diagram whose boundary label 
is the same as that of $\Delta$ and which has area 
$$ 
\area(\Delta) + \area(\Delta_1) + \area(\Delta_2) - |w| - |w'|. 
$$ 
Since $\Delta$ 
is assumed to be least-area, this implies that  
$ \area(\Delta_1) + \area(\Delta_2) \ge  |w| + |w'|.$

 Let $B_0$ be an upper bound on the area of 
all least-area van Kampen diagrams portraying equalities of the  
form  $txt^{-1}\phi^{-1}(x)^{-1}=1$ with $|x|\le L$. 
(It suffices to take $B_0=LL_{inv}$, where $L_{inv}$ is the maximum 
of the lengths of the reduced words $\phi^{-1}(a_i)$.) By definition,  
$ \area(\Delta_1) + \area(\Delta_2)\le 2B_0$, and hence 
$|w| + |w'|\le 2B_0$. Thus for (2) it suffices to let $B=2B_0 + 1$. 
 
The length of the singularity in the above argument 
is less than the sum of the lengths of  the naive 
expansions of $\phi(w)$ and $\phi(w')$. Since  $|w|+|w'|\le B$,  
the singularity has length less than $LB$.

It remains to bound the number of non-degenerate
singularities in $\Delta$. To this
end, we consider the subcomplex  $\Gamma\subset\Delta$ formed by the union of
the tops of all folded corridors. Arguing as in (1), we see that the
graph $\Gamma$ contains no non-trivial loops, i.e. it is a forest. Let 
$V$ denote the set of vertices in $\Gamma$ that have valence at least
3 or else lie on $\partial \Delta$. (Thus $V$ is the set of 
degenerate singularities, endpoints of non-degenerate singularities,
and endpoints of the tops of corridors.)
Let $E$ be the set of  connected components of
$\Gamma\smallsetminus V$.

$|V|-|E|$ is the number $\pi_0$  of connected components of the forest $\Gamma$.
The valence 1 vertices  $V^1\subset\Gamma$ are a subset of the endpoints of
the tops of corridors, so there are less than $|\partial\Delta|$
of them. One can calculate $|E|$ as half the sum of the valences of
the vertices $v\in V$, so $3(|V|-|V^1|) +|V^1| \le 2|E|$. 
Hence
$$
|E| = |V| - \pi_0 \le \frac 2 3 \big{(}|E| + |V^1|\big{)} -\pi_0
< \frac 2 3 \big{(}|E| + |\partial\Delta|\big{)}.
$$
Therefore $|E| < 2|\partial\Delta|$. 

Each non-degenerate singularity determines an element of $E$, so
the (crude) estimate in (3) is established. 
\end{proof}

\index{Bounded Cancellation Lemma}
\begin{lemma}[Bounded Cancellation Lemma] \label{BCL} There is a constant $B$, 
depending only on $\phi$, such that if 
$I$ is an interval consisting of $|I|$ edges 
on the bottom of a (folded) corridor $S$ in a least-area diagram over $\P$, 
and every edge of $I$ dies in $S$, then $|I| < B$. 
\end{lemma} 
\smallskip

\medskip

\begin{figure}[htbp] 
\begin{center} 
  
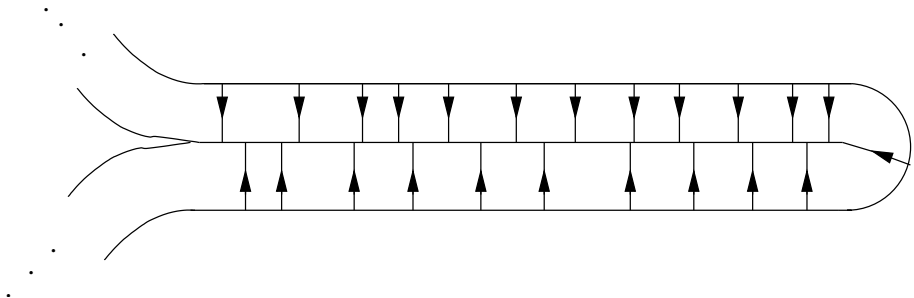 
  
\caption{Bounded Cancellation Lemma} 
\label{figure:BCL} 
\end{center} 
\end{figure} 
 
\begin{proof} The argument is entirely similar to that given for part (2) 
of the previous proposition. 
\end{proof}

The above lemma is a reformulation of the 
 Bounded Cancellation Lemma from \cite{Cooper}, 
which Cooper attributes to Thurston.

\begin{remark} {\em `Singularities are only 1 pixel large.'} 
The reader may find it useful to keep 
in mind the following picture: think of   
a least-area van Kampen diagram rendered on a computer 
screen and assume that the length of the boundary of  
the diagram is  large, so large that the constant $B$ 
in Proposition \ref{SingularityProp} has to be scaled to something less 
than 1 pixel in order to fit the picture on to the  computer's 
screen. 
In the resulting 
image one sees blocks of $t$-corridors as shown in Figure \ref{figure:singflow} 
below, and the singularities take on the appearance of classical $k$-prong
 singularities in the \index{corridor!time-flow of} time-flow of $t$-corridors. 
\end{remark}

\medskip

\begin{figure}[htbp] 
\begin{center} 
  
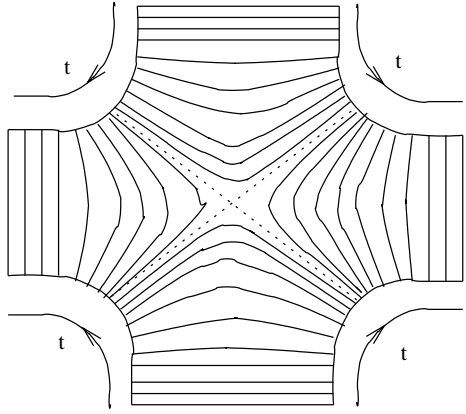 
  
\caption{Schematic depiction of a singularity} 
\label{figure:singflow} 
\end{center} 
\end{figure} 
 
\medskip

\section{Past, Future and Colour} \label{time} 
 
Our investigations thus far have led us to regard van Kampen diagrams 
over $\P$ as flows of corridors  
(at least schematically). We require some more vocabulary to pursue 
this approach. 
 
We continue to work with a fixed disc diagram $\Delta$ over $\P$.

\begin{definition}[Ancestors and Colour]\label{defMu} 
 Each edge $\varepsilon_1$ on the bottom of a corridor either 
lies in the boundary of $\Delta$, or else lies in the top of 
a unique 2-cell, the bottom of which we denote $\e_0$. We consider the 
partial ordering on the set $\mathcal E$ of edges from the bottom of all corridors 
generated by setting $\varepsilon_0 < \varepsilon_1$ whenever edges are related 
in this way.

If $\e'<\e$ then we 
call $\e'$ an \index{ancestor!of an edge}{\em ancestor} of $\e$. The  \index{past!of an edge}{\em past} of $\e$ 
is the set of its ancestors, and the \index{future!of an edge}{\em future} of $\e$ is  
the set of edges $\e''$ such that $\e<\e''$. 
 
Two edges are defined to be of the same \index{colour!of an edge}{\em colour} if  
they have a common ancestor. Since every edge has 
a unique ancestor on the boundary, colours are in 
bijection with a subset\footnote{namely,  
those edges of $\partial\Delta$ that lie on the bottom of 
some 2-cell} of the edges in $\partial\Delta$ whose 
label is not $t$; in particular there are less than 
$|\partial\Delta|$ colours.

Each 2-cell in $\Delta$ has a unique edge 
in the bottom of a corridor. Thus 
we may also regard $\le$ as a partial 
ordering on the 2-cells 
of $\Delta$ and define the past, future and colour 
of a 2-cell. 
 
We define the past (resp. future) of a {\em corridor} 
 to be the union 
of the pasts (resp. futures) of its closed 2-cells. 
\end{definition} 
 
\begin{remark}\label{tree} 
Each $e\in\E$ and each 2-cell has at most one immediate 
ancestor (i.e. one that is maximal among its ancestors). 
Consider the graph $\mathcal F$ with vertex set $\E$ that has an edge 
connecting a pair of vertices if and only if 
one is the immediate ancestor of the other. Note \index{family forest}
that $\mathcal F$ is a forest 
(union of trees). 
 
The \index{van Kampen diagram!colours in}{\em colours} in the diagram correspond to the 
 connected components 
(trees) of this forest.  
 
 There is a natural embedding 
of $\mathcal F\hookrightarrow\Delta$: choose a point (`centre') 
 in the interior of each 2-cell 
and connect it to the centre of its immediate ancestor by an 
arc that passes through their common edge. 
 
\end{remark} 
 
If the future of a corridor $S'$ intersects a corridor $S$ then 
the intersection is connected: 
 
\begin{lemma}[Connected Pasts] \label{Connected}  
If a pair of 2-cells $\alpha$ and $\beta$ in a 
corridor $S$  have ancestors $\alpha'$ and $\beta'$ in a corridor $S'$, then every 
$2$-cell $\gamma$ that lies between $\alpha$ and $\beta$ in $S$  has 
an ancestor $\gamma'$ that lies between 
 $\alpha'$ and $\beta'$  in $S'$. 
 \end{lemma}

\begin{figure}[htbp] 
\begin{center} 
  
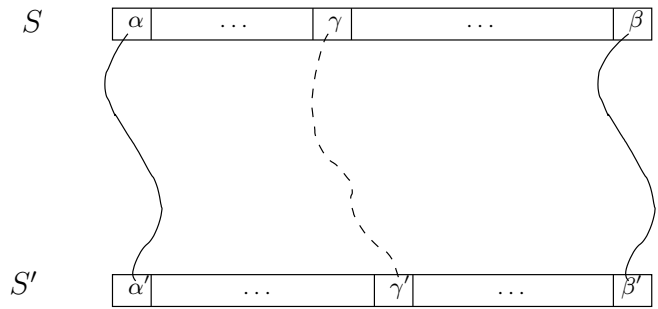 
  
\caption{The `loop' picture} 
\label{figure:loop} 
\end{center} 
\end{figure} 
 
\medskip 
 
\begin{proof} 
Connect the centres of  $\alpha$ and $\beta$ by an arc in the 
interior of $S$ that intersects only those 2-cells lying 
between $\alpha$ and $\beta$, and connect the centres of $\alpha'$ 
and $\beta'$ by 
a similar arc in the interior of $S'$. Along with these two 
arcs, we consider the embedded arcs 
connecting $\alpha$ to $\alpha'$ and $\beta$ to $\beta'$ in the forest 
$\mathcal F$ 
described in Remark \ref{tree}. 
These four arcs together form a  
loop, and the disc that this loop encloses does not intersect 
the boundary of $\Delta$. (Recall that $\Delta$ is a disc.) 
 
Consider the tree from $\mathcal F$ that contains $\gamma$. 
 We may assume that 
the arc in this tree 
that connects $\gamma$ to its ancestor on the boundary does 
not intersect the arc we chose in $S$. It must therefore intersect  
our loop either 
in $S'$, yielding the desired ancestor  $\gamma'$ in $S'$, or 
else in one of the arcs connecting $\alpha$ to $\alpha'$, or 
$\beta$ to $\beta'$. If the latter alternative pertains, $\alpha'$ or 
$\beta'$ is an ancestor of $\gamma$, and we are done. 
\end{proof} 
 
We highlight the degenerate case where the 2-cells $\alpha'$ and $\beta'$ 
are equal and have their bottom   on $\partial\Delta$: 
 
\begin{corollary}\label{muConn} 
Within a corridor, the 2-cells of each colour form a connected region.  
\end{corollary}

\section{Strategy, Strata and Conditioning}  
\label{StrategySection}

Everything that has been said up to  this point has 
 been true for  mapping tori of arbitrary automorphisms of 
 finitely generated free groups. {\em Henceforth, for the remainder
 of Part \ref{Part:BG1}, 
  we assume that the automorphism $\phi$ is positive.} 
 
A van Kampen diagram whose boundary cycle  has length $n$ contains at 
most $n/2$ corridors. Thus  Theorem \ref{MainThmPos} is an immediate consequence of: 
 
\index{corridor!length of}
\begin{theorem}\label{BoundS} 
 There is a constant $K$ depending only on $\phi$ 
such that each  corridor in a least-area diagram $\Delta$ over $\P$ 
has length at most $K\,|\partial\Delta|$. 
\end{theorem}

In order to establish the desired bound on the 
length of corridors, we must analyse how 
corridors grow as they flow into the future, and 
assess what cancellation can take place to inhibit this 
growth. In the remainder of this section we shall 
condition the automorphism to simplify 
the discussion of growth. 
 
\begin{remark}  
 The mapping torus $F\rtimes_{\phi^k}\mathbb Z$ is isomorphic to a  
 subgroup of finite index in  $F\rtimes_{\phi}\mathbb Z$, namely 
  $F\rtimes_{\phi}k\mathbb Z$. Thus, since the Dehn functions of 
  commensurable 
  groups are  
  $\simeq$ equivalent, we are free to replace $\phi$ by a convenient 
  positive power in our proof of the Main Theorem. 
  \end{remark}

\subsection{Strata} 
 
 In the following discussion we shall write $x$ to denote an 
arbitrary choice of letter from our basis $\{a_1,\dots,a_m\}$ 
for $F$.

Naturally associated to any positive automorphism one has 
{\em supports} and  \index{strata}
 {\em strata}.  The support 
  $\supp(x)$ associated to $x$ is  
   the set of all letters which appear in the freely 
   reduced word $\phi^j(x)$ for some $j \geq 0$.  
   The stratum $\S(x)\subset\supp(x)$ associated to $x$ consists 
   of those $y\in\supp(x)$ such that $\supp(x)=\supp(y)$. 
   
   Note that $y\in \supp(x)$ implies $\supp(y)\subseteq\supp(x)$, 
   and $y\in\S(x)$ implies $\S(y)=\S(x)$.

   There are two kinds of strata. 
     The first are  \index{strata!parabolic}{\em parabolic\footnote{Bestvina 
     {\em et al.} \cite{BFH} 
      use the terminology {\em non-exponentially-growing} strata} strata},  
      which are those of the 
     form $\S(x)$ with $x\notin \supp(y)$ 
      for all $y\in\supp(x)\ssm\{x\}$. 
     The second kind  
 are  \index{strata!exponential}{\em exponential strata}, where one has $\S(x)=\S(y)$ for 
 some distinct 
 $x$ and $y$. The letter $x$ is defined to be \index{letter!parabolic}{\em parabolic} 
 or  \index{letter!exponential}{\em  exponential} according to the type of $\S(x)$.

        If $x$ is exponential then $|\phi^j(x)|$ grows exponentially with $j$.  If 
 all the edges of $\supp(x)$ are 
  parabolic then $|\phi^j(x)|$ grows polynomially 
  with $j$.  However, it may also happen that $x$ is a parabolic letter 
  but $|\phi^j(x)|$ 
  grows exponentially; this  
 will be the case if  $\supp(x)$ contains 
exponential letters.

\begin{example}  
 Define $\phi: F_3\to F_3$ by $a_1\mapsto a_1^2a_2,\ a_2\mapsto 
  a_1a_2,\ 
 a_3\mapsto a_1a_2a_3$. Then $\S(a_1)=\S(a_2) =  
 \{a_1, a_2\}$ is an 
 exponential stratum, while $\S(a_3)=\{a_3\}$ 
 is a parabolic stratum with $\supp(a_3)=\{a_1,a_2,a_3\}$. 
 \end{example} 
  
 \begin{remark}\label{induct} 
  The relation {\rm $[y< x$ if  
 $\S(y)\subset\supp(x)\ssm \S(x)]$} generates a partial 
 ordering on the letters $\{a_1,\dots,a_m\}$. For each 
 $x$, the subgroup of $F$ generated by $\pre(x)=\{y\mid y<x\}$ 
 is $\phi$-invariant. Let  $F\lfloor x \rfloor $ denote 
  the quotient of $\langle\supp(x)\rangle$ 
 by the normal closure of $\pre(x)\subset\supp(x)$, and 
 let $F\lceil x\rceil $  denote 
  the quotient of $F$ 
 by the normal closure of $\pre(x)\subset F$. Note 
that  $F\lfloor x \rfloor $  
is a free group with basis (the images of) the letters in $\S(x)$, and $F\lceil x\rceil $ is the free 
 group with basis $\{a_1,\dots,a_m\}\ssm\pre(x)$. 
 
The automorphisms of $\pre(x),\  F\lfloor x \rfloor$ 
and $F\lceil x\rceil $ induced by $\phi$ are positive 
with respect to the obvious bases, and their strata 
are images of the strata of $\phi$. 
 \end{remark} 
  
 \subsection{Conditioning the automorphism} 
  
 In the following proposition, the strata considered are those 
  of $\phi^k$. 
 (These may be smaller than the strata of $\phi$; consider 
 the periodic case for example.) 
  
  \begin{proposition}\label{power} There exists a positive 
   integer $k$ 
  such that $\phi_0:=\phi^k$ has the following properties: 
\begin{enumerate} 
\item[1.]       Each letter $x$ appears in its own image under $\phi_0$. 
\item[2.]       Each exponential letter $x$ appears 
 at least $3$ times in its own image under $\phi_0$. 
\item[3.]       For all $x$, each letter $y\in\supp(x)$ appears 
 in $\phi_0(x)$. 
\item[4.]       For all $x$ and all $j \geq 1$, the 
leftmost  and rightmost letters of $\phi_0^j(x)$ 
are the same as those of $\phi_0(x)$. 
\item[5.]  For all $x$, all $j\geq 1$ 
 and all strata $\S\subseteq\supp(x)$, 
  the leftmost  (respectively, 
  rightmost) letter 
 from $\S$ in the reduced word $\phi_0^j(x)$ is the same as 
  the leftmost (resp. 
 rightmost) 
 letter from $\S$ in $\phi_0(x)$. 
\end{enumerate} 
\end{proposition} 
 
\begin{proof} Items (1) to (3) can be seen as simple facts about 
positive integer matrices, read-off from the action of $\phi$ 
on the abelianization of $F$. 
(By definition $a_j\in \S(a_i)$ if and only if the $(i,j)$ 
 entry of some power 
of the matrix describing this action is non-zero.)

Assume that $\phi_1$ is a power of $\phi$ 
that satisfies (1) to (3). Note that (3) implies 
that the strata of $\phi_1$ coincide with those of any proper 
power of it.  
 
Replacing $\phi_1$ by a positive power if necessary, we may 
assume that if $\phi_1^j(x)$ begins with the letter $x$, for 
any $j\ge 1$, then $\phi_1(x)$ begins with $x$. This ensures  
that {\em{$[y\preceq_L x$ if some $\phi^j(x)$ begins with $y]$}} 
 is a partial 
ordering, for if $\phi_1^{j_k}(x_k)$ begins with $x_{k+1}$ for 
$k=1,\dots,r$ and if $x_{r+1}=x_1$,  
then $\phi_1^{\Sigma j_k}(x_1)=x_1$ 
and hence $x_1=x_2=\dots = x_r$. 

If $\phi_1(x)$ begins with $z$ then $z\preceq_L x$, so 
by raising $\phi_1$ to a suitable power 
we can ensure for all $x$ 
that $\phi_1(x)$ begins with a 
letter that is $\preceq_L$-minimal. The $\preceq_L$-minimal 
letters $y$ are precisely those such that $\phi_1(y)$ begins with $y$. 
An entirely similar argument applies to the relation 
 {\em{$[y\preceq_R x$ if some $\phi^j(x)$ ends with $y]$}}. 
 This proves (4). 
 
Now assume that $\phi_0$ satisfies (1) to (4). The assertion  
in (5) concerning leftmost letters from $\S$ is clear 
for those $x$ where  $\phi_0(x)$ begins with $x$. If $\phi_0(x)$ begins
with $y\neq x$, then either $\S\subset\supp(y)$ 
or else the occurrences 
of letters from $\S$ in $\phi_0^j(x)$ are 
in 1-1 correspondence with the occurrences in the  
image of $\phi_0^j(x)$ in $F\lceil y\rceil $.  (Notation of
Remark \ref{induct}.) In the latter case,
arguing by induction on the size of $\pre(y)$ we 
may assume that  the induced automorphism   $\lceil \phi_0\rceil_y 
:F\lceil y\rceil\to F\lceil y\rceil $ has 
the property asserted in (5); the desired conclusion  for $\phi_0^j(x)$
is then  tautologous. In the former case, if 
we replace $\phi_0$ by $\phi_0^2$ then the conclusion 
becomes as immediate as it was when $\phi_0(x)$ began with $x$.  
 
An entirely similar argument applies to rightmost letters. 
\end{proof} 
  
\begin{remark} 
Although we shall have no need of it here, it seems worth 
recording that item (5) of the above proposition 
remains true if one replaces strata $\Sigma \subset \supp(x)$ 
by supports $\supp(y)\subset\supp(x)$. 
\end{remark} 
 
\smallskip

 \begin{quote}{\em We now fix an automorphism $\phi=\phi_0$ and assume that 
 is satisfies conditions (1)-(5) above. 
 All of the constants discussed in the sequel 
 will be calculated with respect to this $\phi$.} 
\end{quote} 
 
\section{Preferred Futures, Fast Letters and Cancellation}  
\label{PrefFutSec} 
 
Having conditioned our automorphism appropriately, we 
are now in a position to analyse the fates of (blocks of) edges 
 as they evolve in time. 
  
 \begin{definition}[Preferred futures]\label{pref-fut} 
For each element $x\in\{a_1,\dots,a_n\}$ of the basis,  
 we choose an occurrence of $x$ in the reduced word 
 $\phi(x)$ to be the (immediate) 
\index{future!preferred}{\em preferred future of $x$}: 
   if $x$ is a 
  parabolic letter, there is only one possible choice; 
  if $x$ is an  
  exponential letter, 
   we choose an occurrence of $x$ that is neither  
  leftmost nor rightmost (recall that we have 
  arranged for $x$ to appear 
 at least three times in $\phi(x)$). More generally, we 
 make a recursive definition of the {\em preferred future 
 of $x$ in $\phi^n(x)$}: 
 this is the occurrence of $x$ in $\phi^n(x)$ that 
 is the preferred future of the 
 preferred future of $x$ in $\phi^{n-1}(x)$.  
  
 The above definition distinguishes an edge $\e_1$  on the top of 
 each 2-cell in our diagram $\Delta$, namely the edge 
 labelled by the preferred future of the label at 
 the bottom $\e_0$ of the 2-cell. We define $\e_1$ to 
 be the (immediate) 
  {\em preferred future} of $\e_0$. As with letters, 
 an obvious recursion then defines a preferred future of $\e_0$ 
 at each step in its future (for as long as it continues 
 to exist). 
  
 Note that $\e_0$ has at most  one preferred 
 future at each time. (It has exactly one until a preferred 
 future dies in a corridor, 
 lies on the boundary, or hits a singularity.) 
  
 If the bottom edge of a 2-cell is $\e_0$, then we define 
 the preferred future of that 2-cell at time $t$ to be the unique 2-cell 
 at time $t$ whose 
 bottom edge is the preferred future of $\e_0$. 
 \end{definition}

\subsection{Left-fast, constant letters, etc.} \label{ss:speed}
In this paragraph, we divide the letters $x\in\{a_1^{\pm 1},\dots,a_m^{\pm 1}\}$ into 
classes according to the growth of the words 
$\phi^j(x), j=1,2,\dots$, and divide the edges of $\Delta$ into 
classes correspondingly. 
\begin{enumerate} 
\item[$\bullet$] 
If $\phi(x) = x$ then $x$ is called a {\em constant 
 letter}.  
 \index{letter!constant, fast, slow}

\item[$\bullet$] If $x$ is a {\em non}-constant letter, then 
  the function $n\mapsto |\phi^n(y)|$ grows 
 like a polynomial of degree $d\in\{1,\dots,m-1\}$ or else as an exponential 
   function of $n$. 
 
\item[$\bullet$] Let $x$ be a non-constant letter. 
 If the distance between the preferred future of $x$ and the   
 beginning of the word $\phi^n(x)$ grows at least quadratically as 
 a function of $n$, we say that 
 $x$ is {\em left-fast}; if this is not the case, 
 we say that $x$ is 
 {\em left-slow}.  {\em Right-fast} and  
 {\em right-slow} are defined similarly. Note that $x$ is 
left-fast (resp. slow) if and only if $x^{-1}$ is right-fast (resp. slow). 
 
\index{letter!para-linear}
\item[$\bullet$] Let $x$ be a non-constant letter. If $\phi(x) = uxv$ (the shown occurrence of $x$ need not be the preferred future), where $u$ consists only of constant letters,   
 then we say that $x$ is {\em \lpl}. (We place no restriction on $v$; in particular 
 it may contain occurrences of $x$.) {\em Right para-linear} is defined 
 similarly. 
\end{enumerate}

\bd 
For \lpl letters, we define the \index{future!para-preferred}{\em (left) para-preferred future} 
 (pp-future) to be the left-most occurrence of $x$ in $\phi(x)$. 
  The (right) pp-future of a \rpl letter is defined similarly, and 
  edges in $\Delta$ inherit these designations from their labels. 
  
  (It is possible that a letter  
   might be both \lpl  and right para-linear, and in such cases the 
   left and right 
    pp-futures need not agree. But when we discuss pp-futures, 
    it will always be clear  
     from the context whether we are favouring the left or the  right.) 
\ed

The following lemma indicates the origin of the 
terminology `left-fast' (cf.~\cite[Lemma 4.2.2]{BFH}).  
(A slight irritation arises from the fact that 
there may exist letters $x$  such that $x$ is not left-fast but  
 $\phi(x)$ contains left-fast letters; this difficulty accounts 
 for a certain clumsiness in the statement of the lemma.)

\begin{lemma}\label{C_0} There exists a constant $C_0$ with the 
following property: if $x\in\{a_1,\dots,a_n\}$ is such that 
 $\phi(x)$ contains a left-fast letter $x'$ 
and if $UVx\in F$ is a reduced word with $V$  positive\footnote{i.e. no 
 inverses $a_j^{-1}$ appear 
in $V$} and $|V|\ge C_0$, 
then for all $j\ge 1$, the preferred 
future of $x'$ is not cancelled 
when one freely reduces $\phi^j(UVx)$. Moreover, 
$|\phi^j(UVx)|\to\infty$  as $j\to\infty$. 
\end{lemma} 
 
\begin{proof} We factorize the 
reduced word $\phi^j(x)$ as $Y_{x,j}x'Z_{x,j}$ to emphasise the 
placement of the preferred future of a fixed left-fast letter 
$x'$ from $\phi(x)$. The fact that $x'$ is left-fast implies that $j\mapsto |Y_{x,j}|$ grows at least quadratically.  
 
Fix $C_0$ sufficiently large to 
ensure that for each of the finitely many possible 
$x\in\{a_1,\dots,a_n\}$, the integer 
$|Y_{x,j}|$ is greater than $Bj$ whenever $j\ge C_0/B$, 
where $B$ is the bounded cancellation constant. 
 
The Bounded Cancellation Lemma assures us that during the 
free reduction of the naive expansion of $\phi(UVx)$, 
at most $B$ letters of the positive word $\phi(Vx)$ will 
be cancelled. At most $B$ further letters will be cancelled 
when the naive expansion of $\phi^2(UVx)$ 
is freely reduced, and so on. Since $V$ and $\phi$ are positive and 
$ |V| \ge C_0$, it follows that $\phi^j(V)$ will 
not be completely cancelled during the free reduction  
of $\phi^j(UVx)$ if $j\le C_0/B$. When $j$ reaches $j_0:=\lceil C_0/B\rceil$ the 
distance 
from the preferred future of $x'$ to the left end of 
the uncancelled segment of $\phi^j(Vx)$ is 
at least  $|Y_{x,j_0}|$, which is greater than $Bj_0$ and hence $C_0$. 
Repeating the argument with $Y_{x,j_0}$ in place of $V$, we conclude that 
the length of the uncancelled segment of $\phi^j(Vx)$ in $\phi^j(UVx)$ 
remains positive and goes to infinity  with $j$. 
\end{proof}

Significant elaborations of the previous argument will be developed in 
Section \ref{ConstantSection}.  
 
\begin{definition}[New edges, cancellation and consumption]\label{new} 
Fix a 2-cell in $\Delta$. One edge in the top of 
the cell is the preferred future of the bottom 
edge; this will be called  \index{edge!old}{\em old} and the 
remaining edges will be called \index{edge!new}{\em new}. (These 
concepts are unambiguous relative to a fixed 2-cell or (folded) corridor, but `old edge' would be 
ambiguous if applied simply to a 1-cell of $\Delta$.) 
 
Two (undirected) edges $\e_1, \e_2$ 
in the naive top of a  
corridor are said to {\em cancel} each other if their images in the 
folded corridor coincide. If $\e_1$ lies to the 
left\footnote{Recall 
that corridors have a left-right orientation.} of 
$\e_2$, we say that $\e_2$ has been \index{edge!cancelled}cancelled {\em 
from the left} and $\e_1$ has been cancelled {\em 
from the right}. 
If $\e_1$ is the preferred future of an edge $\e$ 
in the bottom of the corridor and $\e_2$ is a new 
edge in the 2-cell whose bottom is $\e'$, then we 
say that $\e'$ has \index{edge!consumed}{\em (immediately) consumed} $\e$ 
{\em from the right}. `Consumed 
from the left' is defined similarly. 
 
Let $e$ and $e'$ be edges in $\bot(S)$ for some 
corridor $S$, with $e$ to the left (resp. right) of $e'$. 
If an edge in the future of $e$ 
cancels a preferred future of $e'$, then we say 
that $e$ {\em eventually consumes} $e'$ {\em from 
the left (resp. right).}  
\end{definition} 
\begin{lemma} \label{NoOldCanc} 
No pair of old edges can cancel each other. 
\end{lemma} 
 
\begin{proof} 
Suppose that 
two old edges in the naive top of a corridor $S$ 
are labelled $x$ and cancel each other.  These 
edges are the  preferred futures of edges on $\bot(S)$ 
that bound an arc $\alpha$ labelled by a reduced word 
 $x^{-1}wx$. 
Consider the freely-reduced factorisation $\phi(x) = uxv$ where 
the visible $x$ is the preferred future. 
The arc in the naive top of $S$ corresponding 
to $\alpha$ is labelled  
 $v^{-1}x^{-1}u^{-1}Wuxv$, where $W$ is the naive 
 expansion of $\phi(w)$. The old edges that we are considering 
 are labelled by the visible occurrences of $x$ in this word and 
 our assumption that these edges cancel means that the subarc 
 labelled $x^{-1}u^{-1}Wux$ becomes a loop (enclosing a 
 zero-area sub-diagram) in the diagram $\Delta$. 
  
 But this is impossible, because $x^{-1}wx$ is freely reduced, 
which means that $W$ is not freely equal to the empty 
 word, and hence neither is $x^{-1}u^{-1}Wux$. 
 \end{proof}

 \begin{corollary}\label{parabolicC} 
  An edge labelled by a 
 parabolic letter $x$ 
 can only be consumed by an edge labelled $y$ with   
 $\supp(x)$ strictly contained in $\supp(y)$. 
 \end{corollary} 
  
\begin{remark} 
A non-constant letter can only be (eventually) consumed from the left (resp. right) by a right-fast  
(resp. left-fast) letter. 
\end{remark}

\begin{remark} The number of old letters in 
 the naive top of a corridor $S$ is $|S|$, so 
 the length of corridors in the future of $S$ 
 will grow relentlessly unless old letters are 
 cancelled by new letters or the corridor hits a 
 boundary or a singularity. 
 \end{remark}

An obvious separation argument provides 
us with another useful observation concerning cancellation:

\begin{lemma}\label{perfect} Let $\e_1,\ \e_2$ and $\e_3$ be three 
(not necessarily adjacent) edges that appear in 
order of increasing subscripts as one reads from 
left to right along the bottom of a corridor. If 
the future of $\e_2$ contains an edge of $\partial\Delta$ 
or of a singularity, then no edge in the future of 
$\e_1$ can cancel with any edge in the future of $\e_3$. 
\end{lemma}

\section{Counting Non-constant Letters} \label{NonConstantSection} 
\index{letter!non-constant}
In this section we fix a corridor $S_0$ in $\Delta$ and 
bound the contribution of non-constant letters to the 
length of $\bot(S_0)$.

\subsection{The first decomposition of $S_0$}\label{decomp} 
 \index{corridor!first decomposition of}
Choose an edge $\e$ on the bottom of $S_0$.  As we follow the 
preferred future of $\e$ forward one of the following (disjoint) events must 
occur: 
 
\begin{enumerate} 
\item[1.] The last 
  preferred future of $\e$ lies on the boundary of  
 $\Delta$. 
  
\item[2.]  The last 
  preferred future of $\e$ lies in a singularity. 
 
\item[3.]   The last 
  preferred future of $\e$ dies in a corridor $S$ (i.e. 
  cancels with another edge from the naive top of $S$). 
\end{enumerate} 
 
We shall bound the length of $S_0$ 
by finding a bound on the number of edges in each of these three 
cases.  
 
We divide Case (3)  into two sub-cases:

\begin{enumerate} 
\item[3a.] 
 The preferred future of $\e$ dies when it is cancelled by an edge 
that is not in the future of $S_0$.  
 
\item[3b.] 
The preferred future of $\e$ 
dies when it is cancelled by an edge 
that is in the future of $S_0$.  
\end{enumerate} 
 
\subsection{Bounding the easy bits} \label{EasyBounding} 
 
Label the sets of edges in $S_0$ which fall into the above classes 
$S_0(1), S_0(2), S_0(3a)$ and $S_0(3b)$ respectively.   
 We shall see that $S_0(3b)$ is by far the most troublesome 
 of these sets. 
    
The first of the bounds in the following lemma is obvious, and the 
second follows immediately from  Proposition \ref{SingularityProp}.

\begin{lemma}\label{bound1and2} 
 $|S_0(1)| \leq \n \text{   \rm{and}   } 
  |S_0(2)| \leq 2B\n$.    
  \end{lemma}

\begin{lemma}\label{bound3a} $|S_0(3a)|\le B\n$. 
\end{lemma} 
 
\begin{proof} The preferred future of each $\e\in S_0(3a)$ 
dies in some corridor in the future of $S_0$. Since 
there are less than $\n /2$ corridors, we will be done 
if we can argue that the preferred future of at most $2B$  
such edges can die in each corridor $S$. 
 
Lemma \ref{Connected} tells us that the future of $S_0$ 
intersects $S$ in a connected region, the 
bottom of which is an interval $I$. The Bounded Cancellation Lemma 
assures us that only the edges within a distance $B$ of the 
ends of $I$ can be consumed in $S$ by an edge from 
outside the interval. And by definition, if a preferred future 
of an edge from  $S_0(3a)$ is to die in $S$, then it must 
be consumed by an edge from outside $I$. 
\end{proof} 
 
We have now reduced Theorem \ref{BoundS} to 
the problem of bounding $S_0(3b)$, 
i.e. of understanding cancellation {\em within} the future of $S_0$. 
This will require a great deal of work. As a first step, 
we further decompose 
$S_0$, mingling the above decomposition based on the fates 
of preferred futures of  
edges with the natural decomposition of $S_0$ into 
colours, as defined in Definition \ref{defMu}.

\subsection{The chromatic decomposition of $S_0$} \label{chromatic} 
\index{corridor!chromatic decomposition of}
We fix a colour $\mu$ and 
write $\mu(S_0)$ to denote the interval of $\bot(S_0)$ 
consisting of edges coloured $\mu$.  
We shall abuse terminology to the 
extent of referring to $\mu(S_0)$ as {\em a colour}, evoking 
the mental picture of the 2-cells in $S_0$ being painted 
with their respective colours. (Recall 
that the 2-cells of $S_0$ are in 1-1 correspondence with 
the edges of $\bot(S_0)$.) 
 
We shall subdivide $\mu(S_0)$ into five subintervals 
according to the fates  
of the preferred futures of edges. To this end,  
we define  $l_{\mu}(S_0)$ to be the rightmost edge in $\mu(S_0)$   
whose immediate future contains  a left-fast edge that is 
ultimately consumed 
from the left by an edge of  $S_0$, and we define 
 $A_1(S_0,\mu)$ to be the set of edges in $\bot(S_0)$ 
 from the left end of $\mu(S_0)$ to $l_{\mu}(S_0)$, inclusive. 
We define $A_2(\mu,S_0)\subset\mu(S_0)$ 
 to consist of the remaining 
 edges in $\mu(S_0)$ whose 
 preferred futures  are  ultimately consumed 
from the left by an edge of  $S_0$. 
 
Similarly, we define $r_{\mu}(S_0)$ to be 
 the leftmost edge $\mu(S_0)$ that has a  right-fast edge in its immediate future 
 that is ultimately consumed 
from the right by an edge of  $S_0$, and we define 
 $A_5(S_0,\mu)$ to be the set of edges in $\bot(S_0)$ 
 from the right end of $\mu(S_0)$ to $r_{\mu}(S_0)$, inclusive. 
We define $A_4(\mu,S_0)\subset\mu(S_0)$ 
 to consist of the remaining 
 edges in $\mu(S_0)$ whose 
 preferred futures  are  ultimately consumed 
from the right by an edge of  $S_0$. 
 
Finally, we define $A_3(S_0,\mu)$ to be the 
 remainder of the edges in $\mu(S_0)$.

\begin{figure}[htbp] 
\begin{center} 
  
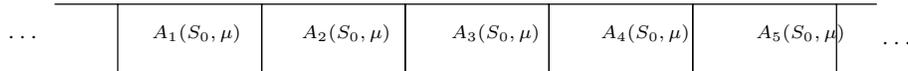 
 
\caption{The second decomposition of $S_0$} 
\label{figure:A1A5} 
\end{center} 
\end{figure} 
 
\medskip

Modulo the fact that any of the $A_i(S_0,\mu)$ 
might be empty, Figure 
10 is an accurate portrayal of  $\mu$: 
the $A_i(S_0,\mu)$ are  connected 
and they occur in 
ascending order of suffix from left to right. 
 
The chromatic decomposition of $S_0$ is connected to the 
decomposition of Subsection \ref{decomp} by the equality 
in the following lemma, which is a tautology. The  
inequality in this lemma is a restatement 
of Lemmas \ref{bound1and2} and \ref{bound3a}. 
 
\begin{lemma} \label{A3Lemma} 
$$ 
 \bigcup_{\mu} A_3(S_0,\mu) = S_0 \ssm S_0(3b)\ \ \  
 \text{  {\rm{and}}  }\ \ \  
 \sum_{\mu}|A_3(S_0,\mu)| \le \left({3B} + 1\right)\n . 
 $$ 
 \end{lemma}    
 
Thus the following lemma is a step towards bounding the 
size of $S_0(3b)$. 
 
\begin{lemma} \label{A1A5Lemma} 
 
$$ 
|A_1(S_0,\mu)|  \leq  C_0 \ \ \  
 \text{  {\rm{and}}  }\ \ \  
|A_5(S_0,\mu)|  \leq  C_0. 
$$ 

\end{lemma} 
 
\begin{proof} 
We prove the result only for $A_1(S_0,\mu)$; 
the proof for $A_5(S_0,\mu)$ is entirely similar.

As in 
Lemma \ref{perfect}, we know that the entire  future of the edges of  
$A_1(S_0,\mu)$ to the left of $l_\mu(S_0)$ must 
eventually be consumed from the left  by edges of $S_0$. This means 
that we are essentially in the setting of Lemma \ref{C_0}, with  
$l_{\mu}(S_0)$ in the role 
of $x$ and $A_1(S_0,\mu)$ in the role of $Vx$.

Thus if the length of $A_1(S_0,\mu)$ were greater than $C_0$,  
then we would conclude that  no  left-fast edge in the  
immediate future of  $l_{\mu}(S_0)$ would be cancelled from the left by 
an edge  of $\bot(S_0)$, contradicting the definition of $l_{\mu}(S_0)$.   
\end{proof} 
 
\begin{corollary} 
$$ 
\sum_{\mu}|A_1(S_0,\mu)| \, \leq \, C_0\n 
 \ \ \  
 \text{  {\rm{and}}  }\ \ \  
\sum_{\mu}|A_5(S_0,\mu)| \, \leq \, C_0\n . 
$$ 
\end{corollary} 
 
\subsection{A further decomposition of $A_2(S_0,\mu)$ 
 and $A_4(S_0,\mu)$} \label{ss:furtherdecomp}

It remains to bound $A_2(S_0,\mu)$ and $A_4(S_0,\mu)$.  
We deal only with $A_4(S_0,\mu)$, the argument for $A_2(S_0,\mu)$ 
being entirely similar. 
 
\begin{figure}[htbp] 
\begin{center} 
  
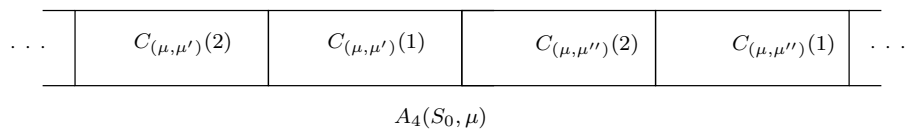 
  
\caption{$C_{(\mu,\mu')}(1)$ and $C_{(\mu,\mu')}(2)$.} 
\label{figure:Cmumu} 
\end{center} 
\end{figure}

First partition $A_4(S_0,\mu)$ into subintervals $C_{(\mu,\mu')}$ 
that consist of edges  that are eventually consumed by edges of a specified 
colour $\mu'$. Then partition $C_{(\mu,\mu')}$ into two subintervals:  
$C_{(\mu,\mu')}(1)$ begins at the   
right of  $C_{(\mu,\mu')}$ 
and ends with the last non-constant edge;  
 $C_{(\mu,\mu')}(2)$ consists of the remaining (constant) edges. 
See Figure \ref{figure:Cmumu}. 
\medskip

In the course of this section we will bound the size of the intervals $C_{(\mu,\mu')}(1)$ and during the following four sections we bound the sum over all pairs $(\mu,\mu')$
of the  sizes of the intervals $C_{(\mu,\mu')}(2)$ to get the desired bound on $|S_0(3b)|$.  In order to control this sum, we have 
to address the question of which colours can be adjacent. 
 
\subsection{Adjacent Colours}

In Corollary \ref{muConn} we saw that in any corridor 
$S$, the edges in $\bot(S)$ of a fixed colour form an interval. 
We say that two distinct colours $\mu$ and $\mu'$ are \index{colours!adjacent}
{\em adjacent} in $S$ if the closed intervals 
$\mu(S)$ and $\mu(S')$ 
have a common endpoint in $\bot(S)$. (Equivalently, 
there is a pair of 2-cells in $S$, one  coloured $\mu$ and 
the other $\mu'$, that share an edge labelled $t$.) 
We write 
 ${\vecZ}$ to denote the set of ordered pairs  
$(\mu,\mu')$ such that  
 $\mu$ and $\mu'$ are adjacent in some corridor $S$ 
with $\mu(S)$ to the left of  $\mu'(S)$  in $\bot(S)$.  
 
\begin{lemma} \label{NoOfAdjacencies} 
$$ 
|{\vecZ}| < 2\n -3 .     
$$ 
\end{lemma} 
 
\begin{proof}  
 We shall express this proof in the language 
of the forest $\mathcal F$ introduced in Remark \ref{tree}. Suppose 
that $\mu$ and 
$\mu'$ are adjacent in $S$. 
In $S$ we can connect the centre 
of some 2-cell coloured $\mu$ to the centre of some 2-cell 
coloured $\mu'$ by an arc contained in the union of 
the pair of 2-cells. The union of this arc and the trees in $\mathcal F$ 
corresponding to the colours $\mu$ and $\mu'$ disconnects 
the disc $\Delta$; each of the other trees in $\mathcal F$ 
is entirely contained in a 
component of the complement, and the 
colours with trees in different components can 
never be adjacent in any corridor. 
 
We can encode adjacencies of colours by a chord diagram: draw 
a round circle with marked points representing the colours of 
$\Delta$ in the cyclic order that they appear in $\partial\Delta$, 
then connect two points by a straight line if the corresponding 
colours are adjacent in some corridor.  
 The final phrase of 
the preceding paragraph tells us that the lines in this 
chord diagram do not intersect in the interior of the disc. 
A simple count shows that since there are less than 
$\n$ colours, there are less than $2\n -3$ lines in this diagram. 
\end{proof}

\subsection{Non-constant letters in 
 $C_{(\mu,\mu')}$ that are not left-fast} 
 \label{NonConstantSubsect} 
 
We stated in the introduction that a careful analysis of 
van Kampen diagrams would allow us to reduce  Theorem \ref{MainThmPos} 
to the 
study of blocks of constant letters.  
In this section we achieve the last step of this reduction. 
 
\begin{lemma} \label{C1Lemma}  
There is a constant $C_1$ depending only on $\phi$ 
with the following property: 
 
Let $S$ be a corridor and let $\mu_1$ and $\mu_2$ be 
 colours that occur in $S$ with $\mu_1$ to the left of $\mu_2$ (but do not 
assume that $\mu_1(S)$ is adjacent to $\mu_2(S)$).  Let 
 $I\subset A_4(S,\mu_1)$ 
  be a sub-interval that satisfies the 
 following conditions 
\begin{enumerate} 
\item[1.] the left-most edge of $I$ is 
 non-constant \mbox{and } 
   
\item[2.]  the preferred future of each edge in $I$ 
 is eventually consumed by an edge of $\mu_2(S)$.\\ 
\end{enumerate} 
\noindent Then $|I| \leq C_1$. 
In particular, $|C_{(\mu,\mu')}(1)| \leq C_1$ 
for all $(\mu,\mu') \in {\vecZ}$. 
 
It suffices to take $C_1 = 2mB^2$, where $m$ is the rank of $F$, and 
$B$ is the constant from the Bounded Cancellation Lemma. 
\end{lemma}

\begin{proof} 
  The region $I$ being considered contains no edge with a right-fast  
letter in the $\phi$-image of its label. 
 Since all exponential letters 
 are both left-fast and right-fast, all non-constant edges in the future of $I$ are  parabolic.   
 
We begin the argument at the stage in time where $\mu_2$ 
starts cancelling $I$.  For 
notational convenience we assume that this time 
is in fact $\height(S)$. (If it is not, then the  
fact that the length of $I$ may 
have increased in passing from $\time(S)$ to 
this time adds greater strength to the 
bound we obtain.) 
  
 We focus on the leftmost 
 edge $\e_0$ of $I$ that is labelled 
 by a non-constant letter $x$ for which $\supp(x)$  is maximal  
 among the supports of all edge-labels 
 from $I$ (with respect to inclusion). 
 Let $y$ be the label on the edge 
 $\e_0'$ of $\mu_2(S)$ that 
 eventually consumes $\e_0$ (oriented as shown in Figure \ref{C0Pic}). 
 Note that $\supp(x)$ is strictly contained in $\supp(y)$, 
 by Corollary \ref{parabolicC}. If $\e_0'$ consumes $\e_0$ immediately, 
 then the Bounded Cancellation Lemma tells us that 
$\e_0$ is a distance less than $B$ from the righthand end 
of  $I$. If not, then we 
 proceed one step into the future\footnote{proceeding one 
step into the future also allows us to assume that there 
are no letters coloured $\mu_1$ to the right of $I$} 
 and  appeal to the 
 conditioning done in Proposition \ref{power}(5) to assume 
 that for all $j\ge 1$, the rightmost letter in $\phi^j(y)$ 
 whose support includes $x$ is $y$. 
 We shall call the edge in the future of $\e_0'$ 
carrying the rightmost $y$  
 the \index{future!highlighted}{\em highlighted} future of $\e_0'$ (perhaps it is not 
 the preferred future).  

\begin{figure}[htbp] 
\begin{center} 
  
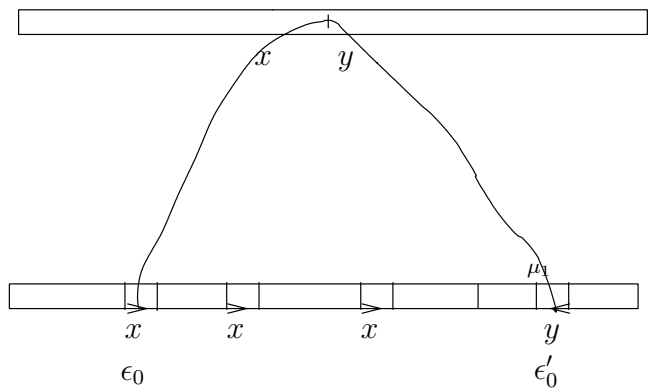 
  
\caption{The edge labelled $\epsilon_0'$ will eventually consume
$\epsilon_0$.}  
\label{C0Pic} 
\end{center} 
\end{figure} 
  
  The first important point to observe is that  
  the maximality of $\supp(x)$ ensures that 
  there will never be any new edges labelled $x$ 
  in the future of $I$ 
  (`new' in the sense of \ref{new}). 
   
 The second important point to note is that  
 the edges labelled $x$ in the future of 
 $\e_0'$ that are to cancel with 
 the futures of the 
 edges labelled $x$ in $I$ must all lie to the 
 left of the highlighted future of $\e_0'$. The point here is 
 that the highlighted future of $\e_0'$ 
 cannot be cancelled by an edge of $I$ (by the maximality 
 of $x$), and in order for it to be cancelled from the 
 other side, all the edges to 
 its right labelled $x$ 
 would have to be cancelled first, which would mean that they too 
 were cancelling with something not in the future of $I$. 
  
 We now come to the key observation of the proof: at 
 each stage $j$ steps into  the  
  future of $S$, the leftmost\footnote{we 
  have already noted that this is to 
  the left of the highlighted future of $\e_0'$} 
   edge $\e_j'$ in the future of $\e_0'$ 
  that is labelled  
  $x$ must be cancelled by an edge from 
  the future of $I$ {\em immediately}, i.e. in the corridor 
  where it appears at  $\height(S)+j$. 
   Indeed if this 
  were not the case, then  $\e_j'$ would develop a preferred 
  future which, being an old  
   edge (in the sense of Definition \ref{new}), could only  
  cancel with a  new edge (Lemma \ref{NoOldCanc}) 
  in the future of $I$. And since 
  we have arranged that there be no new edges labelled $x$, 
  the preferred future of $\e_j'$ would never cancel with 
  an edge in the future of $I$. But this cannot be, because 
  the continuing existence of a preferred future for $\e_j'$ would prevent 
  anything to its {\em right} consuming  an 
  edge in the future of $I$, and the penultimate sentence in the 
  third paragraph of this proof implies that no new 
edges labelled $x$ will ever appear to its {\em left} in the future of $\e_0'$. 
Thus if $\e_j'$ is not 
  cancelled immediately then we have a  contradiction 
  to the fact that $\e_0'$ must 
 eventually  consume $\e_0$.

  We have just proved that at $\height(S)+j$ the 
  edge $\e_j'$ must cancel with the preferred 
  future of an edge $\e_j$ in $I$ that is labelled $x$.  
  According to 
  the Bounded Cancellation Lemma, the preferred 
  future of $\e_j$ at $(\height(S)+j-1)$ must lie within 
  a distance $B$ of the right end of the future of $I$.  
  Since there is  no cancellation within the 
  future $I$, an iteration of this argument shows that 
   for as 
  long as there exist edges labelled $x$ in the future 
  of $I$, each successive pair of these edges is separated 
  by less than $B+|\phi(y)|\le 2B$ edges at each moment in time, 
  and the rightmost must be within a distance $B$ of the 
  right end of the future of $I$. 
   
  But since $\phi(x)$ contains at least 
  one letter other than  the preferred future of $x$, 
  it follows that there cannot be a pair 
  of  edges of $I$ labelled $x$ that remain unconsumed 
  at  $\height(S)+2B$, for otherwise 
  they would have grown a distance more than 
  $2B$  apart, contradicting the 
  conclusion of the previous paragraph. And proceeding 
  one more step into the future, the last edge labelled $x$ 
  must be consumed.  
   
  Since at most 
  $B$ letters of $I$ are cancelled at the right  
   at each stage in its future, all of the edges of $I$ labelled $x$ 
  are within a distance less than $2B^2$ of the right end of $I$, 
  and they are all consumed when $I$ has flowed $2B$ steps 
  into the future. 
  If no non-constant edges remain in the future of $I$ 
  at this stage, then we know 
  that $|I|\le 4B^2$. 
   
  If there do remain non-constant edges, we take the maximal interval of the 
   future of $I$ at   $\height(S)+2B$ whose leftmost 
   edge is non-constant, and we repeat the argument. (This 
   interval is obtained from the complete future of $I$ by 
   removing a possibly-empty collection of constant edges 
   at its left extremity.) 
    
  We proceed in this manner. The interval that 
  we begin with at each iteration has strictly fewer 
  strata than the previous one 
  and therefore the procedure 
  stops before $m=\text{\rm{rank}}(F)$ iterations. At 
  the time when it stops (at most $\height(S)+2mB$), the future 
  of $I$ has been cancelled entirely, except possibly for 
  a block of 
   constant edges at its left extremity.  
   With one final appeal to the Bounded Cancellation 
   Lemma, we deduce that  $|I|\le 2mB^2$. 
  \end{proof}

  \begin{corollary}\label{C1Corollary} 
\[ \sum_{(\mu,\mu') \in {\vecZ}}|C_{(\mu,\mu')}(1)| < 2C_1\n . \] 
\end{corollary} 
 
\begin{proof} This follows immediately 
 from Lemmas \ref{NoOfAdjacencies} and \ref{C1Lemma}. 
\end{proof}

\section{The Bound on $\sum\limits_{\mu \in S_0}|A_4(S_0,\mu)|$ and $\sum\limits_{\mu \in S_0}|A_2(S_0,\mu)|$}\label{A4sec} 

The sum of our previous arguments has reduced us to the nub of the
difficulties that one faces in trying to prove the Theorem \ref{MainThmPos},
namely the possible existence of large blocks of constant letters in
the words labelling the bottoms of corridors.
Now we must obtain a bound on
 $$\sum\limits_{(\mu,\mu') \in {\vecZ}} |C_{(\mu,\mu')}(2)| 
$$ 
that will enable us to bound 
$\sum\limits_{\mu \in S_0}|A_4(S_0,\mu)|$ and\footnote{In 
practice 
 we  only need concern ourselves with $A_4$, the arguments for $A_2$ 
being entirely similar}   $\sum\limits_{\mu \in 
S_0}|A_2(S_0,\mu)|$ by a linear function of $\n$.  These are the final 
estimates required to complete the proof of  Theorem  \ref{MainThmPos} --- see 
Section \ref{summary} for a r\'esum\'e of the proof.  

 The regions $C_{(\mu,\mu')}(2)$ are static, in the sense 
that they do not change under iteration by $\phi$, so  
the considerations of future growth 
that helped us so much in previous 
sections cannot be brought to bear directly. Rather, we must 
analyse the complete history of blocks of constant letters, 
understand how large blocks come into existence, and 
 use global considerations to limit the sum of 
 the sizes of all such blocks. 
  
Because of the global 
nature of the arguments, 
 we shall not obtain bounds on the sizes of the individual sets $C_{(\mu,\mu')}(2)$.
 Instead, 
we shall identify an associated block of 
constant letters
elsewhere in 
the diagram  (a ``team") that is amenable to a delicate string of 
balancing arguments that facilitates a bound on a union of 
associated regions $C_{(\mu,\mu')}(2)$. 

Our strategy is motivated by the following considerations. 
Believing Theorem \ref{BoundS} to be true, we seek 
payment from the global geometry of $\Delta$ to compensate 
us for having to handle the troublesome blocks of constant 
edges $C_{(\mu,\mu')}(2)$;  the currencies of payment are \index{colours!consumed}
{\em consumed colours} 
and dedicated subsets of 
edges on $\partial\Delta$ --- since $\Delta$ can have at 
most $\n$ of each, if we prove that adequate payment is available
then our troubles will be bounded and Theorem  \ref{MainThmPos}
will follow.  
The chosen currencies are apposite
because, as we shall  see in Section \ref{ConstantSection}, 
a large block of edges labelled by constant letters can only 
come into existence if  a colour (or colours) associated to a 
 component of this block in the past was consumed completely, 
  or else the boundary of $\Delta$ intruded into the past of 
   the block (or else something nearby) causing smaller regions of constant edges to elide.

In the remainder of this section we shall explain how various estimates on the behaviour of 
blocks of constant letters in $\Delta$ can be combined to obtain 
the bounds that we require on  
$\sum\limits_{\mu \in S_0}|A_4(S_0,\mu)|$ and  $\sum\limits_{\mu \in 
S_0}|A_2(S_0,\mu)|$. We hope that this explanation will provide the 
diligent reader with a useful road map and sufficient motivation to 
sustain them through the many technicalities needed to establish the 
estimates in subsequent sections.

In the following proposition, $L$ is the maximum length of the 
images $\phi(x)$ of the basis elements of $F$, while
$\ttt$ is the constant from the Pincer Lemma \ref{PincerLemma}, and
$C_1$ is the upper bound on the lengths of the intervals
$C_{(\mu,\mu')}(1)$ from Lemma \ref{C1Lemma}, $T_0$ comes from the Two
Colour Lemma \ref{TwoColourLemma} and $C_4$ comes from Lemma
\ref{G34pics}. The constant $\ll$ is
defined above Definition \ref{NestingDef}, and $B$ is the Bounded Cancellation
constant from Lemma \ref{BCL}.
\smallskip

\noindent{\bf{The Constant $K_1$ is defined to be}}
$$
 \AFourC.  
$$ 

\begin{proposition}\label{SummaryLemma} 
\[      \sum_{\mu \in S_0}|A_4(S_0,\mu)| \leq K_1  \n . 
 \]  
\end{proposition}

\subsection{Dramatis Personae} 
The ``proof" that we are about to present is essentially a scheme for
reducing 
the proposition to a series of technical lemmas that will be proved
in Sections \ref{teamSec} and \ref{BonusScheme}. These lemmas are
phrased in the language associated to \index{teams}{\em teams}, the precise
definition of which will
also be given  in  Section \ref{teamSec}. 
Many of the proofs involve global cancellation arguments
based on the {\em Pincer Lemma}, which will be proved in the
next section. 
Intuitively speaking, a {\em team}  
(typically denoted $\T$) is a  contiguous 
region of $\|\T\|$ constant letters all of which 
are to be consumed by a  fixed \lpl edge (the {\em reaper}). Notwithstanding 
this intuition, it is preferable for 
technical reasons  to define a team to be a set of pairs of colours $(\mu,\mu')\in\vecZ$, 
where $\mu'$ is fixed and the different \index{teams!members}{\em members} of the team correspond to  
different values of $\mu$. We write $(\mu,\mu')\in\T$ to denote 
membership. Teams also have {\em virtual members}, denoted $(\mu,\mu') \vin \T$ (see 
Definition \ref{Virtual}). There are less than $2\n$ teams (Lemma \ref{allIn}).
 
Each pair $(\mu,\mu')$ with $C_{(\mu,\mu')}(2)$ non-empty 
is either a member or a  virtual member  of a team (Lemma \ref{allIn}).  
There are {\em short} teams (Definition \ref{newTeams}) and long teams, 
 of which some are {\em distinguished} (Lemma \ref{Aget2Lemma}). 
There are four types of  \index{teams!genesis of}{\em genesis} of a team, (G1), (G2), (G3) and
(G4) (see Subsection \ref{genesis}).  Teams of genesis (G3) have
associated to them a pincer $\Pin_{\T}$ (Definition \ref{pl}) yielding
an auxiliary set of colours
$\subT$. There is also a set of colours $\chi_P(\T)$ associated to the
time before the pincer $\Pin_{\T}$ comes into play.  For long,
undistinguished teams, we also need to consider certain sets
$\CT$ and $\chi_{\delta}(\T)$ of colours consumed in the past of $\T$ (see
the proof of Lemma \ref{Aget2Lemma}). Such teams may
also have  three sets of edges in $\partial\Delta$ associated to  
 them: $\partial^\T$, $\down_1(\T)$ and  
$\down_2(\T)$. An important feature of the definitions of 
 $\partial^\T$ and $\down_1(\T)$ is that the sets associated 
  to different teams are disjoint.   This disjointness is crucial 
   in  the following proof, where we use the fact that the sum 
    of their cardinalities is at most $\n$. Similarly, the disjointness of
the sets $\chi_c(\T)$ is used to estimate the sum of their cardinalities by
$\n$ and likewise for $\chi_{\delta}(\T)$ and $\chi_P(\T)$.

It is not necessarily true that the sets $\down_2(\T)$ are disjoint
 for different teams, but we shall explain how to account for the
 amount of `double-counting' that can occur (see Lemma
 \ref{Aget2Lemma}).
 
Associated to every team one has  the time \index{teams!times $t_i(\T)$} $t_1(\T)$ at 
which the reaper starts consuming the team (see Subsection
\ref{t1}). Teams of genesis (G3) also have two
earlier times $t_2(\T)$ and $t_3(\T)$ associated to them as well as an
auxiliary set of edges $\QT$, the definitions of which
are somewhat  technical (see Definition \ref{PincerDef} {\em et seq.}). 

In Section \ref{BonusScheme} we describe a {\em bonus scheme} that
assigns a set  of extra edges, $\bonus(\T)$ to each team.  These
bonuses are assigned so as  to ensure that $|\bonus (\T)|+\|\T\|$
 dominates  the sum of 
 the cardinalities of the sets $\cmm$ 
associated to the  members and virtual members of $\T$. 
 
\smallskip 
\noindent{\bf Proof of Proposition \ref{SummaryLemma}.}

Recall that   $A_4(S_0,\mu)$ is partitioned  into disjoint regions $C_{(\mu,\mu')}$ 
which in turn are partitioned into $C_{(\mu,\mu')}(1)$ and 
$C_{(\mu,\mu')}(2)$.  
 
Given any $\mu_1$ and $\mu_2$, at most one  ordering of $\{\mu_1,\mu_2\}$ can 
arise in  $S_0$. Thus Lemma \ref{NoOfAdjacencies} 
implies that there are less than $2\n$ pairs $(\mu,\mu')\in\vecZ$ with 
$C_{(\mu,\mu')}\subset\bot(S_0)$  non-empty. 
It follows immediately from this observation and Lemma \ref{C1Lemma} that 
\[      \sum_{(\mu,\mu') \in {\mathcal Z}}|C_{(\mu,\mu')}(1)| \leq 2C_1\n . 
\]

Lemma  \ref{Aget2Lemma} accounts for the set of  distinguished 
 long 
teams $\dlong$: 
\[      \sum_{{\mathcal T} \in \dlong}\sum_{(\mu,\mu') \in  
{\mathcal T}}|C_{(\mu,\mu')}(2)| \leq 6B\n(T_1+T_0).   \]  
For all other teams $\T$ we rely on Lemma \ref{C1toTeamLength} which
states
\begin{equation}\label{goodEq} 
\sum_{(\mu,\mu') \in \T \mbox{ \tiny or } (\mu,\mu') \vin \T} |C_{(\mu,\mu')}(2)| \le
\|\T\|  + |\bonus(\T)| + B. 
\end{equation} 
We next consider the {\em genesis} of teams. All teams of genesis (G4)
are short (Lemma \ref{G4lemma}). And by Definition \ref{newTeams} for
the short teams $\T\in\Sigma$ we have
\[      \sum_{{\T} \in \Sigma}\sum_{(\mu,\mu') \in  
{\mathcal T}}|C_{(\mu,\mu')}(2)| \leq 2\ll\n + \sum_{\T\in\Sigma}
\big(|\bonus(\T)| + B\big).   \]

Lemma \ref{TeamAgeLemma} tells us that for teams of genesis (G1) and
(G2) we have 
\[      \|\T\| \leq 2LC_4|\down_1(\T)| + |\partial^{\T}|,  \]
whilst for teams of genesis (G3) we have
\[      \|\T\| \leq 2LC_4\big{(}|\down_1(\T)| +|\QT|\big{)}+
T_0\big{(}|\chi_P(\T)| + 1\big{)} +|\partial^{\T}| + \ll.       \]

Let $\Gthree$ denote the set of teams of genesis (G3) with $\QT$
non-empty. In Definition \ref{down2} we break $\QT$ into pieces so
that
$$
|\QT|\ =  t_3(\T) - t_2(\T)  
+ |\down_2(\T)|.
$$  
Making crucial use of the Pincer Lemma, in Corollary \ref{t1-t2Corr} we prove that
$$
\sum\limits_{\T \in \Gthree} 
 t_3(\T) - t_2(\T)  \ \le \ 3\ttt \, \n, 
$$ 
and in   Corollary \ref{downbound} we prove that
$$ 
\sum\limits_{\T \in \Gthree}|\down_2(\T)|
\leq (2 + 3\ttt + 5T_0)\n.
$$
This completes the estimate on $|\QT|$ and hence $\|\T\|$. 

Section 10 is dedicated to the proof of Proposition  \ref{BonusBound}, which states
\[      \sum_{\text{\small{teams}}}|\bonus(\T)| \leq \big( \Bb \big)\n. \]

Adding all of these estimates and recalling that there are less than $2\n$ teams, we deduce:
\[      \sum_{\mu \in S_0}|A_4(S_0,\mu)| \leq K_1 \n ,    \]
where $K_1$ is
$$ 
 \AFourC.
$$
Thus the proposition is proved.
\hfill$\square$
\smallskip 
\begin{remark} The stated value of the constant $K_1$ 
  is an artifact of our proof: we
have simplified the estimates at each stage for the sake
 of clarity rather than trying to optimise the
constants involved. Nevertheless, we have made some effort
 to make the arguments constructive
so as to prove that there exists an algorithm to calculate  
the Dehn function of $F\rtimes_\phi\mathbb Z$ directly from $\phi$.  
\end{remark} 
By a precisely analogous argument, we also have:
\bp \label{A2Prop}
\[      \sum_{\mu \in S_0}|A_2(S_0,\mu)| \leq K_1 \n ,
 \] 
where $K_1$ is the constant defined prior to Proposition \ref{SummaryLemma}.
\end{proposition}

\section{The Pleasingly Rapid Consumption of Colours} \label{ConstantSection}

This section contains the  cancellation lemmas 
that we need to control the manner in which colours are consumed. 
The key result in this direction is the {\em Pincer Lemma} (Theorem \ref{PincerLemma}). 
   
\subsection{The Buffer Lemma} 
  
\index{Buffer Lemma}
\begin{lemma}\label{BufferLemma}  
Let $I\subset\bot(S)$ be an interval of edges labelled by constant letters, and 
suppose that the colours $\mu_1(S)$ and $\mu_2(S)$ lie either side 
of $I$, adjacent to it.  
 Provided that the whole of $I$ does not die in $S$,  no 
non-constant edge  coloured  $\mu_1$ will ever cancel with  
a non-constant edge coloured $\mu_2$. 
\end{lemma}  
 
\begin{proof} Suppose that the future of $I$ in $\top(S)$ 
is a non-empty interval labelled $w_0$. If $\mu_1(S)$ is to the left of $I$, 
then reading from the left beginning with the last non-constant 
edge coloured $\mu_1$, on the naive top of $S$ we have an interval labelled 
$x w_1 y$, 
where $y$ is a non-constant letter coloured $\mu_2$ and  
$w_1$ contains $w_0$ and perhaps some constant letters 
from $\mu_1$ and $\mu_2$.  
 
Our conditioning of $\phi$ (Proposition \ref{power}) ensures that, for all non-constant letters $z$, the rightmost non-constant letter in $\phi^j(z)$ is the same for all $j \geq 1$.  Therefore, in order for there to ever be cancellation between non-constant letters coloured $\mu_1$ and $\mu_2$, we must have $x = y^{-1}$.
Thus on $\top(S)$ there is an interval labelled  $xwx^{-1}$, where $w$ is the  
(non-empty) free-reduction 
of $w_1$. 
 
At times greater than $\height(S)$, the future of the interval that 
we are considering will continue to have a core subarc labelled $xw_jx^{-1}$,  
where $w_j$ is a conjugate of $w$ by a (possibly-empty) 
 word in constant letters (unless the interval hits a singularity or 
the boundary). In particular, no non-constant letters from $\mu_1$ 
and $\mu_2$ can ever cancel each other.  
\end{proof}

In the light of the Bounded Cancellation Lemma we deduce:

\begin{corollary} \label{BufferCorollary} 
Let $I\subset\bot(S)$ be an interval of edges labelled by constant letters, and 
suppose that the colours $\mu_1(S)$ and $\mu_2(S)$ lie either side 
of $I$, adjacent to it. 
If $|I|\ge B$ 
then there is never any cancellation between non-constant letters in $\mu_1$ and $\mu_2$. 
\end{corollary}

\subsection{The Two Colour Lemma} 
  
\index{neutering}
\bd \label{Neuters} 
Suppose that $U$ and $V$ are positive words\footnote{i.e. none of their letters are 
 inverses $a_j^{-1}$} 
 and that for some $k>0$ the only  negative exponents occurring   in $\phi^k(UV^{-1})$  
are on constant letters.  Then we say that $U$   
{\em $\phi$-neuters   $V^{-1}$ in at most $k$ steps}.  
\ed 
 
We shall also apply the term $\phi$-neuters to describe  the 
cancellation between  colours $\mu(S), \mu'(S) \subseteq \bot(S)$ that are adjacent in  corridors of van Kampen diagrams,  
and the following lemma remains valid in that context. 
 
\index{Two Colour Lemma}
\begin{prop}[Two Colour Lemma] \label{TwoColourLemma} 
There exists a constant $\tz$ depending only on $\phi$ so that 
 for all positive words $U$ and $V$, if  
$U$ $\phi$-neuters  $V^{-1}$ then it does so in at most $\tz$ steps.  
\end{prop} 
 
\begin{proof} We express 
 $V^{-1}$ as a product of  three subwords:  reading from the left of
$V^{-1}$, the first subword ends with the last letter $y$ such that
$\phi(y)$ contains a left-fast letter; the second subword follows the
first and ends with the last non-constant letter in $V^{-1}$; the
remainder of $V^{-1}$ consists entirely of constant letters.
 
Lemma \ref{C_0} tells us that the length of the first subword is less
than $C_0$, and the proof of Lemma \ref{C1Lemma} provides a  bound of
$C_1$ on the length of the second subword.
 
Now consider the freely reduced form of $\phi^k(UV^{-1})$, and let $v_k$ denote its subword 
that begins with the first letter of negative exponent and ends with the final non-constant 
letter. The argument just applied to $V^{-1}$ shows that $v_k$ has length less than $C_0+C_1$ 
for all $k\ge 0$. 
 
Suppose that $U$ $\phi$-neuters $V^{-1}$ in exactly $N$ steps, let $\alpha_{N-1}$ 
be the letter of $\phi^{N-1}(UV^{-1})$ that consumes the last  letter of $v_{N-1}$, and 
let $\alpha_k$ be the ancestor of $\alpha_{N-1}$ in $\phi^k(UV^{-1})$. Write 
$\phi^k(UV^{-1}) = w_k\alpha_k u_k v_k w_k'$. 
 
Lemma \ref{C_0} shows that $|u_k| < C_0$ for all $k< N$, and we have just argued  
that $|v_k| < C_0+C_1$. Thus we obtain a bound  (independent of $U$ and $V$) on 
the number of words $\alpha_k u_k v_k$ that arise as $k$ varies --- call this 
number $\tz$. If $N$ were greater 
than $\tz$, then some configuration $\alpha_k u_k v_k$ with $v_k$ non-empty would recur. But 
this is nonsense, because once there is this repetition, the words $v_k$ will continue to repeat, and thus $V^{-1}$ will never be $\phi$-neutered, contrary to assumption.  
\end{proof} 

\begin{corollary} There exists a constant $\tz'$, depending only on $\phi$, with the following 
property: if $U$ and $V$ are positive words, $V$ begins with a non-constant letter 
and $\phi^k(UV^{-1})$ is positive for some $k>0$, then the least such $k$ 
is less than $\tz'$. 
\end{corollary}  
 
\begin{proof} The preceding lemma provides an upper bound on the least integer $N$ such that 
$\phi^N(UV^{-1})$ contains no non-constant  letters with negative exponent. Up to 
this point, the rightmost non-constant letter in $\phi^k(UV^{-1})$ may have been spawning 
constant letters to its right, and thus  $\phi^k(UV^{-1})$ may have a terminal segment consisting 
of constant letters. Since the rightmost non-constant letter of $\phi^k(V^{-1})$ does 
not vary with $k$ when $k<N$ (by Proposition \ref{power}),  the length of this segment  
grows at  a constant rate ($<L$) during each application of $\phi$. Similarly, its length 
changes at a constant rate after time $N$, decreasing until it is eventually cancelled. 
 
Since $N\le \tz$, this segment of constant letters has length less than $L\tz$ 
at time $N$, and hence is cancelled entirely before time $T_0(L+1)$. 
\end{proof}

\subsection{The disappearance of colours: Pincers and implosions} 
 
\begin{figure}[htbp] 
\begin{center} 

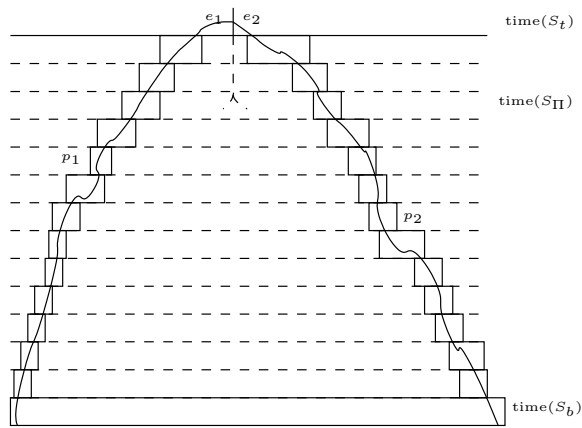 

\caption{A pincer.} 
\label{PincerPic} 
\end{center} 
\end{figure}

In this subsection we turn our attention to the detailed study of how non-adjacent colours 
along a corridor in $\Delta$ can come together solely as a result of the mutual 
annihilation of the intervening colours. Such an event determines a {\em pincer} (Figure \ref{PincerPic}), which is defined as follows.

\begin{definition}\label{pincerDef} 
Consider a pair of paths $p_1, p_2$ in $\mathcal F \subseteq \Delta$  
tracing the histories of $2$ non-constant edges $e_1, e_2$ that cancel in a corridor $S_t$.   
Let $\mu_i$ denote the colour of the 2-cells along  $p_i$. 
Suppose that at time $\tau_0$ these paths lie in a common corridor $S_b$.
Under these circumstances, we define \index{pincer}
the {\em pincer}  
$\Pin=\Pin( p_1, p_2, \tau_0 )$ to be the subdiagram of $\Delta$ enclosed by the chains of  
$2$-cells along $p_1$ and $p_2$, and the chain of $2$-cells connecting them in $S_b$. 
  
When it creates a desirable emphasis, we shall write $S_b(\Pin)$ and $S_t(\Pin)$ 
in place of $S_b$ and $S_t$. 
  
We define $S_\Pin$ to be the earliest corridor of the pincer in which $\mu_1(S_\Pin)$ and 
$\mu_2(S_\Pin)$ are adjacent. \index{pincer!colours $\chi(\Pi))$}
We define  $\tilde{\chi}(\Pin)$ to be the set of colours $\mu\notin\{\mu_1,\mu_2\}$ such 
that there is a 2-cell in $\Pin$ coloured $\mu$. And we define \index{pincer!life of}
$$\life(\Pin) = \time(S_\Pi)-\time(S_b). 
$$ 
\end{definition} 
 
\index{Pincer Lemma}
\begin{proposition}[Unnested Pincer Lemma]\label{prePincerLemma} 
There  exists a constant $\hat{T_1}$, depending only on $\phi$, such that for any pincer $\Pin$ 
\[      \life(\Pin) \leq \hat{T_1}(1+ |\tilde{\chi}(\Pin)|).        \] 
\end{proposition} 
 
Fix a pincer $\Pin$ and assume $\life(\Pin)\neq 0$.
 The idea of the proof of Proposition \ref{prePincerLemma} 
 is as follows:  we shall identify a constant $\hat{T_1}$ and argue 
 that if none of the colours $\mu\in\tilde{\chi}(\Pin)$ were consumed entirely  
  by  $\time(S_b) + \hat{T_1}$, 
the situation reached would be so stable  that no colours could be consumed in $\Pin$ at 
subsequent times,  
 contradicting the fact that all but $\mu_1$ and $\mu_2$ must be consumed by $\time(S_\Pi)$. 
 
With this approach in mind, we make the following definition: 
 
\begin{definition} Let $p$ be a positive integer.  
 A \index{implosive array}{\em $p$-implosive array} of colours  in a corridor $S$ 
is an ordered tuple $A(S)=[\nu_0(S),\dots,\nu_r(S)]$, with $r>1$, 
 such that: 
\begin{enumerate}
\item each  pair of colours $\{\nu_j, 
\nu_{j+1}\}$ is \index{colours!essentially adjacent}{\em essentially adjacent} in  $S$, meaning that there are no
non-constant edges of any other colour separating $\nu_j(S)$ from $\nu_{j+1}(S)$;
\item in each of the corridors $S=S^1,S^2,\dots, S^{p}$ in the future of $S$, 
 every $\nu_j(S^i)$ contains a non-constant edge;
\item in $S^p$,  {\em either} a non-constant edge coloured $\nu_0$ cancels a non-constant edge coloured $\nu_r$ 
(and hence the colours $\nu_j$ with $j=1,\dots,r-1$ are consumed 
entirely), {\em or else} all of the non-constant letters in $\nu_j(S^p)$, for $j=1,\dots, r-1$,
are cancelled in $S^p$ by edges from one of the colours of the array, while
$\nu_0(S^p)$ and $\nu_r(S^p)$ contain non-constant letters that survive in
the free-reduction of the naive future of the interval $\nu_0(S^p)\dots\nu_r(S^p)\subset\bot(S^p)$
(but  may nevertheless be cancelled in $S^p$ by edges from colours external to the array). 
\end{enumerate}
Arrays satisfying the first of the conditions in (3) are said to be of Type I, and those
satisfying the second condition are said to be of Type II. (These types are not
mutually exclusive.)

The {\em residual block} of an array of Type II is the interval of constant edges  between 
the rightmost non-constant letter of $\nu_0$ and 
the leftmost non-constant letter of $\nu_r$ in 
the free reduction of the naive future of  $\nu_0(S^p)\dots\nu_r(S^p)$.
The {\em enduring block} of the array  is the set of constant edges in 
$ \bot(S)$ that have a future  in the residual block.

 Note that there may exist {\em unnamed colours} between
$\nu_j(S)$ and $\nu_{j+1}(S)$ consisting entirely of constant edges.
\end{definition} 
 
\begin{remarks}\label{subarray} Let $[\nu_0(S),\dots,\nu_r(S)]$ be a $p$-implosive array.

\smallskip
(1) Any implosive subarray of $[\nu_0(S),\dots,\nu_r(S)]$ is $p$-implosive (same $p$).

(2) If an edge of $\nu_i$ cancels 
with an edge of $\nu_j$ and $j-i>1$, then this cancellation 
can only take place in $S^p$. If the edges cancelling are non-constant, 
 then the subarray $[\nu_i(S),\dots, \nu_j(S)]$ is   $p$-implosive of Type I.

(3) Given $x,y,w\in F$,
if the freely reduced words representing $x, y$ and $\phi(xwy)$ consist only 
of constant letters, then
so does the reduced form of $w$, since the subgroup generated by the constant
letters is invariant under $\phi^{\pm 1}$. It follows that the residual block
of any array of Type II contains edges from at most two of the colours $\nu_j$, and if
there are two colours they must be essentially adjacent, i.e. $\nu_j(S^p), \nu_{j+1}(S^p)$.

(4) For the same reason, the enduring block of an implosive
array of Type II is an interval involving at most two of the $\nu_j$, and if
there are two such colours then they must be essentially adjacent.
\end{remarks}

\begin{lemma}\label{haveImp}  The ordered list of colours along each corridor before 
$\time(S_\Pi)$ in a pincer $\Pi$ must contain an implosive array.
\end{lemma}

\begin{proof} At the top of the pincer there is cancellation between non-constant
edges. Lemma \ref{BufferLemma} tells us that before $\time(S_\Pi)$ the colours
of these edges must have been separated by a non-constant letter of a different
colour, hence the list of non-constant colours along the bottom of $S_\Pi$ is a
1-implosive array. This same list of colours defines an implosive array at
each earlier time in the pincer until, going backwards in time, further non-constant
colours appear. Suppose $\mu$ has non-constant letters in $\Pi$ at
time $t$ but not time $t+1$. Let $\nu_0$ be  the first colour to the  left of $\mu$ that
contains non-constant letters at time $t+1$, and let $\nu_r$ be the first such colour
to the right. If $S_t$ is the corridor at time $t$, then the list of essentially-adjacent
non-constant colours $[\nu_0(S_t),\dots,\mu(S_t),\dots,\nu_r(S_t)]$ is a
1-implosive array. Furthermore, the array $[\nu_0(S_{t'}),\dots,\mu(S_{t'}),\dots,\nu_r(S_{t'})]$
is a $(t'-t+1)$-implosive array for each earlier time $t'$ until (going backwards in
time) either further non-constant colours appear or else we reach the bottom of the
pincer.
\end{proof}

 If, further to the above lemma, we can argue that there is a constant $\hat{T_1}$ such
that each corridor before $\time(S_\Pi)$ contains a
$p$-implosive array with $p\le\hat{T_1}$, then we will know that at least one of the colours 
from $\tilde\chi(\P)$ is {\em essentially consumed} (i.e. comes to consist of constant
edges only)
 during each interval of $\hat{T_1}$ units in time during 
the lifetime of the pincer. Thus Proposition \ref{prePincerLemma} 
 is an immediate consequence of the following result, which will be proved
in (\ref{RIP}).

\begin{proposition}[Regular Implosions]\label{implosion} 
There is a constant $\hat{T_1}$ depending only on $\phi$ 
such that every \index{implosive array}
 implosive array in any minimal area diagram $\Delta$ is $p$-implosive 
for some $p\le\hat{T_1}$. 
\end{proposition} 
 
The first restriction to note concerning implosive arrays is this: 
 
\begin{lemma} \label{OnlyBColours} 
If $[\nu_0(S),\dots,\nu_r(S)]$ is implosive 
of Type I, then $r\le B$. If
it is implosive of Type II, then $r< 2B$.
\end{lemma} 
 
\begin{proof} In Type I arrays,  the interval 
 $\nu_1(S^{p})\dots\nu_{r-1}(S^p)\subset\bot(S^p)$
  is to die in $S^p$, so $r-1<B$ by
  the Bounded Cancellation Lemma.
  For Type II arrays, one applies
the same argument to the intervals 
 joining $\nu_0(S^p)$ and $\nu_r(S^p)$ to the residual block of constant letters. 
\end{proof} 
 
\begin{remark} \label{shortisenough} 
In the light of Lemma \ref{OnlyBColours}, an 
obvious finiteness argument would provide the bound required for 
Lemma \ref{implosion} if we were willing 
to restrict ourselves to implosive arrays  with
 a uniform bound on their
length.

Motivated by this observation, we seek to prove  that every implosive array contains an
implosive sub-array that is uniformly {\em short}.
\end{remark} 
 
\smallskip 
 
In order to identify a suitable notion of {\em short}, 
 we need to consider a further decomposition 
of the colours $\nu_j(S_b)$ in a $p$-implosive array $[\nu_0(S_b),\dots,\nu_r(S_b)]$. 
 
Previously (Subsection \ref{chromatic}) we partitioned each colour 
$\nu_j(S_b)$ into five intervals $A_1(S_b,\nu_j),\dots, 
A_5(S_b,\nu_j)$ and then further decomposed $A_4$ into subintervals 
$C_{(\nu_j,\nu')}(1)$ and $C_{(\nu_j,\nu')}(2)$ according to the 
colours of the edges that were going to consume these subintervals in 
the future. There is a corresponding decomposition of $A_2$ into 
intervals which we denote $C^2_{(\nu_j,\nu')}(1)$ and 
$C^2_{(\nu_j,\nu')}(2)$ (where $\nu'$ is now to the left of $\nu_j$ in 
$S_b$).  
 
\begin{figure}[htbp] 
\begin{center} 

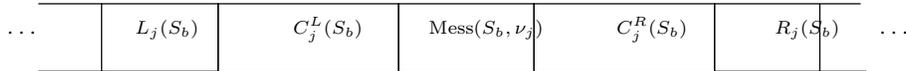 

\caption{The decomposition of the colour $\nu_j$} 
\label{PincerDecomp} 
\end{center} 
\end{figure}

Adapting to our new focus, we now define $R_j(S_b)=A_5(\nu_1,S_b)\cup 
C_{(\nu_j,\nu_{j+1})}(1)$, and $L_j(S_b)=A_1(\nu_1,S_b)\cup 
C^2_{(\nu_j,\nu_{j+1})}(1)$. We also define  $C_j^R(S_b)$ to be 
$C_{(\nu_j,\nu_{j-1})}(2)$ minus any edges from the excluded block, and $C_j^L(S_b)$ to 
be $C^2_{(\nu_j,\nu_{j-1})}(2)$ minus any edges from the excluded block.
Thus we obtain a decomposition of 
$\nu_j(S_b)$ into five intervals (see Figure \ref{PincerDecomp})  
$$ 
 L_j(S_b),\ C_j^L(S_b), \ \mess(S_b,\nu_j),\    C_j^R(S_b),\ R_j(S_b) 
$$ 
where $\mess(S_b,\nu_j)$ contains the edges 
whose preferred future dies at the time of  implosion together
with edges from the excluded block\footnote{At this point the reader may 
find it helpful to recall that only arrays of Type II have excluded
blocks, and such a block is either contained in a single colour,
or in adjacent colours $\nu_j(S_b)\cup\nu_{j+1}(S_b)$ with
the intervening intervals $R_j(S_b) \dots L_{j+1}(S_b)$ empty.}.

The terminal colours in our array, $\nu_0$ and $\nu_r$, play a special 
role. This is reflected in the fact that we shall only need to consider 
the segment of $\nu_0$ from its right end 
 up to and including the edge one to the left of $\mess(S_b,\nu_0)$. And 
 in $\nu_r$ we shall only need to consider the segment from its left end 
up to and including the edge one to the right of $\mess(S_b,\nu_r)$. 
 We write $\mathcal L(\nu_j,S_b)$ and $\mathcal R(\nu_j,S_b)$, 
  respectively, to denote these sub-intervals of $\nu_j(S_b)$. 
  
 \begin{definition} The length of $A(S)=[\nu_0(S),\dots,\nu_r(S)]$, written $\|A(S)\|$,
 is the number of edges in the interval
$\mathcal L(\nu_0,S)\dots\mathcal R(\nu_r,S)\subset \bot(S)$. (Note that $\|A(S)\|$ takes account of the
unnamed colours.)
\end{definition}

In keeping with the notation in the definition of $p$-implosive, we shall 
write $S^t$ for the corridor $t$ steps into the future of $S_b$; in particular  $S^0=S_b$ and
each $\nu_j$ with $j=1,\dots,r-1$ essentially vanishes in $S^p$. 
 
By definition, no preferred future of any 
 edge in  $\mess(\nu_j,S_b)$ is cancelled 
before $S^p$. Hence these intervals do not shrink in length before 
that time, and as in the proof of Lemma \ref{OnlyBColours} we can use
the Bounded Cancellation Lemma to bound the sum of their
lengths:

\begin{lemma} \label{2Bcols}
After excluding the edges of the enduring block, the sum of the lengths of the 
intervals  $\mess(\nu_j,S_b)$ is   at most $2B$. 
\end{lemma}

Combining this estimate with the bounds from Lemmas \ref{C_0} 
and \ref{C1Lemma}, we deduce that for $j=1,\dots,r-1$  
 
\[   |\nu_j(S_b)| \le    |C^L_j(S_b)| + |C^R_j(S_b)|  + 2C_0 + 2C_1 + 2B +\mathcal E_j,  \] 
where $\mathcal E_j$ is the number of edges from the excluded block coloured $\nu_j$.

Similarly, 
$$ 
|\L(\nu_0,S_b)|\le   
|C^R_0(S_b)| + C_0 + C_1 + B +\mathcal E_0  
$$ 
and  
$$ 
|\R(\nu_{r},S_b)|\le  
|C^L_{r}(S_b)| +C_0 + C_1 + B + \mathcal E_r. 
$$ 
This motivates us to define an array of colours $[\nu_0(S),\dots,\nu_r(S)]$ 
to be {\em very short} if for $j=1,\dots,r-1$ we have 
\[      |\nu_j(S)|  \leq 2C_0 + 2C_1 + 5B  + 1,    \]  
and 
\[      |\L(\nu_0,S)|  \leq C_0 + C_1 + 5B  + 1,    \] 
 and  \[      |\R(\nu_r,S)|  \leq C_0 + C_1 + 5B + 1,        \] 
and for $j=0,\dots,r-1$ the interval formed by the unnamed colours between $\nu_j(S)$
and $\nu_{j+1}(S)$ has total length at most $B$.

An implosive array is said to be {\em short} if it satisfies the weaker 
inequalities obtained by increasing  each of these bounds by $2B{\tz}$.   
 
\begin{lemma} \label{vshort} Let  $A=[\nu_0(S^0),\dots,\nu_r(S^0)]$ be  a 
$p$-implosive array 
with $p\ge \tz$. 
\begin{enumerate}
\item If $[\nu_0(S^\tz),\dots,\nu_r(S^\tz)]$ is very short, 
then $A$ is short. 
\item If $A$ is short, then $\|A\|\le 2B(2C_0+2C_1 + 5B +1+2BT_0) +2B^2(1+2T_0).$
\end{enumerate}
\end{lemma} 
 
\begin{proof} Item (1) is an immediate consequence of the Bounded Cancellation 
Lemma \ref{BCL}. The (crude) bound in (2) is an immediate consequence of Lemma \ref{2Bcols}
and the inequalities in the definition of {\em short}; the first summand is an estimate
on the sum of the lengths of the named colours, and the second summand accounts for
the unnamed colours.
\end{proof} 
   
The following 
lemma is the key step in  the proof of Proposition \ref{prePincerLemma}. 
 
\begin{lemma} \label{shorty} 
If $A(S^0)=[\nu_0(S^0),\dots,\nu_r(S^0)]$
 is a $p$-implosive array,
then at least one of the following statements is true:
\begin{enumerate}
\item $p\le 2T_0$; 
\item $A(S^0)$ is short;
\item  $p > 2T_0$ and the array $A(S^{T_0})$ contains a very short
implosive sub-array $[\nu_k(S^\tz),\dots,\nu_l(S^\tz)]$.
\end{enumerate}
\end{lemma} 
 
\begin{proof} Assume $p > 2T_0$ and that
 $[\nu_0(S^0),\dots,\nu_r(S^0)]$ is not short. We claim that there is a 
 block of at least $B+1$ constant letters in the interval determined by the
array
 $\mathcal L(\nu_0,S^{T_0})\dots\mathcal L(\nu_r,S^{T_0})$.
Indeed, by definition, if an array
 is not short then either one of the $\mathcal E_j$ has
length at least $B+1$, or one of the
 blocks of unnamed colours has length at
least $B(2T_0+1)+1$, or
else at least one of the intervals of 
constant letters $C^L_j(S^{0})$ 
or $C^R_j(S^{0})$ has length at least $B(T_0+1)+1$. 
In the first case, since  $\mathcal E_j$  is in
the excluded block, none of its edges are cancelled
 before the moment of implosion, and
hence it contributes a block of at least $ B+1$ constant
 letters to $A(S^{T_0})$; in the second
case, the Bounded Cancellation 
Lemma assures us that the length of the appropriate block of unnamed colours  can decrease by 
at most $2B$ at each step before the implosion of the array, 
and hence it still contributes
a block of  at least $ B+1$ constant edges to $A(S^{T_0})$;
and similarly, in the third case,
 $C^{\ast}_j(S^{0})$  can decrease by 
at most $B$ at each step before the implosion of the array.

Let $\beta$ be a block of at least $ B+1$ constant edges in $A(S^{T_0})$ with non-constant
edges $e_l$ and $e_{\rho}$ immediately to its left and right, respectively.

The Buffer Lemma \ref{BufferLemma} assures us that the
non-constant edges in the future of $e_l$ will never interact with the non-constant edges in the
future of $e_{\rho}$. Thus at least one of $e_l$ or $e_{\rho}$ must be {\em stabbed in the back}, i.e. 
its entire non-constant future must be consumed by edges on its own side of $\beta$. Suppose,
for ease of notation, that it is $e_l$ and let $\nu_i$ be the colour of $e_l$. We claim that if $\nu_k$
is the colour of the letter that ultimately consumes $e_l$, then $k\le i-2$.

We shall derive a contradiction
from the assumption that the edge which ultimately
 consumes  $e_l$ is  coloured $\nu_{i-1}$.
There are two cases to
consider according to whether $e_{\rho}$ is also coloured $\nu_i$. If it is, then we consider the
word $V$ labelling the arc of $\bot(S^0)$
 from the left end of $\nu_i(S^0)$ to the past of
$e_l$; the consumption
of the non-constant future of 
$e_l$ completes the $\phi$-neutering  of $V$
 by the word labelling $\nu_{i-1}(S^0)$, 
in particular this neutering will have taken more than $T_0$ steps in time, contradicting
the Two Colour Lemma \ref{TwoColourLemma}. If $e_{\rho}$ is not coloured $\nu_i$, then the
consumption of the non-constant future of 
$e_l$ results in a new essential adjacency of colours and hence can only be complete
 at the moment of
implosion, i.e. $\time(S^p)$. But this consumption constitutes the neutering of $\nu_i(S^{T_0})$
by $\nu_{i-1}(S^{T_0})$, and according to the Two Colour Lemma this neutering 
must be accomplished in at most $T_0$ units of time. Thus $p\le 2T_0$,
contrary to our hypothesis. 

\def\kill{\!\!\!\searrow\!}

Thus we have proved that the edge which ultimately consumes $e_l$ is coloured
$\nu_k$ where $k\le i-2$. Under these circumstances (or the 
symmetric situation with $e_{\rho}$ in place of $e_l$) we say that {\em $\nu_k$
\index{neutering} neuters $\nu_i$ from behind} and write  $\nu_k \kill \nu_i$.

\begin{figure}[htbp] 
\begin{center} 
  
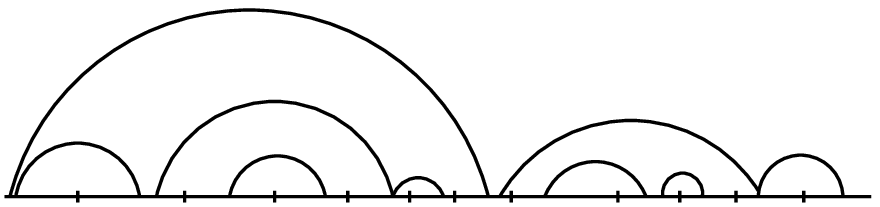 
  
\caption{The nesting associated to $\kill$} 
\label{figure:Arc} 
\end{center} 
\end{figure} 
 
\medskip 
 
There is a natural {\em nesting} among the $\kill$-related pairs of colours from the
array: 
 $(\nu_{k_1},\nu_{j_1}) < (\nu_{k_2},\nu_{j_2})$ if $\nu_{k_1}$ and 
$\nu_{j_1}$   
both lie between $\nu_{k_2}$ and $\nu_{j_2}$ in $S^0$. See Figure 
\ref{figure:Arc}. 
 
We focus our attention on an innermost (i.e. minimal) 
pair with  $\nu_k\kill\nu_i$.  By definition $|k-i|\ge 2$. If there were
a block of at least $ B+1$ constant letters between the closest non-constant
letters of $\nu_k(S^{T_0})$ and $\nu_i(S^{T_0})$, 
then the preceding argument 
would yield a neutering from behind that contradicted the innermost nature
of $\nu_k\kill\nu_i$. Thus  
$[\nu_j(S^\tz),\dots,\nu_k(S^\tz)]$ is a very short array, and we are done.  
\end{proof}

\begin{RIproof}\label{RIP}{\em Proof of Regular Implosions (Prop.\ref{implosion}):} Given the bound in Lemma \ref{vshort}(2), an obvious finiteness
argument provides a constant $\tau$ such that every short implosive array is $p$-implosive
with $p\le \tau$. And the same bound applies to implosive arrays that contain a short
sub-array (Remark \ref{subarray}(1)). So in the light of Lemmas  \ref{shorty} and   \ref{vshort}(1),
it suffices to let $\hat{T_1} = \max\{2T_0, \tau\}$.\hfill $\square$
\end{RIproof}
 
\subsection{Super-Buffers}
\index{super-buffers}

In this subsection we prove an important cancellation lemma based on
Proposition \ref{prePincerLemma}, this lemma involves the following
constant.

\begin{definition} \label{T1'Lemma}
We fix an integer $T_1'$ such  that one gets repetitions  in all $T_1'$-long subsequences of $5$-tuples of reduced words
\[	U_k:=\Big( u_{k,1}, u_{k,2}, u_{k,3}, u_{k,4}, u_{k,5}	 \Big) \ \ \ \ k=1,2,\dots \]
with $|u_{k,1}| $ and $|u_{k,1}| $ at most $ C_0+ C_1 + 2B +1$, while
 $|u_2^k|$ and $ |u_4^k|$ are at most $ C_0 + C_1$, and $|u_3^k| \leq 4B+1$. That is,
for some $t_1\le t_2\leq T_1'$ and 
\[	\Big( u_{t_1,1}, u_{t_1,2}, u_{t_1,3}, u_{t_1,4}, u_{t_1,5}	 \Big) = \Big( u_{t_2,1},  u_{t_2,2}, u_{t_2,3}, u_{t_2,4}, u_{t_2,5}	 \Big) .	\]
\end{definition}

\newtheorem{stipulation}[theorem]{Stipulation}

\begin{stipulation} Assume $T_1' \ge \hat{T_1}$.
\end{stipulation}
 
The cancellation lemma we need is most easily phrased in terms of 
colours of subwords, which we define as follows, keeping firmly in mind
the example of a stack of partial corridors excised from the interior of a van Kampen
diagram, retaining their memory of the colours to which the edges belong.

We have a word $W$ with a  decomposition into preferred subwords
$V = V_1 V_2
\cdots V_k$, where each $V_i$ is either positive or negative;
we think of these subwords as having colours $\mu_1, \ldots \mu_k$.
Take the freely reduced words $\phi(V_i)$, concatenate them, then
cancel to form a freely reduced word. There is some freedom in the
choice of cancellation scheme, as in the folding of corridors, but we fix
a choice, thus assigning to each letter of the freely reduced form of
$\phi(V)$ the colour  $\mu_i$ of its ancestor. We repeat this process,
thus assigning colours to the letters in the reduced form of $\phi^k(V)$ for
each integer $k>0$.

The process that we have just described is an algebraic description of
a choice of  minimal area van Kampen diagram for $t^{-k}Vt^k\phi^k(V)^{-1}$.
Thus the following lemma is a comment on the form of  such
diagrams.

\begin{proposition} \label{NoDoubleNeuter}  
Let $V=V_1V_2V_3$ be a concatenation of words (coloured $\nu_1, \nu_2, \nu_3$)
each of which is either positive or negative.
If $W$ is a subword of the reduced
form of  $\phi^{T_1'}(V)$ and   $W$ has a non-constant
letter coloured $\nu_i$ for each $i\in\{1,2,3\}$, 
then for all $k \geq 0$ there are
non-constant letters  in $\phi^k(W)$ coloured $\nu_2$.
\end{proposition}

\begin{proof} Let $\nu_i(W)$ denote the subword of $W$
coloured $\nu_i$, and let $\nu_i^j$ denote  the maximal subword coloured $\nu_i$ in
(the reduced word representing) $\phi^i(V_1V_2V_3)$ .
Note that $\nu_2(W)=\nu_2^{T_1'}$, and more generally $\nu_2^{T_1'+j}$
is the maximal word  in $\phi^j(W)$ coloured $\nu_2$.

Fix $k>T_1'$ and consider the diagram formed by the stack of corridors
described prior to the proposition. The bottom
of the first corridor is labelled $V$, and we regard it as being divided into
three coloured intervals according to the decomposition $V_1V_2V_3$.
Since $\nu_2(W)$ contains non-constant letters and $T_1'>\hat{T_1}$, 
the array formed by these colours is not implosive (Proposition \ref{prePincerLemma}),
and hence  
$\nu_1(W)$ and $\nu_3(W)$ will never essentially consume $\nu_2(W)$.
However, the proposition is not yet proved because there remains
the possibility
that  $\nu_2$ may essentially vanish because it  neuters $\nu_1(W)$, say, and is then
neutered by $\nu_3(W)$. We proceed under this assumption, seeking
a contradiction. (The case where the roles of 
$\nu_1$ and $\nu_3$ are reversed is entirely similar.)

For each $1 \leq i \leq T_1'$, we have
$\phi^i(V_1V_2V_3)=\nu_1^i, \nu_2^i$ and $\nu_3^i$.
Write $\nu_2^i\equiv V^i(1)  V^i(2)  
V^i(3)$, where $V^i(1)$ ends with last  non-constant letter in
$\nu_2^i$ whose entire non-constant future is eventually consumed by
letters coloured $\nu_1$, and $V^i(3)$ begins with the leftmost
non-constant letter whose entire non-constant future is  
eventually consumed by letters coloured $\nu_3$.
Lemmas \ref{C_0} and \ref{C1Lemma} tell us that $V^i(1)$ and $V^i(3)$ have
length at most
$C_0 + C_1$.  

\noindent{\em Claim:} $V^i(2)$  contains exactly one non-constant edge
and has length no more than $4B+1$.

We are assuming that $\nu_2(W)$ neuters $\nu_1(W)$. Consider the
(non-constant) edge $\e_i$ in
$\nu_2^i$ that will eventually consume the final non-constant edge in
$\nu_1(W)$. Note that $\e_i$ is  the
leftmost non-constant edge in $V^i(2)$. Moreover, we are assuming 
that  $\nu_3(W)$ ultimately neuters
$\nu_2(W)$, so in particular  it consumes the entire future of
any edge to the right of $\e_i$, which
forces $\e_i$ to be the rightmost non-constant edge in $V^i(2)$. The Buffer
Lemma tells us that $\e_i$ must lie within $2B$ of both ends of $V^i(2)$,
and hence the claim is proved.

Looking to the left
of $V^i(1)$, we now consider the  subword    $L^i$ of $\nu_1^i$  
that begins with the leftmost  non-constant edge in the future of
which there is a non-constant letter that cancels with a letter
coloured $\nu_2$. And looking to the right of  $V^i(3)$, we consider
the subword that ends with the rightmost non-constant letter  in the future of
which there is a non-constant letter that cancels with a letter
coloured $\nu_2$. 
any of whose non-constant future cancels
with an edge painted $\nu_2$.  As in previous arguments, The Buffer Lemma and
Lemmas \ref{C_0}, \ref{C1Lemma} tell is that  $|R^i|, |L^i| \le C_0 +
C_1 + 2B +1$, for all $i$.

We have already bounded the lengths of $V^i(1), V^i(2)$ and $V^i(3)$
by $C_0+C_1, 4B+1$ and $C_0+C_1$, respectively. Thus we are
now in a position to invoke the repetitive behaviour described in Definition
\ref{T1'Lemma}:
for some positive integers $i $ and $t$ with $i+t\le T_1'$, we get a repetition 
\[	\Big( R^i, V^i(1), V^i(2), V^i(3), L^i \Big) =  
\Big( R^{i+t}, V^{i+t}(1), V^{i+t}(2), V^{i+t}(3), L^{i+t} \Big).
\]
For as long as we are assured of the continuing presence
of $\nu_1^{i+s}$ and $\nu_3^{i+s}$,
 the fate of $\nu_2^i=V^i(1)V^i(2)V^i(3)$ under $s$ iterations of $\phi$ depends
only on $(R^i, V^i(1), V^i(2), V^i(3), L^i)$. Thus
$$
\Big( V^j(1), V^j(2), V^j(3) \Big) = \Big( V^{j+t}(1),
V^{j+t}(2), V^{j+t}(3) \Big)
$$
for all $j\ge i$ within the time scale of this assurance. However this leads us
to an absurd conclusion, because once $ \nu_1$ has become constant, 
at all subsequent time,
the surviving word coloured  $\nu_2$ contains as a proper subword, the
$\nu_2$ word that existed at the corresponding times in the  cycles (of
period $t$) before $T_1'$, and in particular they can never essentially
vanish, contrary to our assumption that $\nu_3$ eventually neuters
$\nu_2$. 
\end{proof}

\subsection{Nesting and the Pincer Lemma}

In subsequent sections we would like to
 bound the life of pincers by arguing
that during the lifetime of a pincer,
 colours must be consumed at a predictable rate (appealing
to Proposition \ref{prePincerLemma}),
noting that there  are only a limited number of colours. However, the bounds
we need will require us to ascribe each consumed colour to a {\em unique}
pincer. Thus we encounter problems whenever one pincer is contained in another.
For reasons that will become apparent in subsequent sections,
in situations where
we must confront this problem, the inner of the two pincers will have a long block
of constant edges along the corridor immediately above its peak. More precisely,
we will find ourselves in the situation described in the following definition. The
appearance of the constant $\ll := 2B(T_0+1)+1$ in the following definition is
explained by the role that this constant played in the course of Lemma \ref{shorty}.

\begin{figure}[htbp] 
\begin{center} 
  
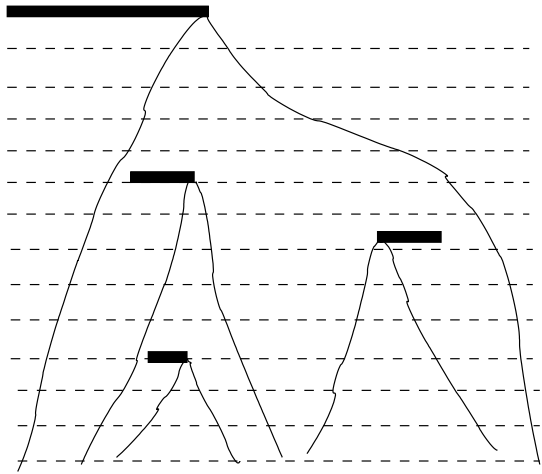 
  
\caption{A depiction of nesting} 
\label{figure:Nest} 
\end{center} 
\end{figure}

\begin{definition} \label{NestingDef}
Consider one pincer $\Pin_1$ contained in another $\Pin_0$.  Suppose
that in the corridor $S \subseteq \Pin_0$ at the top of $\Pin_1$
(where its boundary paths $p_1(\Pin_1)$ and $p_2(\Pin_1)$ come
together) the future in $\top(S)$ of at least one of the edges
containing $p_1(\Pin_1) \cap \bot(S)$ or $p_2(\Pin_1) \cap \bot(S)$
contains no non-constant edges, and this future\footnote{We allow this
future to be empty, in which case ``contained in" means that the immediate past
of the long block of constant edges is not separated from $\Pin_1$ by any
edge that has a future in $\top(S)$.}
lies in an interval of at least $\ll$ constant edges contained in
$\Pin_0$.  Then we say that $\Pin_1$ is \index{pincer!nested}{\em nested in} $\Pin_0$. (in
Figure \ref{figure:Nest}, the $\ll$-long block of constant edges are shown in
black.) We say that $\Pin_1$ is \index{pincer! left(right)-loaded}{\em left-loaded} or {\em right-loaded}
according to the direction in which  the $\ll$-long block of constant edges
extends from the peak of $\Pin_1$.
\end{definition} 

\begin{remark}  A nested pincer cannot be both left-loaded and right-loaded (cf.
Remark \ref{subarray}(3)).

If $\Pin_1$ is left-loaded, then  the future of  
$p_1(\Pin_1) \cap \bot(S)$ contains no non-constant edges.
It may happen that the future of $p_2(\Pin_1)$ also contains no 
non-constant edges; in this case the colour $\mu$ of $p_2(\Pin_1)$
essentially vanishes in $S$ due to   cancellation between non-constant edges 
of $\mu$ and some colour to its right. Symmetric considerations apply
to right-loaded pincers.
\end{remark}

\begin{definition} \label{chiP}\index{pincer!colours $\chi(P)$}
For a pincer $\Pin_0$, let $\{ \Pin_i \}_{i \in I}$ be the set of all
pincers nested in $\Pin_0$.  Then define
\[	\chi(\Pin_0) = \tilde{\chi}(\Pin_0) \smallsetminus \bigcup_{i \in
I} \tilde{\chi}(\Pin_i).	\]
\ed

\bl \label{nestLife}
If the pincer $\Pin_1$ is nested in $\Pin_0$ then 
$\time(S_t(\Pin_1)) < \time(S_{\Pin_0}).	$
\end{lemma}
\begin{proof}
The presence of the hypothesised block of constant letters in
$\top(S_t(\Pin_1))$ makes this an immediate consequence of the Buffer
Lemma \ref{BufferLemma}.
\end{proof}

Define $T_1 :=  T_1' + 2T_0$.
The following theorem is the main result of this section.

\index{Pincer Lemma}
\begin{theorem}[Pincer Lemma] \label{PincerLemma} 
For any pincer $\Pin$
\[	\life(\Pin) \leq T_1(1 + |\chi(\Pin)|).	\]
\end{theorem}

\begin{proof}
The heart of our proof of Proposition \ref{prePincerLemma} was that
in each block of $\hat T_1$  steps in time between  $\time(S_b)$ and
$\time(S_\Pin)$ at least one colour essentially disappears. Our proof
of the present theorem is an elaboration of that argument: we must
argue for the essential disappearance of a  colour that is not
contained in any of pincers nested in $\Pin$.  Thus we concentrate
on that region of the pincer $\Pin$ that is exterior to the set of 
{\em co-level\footnote{i.e. those that are maximal
with respect to inclusion among the pincers nested in $\Pin$} 1} pincers nested in it;
let $\{ \Pin_j\}, \ j=1,\dots,{J}$ be the set of such, indexed in order of appearance from
left to right.

For $j=1,\dots,J-1$, let $\Sigma_j$ denote the set of  colours  along the bottom of $\Pin$
that have a non-constant edge strictly
between $\Pin_j$ and $\Pin_{j+1}$; if $\Pin_j$ is left-loaded, then we include
the colour of $p_2(\Pin_j)$ in $\Sigma_j$, and if $\Pin_j$ is right-loaded, then we include
the colour of $p_1(\Pin_j)$ in $\Sigma_{j-1}$. Likewise, we define $\Sigma_0$ to be
the set of non-constant colours that lie to the left of $\Pin_1$ together with the
colour of $p_1(\Pin)$, and we define
$\Sigma_{J}$ to be
the set of non-constant colours that lie to the right of $\Pin_J$ together with the
colour of $p_2(\Pin)$. 

In order to prove the theorem, we derive a contradiction from the assumption that
in the first $T_1$ units of time in the life of $\Pin$ no colours in the union of the
$\Sigma_j$ essentially vanish. (There is no loss of generality in starting at the
bottom of the pincer, since given any other starting time, one can discard the
pincer below that level.) We label the corridors, beginning at the bottom of $\Pin$
and proceeding in time as $S^0,S^1,\dots$

We focus on a single $\Sigma_j$, and write its colours in order as    $\nu_1, \ldots , \nu_r$.
We analyse how the colours in $\Sigma_j$ come to vanish.  
The first important observation is that  $2 \le i \le r-1$,
 it is not possible for the colour $\nu_i$ to essentially vanish (at any time)
due to cancellation merely between the colours in $\Sigma_j$.  
For if this happened,  there would be an implosive array in $S^0$
containing $\nu_i(S^0)$  and so, by Proposition \ref{prePincerLemma}, $\nu_i$ would vanish before $S^{T_1}$, contrary to our assumption.

There remains the possibility that $\nu_2$ may neuter  $\nu_1$ (after
$S^{T_1}$).  This can happen in two ways.  The first is that $\Pin_{j-1}$ is left-loaded: in
this case
the neutering happens within time $T_0$ of the top of $\Pin_{j-1}$ (by Two Colour Lemma),
and we are then in a stable situation in the sense that $\nu_3$ cannot subsequently neuter $\nu_2$,
by Proposition \ref{NoDoubleNeuter}.  Now suppose that $\Pin_{j-1}$ is right-loaded.
Consider the earliest time $t_0$ at which there is a block of  at least $B+1$ constant edges in the
past of the $\lambda_0$-long block associated to $\Pin_{j-1}$. If $\nu_2$ is to neuter
$\nu_1$, then it must do so within   $T_0$ steps of this time. Indeed, within $T_0$ steps,
if the non-constant edges of $\nu_1$ to the right of the block have not been consumed
by $\nu_2$, 
then they will never be consumed by a colour from $\Sigma_j$. 

There is a further event that we must account for, which is closely related to
neutering: it may
happen that $\nu_1$ is the colour of $p_2(\Pin_{j-1})$ and that $\nu_2$ consumes
all of the non-constant edges to the right of the block of constant edges discussed above;
this is not a neutering but nevertheless the Two Colour Lemma applies. We would like
to apply Proposition \ref{NoDoubleNeuter} in this situation to conclude that
$\nu_3$ cannot subsequently neuter $\nu_2$.
This is legitimate provided
$t_0\ge\time (S^{T_1'})$.  If $t_0< \time (S^{T_1'})$, then we still know that $\nu_3$
cannot  neuter $\nu_2$ before $S^{T_1}$, because by hypothesis no colour from
$\Sigma_j$ essentially vanishes before this time. On the other hand, the Two Colour Lemma
tells us that if $\nu_3$ is to neuter $\nu_2$, then it must do so within $T_0$ steps
from $t_0$, and $t_0+T_0\le  \time (S^{T_1})$. Thus, once again, we conclude that
$\nu_3$ can never neuter $\nu_2$.

Entirely similar arguments show that it cannot happen that $\nu_r$ is neutered
by $\nu_{r-1}$ and that subsequently $\nu_{r-2}$ neuters $\nu_{r-1}$.

We have established the existence of a stable situation: proceeding past the point where
the restricted amount of possible neutering within $\Sigma_j$
has occurred, we may assume that the next
essential disappearance of a colour from $\Sigma_j$ can only occur as a result of
cancellation with a colour from some $\Sigma_i$ with $i\neq j$. Such further cancellation
must occur, of course, because all but two\footnote{Degenerate cases with few
colours are covered by the Two Colour Lemma and the Buffer Lemma.} of the
 colours in $\bigcup_j\Sigma_j$ must be consumed within $\Pin$.

Passing to innermost pair of interacting $\Sigma_k$
we may assume $i=j-1$ (cf. proof of Lemma \ref{shorty}). Thus our proof will be
complete if we can argue that cancellation between non-constant edges
from $\Sigma_{j-1}$ and $\Sigma_j$ is impossible. We have
argued that the colours which are to cancel will be essentially adjacent within
time $T_0$ of 
the top of $\Pin_{j-1}$. On the other hand, there is a block of $\ll$ constant
edges separating $\Sigma_{j-1}$-nonconstant edges and $\Sigma_{j}$-nonconstant edges
at the top of $\Pin_{j-1}$. Since $\ll > 2B(T_0+1)$ at least $B+1$ of these constant edges
remain $T_0$ steps later. The Buffer Lemma now obstructs the supposed
cancellation between non-constant edges in $\Sigma_{j-1}$ and $\Sigma_j$. 
\end{proof}

\section{Teams and their Associates}\label{teamSec} 
\index{teams}
We begin the process of grouping pairs of colours $(\mu,\mu')$ into 
teams.  
 
\subsection{Pre-teams} \label{t1} 

The whole of $C_{(\mu,\mu')}(2)$ will ultimately be  consumed by  
a single edge $\e_0\in\mu'(S_0)$. 
We consider the time $t_0$ at which the future of $\e_0$ 
starts consuming the future of  $C_{(\mu,\mu')}(2)$.  
If $|C_{(\mu,\mu')}(2)|> 2B$, then this consumption will not be completed in 
three steps of time  (Lemma \ref{BCL}). We claim that in  this circumstance, the 
leftmost $\mu'$-coloured edge after the first two steps of the
cancellation must be left para-linear. Indeed it is not left-constant since it must consume edges in 
the future of $C_{(\mu,\mu')}(2)$, and since no non-constant $\mu'$-edges 
are cancelled by $\mu$ in passing from the first to the second stage of  
cancellation, the leftmost non-constant $\mu'$-label must remain the same (Proposition 
\ref{power}). We denote this \lpl edge at time $t_0+2$ by $\e^\mu$. 
 
Let $\e_\mu$ be the rightmost edge in the future of $C_{(\mu,\mu')}(2)$  at time $t_0$. 
We trace the ancestry of $\e_\mu$ and $\e^\mu$ in the trees of \index{family forest}
$\F\subset\Delta$ corresponding 
to the colours $\mu$ and $\mu'$ (as defined in \ref{tree}). 
 We go back to the last point in time \index{teams!times $t_i(\T)$}
$\ptmm$ at which  
both ancestors  lay in a common corridor 
 {\em and} the interval on the bottom of this corridor between the pasts of 
  $\e_{\mu}$ and $\e^\mu$  
is comprised entirely of constant edges whose future is eventually  
 consumed by the ancestor of 
$\e^\mu$ at this time. We denote this corridor  $S_{\uparrow}$. 
 
\begin{definition}\label{preteam}  The ancestor of $\e^\mu$ at time $\ptmm$ is called the 
\index{reaper}{\em reaper} and is denoted $\prmm$.  The set of edges in $\bot(S_{\uparrow})$ 
 which are eventually consumed by $\prmm$ is denoted $\pEmm$.  
  This is a contiguous set of edges. 
The \index{pre-team}{\em pre-team} $\pTmm$ is defined to 
be the set of pairs $(\mu_1,\mu')$ such that $\pEmm$ contains 
edges coloured $\mu_1$. The number of edges in $\pEmm$ is denoted 
$\|\hat \T\|$. 
\end{definition} 
 
In a little while we shall define {\em teams} to be   
pre-teams satisfying a certain maximality condition (see Definition \ref{newTeams}).

\begin{remark} \label{rem:virt}
 If $\ptmm<\time(S_0)$ then 
near the right-hand end of $\pEmm$ one may have  an interval of 
colours  $\nu$ such that $\nu(S_0)$ is empty.  
\end{remark} 
 
\smallskip 
 
In the proof of Proposition \ref{SummaryLemma} we saw that it would be
desirable if (whatever our final 
definition of {\em team} and $\bonus$ may be) the following inequality
(\ref{goodEq}) should hold for all teams: 
\begin{equation} \label{preTeamInequality} 
\sum_{(\mu,\mu') \in \T \mbox{ \tiny or } (\mu,\mu') \vin \T} |C_{(\mu,\mu')}(2)| \le 
\|\T\|  + |\bonus(\T)| + B. 
\end{equation}

The following lemma shows that, even without introducing a  bonus scheme 
or virtual members, the 
desired inequality is straightforward for pre-teams with $\ptmm \geq \time(S_0)$. 
 
\begin{lemma} \label{t1high} 
 If $\ptmm \geq \height(S_0)$ then $\pTmm$ satisfies 
\[ \sum_{(\mu,\mu') \in \pTmm}|\cmm| \leq \|\pTmm\| + B.      \] 
\end{lemma} 
 
\begin{proof} By definition   
$\mu'(S_0)$ does not start consuming any  of the
 $C_{(\mu_1,\mu')}(2)$ with $(\mu_1,\mu')\in\pT$ before $\ptmm$  
(apart from a possible nibbling of  length $< B$ from the rightmost team 
member  at time $\ptmm -1$).  
Since each $C_{(\mu_1,\mu')}(2)$ consists only of edges consumed 
by $\mu'(S_0)$, the future  of each $C_{(\mu_1,\mu')}(2)$ at time $\ptmm$ 
will have the same length as $C_{(\mu,\mu')}(2)$ 
(except that the rightmost may have lost these $< B$ edges).  
And these futures are contained in $\pEmm$. 
\end{proof}

The case where $\ptmm <\height(S_0)$ is more troublesome. As $\pEmm$ flows forwards in 
time, the number of constant letters in the future of $\pEmm$ that are  
consumed by $\prmm$ between $\ptmm$ and $\time (S_0)$ may be outweighed by the number of 
constant letters generated to the left of the future of $\pEmm$ that will 
ultimately be consumed by $\prmm$. 
 
It is to circumvent the failure of inequality 
(\ref{preTeamInequality}) in this setting that we are 
obliged to instigate the bonus scheme described in Section \ref{BonusScheme}.

\subsection{The Genesis of pre-teams} \label{genesis} 
 
\index{teams!genesis of}
We fix $\pTmm$ with $\ptmm < \time(S_0)$ and consider 
the various events that occur  at  $\ptmm$ to prevent 
us pushing the pre-team back one step in time.  We write $S_\omega$ 
to denote the corridor at time $\ptmm$ containing $\pEmm$. 
 
\begin{figure}[htbp] 
\begin{center} 
  
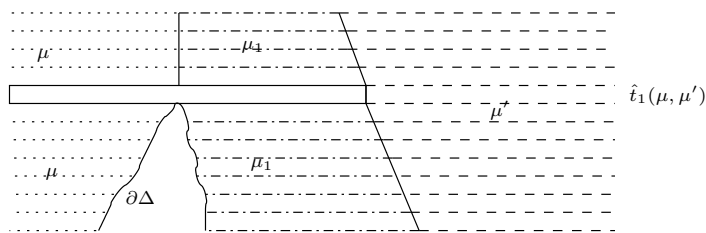 
  
\caption{A team of genesis (G1)}
\label{G1Pic} 
\end{center} 
\end{figure}

\begin{figure}[htbp] 
\begin{center} 
  
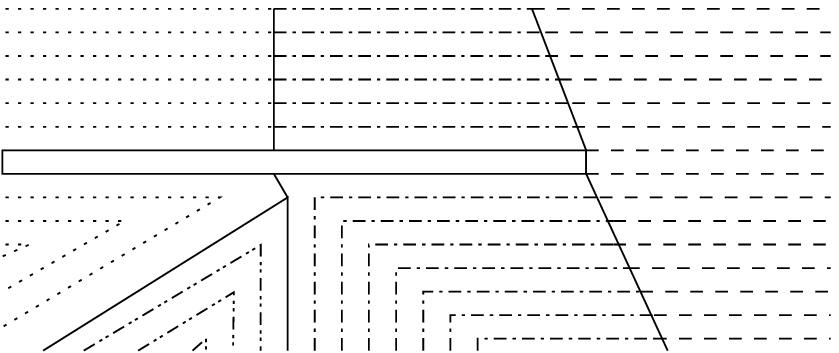 
  
\caption{A team of genesis (G2)}
\label{G2Pic} 
\end{center} 
\end{figure} 

\begin{figure}[htbp] 
\begin{center} 
  
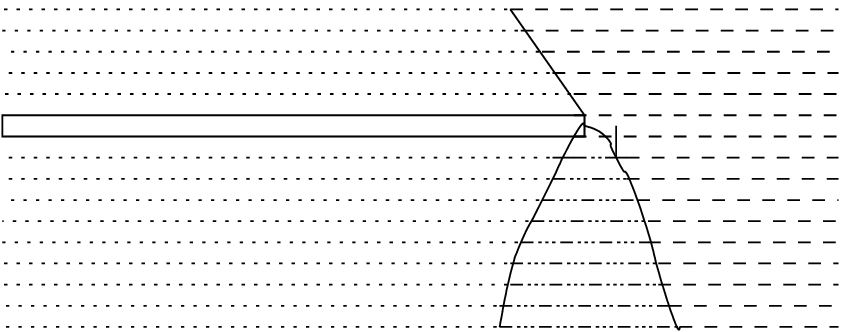 
  
\caption{A team of genesis (G3)}
\label{G3Pic} 
\end{center} 
\end{figure} 

\begin{figure}[htbp] 
\begin{center} 
  
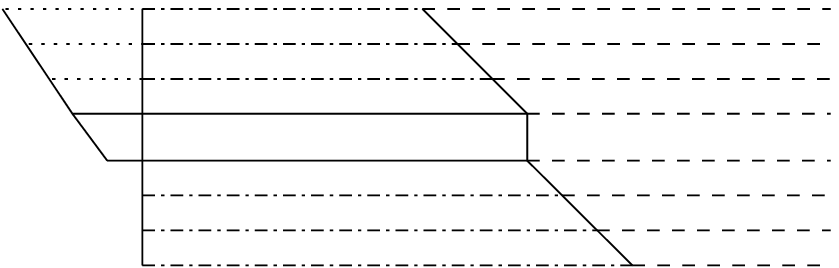 
  
\caption{A team of genesis (G4)}
\label{G4Pic} 
\end{center} 
\end{figure}

\def\CC{C}
 
There are four types of events: 
 
\begin{enumerate} 
\item[(G1)] The immediate past  of $\CC_{(\mu,\mu')}(S_\omega) 
$ is separated from the 
past of $\prmm$ by an intrusion of $\partial\Delta$ (Figure \ref{G1Pic}). 
\item[(G2)] We are not in case (G1), but the immediate past  of
$\CC_{(\mu,\mu')}(S_\omega)$ is separated from the 
past of $\prmm$ because of a singularity (Figure \ref{G2Pic}). 
\item[(G3)] The immediate past of $C_{(\mu,\mu')}(S_\omega)$ is  still
in the same corridor as the past of $\prmm$, but it is  separated from
it by a non-constant letter (Figure \ref{G3Pic}). 
\item[(G4)] We are not in any of the above cases,  
but the immediate past of the rightmost letter in  
 $C_{(\mu,\mu')}(S_\omega)$ is not constant (Figure \ref{G4Pic}). 
\end{enumerate} 
 
\smallskip  
 
\def\ST{S_\T} 
\def\STmm{\S_{\T(\mu,\mu')}}

The following lemma explains why Figures \ref{G3Pic} and \ref{G4Pic} are an
accurate portrayal  of cases (G3) and (G4). 
 
Let $L_{inv}$ be the maximum length of $\phi^{-1}(x)$ over generators
$x$ of $F$, and $C_4 = L_{inv}.L$.  
 
\begin{lemma}\label{G34pics} 
If $I$ is an interval on $\top (S)$ labelled by a word $w$ in constant
letters
then the reduced word labelling the past of $I$ in $\bot(S)$ is of the 
form $u\alpha v$, where $\alpha$ is a word in constant letters and 
$|u|$ and $|v|$ are less than $C_4$. Moreover, if the past of the leftmost  
(resp. rightmost) letter 
in $w$ is constant, then $u$ (resp. $v$)  is empty. 

In particular, $|I| \leq |\alpha| + 2LC_4$.
\end{lemma}

\begin{figure}[htbp] 
\begin{center} 
  
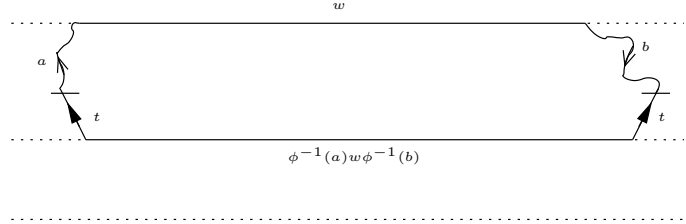 
  
\caption{The proof of Lemma 1.9.4}
\label{bottom-constant} 
\end{center} 
\end{figure} 

\begin{proof} See Figure \ref{bottom-constant}.  Follow the path from
the left end of $I$ to $\bot(S)$.  This passes through a (possibly
empty) path $a^{-1}$, followed by an edge labelled $t^{-1}$, where the
length of $a$ is less than $L$ (since it can be chosen to be on the
top of a $2$-cell which has an edge in $I$).
Similarly, at the right end of $I$ we have a path labelled
$bt^{-1}$, where the length of $b$ is less than $L$. The path along
$\bot(S)$ joining the two endpoints of these paths is labelled by the
reduced word freely equal in $F$ to $\phi^{-1}(awb) =
\phi^{-1}(a)w\phi^{-1}(b)$.  The only non-constant edges in this word
come from $\phi^{-1}(a)$ and $\phi^{-1}(b)$, which have lengths at
most $L.L_{inv}$.  This proves the assertion in the first sentence.
 
The assertion in the second sentence follows from the observation that
if $x,\, y$ and $\phi(x\beta y)$ consist only of constant letters,
then so does the reduced form of $\beta$, and the assertion in the
final sentence follows immediately from the first.
\end{proof} 
 
\begin{remark}\label{decreell} 
 It is convenient to assume that $LC_4 < \ll$. (In the unlikely 
event that this is not the case, we simply increase $\ll$.) 
\end{remark}

We are finally in a position to make an appropriate definition of a team.

\begin{definition}\label{newTeams} \label{shortDef}
All pre-teams $\pTmm$ with $\hat 
t_1(\mu_1,\mu')\ge\time(S_0)$ 
are defined to be teams, but the qualification criteria for pre-teams with 
 $\hat t_1(\mu_1,\mu')<\time(S_0)$ are 
more selective. 
 
If the genesis of  $\pTmm$ is of type  (G1) or (G2), then 
 the rightmost component of the pre-team may  form a pre-team 
at times before $\ptmm$.  In particular,  it may happen  
that $(\mu_1,\mu')\in\pTmm$ but $\ptmm > \hat t_1(\mu_1,\mu')$ and hence 
$(\mu,\mu')\not\in\pT(\mu_1,\mu')$. To avoid double 
counting in our estimates on $\|\T\|$ we disqualify the  
(intuitively smaller) pre-team $\pT(\mu_1,\mu')$ in these settings.  
 
If the genesis of $\pTmm$ is of type (G4), then again it may happen 
that what remains to the right of $\pTmm$ at some time  before $\ptmm$ is a pre-team. 
In this case, we disqualify the (intuitively larger) pre-team  $\pTmm$.   
  
The pre-teams that remain after these disqualifications 
are now defined to be {\em teams}.  
 
A typical team will be denoted $\T$ 
and all hats will be dropped from the notation for their associated objects 
(e.g. we write $\Emm$ instead of $\pEmm$).  
 
A team is said to be {\em short} if $\|\T\|\le \ll$
or $\sum\limits_{(\mu,\mu')\in\T} |C_{(\mu,\mu')}(2)| \le \ll$. Let
$\S$ denote the set of short teams.
\end{definition} 
 
\begin{lemma} \label{G4lemma} Teams of genesis (G4) are short. 
\end{lemma} 
 
\begin{proof} Lemma \ref{G34pics} implies that $\ET$ is in the immediate 
future of an interval of length at most $C_4$. And we have decreed (Remark 
\ref{decreell}) that $LC_4< \ll$. 
\end{proof} 
 
We wish our ultimate definition of a team to be such that every  pair  $(\mu,\mu')$ 
with $C_{(\mu,\mu')}(2)$ non-empty is assigned to a team. The above definition 
fails to achieve this because of two phenomena: first, a pre-team 
 $\pTmm$ with genesis of type (G4) may
have been disqualified, leaving $(\mu,\mu')$ teamless; second, in our initial discussion of  
pre-teams (the first paragraph of Section \ref{t1}) we excluded pairs $(\mu,\mu')$ 
with $|C_{(\mu,\mu')}(2)|\le 2B$. The following definitions remove these difficulties. 
 
\bd[Virtual team members] \label{Virtual} 
If a pre-team $\pTmm$ of type (G4) is disqualified under the terms of Definition \ref{newTeams} 
and the smaller team necessitating disqualification is $\pT(\mu_1,\mu')$,  
then we define $(\mu,\mu')\vin\pT(\mu_1,\mu')$ and $\pTmm\subset_v\pT(\mu_1,\mu')$. 
We extend the relation $\subset_v$ to be transitive and extend $\vin$ correspondingly. 
If $(\mu,\mu')\vin\T$ then $(\mu_2,\mu')$ is said to be a 
\index{teams!virtual members}{\em virtual member} of  
the team $\T$. 
\ed 
 
\bd If  $(\mu,\mu')$ is such that  $1\le |C_{(\mu,\mu')}(2)|\leq 2B$ and 
$(\mu,\mu')$ is neither a member nor a virtual member of any previously 
defined team, then we define $\T_{(\mu,\mu')}:=\{(\mu,\mu')\}$ to be a
(short) team with $\|\T_{(\mu,\mu')}\|=|C_{(\mu,\mu')}(2)|$.  
\ed 
 
\begin{lemma}\label{allIn}  
Every   $(\mu,\mu')\in\vecZ$ with $C_{(\mu,\mu')}(2)$ non-empty is a member 
or a virtual member of exactly one team, and there are less than $2\n$ teams. 
\end{lemma} 
 
\begin{proof} The first assertion is an immediate consequence of the preceding 
three definitions, and the second  follows 
 from the fact that $|\vecZ| < 2\n$. 
\end{proof}

\subsection{Pincers associated to teams of Genesis (G3)} \label{ss:PinG3}
 
In this subsection we describe the pincer $\Pin_\T$ canonically
associated to each team of genesis $(G3)$. 
The definition of $\Pin_\T$
involves the following concept which will prove important also for teams
of other genesis.

\begin{definition}\label{narrowPast}
We define the \index{teams!narrow past of}{\em narrow past} of a team 
$\T$ to be the set of constant edges that have a future in $\ET$. The narrow 
past may have several components at each time, the set of which  
are ordered left to right according to the ordering in $\ET$ of their futures.  We call these components {\em sections}. 
\end{definition}

{\center{\em{For the remainder of this subsection we 
consider only long teams of genesis (G3). }}}

\index{pincer!of a team|(}
\begin{definition}[The Pincer $\tilde\Pin_\T$] \label{pl}\label{t2}
The paths labelled $\hat p_l$ and $\hat p_r$ in Figure \ref{G3Pic}
 determine a pincer and are defined as follows. Let $\xT$ be the
 leftmost non-constant edge to the right of $\mu$ in the immediate
 past  of  $\T$, and let  $x_1(\T)$  be the edge that consumes it.   
Define $\tplT$ to be the path in $\F$ that traces the history of $\xT$
to the boundary,  and let $\tprT$ be the path that traces the history
of  $x_1(\T)$.
(Note that  $x_1(\T)$ is left-fast.)

Define $\tilde t_2(\T)$ to be the earliest time at which the  
paths $\tplT$ and $\tprT$ lie in the same corridor.
The segments of the paths $\tplT$ and $\tprT$ after this time, together
with the path joining them along the bottom of the
corridor at  time $\tilde t_2(\T)$ form a pincer. We denote this pincer
$\tilde\Pin_{\T}$. 
\end{definition}

The Pincer Lemma argues for the regular disappearance of colours
within a pincer during those times when more than two colours continue
to survive along the corridors of $\tilde\Pin_\T$. 
However, when there are only two colours the situation  is
more complicated.  

We claim that the following situation cannot arise:
$\time(S_{\hat\Pin_\T}) \leq \tone - T_0$, the path
$\tplT$ and the entire narrow past of $\T$ are in the same corridor at
time $\tone - T_0$, and at this time they 
are separated only by constant edges. For if this were the case,
then the colour of $\tprT$ would $\phi$-neuter the colour of $\tplT$
but would take more than $T_0$ steps to do so, contradicting
the Two Colour Lemma.  Thus at least one of the three hypotheses in
the first sentence
of this paragraph is false; we consider the three possibilities. The
troublesome case (3) leads to a cascade of pincers as depicted in
Figure \ref{cascade}.

\begin{definition}[The Pincer $\Pin_{\T}$ and times $t_2(\T)$ and
$t_3(\T)$] \label{PincerDef} 
\index{teams!times $t_i(\T)$}
\ 

\begin{enumerate}
\item {\em Some section of the narrow past of $\T$ is not in the same corridor as
 $\tplT$ at time $\tone - T_0$:} In this case\footnote{this includes the
possibility that $\tplT$ does not exist at time  $\tone - T_0$}
we define $t_2(\T)=t_3(\T)$ to be the earliest time at which the entire
narrow past of $\T$ lies in the same corridor as $\tplT$ and has length at least
$\ll$.
\item {\em Not case (1), there are no non-constant edges between $\tplT$
and the narrow past of $\T$ at time $\tone - T_0$:} In this case
$\time(S_{\tilde\Pin_{\T}}) > \tone - T_0$.  We define
$\Pin_{\T} = \tilde\Pin_{\T}$ and  $t_3(\T) =
\time(S_{\Pin_{\T}})$. If the narrow past of $\T$ at time $\tone - T_0$ 
has length less than $\ll$, we define
$t_2(\T) =  t_3(\T)$, and otherwise $t_2(\T) =\tilde t_2(\T)$. 
\item 
{\em Not in case (1) or case (2):}
In this case there is at least one non-constant edge between the
narrow past of $\T$ and $\tplT$ at 
$\tone - T_0$. We  pass to the latest time at which there is such an
intervening
non-constant edge and consider  the path $\tilde p_l^\prime(\T)$
that traces the history of the
leftmost intervening non-constant edge $x'(\T)$ and the path $\tilde
p_r^\prime(\T)$
that traces the history of the edge $x_1^\prime(\T)$ that cancels
with $x'(\T)$.
We define  $\tilde t_2'(\T)$ to be the earliest time at which the  
paths $\tilde p_l^\prime(\T)$ and $\tilde p_r^\prime(\T)$  lie in the
same corridor
and consider  the pincer formed by the  segments of the paths $\tilde
p_l^\prime(\T)$ and $\tilde p_r^\prime(\T)$
after  time $\tilde t_2'(\T)$ together
with the path joining them along the bottom of the
corridor at  time $\tilde t_2^\prime(\T)$.

We now repeat our previous analysis with the primed objects $\tilde p_l^\prime(\T), \tilde t_2^\prime(\T)$ {\em etc.} in
place of $ \tilde p_l(\T), \tilde t_2(\T)$ {\em etc.}, checking whether we now fall into case (1) or (2);
if we do not then we pass to $ \tilde p_l''(\T), \tilde t_2''(\T)$ {\em etc.}, and iterate the analysis until
we do indeed fall into case (1) or (2), at which point we acquire the desired definitions of 
$\Pin_\T,\, t_2(\T),\,  t_3(\T)$. \index{pincer!of a team|)}
\end{enumerate}

Define $p_l(\T)$ (resp. $\prT$) to be the left (resp. right)
boundary path of the pincer $\Pin_\T$ extended backwards
in time through $\F$ to $\partial\Delta$. Define $p_l^+(\T)$ to be the
sequence of non-constant edges (one at each time) lying immediately to
the right of the narrow past of $\T$ from the top of $\Pin_{\T}$ to
time $\tone$. (These are edges of the leftmost of the primed $\tplT$
considered in case (3).)
\ed

\begin{definition} \label{def:chiP}
Let $\T$ be a long team of genesis (G3).  Let \index{teams!colours $\chi_*(\T)$} $\chi_P(\T)$ be the set
of colours containing the paths $\tilde p_l(\T), \tilde p_l'(\T),\tilde p_l''(\T),\dots$ that
arise in (iterated applications of) case (3) of Definition \ref{PincerDef} but
do not become $p_l(\T)$.
\ed 

\begin{figure}[htbp] 
\begin{center} 
  
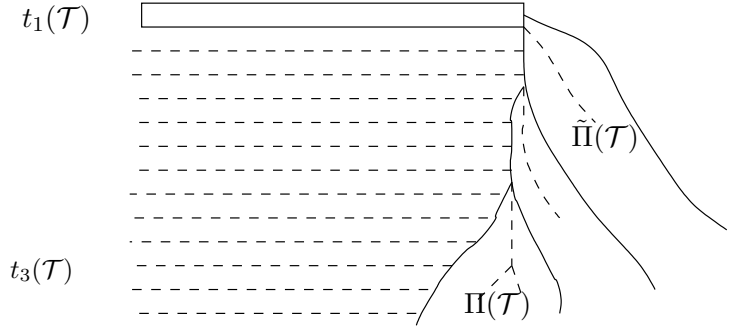 
  
\caption{The cascade of pincers.}
\label{cascade} 
\end{center} 
\end{figure}

The preceding definitions are framed so as to make the following important
facts self-evident.

\begin{lemma} 
 \label{t1t3forTeam}\label{disj1}

\ 

\begin{enumerate}
\item
If $\T$ is a long team of genesis (G3),  
\[      t_1(\T) - t_3(\T) \leq T_0(|\chi_P(\T)| + 1).   \] 
\item
If $\T_1$ and $\T_2$ are disjoint then  $\chi_P(\T_1)\cap \chi_P(\T_2)=\emptyset$.
\end{enumerate}
 \end{lemma}

\subsection{The length of teams} \label{TeamLemmas}

\begin{definition}\label{down1}  \index{teams!$\down_i(\T)$}
Define $\down_1({\mathcal T})\subset\partial\Delta$ to consist of those 
edges $e$ that are labelled $t$ and satisfy one of the following conditions: 
\begin{enumerate} 
\item[1.] $e$ is at the left end of a corridor containing a section of the narrow 
past of $\T$ that is not leftmost at that time; 
\item[2.]  $e$ is at the right end of a corridor containing a section of the narrow 
past of $\T$ that is not rightmost at that time;  
\item[3.]  $e$ is at the right end of a corridor which 
contains the rightmost section of the narrow past of $\T$ at that time but which does 
not intersect $\plT$.\\ 
\end{enumerate} 
\end{definition} 
 
All of the edges shown on the boundary in
Figure \ref{figure:TeamAge} are contained in $\down_1(\T)$. 

\begin{definition} \label{def:partialT} \index{teams!$\partial^\T$}
 Define $\partial^\T\subset\partial\Delta$ to be the
set of (necessarily constant) edges that have a preferred future in
$\ET$.
\end{definition}

We record an obvious disjointness property of the sets defined above.

\begin{lemma}\label{disj2}

\ 

\begin{enumerate}
\item For distinct teams $\T_1$ and $\T_2$, $\partial^{\T_1}$ and
$\partial_{\T_2}$ are disjoint.
\item For distinct teams $\T_1$ and $\T_2$, $\down_1(\T_1)$ and
$\down_1(\T_2)$ are disjoint. 
\end{enumerate}
\end{lemma}

\bd\label{QT}\index{teams!$\QT$}
Suppose that $\T$ is a team of genesis (G3).  We define  $\QT$ be the
set of edges $\e$ with the following properties:
$\plT$ passes through $\e$  before time $t_3(\T)$, and the corridor $S$
with $\e\in\bot(S)$ contains the entire narrow past of $\T$ and
this narrow past has length at least $\ll$.
\ed

The following lemma gives us a bound on $|\ET|$, which will reduce 
our task to that of bounding $|\QT|$ for teams of genesis (G3). 
 
\begin{lemma}  \label{TeamAgeLemma}  

\ 
 
\begin{enumerate} 
\item[1.] If the genesis 
of $\T$ is of type (G1) or (G2), then 
$$
\|\T\| \leq 2LC_4\,|\down_1({\mathcal T})| +  |\partial^\T|  .
$$ 
\item[2.] If the genesis of $\T$ is of type (G3), then 
$$
\|\T\| \leq 2LC_4\,|\down_1({\mathcal T})| +|\partial^\T| + 2LC_4\,|\QT| +
2LC_4T_0\big(|\chi_P(\T)| +1\big) + \ll .
$$  
\end{enumerate} 
\end{lemma} 
 
\medskip 
 
\begin{figure}[htbp] 
\begin{center} 
  
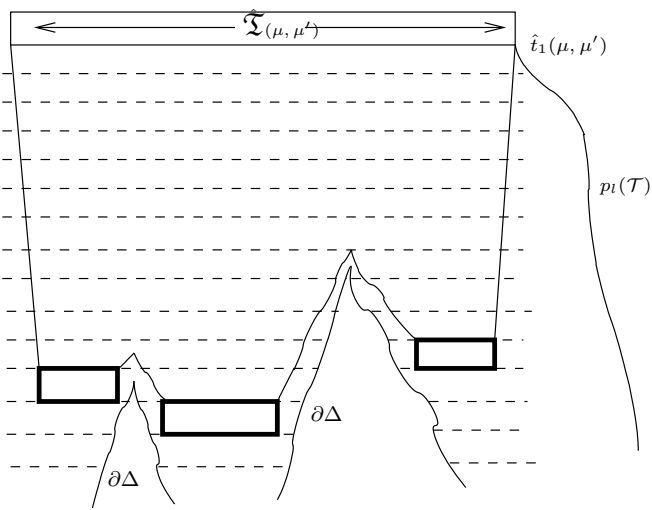 
 
\caption{Bounding the size of a team in terms of $|\down_1|$ and $|p_l|$} 
\label{figure:TeamAge} 
\end{center} 
\end{figure} 
 
\begin{proof}     The first thing to observe is that at any 
stage in the past of $\ET$ the set of letters lying in a 
single corridor form a connected region. As in Lemma \ref{G34pics}, 
this is simply a matter of noting that if $\phi(aub)=w$ where $w, a$ and $b$ 
consist only of constant letters, then $u$ must equal a word in constant letters. 
 
Consider the   past of $\ET$ at a time $t$. Write $k_t$ for the number of 
corridors that contain a non-trivial component of this past. 
The total 
increase in length of these components  when one goes forward to time $t+1$ 
is bounded by $2LC_4k_t$, since the connectedness of the past 
implies that  the only growth that can happen for existing components occurs 
at their extremities, where a block of at most $LC_4$ constant letters
may be added. This follows from Lemma \ref{G34pics}.  Also at time $t+1$,
constant letters from $\partial\Delta$ may join the past of $\ET$, and
there may be new components of constant letters (each of length less
than $2LC_4$) whose ancestors at time $t$ were non-constant
letters. Thus we have three possible causes of increase. The first and
third account for growth of at most $2LC_4k_{t+1}$ and the second
(boundary) contribution is the number of elements of $\partial^\T$
that occur at time $t+1$.  
If the genesis of $\T$ is of type (G1) or (G2), then at least $k_{t+1}$ 
edges of $\down_1(\T)$ occur at time $t$, compensating us for the growth summand 
$2LC_4k_{t+1}$. If the genesis of $\T$ is of type (G3) then we still have the 
above compensation {\em except} at those times where no edges of $\down_1(\T)$ occur. 
At these latter times the whole of the narrow past of $\T$ 
lies in a single corridor through which $\plT$ passes. Since the
narrow past lies
in a single corridor, it is connected and  grows at most $2LC_4$ when moving 
forward one unit 
of time (unless added to by $\partial^\T$).  

The summands $2LC_4\,|\QT|$ and  $2LC_4T_0\big(|\chi_P(\T)| +1\big) $ in item (2) of the lemma
account for the growth of the narrow past in the intervals of time
below  $t_3(\T)$, and from $t_3(\T)$ to $t_1(\T)$, respectively. The additional summand $\ll$
allows us to desist from our estimating if the narrow past of $\T$ ever shrinks to have
length less than $\ll$. 
\end{proof}

\subsection{Bounding the size of $\QT$}\label{Proofs} 

For the remainder of this section we concentrate exclusively on long teams of genesis (G3)
with $\QT$ non-empty. We denote the set of such teams by $\Gthree$.  
Our goal is to bound $|\QT|$. (In the light of our
previous results, this will complete the required analysis of the
length of teams.)  

Recall from  Definition \ref{PincerDef} that for teams of genesis 
(G3), the paths $\plT$ and $\prT$ and the chain 
of 2-cells joining them in the corridor at time $\ttwo$ form a pincer denoted $\Pin_{\T}$. 
The set $\subT$ was defined in Definition \ref{chiP}.

An important
feature of  teams in $\Gthree$ is:

\begin{lemma}\label{GotBlock}
If $\T \in \Gthree$ then there exists a block of  at least $\ll$
constant edges   immediately adjacent to $\Pin_{\T}$ at each time from
$t_3(\T)$ to the top of $\Pin_{\T}$, and adjacent to $p_l^+(\T)$ from
then until $\tone$. (At time $\tone$ this block contains $\ET$.)
\end{lemma}

\begin{proof} The hypothesis
that $\QT$ is non-empty means that the narrow past
of $\T$ at  some time before $t_3(\T)$ has length at
least $\ll$ and is contained in the same corridor as 
$\plT$ (see Definition \ref{QT}).  The definition of $t_3(\T)$ implies that
the  narrow past of $\T$ is contained in a block of constant letters 
immediately adjacent to $\plT$ or $p_l^+(\T)$ from time $t_3(\T)$
until $\tone$. Since the  length of the narrow past of $\T$ does not
decrease before $\tone$, these blocks of constant letters must have
length at least $\ll$.
\end{proof}

The following is an immediate consequence of the Pincer Lemma. 

\begin{lemma} \label{t1-t2Lemma}  
For all $\T \in \Gthree$, 
$$ 
  t_3(\T) - \ttwo = \life(\Pin_{\T}) \leq \ttt (|\subT|+1) .   
$$
\end{lemma}

\begin{lemma} \label{nesters}
If $\T_1, \T_2 \in \Gthree$ are distinct teams then $\chi(\Pin_{T_1})
\cap \chi(\Pin_{T_2})  = \emptyset$.
\end{lemma}  
 
\begin{proof} The pincers $\Pin_{\T_i}$ are either disjoint or else
one is contained in the
other. In the latter case, say  $\Pin_{\T_1}\subset\Pin_{\T_2}$,  
the existence of the  block of $\ll$ constant edges established in
Lemma \ref{GotBlock} means that $\Pin_{\T_1}$ is  
actually nested in $\T_2$ in the sense of Definition \ref{chiP}. Thus 
$\chi(\Pin_{\T_1}) \cap \chi(\Pin_{\T_2})  = \emptyset$ (by 
Definition \ref{chiP}). 
\end{proof}
 
\begin{corollary} \label{t1-t2Corr}
$\sum\limits_{\T \in \Gthree} t_3(\T) - t_2(\T) \leq 3\ttt \n$.
\end{corollary}

It remains to bound the number of edges in $\QT$ which occur before $\ttwo$;
this is  cardinality of the following set.
 
\begin{definition}\label{down2} 
For $\T\in\Gthree$ we define  \index{teams!$\down_i(\T)$}
$\down_2(\T)$ to be the set of edges in $\partial\Delta$ that 
lie at the righthand end of a corridor containing an edge in $\QT$
before time $\ttwo$.
\end{definition}

The remainder of this section is dedicated to obtaining a bound on
$$
\sum\limits_{\T\in\Gthree}|\down_2(\T)|,
$$ 
(see Corollary \ref{downbound}).
 
At this stage our task of bounding $\|\T\|$ would be complete if
the the sets $\down_2(\T)$ associated to distinct teams 
were disjoint --- unfortunately they need not be, because of the possible 
nesting of teams as shown in Figures \ref{figure:Nest}
and \ref{figure:doublecount}. Thus we shall 
be obliged to seek further pay-off for our troubles. To this end we 
shall identify two sets of consumed colours $\chi_c(\T)$ and 
$\chi_{\delta}(\T)$ that arise from 
the nesting of teams. 
 
In order to analyse  the 
effect of nesting we need the following vocabulary.

There is an obvious left-to-right ordering of those paths in the forest 
$\F$ which begin on the arc of $\partial\Delta\ssm\partial S_0$ that commences 
at the initial vertex of the left end of $S_0$. (First one orders the trees, then 
the relative order between paths in a tree is determined by the manner in 
which they diverge; the only paths which are not ordered relative to each 
other are those where one is an initial segment of the other, and this 
ambiguity will not concern us.) 

\smallskip

\noindent{\bf Notation:} We write $\Gthree'$ for the set of teams $\T
\in \Gthree$ such that $\down_2(\T) \neq \emptyset$.

\smallskip

We shall need the following obvious separation property.

\begin{lemma}\label{separate} 
Consider $\T \in \Gthree'$.  If a
path $p$ in $\F$ is to the left of $\plT$ and a path $q$ is the right
of $\prT$,
then there is no corridor connecting $p$ to $q$ at any time $t<\ttwo$. 
\end{lemma} 

\begin{proof}  The hypothesis $\down_2 (\T) \neq \emptyset$ implies
that before $\ttwo$ the paths $\plT$ and $\prT$ are not in the same
corridor.
\end{proof}

\begin{definition} \label{depthDef} 
$\T_1\in\Gthree'$ is said to be {\em below} $\T_2\in\Gthree'$ if
$p_l(\T_2)$ and $p_r(\T_2)$ both lie  between $p_l(\T_1)$ and
$p_r(\T_1)$ in the left-right ordering described above.

$\T_1$ is said to be {\em to the left} of $\T_2$ if both  $p_l(\T_2)$ and $p_r(\T_2)$ 
lie to the right of $p_r(\T_1)$. 

 \index{teams!depth of}
 
We say that $\T$ is at {\em depth} $0$ if there are no teams above it. 
Then, inductively, we say that a team is at depth $d+1$ if $d$ is the maximum 
depth of those teams above $\T$.   
 
A {\em final depth} team is one with no teams below it. 
 
Note that there is a complete left-to-right ordering of teams $\T\in\Gthree$ at any given depth. 
\end{definition}

\begin{lemma}\label{trapS_0} If there is a team from $\Gthree'$ below
$\T\in\Gthree'$, then $\tone \ge \time (S_0)\ge \ttwo$. 
\end{lemma}

\begin{proof} The first thing to note is that  if
 $\time (S_0)$ were less than $\ttwo$, then  the narrow
past of $\T$  at time $t_2(\T)$ must contain at least
$\ll$ edges. This is because the length of the narrow past of $\T$ cannot
decrease before $\tone$, and
at $\time(S_0)$ the narrow past
is the union of the intervals $C_{(\mu,\mu')}(2)$ with $(\mu,\mu')\in\T$,
which has length at least $\ll$ since $\T$ is assumed not to be short.

Thus if $\time(S_0)<\ttwo$ then  we are in the non-degenerate
situation of Definition \ref{PincerDef} and the defining property of $\ttwo$
means that  before time $\ttwo$  no edge to the right of $\prT$ lies
in the same
corridor as all the colours of $\T$ (cf. Lemma \ref{separate}).  In
particular this is true of
the past of the reaper of $\T$ (assuming that it has a past at time
$\ttwo$). On the
other hand,   the reaper of $\T$ has a past in $S_0$ (by the very
definition of a team), as do all of the colours of $\T$. And since they
lie in a common corridor at $\time(S_0)$, they must also do so
at all times up to $\tone$. This contradiction implies
that in fact $\time(S_0)\ge\ttwo$. 
 
Consider Figure \ref{figure:Nest}.  
Suppose that $\T'\in\Gthree'$ is below $\T$. The proof of Lemma \ref{GotBlock}
tells us that there is a block of constant edges extending from the
top of $\Pin_{\T'}$ containing the narrow past of $\T'$, and there is
a similarly long block
extending from the path $p_l^+(\T)$ at each subsequent time until
$t_1(\T')$. Thereafter the future of the block is contained in the
block of constant edges that evolves into the union of the
$C_{(\mu,\mu')}(2) \subseteq \bot(S_0)$ with $(\mu,\mu') \in \T'$,
which is long by hypothesis.
 
 At no time can this evolving block extend across $\plT$ 
because by  definition the edges along $\plT$ are labelled by non-constant
letters.  Thus the evolving
block is trapped to the right of $\plT$ and to the left of 
$\prT$. In particular, it must vanish entirely before the time
at the top of the pincer $\Pin_\T$, which is no later than
$\tone$ and therefore $\tone \ge \time (S_0)$. 
\end{proof} 
 
The following is the main result of this section.

\begin{lemma}   \label{Aget2Lemma} \label{chiT} There exist sets of
colours \index{teams!colours $\chi_*(\T)$} 
$\chi_c(\T)$ and 
$\chi_{\delta}(\T)$ associated to each team
$\T\in\Gthree'$ such that the sets associated to
distinct teams are disjoint and the following inequalities hold. 

For each fixed team $\T_0 \in \Gthree'$ (of depth $d$ say),  
the teams  of depth $d+1$ that lie below   $\T_0$ 
may be described as follows: 
\begin{enumerate} 
\item[$\bullet$] There is at most one \index{teams!distinguished}
{\em distinguished team} $\T_1$,
and
$$ 
\|\T_1\|\le 2B\Big(\ttt (1+ |\chi(\Pin_{\T_0})|) + T_0(|\chi_P(\T_0)| +
1)\Big).
$$
\item[$\bullet$] There are some number of final-depth teams. 
\item[$\bullet$] For each of the remaining teams $\T$ we have 
  $$|\down_2(\T_0) \cap\down_2(\T)| \le \ttt  \Big( 1 
+ |\chi_c(\T)| \Big) + T_0 \Big( |\chi_{\delta}(\T)| + 2 \Big).
$$ 
\end{enumerate}
\end{lemma}  

\begin{proof}  The first thing to note is that if two teams $\T, \T'
\in \Gthree'$ are at the same depth, then $\down_2(\T)$ and
$\down_2(\T')$ are disjoint. Indeed if $\T$ is to the left of $\T'$,
then at times before $\ttwo$ the paths $\plT$ and $p_l(\T')$ never lie in
the same corridor.  Let $\T \in \Gthree'$ be a team of level $d+1$
that is below $\T_0$ and consider the edge $e$  at the right end of a
corridor earlier than $\ttwo$ that contains an edge in $\QT$. We are
concerned with the fact that this edge may be in $\down_2(\T_0)$.  In
this situation we say that $\T_0$ and $\T$ {\em double count} $e$.

\begin{figure}[htbp] 
\begin{center} 
  
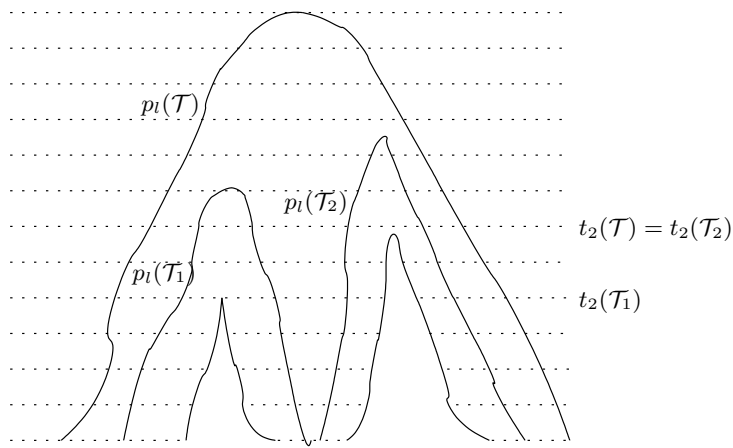 
  
\caption{A depiction of double-counting} 
\label{figure:doublecount} 
\end{center} 
\end{figure} 
 
Let $\T_1, \dots ,\T_r$ be the teams in $\Gthree'$ of depth $d+1$
which double-count with $\T_0$, ordered from left to right, with the
final-depth teams deleted.  We define $\chi_c(\T)$ to be empty for
teams not on this list. $\T_1$ will be the distinguished team.

Since there is no double-counting between teams of the same level, 
 the sets of times  at which 
$\T_1, \dots,\T_r$ double-count with $\T_0$ must be disjoint. Indeed if 
 $i < j$ then the set of times at which $\T_i$ double-counts 
  with $\T_0$ is earlier than the 
set of times at which $\T_j$ 
double-counts with ${\T_0}$ (Lemma \ref{separate}). Moreover, 
 the times for each $\T_i$ form an interval, which we denote $\I_i$.  
 
We assume $r\ge 2$ and describe the construction of 
the sets $\chi_c(\T_i)$ and $\chi_{\delta}(\T_i)$ that account for
double-counting.

The first thing to note is that each $\I_i$ must be later than $t_2(\T_1)$, 
 by Lemma \ref{separate}. 
The second thing to note is that   the entire interval of time $\I_i$ 
must also be earlier than $t_1(\T_1)$. Indeed if  some double-counting by 
$\T_i$  and $\T_0$  were to occur  after $t_1(\T_1)$, then we would
have $t_2(\T_k) > t_1(\T_1)$. But then
$\time (S_0) > t_1(\T_1)$, so Lemma \ref{trapS_0} would imply  that
there was no team below $\T_1$, contrary to hypothesis.  

We separately consider the intervals $\I_i \cap [ t_2(\T_1),t_3(\T_1)
]$ and  $\I_i \cap [ t_3(\T_1),t_1(\T_1) ]$, whose union is all of
$\I_i$.

For that part of $\I_i$ before $t_3(\T_1)$, the proofs of
the Pincer Lemma (Theorem \ref{PincerLemma}) and Proposition 
\ref{prePincerLemma} tell us that colours in $\chi(\Pin_{\T_1})$
will be consumed at the rate of at
least one per $\ttt$ units of time. Define
$\chi_c(\T_i)$ to be this set of consumed colours. We have
$$
\Big|\, \I_i \cap [ t_2(\T_1),t_3(\T_1) ]\, \Big| \leq \ttt (1 +
|\chi_c(\T_i)|) .
$$

Now consider  $\I_i \cap [ t_3(\T_1),t_1(\T_1)]$. Define
$\chi_{\delta}(\T_i)$ as follows. The discussion in Definition
\ref{PincerDef} shows that in any period of time of length $T_0$ in
the interval $[t_3(\T_1),t_1(\T_1) ]$ at least one colour in
$\chi_P(\T_1)$ disappears.  Let $\chi_{\delta}(\T_i)$ be
the set of colours in $\chi_P(\T_1)$ which disappear during $\I_i
\cap [t_3(\T_1), t_1(\T_1)]$ (these disappearances correspond to the
discontinuities in the `path' $p_l^+(\T_1)$).  By construction, we
then have\footnote{There is a 2 rather than the familiar 1 on the
right to account for the colour containing
$p_l(\T_1)$, which is not included in
$\chi_P(\T_1)$; there might be up to $T_0$ corridors between
$t_3(\T_1)$ and the top of $\Pin_{\T_1}$.}

\[	\Big| \, \I_i \cap  [ t_3(\T_1),t_1(\T_1)]\, \Big| \le
T_0(|\chi_{\delta}(\T_i)| + 2),	\]
and combining these estimates we have
$$
|\I_i| \le \ttt  \Big( 1 + |\chi_c(\T_i)| \Big) + T_0
\Big(  |\chi_{\delta}(\T_i)| +  2\Big) , 
$$
as required. 
Since the intervals $\I_i$ are disjoint, 
the sets $\chi_c(\T_i),\, i=2,\dots,r$ are mutually disjoint. 
And by construction, these sets are also disjoint from the sets 
associated to teams other than 
the $\T_i$ under consideration (i.e. those under other depth $d$  
teams, or those of different 
depths). The same considerations hold for the sets
$\chi_{\delta}(\T_i),\, i=2,\ldots, r$.

In Figure \ref{figure:double}, the shaded region is where we recorded
the regular disappearance of the colours forming $\chi_c(\T_i)$,
whilst in Figure \ref{figure:doubletwo}, the shaded region is where we
recorded the regular disappearance of the colours forming
$\chi_{\delta}(\T_i)$.

\begin{figure}[htbp] 
\begin{center} 
  
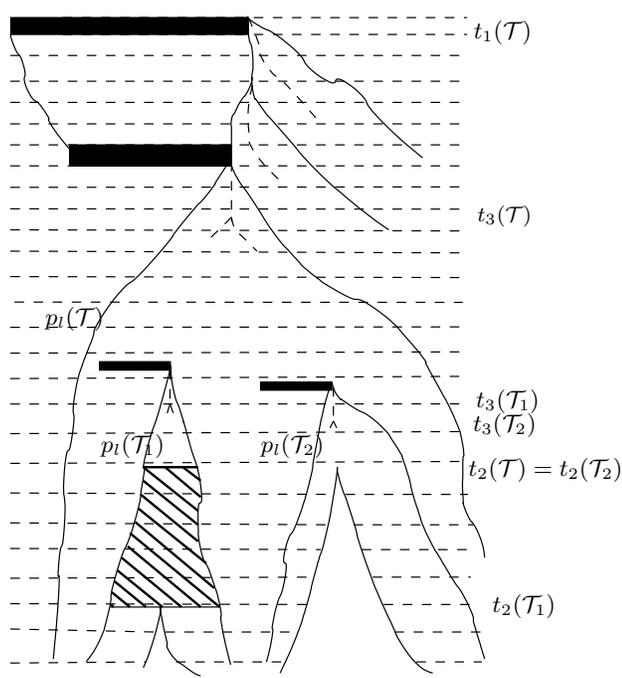 
  
\caption{Finding the colours $\chi_c(\T_i)$} 
\label{figure:double} 
\end{center} 
\end{figure}

\begin{figure}[htbp]
\begin{center}

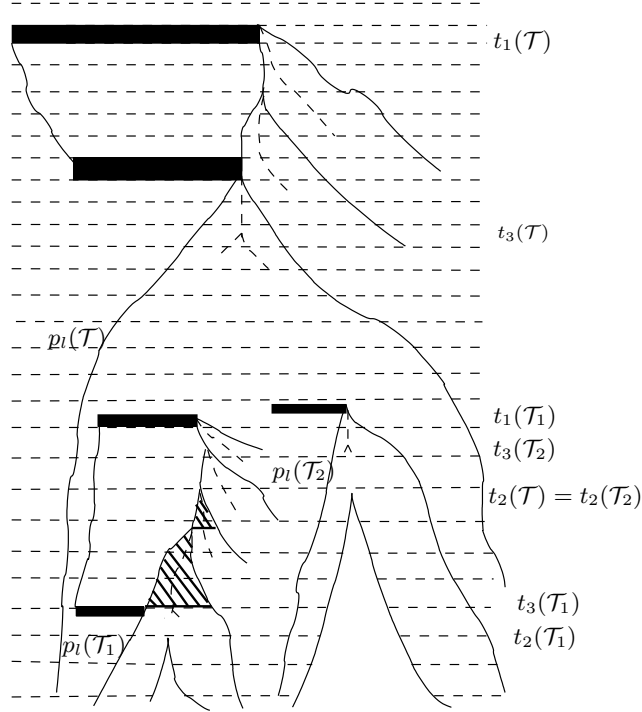

\caption{Finding the colours $\chi_{\delta}(\T_i)$}
\label{figure:doubletwo}
\end{center}
\end{figure}
 
It remains to establish the inequality 
$$ 
\|\T_1\|\le 2B\Big( \ttt (|\chi(\Pin_{\T_0})| + 1) +
(|\chi_P(\T_0)|+1) \Big). 
$$ 
We first note (as in the proof of Lemma \ref{trapS_0})  that 
$\ET_1$ is trapped between $\plT$ and $\prT$, so it must be consumed
entirely between  the times $t_1(\T_1)$ and  $t_1(\T_0)$. But by the
Bounded Cancellation Lemma, the length of  the future of $\ET_1$ can
decrease by at most $2B$ at each step in time. Therefore  $\|\T_1\|\le
2B (t_1(\T_0) - t_1(\T_1))$.
 
$\T_1$ is assumed not be final-depth, so from Lemma \ref{trapS_0} we have 
$t_2(\T_0) \le \time (S_0) \le t_1(\T_1)$. By combining these  
 inequalities with Lemmas \ref{t1-t2Lemma} and \ref{t1t3forTeam}
we obtain: 
\begin{eqnarray*} 
\|\T_1\| &\le& 2B\, \Big(t_1(\T_0) - t_1(\T_1)\Big)\\  
& \leq & 2B\, \Big(t_1(\T_0) - \time (S_0)\Big)\\  
& \leq & 2B\, \Big(t_1(\T_0) - t_2(\T_0)\Big)\\  
& \leq & 2B \Big[ \ttt \Big(1+|\chi(\Pin_{\T_0})|\Big) + T_0
\Big( |\chi_P(\T_0)| + 1 \Big) \Big]. 
\end{eqnarray*} 
\end{proof} 
 
\begin{corollary} \label{Firstdown2sum} Summing over the set of teams $\T \in \Gthree'$ that
are not distinguished, we get 
$$ 
\sum_\T \Big|\down_2 (\T)\Big| \le2\,\Big| \bigcup_\T \down_2(\T)\Big|
+\sum_{\T}  \ttt  \Big( 1 
+ |\chi_c(\T)| \Big) + \sum_{\T} T_0 \Big( |\chi_{\delta}(\T)| + 2
\Big) .
$$ 
\end{corollary} 
 
\begin{proof} Suppose $\T \in \Gthree'$ of depth $d+1$ is not
final-depth and not distinguished, and that $\T$ double-counts with
some $\T_0$ of depth $d$ above it.  Then, by Lemma \ref{Aget2Lemma},
we have
\begin{eqnarray*}
|\down_2(\T)| & = & |\down_2(\T) \smallsetminus \down_2(\T_0)| +
|\down_2(\T) \cap \down_2(\T_0)| \\
& \leq & |\down_2(\T) \smallsetminus \down_2(\T_0)| + \ttt (1 +
|\chi_c(\T)|) + T_0 (2 + |\chi_{\delta}(\T)|).
\end{eqnarray*}
Suppose that $\T' \in \Gthree'$ is a team of depth $k < d$ and that
$\T'$ is above $\T$. If $\T$ double-counts with $\T'$ at time $t$,
then $\T$ double-counts with $\T_0$ at time $t$, by Lemma
\ref{separate}.  Therefore, the set of edges that $\T$ double-counts
with any team of lesser depth is exactly $\down_2(\T) \cap
\down_2(\T_0)$.

Thus we have accounted for all double-counting other than than
involving final depth teams.  The factor $2$ in the statement of the
corollary accounts for this.
\end{proof} 

And summing over the same set of teams again, we obtain: 
\begin{corollary} \label{downbound} 
$$ 
\sum_\T |\down_2(\T)| \ \le \ \n (2 + 3\ttt + 5T_0). 
$$
\end{corollary} 
\begin{proof} The sets of colours $\chi_c(\T)$ and $\chi_{\delta}(\T)$
are disjoint. And the
union of the sets $\down_2(\T)$  is a subset of $\partial\Delta$. The
set of all colours and the  set of edges in $\partial\Delta$ each have
cardinality at most $\n$.   And the number of teams is less than $2\n$
(Lemma \ref{allIn}).
\end{proof}

\section{The Bonus Scheme} \label{BonusScheme}

We have defined teams and obtained a global bound on 
$\sum\|\T\|$. 
If $\cmm$ is non-empty then $(\mu,\mu')$ is a member or  
 virtual member of a  unique team. 
If this team is such that $\tone \ge \time (S_0)$, then no member of 
the team is virtual  and we have the inequality
$$\|\T\|>\sum\limits_{(\mu,\mu') \in \T}|\cmm| - B$$
 established in Lemma \ref{t1high}.
 We indicated following this lemma how this inequality
 might fail in the case where
  $\tone < \time (S_0)$. In this section we take up this
  matter in detail
   and introduce a  \index{bonus scheme}{\em bonus scheme} that
 assigns additional edges to teams in order to compensate for the possible failure
 of the above inequality when   $\tone < \time (S_0)$.
 
By definition, at time $\tone$ the reaper $\r=\rT$ lies
 immediately to the right of $\ET$.  The edges of 
$\ET$ not consumed from the right by $\r$ by  $\time (S_0)$
 have a preferred future in $S_0$ 
that lies in $\cmm$ for some member $(\mu,\mu')\in\T$. 
However, not all of the edges of 
$\cmm$ need arise in this way:
 some may not have  a constant ancestor at time $\tone$.
And  if $(\mu,\mu')$ is only a virtual member of $\T$,
 then no edge of $\cmm$ lies in the 
future of $\ET$. The  {\em bonus} edges in $\cmm$
are a certain subset of those  that do not have a constant 
ancestor at time $\tone$. They are defined
as follows.

\begin{definition} \label{def:swol}
 Let $\T$ be a team with $\tone < \time(S_0)$
and consider a time $t$ with $\tone < t < \time(S_0)$.

The \index{future!swollen}
{\em swollen future} of $\T$ at time $t$ is the interval 
of constant edges beginning immediately to the left of the pp-future of $\rT$. 
 
Let $e$ be a non-constant edge that lies immediately to the left of the 
swollen future of $\T$ but whose ancestor is not a 
right para-linear edge in this position. If $e$ is a right para-linear and
 the (constant) rate 
at which $e$ adds letters to the swollen future of $\ET$ is greater 
than the (constant) rate at  which the future of the reaper cancels letters 
in the future of $\ET$, then we define $e$ to be 
a \index{rascal}{\em rascal}; if $e$ is right-fast then we define it to be a 
\index{terror}{\em terror}.
In both cases,
 we define the
{\em bonus provided by $e$} to be the set of edges in the swollen future
of $\T$ in $S_0$ that have $e$ as their most recent non-constant
ancestor, and are eventually consumed by $\rT$. 

The set $\bonusT$ is the union of the bonuses provided to $\T$ by all
rascals and terrors.
\end{definition} 
 
\begin{lemma}  \label{C1toTeamLength} 
For any team $\T$, 
\[      \sum_{(\mu,\mu') \in \T \mbox{ \tiny{or} } (\mu,\mu')
 \vin T}|\cmm| \leq \|\T\| + |\bonus(\T)| + B.    \] 
\end{lemma}

\begin{proof}
If $\tone \geq \time(S_0)$, this follows immediately from Lemma
\ref{t1high}.  If $\tone < \time(S_0)$ then at each step in time
between $\tone$ and $\time(S_0)$ the only possible cause of growth in
the length of the swollen future of the team is the possible action of
a rascal or terror if such is present at that time.  (There is no
interaction of the swollen future with the boundary or singularities,
because of the exclusions in the second paragraph of Definition
\ref{newTeams}.)

The swollen future has length $\|\T\|$ at time $\tone$ and length at
least  $\sum |C_{(\mu,\mu')}(2)|$ at $\time(S_0)$. By definition,
$|\bonus(\T)|$ is a bound on the growth in length between these
times. (The summand $B$ is thus unnecessary in the case  $\tone <
\time(S_0)$.)
\end{proof}
 
 The following lemma shows that our main task in this
 section will be to analyse the behaviour of rascals.

\begin{lemma} \label{lem:TerrorCont}
The sum of the lengths of the bonuses 
provided to all teams by terrors is less than $2L\n$.
\end{lemma}  

\begin{proof} Since it is right-fast, a terror will be separated
from the team to which it is associated after one unit of
time, and hence the bonus that it provides is less than $L$.
There is at most one terror for each possible adjacency of colours
and hence the total contributions of all terrors is less than
$2L\n$.
\end{proof}
 
The typical pattern of influence of rascals on a team 
is shown in Figure \ref{belowS0};
there may be several times at which rascals
appear at the left of $\T$ and provide a  
bonus for the team before being consumed from the left (or otherwise detached 
from the team).  
 
\begin{figure}[htbp] 
\begin{center} 
  
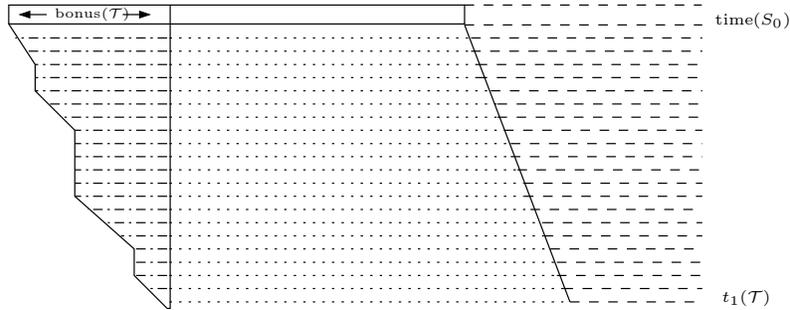 
  
\caption{The generic situation below $\time(S_0)$.} 
\label{belowS0} 
\end{center} 
\end{figure}

\begin{definition}[Rascals' Pincers] \label{RascalPin}
We fix a team $\T$ with $\tone < \time(S_0)$
 and consider the interval of time $[\tau_0(e),\tau_1(e)]$,
where $\tau_0(e)$ 
is the time at which a rascal $e$ appears at the left end of the swollen 
future of $\T$, and $\tau_1(e)$ is the time at which its future is no 
longer to the immediate left of the future of the swollen future of $\T$.  
 
In the case where the pp-future $\hat e$ of $e$ at time  $\tau_1(e)$
is cancelled from the left by an edge $e'$, we define
$\tau_2(e)$ to be the earliest time when the pasts of $\hat e$ and $e'$ are
in the same corridor. The path in $\F$ that traces the pp-future of
$e$ up to  $\tau_1(e)$ is denoted $p_e$ and the path following through the
ancestors of  $e'$ from
 $\tau_2(e)$ to $\tau_1(e)$ is denoted  $p_e'$.
 The pincer\footnote{to lighten the terminology, here we allow the
 degenerate case where the ``pincer" has
  no colours  other than those of $e$
 and $e'$}
  formed by $p_e$ and $p_{e}'$  with base at time  $\tau_2(e)$
is denoted  $\Pin_e$. 
\end{definition}
 
\begin{lemma} \label{lowRascals} 
The total of all bonuses provided to all teams by  rascals $e$ with 
 $\tau_1(e)\le\time(S_0)$ is less than $(3\ttt + 2T_0 +1)L\n$.  
\end{lemma} 
 
\begin{proof} Consider a rascal $e$. We defer the case where $e$ hits
a singularity or the  boundary. If this does not happen, the
pp-future $\hat e$ of $e$ at time $\tau_1(e)$ is cancelled from the left
by an edge $e'$ (which  is right-fast since $e$ is not
constant). We consider the pincer $\Pin_e$ defined above.
The presence of the swollen future of $\T$ at the top of the
pincer allows us to apply the  Two Colour Lemma to conclude that
$\tau_1(e) - T_0 \geq \time(S_{\Pi_e})$ (in the degenerate case
discussed in the footnote, $ \time(S_{\Pi_e})$ is replaced by
$\tau_2(e)$). And the
Pincer Lemma tells us that 
\[      \tau_1(e) - \tau_2(e) \leq \ttt \Big( 1 + |\chi(\Pin_e)| \Big)
+ T_0.       \]
In fact, we could use $\tilde \chi(\Pin_e)$ instead of
$\chi(\Pin_e)$ in this estimate because there cannot be any nesting
amongst the pincers $\Pin_e$ with $\tau_1(e)\le \time(S_0)$,
because nesting would imply that the swollen future of $\T$, which is
immediately to the right of the lower rascal, would be trapped beneath
 the upper pincer,
  contradicting the fact that the team has a non-empty future in $S_0$. 
 
In the case where $e$ hits the boundary or is separated from the team
 by a singularity (at time $\tau_1(e)$) we 
define  $\tau_2(e)=\tau_1(e)$. No matter what the
fate of $e$, we define
 $\partial^e$ to be the set of edges in
 $\partial\Delta$ at the left ends of corridors  
containing the future of $e$ between 
$\tau_0(e)$ and $\tau_2(e)$.
The sets $\partial^e$ assigned to different rascals are disjoint,
so summing over all rascals with $\tau_1(e)\le \time(S_0)$
we have
\begin{eqnarray*}
\sum_e \Big(\tau_1(e) - \tau_0(e)\Big) & = & \sum_e (\tau_1(e) - \tau_2(e)) +
 (\tau_2(e)-\tau_0(e))\\
& \le & \sum_e \ttt\Big(1 + |\chi (\Pin_e)|\Big) + T_0 + |\partial^e|.
\end{eqnarray*}
Since the sets $\chi(\Pin_e)$ and $\partial^e$ are disjoint,
the terms $\ttt |\chi (\Pin_e)|$ and $ |\partial^e|$ 
contribute less than $(\ttt + 1)\n$ to this sum. And since the number of
 rascals is bounded by the number of possible adjacencies of colours, the remaining
terms contribute at most $(\ttt + T_0)2 \n$. Thus
$$
\sum_e  \Big(\tau_1(e) - \tau_0(e)\Big)\ \le\ (3\ttt + 2T_0 + 1) \n .
$$
The bonus produced by each rascal in each unit of time is less than $L$, so  
the lemma  is proved.
\end{proof} 

It remains to consider the size of the bonuses provided by rascals $e$ with 
 $\tau_1(e)>\time (S_0)$.

The bonuses that are not accounted for in Lemma \ref{lowRascals}
reside in blocks  of constant edges along $\bot(S_0)$ each of which is the swollen 
future of some team, with
 a \rpl letter at its left-hand end (the pp-future of
  a rascal) and a \lpl letter at its left-hand end (the pp-future of
  the team's reaper).  
 
 \begin{definition}\index{rascal!left(right)-biased}
A {\em left-biased} rascal $e$ is one with
 $\tau_1(e) > \time(S_0)$ that satisfies the following properties: 
\begin{enumerate} 
\item[1.] the pp-future of the rascal  is (ultimately) 
consumed from the left by an edge of $S_0$,  
\item[2.] the swollen future of  $\T$ at 
time $\tau_1(e)$ has length at least $\ll$ and  
the pp-future of the reaper $\rT$ is still immediately to its right.
\end{enumerate} 
 \end{definition}
 
\def\life{\text{\rm{life}}}  
\def\B{\text{\euf{B}}}

\begin{definition} Let \index{blocks $\B$}$\B\subset\bot(S_0)$ be an interval of constant edges with a  
right para-linear letter at its left-hand end and a left-linear letter $\r$ at its right-hand end. We 
say that $\B$ is {\em right biased} if $\r$ is ultimately consumed by an edge (to its right) 
in $S_0$. We define $\life(\B)$ to be the difference  between $\time(S_0)$ and  the time at which 
the \lpl letter $\rho$ is consumed. 
And we define the {\em effective volume} of $\B$ to be the number of edges in $\B$ 
that are ultimately consumed by $\rho$. 
\end{definition} 
 
We have the following \index{tautologous tetrad}
tautologous tetrad of possibilities covering the swollen teams whose bonuses are 
not entirely accounted for by Lemma \ref{lowRascals}. 
 
\begin{lemma} \label{tetrad} Let $\B\subset\bot(S_0)$ be an interval of constant edges that is 
the swollen future of a team with a rascal 
at its left-hand end and a \lpl letter $\r$ at its right-hand end.  
Then at least one of the following holds: 
\begin{enumerate}  
\item[\rm{(i)}] the length of $\B$ is at most $\ll$;
\item[\rm{(ii)}] $\B$ is the swollen future of a team with a left-biased rascal; 
\item[\rm{(iii)}] $\B$ is right-biased;
\item[\rm{(iv)}] neither of the non-constant letters at the ends of $\B$
 is ultimately consumed by an edge of $S_0$. 
\end{enumerate} 
\end{lemma} 

We note here that when the length of $\B$ is at most $\ll$ then we
have a short team, and we have already accounted for short teams.
The following three lemmas correspond to eventualities (ii) to (iv).
  
\begin{lemma} \label{highRascals} 
The sum of the bonuses provided to all teams by left-biased rascals is
less than $(2L + 6L\ttt + 4LT_0 + 2\ll + 6B\ttt + 4BT_0) \n$.
\end{lemma}

\begin{proof} The proof of this result is similar to the work done in
the previous section.  We have a pincer $\Pin_e$ associated to the
rascal $e$.  Since we are only concerned with the times when the
rascal is immediately adjacent to a block of constant letters, it must
be that at time $\tau_1(e) - T_0$ either we are below $\tau_0(e)$ or
$\time(S_{\Pin_e})$ (cf. Definition \ref{PincerDef}). Therefore the
following is an immediate consequence of the Pincer Lemma.
\[	\tau_1(e) - \tau_2(e) \leq \ttt (1 + |\chi(\Pin_e)|) + T_0 .
\]
It now suffices to bound the amount of time for which $e$ is adjacent
to the narrow past of $\B$ before $\tau_2(e)$.  We define $\tau_0'(e)$
to be the latest time when the rascal $e$ has contributed less than
$\ll$ edges to $\bonus(\T)$.  Then the bonus provided by $e$ is at
most $L(\tau_1(e) - \tau_0'(e)) + \ll$. As in the previous section, we
define $\down_2(e)$ to be those edges on the left end of corridors
containing $e$ at times before $\tau_2(e)$ but after $\tau_0'(e)$.
Just as in Lemma \ref{Aget2Lemma} and the corollaries immediately
following it, we then have a notion of {\em depth} of
rascals describing the nesting of the pincers $\Pin_e$\footnote{One
extends the paths $p_e$ and $p_e'$ of Definition 
\ref{RascalPin} back in time to $\partial\Delta$ so as to define the
order defining depth.}. We also have
{\em distinguished} rascals (corresponding to the distinguished teams
in Lemma \ref{Aget2Lemma}), and proceeding as in the proof of Lemma
\ref{Aget2Lemma} we get the following estimates:

if $e_1$ is a distinguished rascal of depth $d+1$ and $e_0$ is the
rascal of depth $d$ above it, then the bonus provided by $e_1$ is at
most $2B\Big( T_1 (1 + |\chi(\Pin_{e_0}))|) + T_0 \Big)$, since all of
the bonus provided by $e_1$ must disappear before $\tau_1(e_0)$;

for other rascals $e$ of depth $d+1$ which are below $e_0$ we have a
set of colours $\chi_c(e)$, disjoint for distinct teams such that
\[	|\down_2(e) \cap \down_2(e_0)| \leq T_1(1 + |\chi_c(e)|) +
T_0.	\]
Therefore, summing over the set of rascals which are not distinguished
we get (cf Corollary \ref{Firstdown2sum})
\[ \sum_e |\down_2(e)| \leq 2\Big| \bigcup_e \down_2(e)\Big| + \sum_e \Big
( \ttt (1 + |\chi_c(e)|) + T_0 \Big).	\]
And summing over the same set of rascals, we get 
\[ \sum_e |\down_2(e)| \leq (2 + 3\ttt + 2T_0)\n .	\]
Therefore, for undistinguished rascals, we have
\begin{eqnarray*}
\sum_e \tau_1(e) - \tau_0'(e) & = & \sum(\tau_1(e) - \tau_2(e)) +
\sum(\tau_2(e) - \tau_0'(e))\\
& \leq & (3\ttt + 2T_0)\n + (2 + 3\ttt + 2T_0)\n,
\end{eqnarray*}
and so the contribution of all left-biased rascals is at most
\[	\Big( (2 + 6\ttt + 4T_0)L + 2\ll + 6B\ttt + 4BT_0 \Big) \n,	\]
as required.
\end{proof}

\begin{lemma} \label{RightBiased}  The sum $\sum \life(\B)$ over those
$\B$ that are right-biased but  do not satisfy conditions (i) or (ii)
of Lemma \ref{tetrad} is at most $(3\ttt B + 2T_0B)\n$.
\end{lemma} 

\begin{figure}[htbp] 
\begin{center} 
  
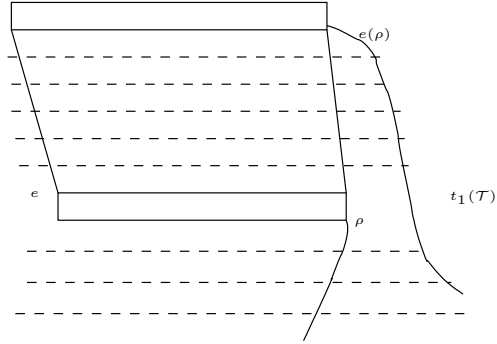 

\caption{A depiction of a right-biased team.} 
\label{r-teams} 
\end{center} 
\end{figure}

\begin{proof}
Once again, as in Lemmas \ref{lowRascals} and \ref{highRascals}, we
obtain compensation for the continuing existence of a non-constant
letter by using the Pincer Lemma to see that colours must be consumed
at a constant rate in order to facilitate the life of $\r$. Thus we
consider the left-fast edge that consumes the pp (i.e. left-most
non-constant) future of $\r$; this edge is denoted $e(\r)$ in Figure
\ref{r-teams}. The Pincer Lemma and the 2 Colour Lemma tell us that if
$\Pin_{e(\r)}$ is the
pincer associated to these paths (with $S_0$ at the bottom) then 
\[      \life(\B) \leq \ttt (1 + |\chi(\Pin_{e(\r)})|) + T_0.   \]
Suppose that $\B$ and $\B'$ are two right-biased blocks with
associated edges $e(\r)$ and $e(\r')$ consuming their reapers. We claim
that the sets $\chi(\Pin_{e(\r)})$ and $\chi(\Pin_{e(\r')})$ are
disjoint.  The key point to observe is that since we are not in case
(ii) of Lemma \ref{tetrad} the length of the swollen future of $\B$
increases from $\time(S_0)$ to the top of $\Pin_{e(\r)}$; since $\B$
had length at least $\ll$, we therefore have a block of more than
$\ll$ of more than $\ll$ constant edges at the top of
$\Pin_{e(\r)}$. Thus the pincers associated to $\B$ and $\B'$ are
either disjoint or nested. Hence $\chi(\Pin_{e(\r)})$ and
$\chi(\Pin_{e(\r')})$ are disjoint. Thus summing over all right-biased
blocks $\B$ we obtain
\[      \sum_{\B \mbox{ right-biased}} \life(\B) \leq (3\ttt B + 2T_0B)
\n, \]
as required.
\end{proof} 
Since any letter consumes less than $L$ constant letters in any unit of time, we conclude: 
\begin{corollary} \label{rightCor} The sum  of the effective volumes of  
all blocks that are right-biased but  do not satisfy conditions (i) and (ii) 
of Lemma \ref{tetrad} is at most $(3L\ttt B + 2LT_0B)\n$.  
\end{corollary} 
 
\begin{lemma} The sum of all blocks that satisfy condition (iv) of  
Lemma \ref{tetrad} is at most $(2B +1)\n$. 
\end{lemma} 
 
\begin{proof} Possibility (iv) involves several subcases: the key
event which halts the growth of the swollen future of $\B$ may be a
collision with $\partial\Delta$ or a singularity;  it may also be that
the key event is that the future of the rascal or reaper adjacent to
$\B$ is cancelled by an edge that is not in the future of $S_0$.  
 
But no matter what these key events may be, since we are in not in
cases (ii) or (iii), associated to the blocks in case (iv) we have the
following set of paths partitioning that part of the diagram $\Delta$
bounded by $S_0$ and the arc of $\partial\Delta$ connecting the
termini of the edges at the ends of $S_0$:

The path $\pi_l$ begins at $\time (S_0)$ and follows the pp-future of
the rascal at the right-end of the future of $\B$ until it hits the
boundary, a singularity, or else is cancelled by an edge $\e_l$ not in
the future of $S_0$; if it hits the boundary, it ends; if it hits a
singularity, $\pi_l$ crosses to the bottom of the corridor $S$ on the
other side of the singularity, and turns left to follow $\bot(S)$ to
the boundary (see Figure \ref{Pi_lOne}); if $\e_l$ cancels with the
pp-future of the rascal, then $\pi_l$ follows the past of $\e_l$
backwards in time to the boundary (see Figure \ref{Pi_lTwo}).

\begin{figure}[htbp] 
\begin{center} 
  
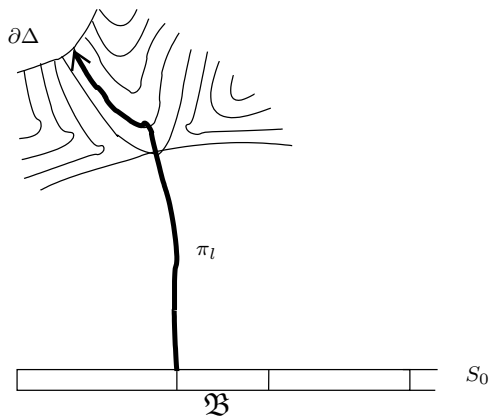 
  
\caption{The path $\pi_l$ hits a singularity.} 
\label{Pi_lOne} 
\end{center} 
\end{figure} 

\begin{figure}[htbp]
\begin{center}

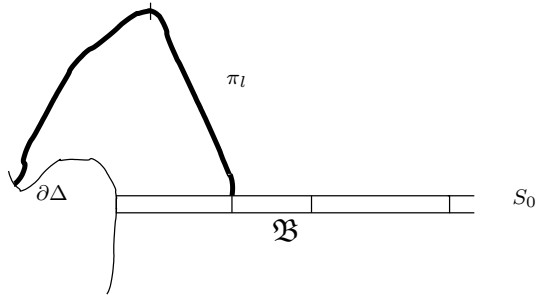

\caption{The path $\pi_l$ in cancelled from outside of the future of
$S_0$.}
\label{Pi_lTwo}
\end{center}
\end{figure}

The path $\pi_r$ describing the fate of $\r$ is defined similarly
(except that it turns right if it hits a singularity). 
 
It is clear from the construction that no two of these paths can cross, thus we have the  partition  
represented schematically in Figure \ref{partition}. 
 
\begin{figure}[htbp]
\begin{center}

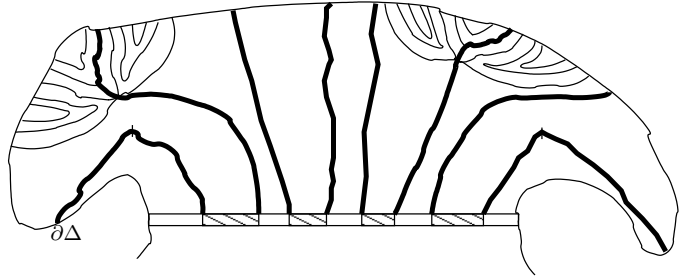

\caption{The schematic partition of $\Delta$ by the paths $\pi_l$ and
$\pi_r$.}
\label{partition}
\end{center}
\end{figure}

\def\bdy{\text{\rm{bdy}}} 
\def\up{\text{\rm{up}}} 
Given a swollen team $\B$ of type (iv), we follow the swollen future of $\B$ until its flow is 
interrupted (at time $\iota(\B)$, say) 
by meeting a singularity, the boundary of $\Delta$, or else its rascal or reaper is cancelled. 
Consider the set of corridors that contain some component of the swollen future of $\B$ 
after $\iota(\B)$. Consider also the set of edges $\bdy(\B) \subseteq
\partial\Delta$ 
that lie in the swollen future of $\B$. We keep account of the set of corridors by recording 
the set of their ends on $\partial\Delta$, except that we ignore an end if we have to cross 
a path $\pi_l$ or $\pi_r$ to reach it. Note that at least one end of each corridor is recorded. 
Let $\up(\B)\subset\partial\Delta$ denote the  set  of ends recorded.  
  
Since the sets $\bdy(\B)$ and $\up(\B)$ are contained in the portion of $\partial\Delta$ accorded 
to $\B$ by the partition formed by the paths $\pi_l$ and $\pi_r$, the sets associated to different 
$\B$ are disjoint. In each unit of time beyond $\iota(\B)$ each
component of the swollen future of $\B$ can shrink by at most $2B$ (by
Lemma \ref{BCL}).  The set $\up(\B)$ measures the sum of the number of
components over all such times, and $|\bdy(\B)|$ is the number of
uncancelled edges. Thus we see that the length of
the swollen future of $\B$ at time $\iota(\B)$ is at most $2B|\up(\B)|
+ |\bdy(\B)|$. Finally, the continued presence of the rascal ensures
that the swollen future of $\B$ grows in each interval of time from
$\time(S_0)$ to $\iota(\B)$.  Thus it follows that the
length of $\B$ is also bounded by this number. So summing over all
$\B$ of type (iv) we have: 
$$ 
\sum |\B| \,\le\, \sum \Big(2B|\up(\B)| + |\bdy(\B)|\Big) \le (2B +1)\n , 
$$ 
as required. 
\end{proof} 
 Summarising the results of this section we have 
 
\bl \label{BonusBound} Summing over all teams that are not short, we
have 
\[      \sum_{\T}|\bonus(\T)| \leq \Big( \Bb \Big) \n .  \]
\end{lemma}

\section{The Proof of Theorem C} \label{summary} 
 
Pulling all of the previous results together, define  
$$
K_1 = \AFourC,
$$
and
$$ 
K = \K. 
$$ 
 
\begin{theorem} \label{S0<=Kn}
$|S_0| \leq K\n$. 
 
\Prf: 
The corridor $S_0$ can be subdivided into distinct colours which form connected regions. 
  Each colour $\mu$ can be partitioned into connected (possibly empty) regions $A_1(S_0,\mu),  
  A_2(S_0,\mu), A_3(S_0,\mu), A_4(S_0,\mu)$ and $A_5(S_0,\mu)$.  By Lemma \ref{A1A5Lemma}, 
  Proposition \ref{SummaryLemma}, Lemma \ref{A3Lemma}, Proposition \ref{A2Prop} and Lemma \ref{A1A5Lemma}, respectively, 
\begin{eqnarray*} 
\sum_{\mu \in S_0}|A_1(S_0,\mu)|  & \leq & C_0\n ,\\ 
\sum_{\mu \in S_0}|A_2(S_0,\mu)|  & \leq &  K_1  \n ,\\ 
\sum_{\mu \in S_0}|A_3(S_0,\mu)|  & \leq & (2B+1)\n ,\\ 
\sum_{\mu \in S_0}|A_4(S_0,\mu)|  & \leq &  K_1 \n,\mbox{ and} \\ 
\sum_{\mu \in S_0}|A_5(S_0,\mu)|  & \leq & C_0\n . 
\end{eqnarray*} 
Summing completes the proof of Theorem \ref{S0<=Kn}.
\et
Since there are at most $\frac{\n}{2}$ corridors in $\Delta$,  
\[      \mbox{Area}(\Delta) \leq \frac{K}{2}\n^2 ,       \] 
which proves the Main Theorem for positive automorphisms, i.e. Theorem \ref{MainThmPos}.

\section{Glossary of Constants}

$B$ -- the Bounded Cancellation constant (Lemmas \ref{BCL} 
and \ref{SingularityProp}). 
 
$C_0$ -- maximum distance a left-fast (right-fast)  
letter can be from the left (right) edge of its colour if it is to be
cancelled from the
left (right) within the future of the corridor.  See Lemma \ref{A1A5Lemma}. 
 
$C_1$ -- an upper bound on the   
lengths of the subintervals 
 $C_{(\mu,\mu')}(1)$ of $A_4(S_0,\mu)$. By definition, $C_{(\mu,\mu')}(1)$ 
 is consumed by $\mu'(S_0)$; it begins at the right end of
$A_4(S_0,\mu)$ and ends   at the last non-constant letter. See Lemma
\ref{C1Lemma}. Note that one can take $C_1=2mB^2$.

$L$ -- the maximum of the lengths of the images  $\phi(a_i)$ of the basis 
elements $a_i$, i.e. the maximum length of $u_1, \ldots, u_m$ in
the presentation $\P$ (see equation \ref{presentation}).

$L_{inv}$ -- the maximum of the lengths of $\phi^{-1}(a_i)$. 
 
$T_0$ -- the constant from the 2-Colour Lemma (Lemma
\ref{TwoColourLemma}). For all positive words $U$ and $V$, if
$U$ neuters  $V^{-1}$ then it does so in at most $T_0$ steps. 

$\hat{T_1}$ -- the constant from the Unnested Pincer Lemma, Theorem \ref{prePincerLemma}.

$T_1'$ -- the constant from Definition \ref{T1'Lemma}.  Recall that we stipulate that $T_1' \ge \hat{T_1}$.

$T_1 := T_1' + 2T_0$ -- $T_1$ is the constant from the Pincer Lemma, Theorem \ref{PincerLemma}. 

$C_4 := LL_{inv}$

$\ll := {\rm max} \{ 2B(T_0 + 1)+1, LC_4 \}$

Finally, $K_1$ is defined to be
\[	\AFourC	,	\]
and $K = 2C_0 + 2K_1 + 2B + 1$.

\part{Train Tracks and the Beaded Decomposition} \label{Part:BG2}

\renewcommand{\thesection}{\ref{Part:BG2}.\arabic{section}}

\setcounter{section}{0}

Part \ref{Part:BG2} of this work is dedicated to the construction and
analysis of a refined \index{topological representative}
topological representative for a suitable iterate  of  an  arbitrary automorphism of a finitely
generated free group. In Part \ref{Part:BG3} we shall use these representatives
to extend the results obtained in Part \ref{Part:BG1} to the
general setting.
Our results rely in a fundamental way on
the theory of {\em improved relative train tracks} developed by Bestvina,
Feighn and Handel in \cite{BFH}.

The  properties of the
topological representative $f:G\to G$ constructed in  \cite{BFH}  allow one
to control the manner in which a path $\sigma$ evolves as
one looks at its iterated images under $f$, and one might naively suppose that this
is the key issue that one must overcome in translating the proof of our Main Theorem from the
positive case (Part \ref{Part:BG1}) to the general case (Part \ref{Part:BG3}). However, upon
closer inspection one discovers this is actually only a fraction  of the story because
when a corridor evolves in the time flow on a van Kampen diagram,
the interaction of the forward iterates of the individual edges is such that the
basic {\em splitting} of paths established in \cite{BFH} may get broken. It is to overcome
this difficulty that we need the notion of \index{hard splitting}{\em hard splitting}.

\begin{unnumb-def} [See Definition \ref{def:HardSplit}]
We say that a decomposition of an edge-path into sub edge-paths 
$\rho = \rho_1 \rho_2$ is a {\em
  hard} $k$-splitting if for {\em any} choice of tightening of $f^k(\rho) =
f^k(\rho_1) f^k(\rho_2)$ there is no cancellation between the image of
$\rho_1$ and the image of $\rho_2$.

A decomposition that is a hard $k$-splitting for all $k \ge 1$ is
called a {\em hard} splitting.  If $\rho_1 \cdot \rho_2$ is a hard
splitting, we write $\rho_1 \odot \rho_2$.
\end{unnumb-def}

In the analysis of van Kampen diagrams that forms the core of the proof of
the Main Theorem, 
the class of ``broken" paths that one must
understand are the residues of the images of a single edge that survive
repeated cancellation during the corridor flow.
 In the language of the topological representative
$f:G\to G$, this amounts to understanding {\em monochromatic
paths}, as defined below. 
Every edge-path $\rho$ in $G$ admits a unique maximal splitting 
into edge-paths (Lemma \ref{l:Hsplit}); 
our main task here in Part  \ref{Part:BG2} is to understand the nature of 
the factors in this splitting and the behaviour
of certain larger units into which they naturally accrete when
$\rho$ is monochromatic.  

To this end, we identify a small number of basic 
units into which the iterated images of 
monochromatic paths split; the key feature of this splitting is that it 
is robust enough   to withstand the
difficulties caused by cancellation in van Kampen diagrams. The basic units
are defined so as to ensure that they enjoy those features of
individual edges that proved important in the positive case (see Part \ref{Part:BG1}). 
We call the units {\em beads}. The vocabulary of beads is as follows.

Let $f:G\to G$ be a topological representative and let $f_\#(\sigma)$
denote the tightening rel endpoints of the image of an edge-path $\sigma$.
Following \cite{BH2}, if $f_\#(\tau)=\tau$  we call $\tau$ a
{\em Nielsen} path.
A path $\rho$ in $G$ is called a \index{growing exceptional path (\gep)}
{\em growing exceptional path} (\gep) if
either $\rho$ or $\bar{\rho}$ is of
the form $E_i\bar{\tau}^k \bar{E_j}$ where $\tau$ is a Nielsen
path, $k \geq 1$, $E_i$ and $E_j$ are parabolic edges,  $f(E_i) = E_i \odot \tau^m$, 
$f(E_j) = E_j \odot\tau^n$, and $n > m > 0$. If it is $\rho$ (resp.~$\overline\rho$) 
that is of this form, then proper initial (resp.~terminal)
sub edge-paths of $\rho$ are called \pep s
({\em pseudo-exceptional paths}). 
\index{pseudo-exceptional path (\pep)}. 

\gep s and \pep s are key objects of study for us in Parts \ref{Part:BG2} and \ref{Part:BG2}. 
 They admit no nontrivial hard splitting, but there is no global
bound on their length.  Therefore, they must be included as basic units in 
the Beaded Decomposition Theorem below.  Also, there is no uniform bound
on the number of iterates required to cancel a \gep \ or \pep \ when it occurs
as a sub-path of the label on a corridor.  This leads to considerable technical
difficulties in Part \ref{Part:BG3}.

Let $f: G \to G$ be an improved relative train track map and $d, J \ge 1$ integers.
Then \index{path!monochromatic}
{\em $d$-monochromatic} paths in $G$ are defined by a simple  recursion:
edges in $G$ are $d$-monochromatic and if $\rho$ is  a
$d$-monochromatic path then every sub edge-path of
$f^d_\#(\rho)$ is $d$-monochromatic.\footnote{See Subsection \ref{ss:Iterate} for a precise definition of 
the map $f_{\#}^d$.} \index{atom}
A {\em $(J,f)$-atom} is a 
$d$-monochromatic edge-path of length at most $J$ that admits no non-vacuous
hard splitting into edge-paths.
 
An edge-path $\rho$ is \index{bead}
{\em $(J,f)$-beaded} if it admits a hard splitting $\rho = \rho_1 \odot \cdots \odot \rho_k$ where each $\rho_i$ is a \gep, a \pep, a $(J,f)$-atom, or
an indivisible Nielsen path of length at most $J$ (where \gep s, \pep s and Nielsen paths are defined with respect to the map $f$).

The following is the most important output of Part \ref{Part:BG2}. 

\index{Beaded Decomposition Theorem}
\renewcommand{\thethmspec}{{\bf{Beaded Decomposition Theorem. \kern-.3em}}}

\medskip
\noindent
\begin{thethmspec}\label{Beaded Decomposition Theorem}
{\em For every $\phi\in{\rm{Out}}(F_n)$,
there exist positive integers $k, d$ and $J$ such that $\phi^{k}$ has an improved
relative train-track representative $f:G\to G$ with the property that every $d$-monochromatic path in $G$ is $(J,f)$-beaded.}
\end{thethmspec}
\medskip

In fact, we do not prove the Beaded Decomposition Theorem {\em{per se}}. Instead, we prove a
more general  statement about futures of arbitrary paths under repeated iteration
and cancellation (Theorem \ref{SuperBDT}). We also need the following:

\begin{addendum} \label{add:SameJ}
If $f$ is replaced by an iterate $f_1 = f_{\#}^l$, then the Beaded Decomposition
Theorem is true for $f_1$ with the same constant as for $f$.
\end{addendum}

This sharpening of the Beaded Decomposition Theorem will prove vital in Part \ref{Part:BG3}:
often, we will need  to replace $f$ by an iterate, but the iterate we choose will depend on $J$, so
 Addendum \ref{add:SameJ} is needed to avoid circularity.  Related to this point, there are
a number of complications concerning how one should interpret beads; these are addressed
in Section \ref{IterateSection}.

As is clear from the preceding   discussion, our main motivation for
developing the Beaded  Decomposition is its application in Part \ref{Part:BG3}. 
The import of Part \ref{Part:BG2} in Part \ref{Part:BG3} has been deliberately distilled
into  this single statement  and Addendum \ref{add:SameJ}
so that a reader who is willing to accept these
as articles of faith can proceed directly from Part \ref{Part:BG1} to Part \ref{Part:BG3}.

We expect that our particular refinement of the train-track technology may prove
useful in other contexts.  This expectation stems from
the general point that the development of refined topological representatives
leads to insights into purely algebraic questions about free-group
automorphisms.  See \cite{BG-Growth} for a concrete illustration of this.\footnote{\cite{BG-Growth} contains results about the growth of words under iterated
automorphisms.  A previous version of Part \ref{Part:BG2} of this book contained an incorrect version of these results.  We thank Gilbert Levitt for bringing this error to our attention.}

\section{Improved Relative Train Track Maps} \label{TrainTracks}

In this section we collect and refine those elements of the train-track
technology that we shall need. Most of the material here is drawn
directly from  \cite{BH2} and \cite{BFH}. 

The philosophy behind
train tracks is to find an {\em efficient} topological representative for
an outer automorphism of $F$.  Precisely what it means for a graph map
to be {\em efficient} is spelled out in this section.

\subsection{Edge-paths and tightening}
Let $G$ be a graph. Following \cite{BFH}, we try to
reserve the term {\em path} for a
map $\sigma : [0,1] \to G$ that is either constant or an immersion (i.e.
{\em tight}).  The reverse path $t\mapsto \sigma(1-t)$ will be denoted $\overline \sigma$.
We
conflate the map $\sigma$ with its monotone reparameterisations (and
even its image, when this does not cause confusion).  Given an arbitrary continuous map
$\rho : [0,1] \to G$, we denote by $[\rho]$ the unique (tight)
{\em path}  homotopic
rel  endpoints to $\rho$. In keeping with the notation of the previous
section, given $f:G\to G$ and a path $\sigma$  in $G$,  we write $f_{\#}(\sigma)$
to denote $[f(\sigma)]$. We are primarily concerned
with {\em edge-paths},  i.e. those paths $\sigma$ for which $\sigma(0)$
and $\sigma(1)$ are vertices.

We consider only maps $f:G\to G$ that send vertices to vertices and  
edges to edge-paths (not necessarily to single edges). If there is 
an isomorphism  $F \cong \pi_1 G$ such that $f$ induces $\mathcal O
\in \text{\rm{Out}}(F)$, then one says that $f$ {\em represents} $\mathcal O$.
\index{topological representative}

\subsection{Replacing $f$ by an Iterate} \label{ss:Iterate}
In order to obtain good topological representatives of outer
automorphisms, one has to replace the given map by a large
iterate. It is important to be clear what one means by {\em iterate}
in this context, since   we wish to consider only
topological  representatives  whose restriction to
each edge is an immersion and this property is not inherited by
(naive) powers of the map.

Thus we deem the phrase\footnote{and obvious
variations on it}
{\em replacing $f$ by an iterate}, to mean that for fixed
$k\in\mathbb N$, we pass from consideration of $f : G \to G$
to consideration of  the map $f_{\#}^k : G \to G$ that sends
each edge $E$ in $G$ to  the tight edge-path $f_{\#}^k(E)$ that is 
homotopic rel endpoints to $f^k(E)$.

\subsection{(Improved) Relative train tracks}

We now describe the properties of Improved Relative Train Track maps,
as constructed in \cite{BH2} and \cite{BFH}.

\bigskip
 \noindent{\bf{Splittings, Turns and Strata.}}
Suppose that $\sigma = \sigma_1\sigma_2$ is a decomposition of a path
 into nontrivial subpaths (we do not assume that
$\sigma_1$ and $\sigma_2$ are edge-paths, even if $\sigma$ is).  We
say that $\sigma = \sigma_1\sigma_2$ is a \index{splitting}
{\em $k$-splitting} if
\[      f_{\#}^k(\sigma) = f_{\#}^k(\sigma_1) f_{\#}^k(\sigma_2)
\]
is a decomposition into sub-paths (i.e. for {\em some} choice of tightening,
there is no folding between the $f^k$-images of $\sigma_1$ and
$\sigma_2$ when $f^k(\sigma_1\sigma_2)$ is tightened). If $\sigma = \sigma_1\sigma_2$
is a $k$-splitting for all $k > 0$ then it is called a {\em splitting\footnote{In the next
  section, we introduce a stronger notion of {\em hard} splittings.}}
and we write $\sigma = \sigma_1 \cdot \sigma_2$.
If one of $\sigma_1$
or $\sigma_2$ is the empty path, the splitting is said to be {\em
  vacuous}.

A {\em turn} in $G$ is an unordered pair of half-edges originating at
a common vertex.  A turn is {\em non-degenerate} if it is defined by
distinct  half-edges, and is {\em degenerate} otherwise.  The
map $f : G \to G$ induces a self-map $Df$ on the set of oriented edges
of $G$ by sending an oriented edge to the first oriented edge in its
$f$-image.   $Df$ induces a map $Tf$ on the set of turns in $G$.

A turn is {\em illegal} with respect to $f : G \to G$ if its image
under some iterate of $Tf$ is degenerate; a turn is {\em legal} if it
is not illegal.

Associated to  $f$ is a {\em filtration} of $G$, 
\[      \emptyset = G_0 \subset G_1 \subset \cdots \subset G_\omega= G,   \]
consisting of $f$-invariant subgraphs of $G$.  We call the sets $H_r
:= \overline{ G_r \ssm G_{r-1}}$ \index{strata}
{\em strata}. To each stratum $H_r$
is associated $M_r$, the {\em transition matrix for $H_r$}; the
$(i,j)^{\rm th}$ entry of $M_r$ is the number of times the $f$-image
of the $j^{\rm th}$ edge crosses the $i^{\rm th}$ edge in either
direction. By choosing a filtration carefully one may ensure that for
each $r$ the matrix $M_r$ is either the zero matrix or is 
irreducible. If
$M_r$ is the zero matrix, then we say that $H_r$ is a \index{strata!zero}{\em zero
stratum}.  Otherwise, $M_r$ has an associated Perron-Frobenius
eigenvalue $\lambda_r \geq 1$, see \cite{Seneta}. If $\lambda_r > 1$
then we say that $H_r$ is an \index{strata!exponential}{\em exponential stratum}; if $\lambda_r
= 1$ then we say that $H_r$ is a 
\index{strata!parabolic}{\em parabolic stratum}\footnote{Bestvina {\em et al.} use the terminology {\em
exponentially-growing} and {\em non-exponentially-growing} for our
exponential and parabolic. This difference in terminology explains the
names of the items in Theorem \ref{MainTrainTrack} below.}. The edges
in strata inherit these adjectives, e.g. ``exponential edge".  

A turn is defined to be  {\em in $H_r$ if both} half-edges 
lie in the stratum $H_r$.  A turn is a {\em
mixed turn in $(G_r,G_{r-1})$} if one edge is in $H_r$ and the other is in
$G_{r-1}$. A path with no illegal turns in $H_r$ is said to be {\em
$r$-legal}.  We may emphasize that certain turns are in $H_r$ by calling them 
{\em $r$-(il)legal turns}.

\begin{definition} \cite[Section 5, p.38]{BH2} \label{TrainTrackDef}
We say that $f : G \to G$ is a \index{train track map!relative}
 {\em relative train track map} if the
following conditions hold for every exponential stratum $H_r$: 
\begin{enumerate}
\item[\rm (RTT-i)] $Df$ maps the set of oriented edges in $H_r$ to
  itself; in particular all mixed turns in $(G_r,G_{r-1})$ are legal. 
\item[\rm (RTT-ii)] If $\alpha$ is a nontrivial path in $G_{r-1}$
  with endpoints in $H_r \cap G_{r-1}$, then $f_{\#}(\alpha)$ is a
  nontrivial path with endpoints in $H_r \cap G_{r-1}$. 
\item[\rm (RTT-iii)] For each legal path $\beta$ in $H_r$,
  $f(\beta)$ is a path that does not contain any illegal turns in
  $H_r$. 
\end{enumerate}
\end{definition}

The following lemma is ``the most important consequence of being a
relative train track map'' \cite[p.530]{BFH}; it follows immediately
from Definition \ref{TrainTrackDef}. 

\begin{lemma} \cite[Lemma 5.8, p.39]{BH2} \label{RTT(i)Lemma}
Suppose that $f : G \to G$ is a relative train track map, that $H_r$
is an exponential stratum and that $\sigma = a_1b_1a_2 \dots b_l$ is
the decomposition of an $r$-legal path $\sigma$ into subpaths $a_j$ 
in $H_r$ and $b_j$ in $G_{r-1}$. (Allow for the possibility
that $a_1$ or $b_l$ is trivial, but assume the other subpaths are
nontrivial.) Then $f_{\#}(\sigma) = f(a_1)f_{\#}(b_1)f(a_2)\dots
f_{\#}(b_l)$ and is $r$-legal.
\end{lemma}

\begin{definition} \label{BasicPaths}\index{basic path}
Suppose that $f : G \to G$ is a topological representative, that the
parabolic stratum $H_i$ consists of a single edge
$E_i$ and that $f(E_i) = E_iu_i$ for some path $u_i$ in $G_{i-1}$.
We say that the paths of the form $E_i\gamma \bar{E_i}$, $E_i\gamma$
and $\gamma \bar{E_i}$, where $\gamma$ is in $G_{i-1}$, are {\em basic paths
of height $i$}.
\end{definition}

\begin{lemma} \cite[Lemma 4.1.4, p.555]{BFH}
\label{BasicPathSplitting}
Suppose that $f : G \to G$ and $E_i$ are as in Definition
\ref{BasicPaths}.  Suppose further that $\sigma$ is a
path or circuit in $G_i$ that intersects $H_i$ nontrivially and that the
endpoints of $\sigma$ are not contained in the interior of $E_i$.
Then $\sigma$ has a splitting each of whose pieces is either a  basic
path of height $i$ or is contained in $G_{i-1}$.
\end{lemma}

\begin{definition} \index{Nielsen path}
A {\em Nielsen path} is a nontrivial path $\sigma$ such that
$f^k_{\#}(\sigma) = \sigma$ for some $k \geq 1$.
\end{definition}
Nielsen paths are called {\em periodic Nielsen paths} in \cite{BFH}, but
Theorem \ref{MainTrainTrack} below allows us to choose an $f$ so  that
any periodic Nielsen path has period $1$ (which is to say that
$f_{\#}(\sigma) = \sigma$), and we shall assume that $f$ satisfies the
properties outlined in Theorem \ref{MainTrainTrack}.  Thus we can
assume that $k=1$ in the above definition.  A Nielsen path is called
{\em indivisible} if it cannot be split as a concatenation of two
non-trivial Nielsen paths.

\begin{definition} [cf. 5.1.3, p.~561 \cite{BFH}] \label{Exceptional}
Suppose that $H_i$ is a single edge $E_i$ and that $f(E_i) =
E_i\tau^l$ for some closed Nielsen path $\tau$ in  \index{path!exceptional}
$G_{i-1}$ and some $l > 0$.  The {\em exceptional paths of height $i$}
are those paths of the form $E_i\tau^k\bar{E}_j$ or
$E_i\bar{\tau}^k\bar{E}_j$ where $k \geq 0$, $j \leq i$, $H_j$ is a
single edge $E_j$ and $f(E_j) = E_j\tau^m$ for some $m > 0$.
\end{definition}

\begin{remark} In \cite{BFH} the authors mistakenly say that $\tau$ is an
indivisible Nielsen path, rather than a primitive Nielsen path (not a proper power).
We omit the modifier entirely.
\end{remark}

In Definition \ref{Exceptional}, the paths do not have a preferred
orientation.  Thus it is
important to note that the paths of the form $E_j \tau^k\bar{E_i}$ and
$E_j \bar{\tau}^k\bar{E_i}$ with $E_i,E_j$ and $\tau$ as above are
also exceptional paths of height $i$.

\subsection{The Theorem of Bestvina, Feighn and Handel}

A matrix is {\em aperiodic} if it has a power in which every entry is
positive.  The map $f$ is {\em eg-aperiodic} if every
exponential stratum has an aperiodic transition matrix.

Theorem 5.1.5 in \cite{BFH} is the main structural theorem for
improved relative train track maps.  We shall use it continually in
what follows, often without explicit mention.  We therefore record
those parts of it which we need. A map $f$ which satisfies the
statements of Theorem \ref{MainTrainTrack} is called an {\em improved
relative train track map}.

\begin{theorem} (cf. Theorem 5.1.5, p.562, \cite{BFH}) \label{Theorem5.1.5}
  \label{MainTrainTrack} 
For every outer automorphism ${\mathcal O} \in {\rm Out}(F)$ there is an
eg-aperiodic relative train track map $f : G \to G$ with filtration
$\emptyset = G_0 \subset G_1 \subset \dots \subset G_\omega = G$ such that
$f$ represents an iterate of ${\mathcal O}$, and $f$ has the following
properties.
\index{train track map!improved relative}
\begin{itemize}
\item Every periodic Nielsen path has period one.
\item For every vertex $v\in G$, $f(v)$ is a fixed point.  If $v$ is
  an endpoint of an edge in a parabolic stratum then $v$ is a fixed
  point.  If $v$ is the endpoint of an edge in an
  exponential stratum $H_i$ and if $v$ is also contained in
  a noncontractible component of $G_{i-1}$, then $v$ is a fixed
  point. 
\item $H_i$ is a zero stratum if and only if it is the union of the
  contractible components of $G_i$. 
\item If $H_i$ is a {\em zero stratum}, then
\begin{enumerate}
\item[\rm z-(i)] $H_{i+1}$ is an exponential stratum.
\item[\rm z-(ii)] $f|H_i$ is an immersion.
\end{enumerate}
\item If $H_i$ is a {\em parabolic stratum}, then
\begin{enumerate}
\item[\rm ne-(i)] $H_i$ is a single edge $E_i$.
\item[\rm ne-(ii)] $f(E_i)$ splits as $E_i \cdot u_i$ for some closed
  path $u_i$ in $G_{i-1}$ whose basepoint is fixed by $f$.
\item[\rm ne-(iii)] If $\sigma$ is a basic path of height
  $i$ that does not split as a concatenation of two basic paths of
  height $i$ or as a concatenation of a basic path of height $i$ with
  a path contained in $G_{i-1}$, then either: {\rm (i)} for some $k$,
  the path $f^k_{\#}(\sigma)$ splits  into pieces, one of which equals
  $E_i$ or $\bar{E_i}$; or {\rm (ii)} $u_i$ is a Nielsen path and, for
  some $k$, the path $f^k_{\#}(\sigma)$ is an exceptional path of
  height $i$. 
\end{enumerate}
\item If $H_i$ is an {\em exponential stratum} then
\begin{enumerate}
\item[\rm eg-(i)] There is at most one indivisible Nielsen path
  $\rho_i$ in $G_i$ that intersects $H_i$ nontrivially.  The
  initial edges of $\rho_i$ and $\bar{\rho_i}$ are distinct (possibly
  partial) edges in $H_i$.
\end{enumerate}
\end{itemize}
\end{theorem}

Suppose that $f : G \to G$ is an improved relative train track map
representing some iterate $\phi^k$ of $\phi \in \On$, and that $\rho$
is a Nielsen path in $G_r$ that intersects $H_r$ nontrivially,
and suppose that $\rho$ is not an edge-path.  Then subdividing the
edges containing the endpoints of $\rho$ at the endpoints, gives a new
graph $G'$, and the map $f' : G' \to G'$ induced by $f$ is an improved
relative train track map representing $\phi^k$.  To ease notation, it
is convenient to assume that this subdivision has been performed.
Under this assumption, all Nielsen paths will be edge-paths, and all
of the paths which we consider in the remainder of Part \ref{Part:BG2} will
also be edge-paths.

\begin{convention} \label{Conv:INPedge}
{\rm Since all Nielsen paths in the remainder of Part \ref{Part:BG2} will be edge-paths, we will use the phrase {\em `indivisible Nielsen path'} to mean a Nielsen edge-path which cannot be decomposed nontrivially as a concatenation of two non-trivial Nielsen {\em edge}-paths.  In particular, a single edge fixed pointwise by $f$ will be considered to be an indivisible Nielsen path.}
\end{convention} 

\medskip

{\em For the remainder of this article, we will concentrate on an improved relative train track map $f : G \to G$ and repeatedly pass to iterates $f_{\#}^k$ in order to better control its cancellation properties.}

\medskip

Recall the following from \cite[Section 4.2, pp.558-559]{BFH}.

\begin{definition}
If $f:G \to G$ is a relative train track map and $H_r$ is an
exponential stratum, then define $P_r$ to be the set of paths $\rho$
in $G_r$ that are such that: 
\begin{enumerate}
\item[(i)] For each $k \geq 1$ the path $f^k_{\#}(\rho)$ contains
  exactly one illegal turn in $H_r$. 
\item[(ii)] For each $k \geq 1$ the initial and terminal (possibly
  partial) edges of $f^k_{\#}(\rho)$ are contained in $H_r$.
\item[(iii)] The number of $H_r$-edges in $f^k_{\#}(\rho)$ is bounded
  independently of $k$.
\end{enumerate}
\end{definition}

\begin{lemma} \cite[Lemma 4.2.5, p.558]{BFH} \label{Pr}
$P_r$ is a finite $f_{\#}$-invariant set.
\end{lemma}

\begin{lemma} \cite[Lemma 4.2.6, p.559]{BFH} \label{ExpSplitting}
Suppose that $f : G \to G$ is a relative train track map, that $H_r$
is an exponential stratum, that $\sigma$ is a
path or circuit in $G_r$ and that, for each $k\geq 0$, the path
$f^k_{\#}(\sigma)$ has the same finite number of illegal turns in
$H_r$.  Then $\sigma$ can be split into subpaths that are either
$r$-legal or elements of $P_r$.
\end{lemma}

\begin{definition} \label{Weightr} \index{path!weight of}
If $\rho$ is a path and $r$ is the least integer such that $\rho$ is in 
$G_r$ then we say that {\em $\rho$ has weight $r$}.
\end{definition}

If $\rho$ has weight $r$ and $H_r$ is exponential, we will say that
$\rho$ is an \index{path!exponential}{\em exponential path}.  We define
\index{path!parabolic} {\em
parabolic paths} similarly.

\begin{lemma} \label{PreNielsen}
Suppose that $\sigma$ is an edge-path and that, for some $k \geq 1$,
$f_{\#}^k(\sigma)$ is a Nielsen path.  Then $f_{\#}(\sigma)$ is a
Nielsen path. 
\end{lemma}

\begin{proof} Suppose that the endpoints of $\sigma$ are $u_1$ and
$v_1$ and that the endpoints of $f_{\#}^k(\sigma)$ are $u_2$ and
$v_2$.  For each vertex $v \in G$, $f(v)$ is fixed by $f$, so
$f(u_1) = u_2$ and $f(v_1) = v_2$.  If $f_{\#}(\sigma) \neq
f_{\#}^k(\sigma)$ then we have two edge-paths with the same
endpoints which eventually get mapped to the same path.  Thus there is
some nontrivial circuit which is killed by $f$, contradicting the fact
that $f$ is a homotopy equivalence.  Therefore $f_{\#}(\sigma) =
f_{\#}^k(\sigma)$ and so is a Nielsen path.
\end{proof}

\medskip

Always, $L$ will denote the maximum of the lengths of the
paths $f(E)$, for $E$ an edge in $G$.

Later, we will pass to further iterates of $f$ in order to find a
particularly nice form.

An analysis of the results in this section allows us to see that there
are three kinds of indivisible Nielsen paths.  The first are those
which are single edges; the second are certain exceptional paths; and the
third lie in the set $P_r$.  We will use this trichotomy frequently
without mention.  The first two cases are where the path is
parabolic-weight, the third where it is exponential-weight.  It is not
possible for Nielsen path to have weight $r$ where $H_r$ is a zero
stratum. 

\begin{observation}
Let $\rho$ be an indivisible Nielsen path of exponential weight $r$.  Then the
first and last edges in $\rho$ are contained in $H_r$.
\end{observation}

Because periodic Nielsen paths have period $1$, the set of Nielsen
paths does not change when $f$ is replaced by a further iterate of
itself.  We will use this fact often.

\begin{lemma} \label{LinearEdges}
Suppose $E$ is an edge such that $|f_{\#}^j(E)|$ grows linearly with
$j$.  Then $f(E) = E \cdot \tau^k$, where $\tau$ is a Nielsen path
that is not a proper power.  The edge-path $\tau$ decomposes into
indivisible Nielsen paths (each of which is itself an edge-path, by Convention
\ref{Conv:INPedge}). 
\end{lemma}
\begin{proof}
The fact that $f(E) = E \cdot \tau^k$, where $\tau$ is a Nielsen path
follows from conditions ne-(ii) and ne-(iii) of Theorem
\ref{MainTrainTrack}.  \footnote{If Theorem \ref{MainTrainTrack},
ne-(iii) held with the Nielsen path $\tau$ in the definition of
exceptional paths being indivisible, we could also insist that
$\tau$ be indivisible here.}
\end{proof}

\begin{lemma} \label{NoTau}
Let $\tau$ be a Nielsen path and $\tau_0$ a proper initial (or
terminal) sub edge-path of $\tau$.  No image $f_{\#}^k(\tau_0)$ contains
$\tau$ as a sub edge-path.
\end{lemma}
\begin{proof} 
It is sufficient to prove the lemma for indivisible Nielsen paths, as
the result for arbitrary Nielsen paths then follows immediately.

If $\tau$ is an indivisible Nielsen path and $\tau_0$ is a proper
non-trivial subpath of $\tau$ then $\tau$ cannot be a single edge.
Therefore, either $\tau$ is either an indivisible Nielsen path of
exponential weight, or an exceptional path.

In case $\tau$ is an indivisible Nielsen path of exponential weight,
suppose the weight is $r$.  By Lemma \ref{ExpSplitting} $\tau$
contains a single illegal turn in $H_r$.  Suppose that $\tau_0$ does not contain
this illegal turn.  Then $\tau_0$ is $r$-legal, and so no iterate of
$\tau_0$ contains an illegal turn in $H_r$.  Therefore no iterate
of $\tau_0$ can contain $\tau$ as a subpath.  

Suppose then that $\tau_0$ {\em does} contain
the $r$-illegal turn in $\tau$.  Then, being a proper subpath of
$\tau$, the path on one side of the illegal turn in $\tau_0$ and its
(tightened) iterates is strictly smaller than the corresponding path in $\tau$.
Once again $\tau$ cannot be contained as a subpath of any iterate of
$\tau_0$. 

Finally, suppose $\tau$ is an exceptional path.  Then $\tau = E_i \rho^k
\bar{E_j}$ where $\rho$ is a Nielsen path and $E_i$ and $E_j$ are
of weight greater than $\rho$.  Any proper sub edge-path $\tau_0$
of $\tau$ contains at most one edge of weight greater than $\rho$. 
T he same is true for any iterate of $\tau_0$, and once again no
iterate of $\tau_0$ contains $\tau$ as a sub-path. 
\end{proof}

\section{Hard Splittings} \label{s:hard}

In this section we introduce a new concept for improved relative train tracks:
{\em hard splittings}.  This plays an important role in the subsequent sections
of Part \ref{Part:BG2}, and also in Part \ref{Part:BG3}.

Recall that a decomposition of a path $\sigma = \sigma_1 \sigma_2$ is
a $k$-splitting if $f_{\#}^k(\sigma) = f_{\#}^k(\sigma_1)
f_{\#}^k(\sigma_2)$; which means that, for {\em some} choice of
tightening, the images of $\sigma_1$ and $\sigma_2$ do not interact
with each other.  This leads to the concept of {\em splittings}.  We
need a more restrictive notion, where the decomposition is preserved
for {\em every} choice of tightening. For this purpose, we make the
following

\begin{definition}\label{def:HardSplit} [Hard splittings] \index{splitting!hard}
We say that a $k$-splitting $\rho = \rho_1 \rho_2$ is a {\em
  hard} $k$-splitting if for {\em any} choice of tightening of $f^k(\rho) =
f^k(\rho_1) f^k(\rho_2)$ there is no cancellation between the image of
$\rho_1$ and the image of $\rho_2$.

A decomposition which is a hard $k$-splitting for all $k \ge 1$ is
called a {\em hard} splitting.  If $\rho_1 \cdot \rho_2$ is a hard
splitting, we write $\rho_1 \odot \rho_2$.

An edge-path is {\em hard-indivisible} (or {\em h-indivisible}) if it admits no non-vacuous hard splitting into edge-paths.
\end{definition}

\begin{remark}
If one works in the universal cover, then $\tilde{\sigma}_1 \tilde{\sigma_2}$ is a $k$-hard
splitting if and only if, inside $\tilde{f}^k(\tilde\sigma_1 \tilde\sigma_2))$, the intersection
$\tilde{f}^k(\tilde\sigma_1) \cap \tilde{f}^k(\tilde\sigma_2)$ is a single point.
\end{remark}

\begin{remark} \label{rem:vacuoussplit}
In the above definition, we allow the possibility that one of the paths in 
the hard splitting is empty.  This is to allow various later statements to be
made more concisely.

For example, the phrase `$\rho$ admits a hard splitting immediately on either side of 
$\sigma$ of $\rho$' (for a path $\rho$ and a sub edge-path $\sigma$) allows the possibility that $\sigma$ is an initial or terminal sub-path of $\rho$.
\end{remark}

\begin{example}
Suppose that $G$ is the graph with a single vertex and edges $E_1,
E_2$ and $E_3$.  Suppose that $f(E_1) = E_1$, $f(E_2) = E_2 E_1$ and
$f(E_3) = E_3 \bar{E_1} \bar{E_2}$.  Then $f$ is an improved
relative train track.  And $E_3E_2 \cdot \bar{E_1}$ is a
$1$-splitting, since
\[      f(E_3E_2\bar{E_1}) = E_3 \bar{E_1}\bar{E_2} E_2
E_1 \bar{E_1},      \]
which tightens to $E_3 \bar{E_1} = f_{\#}(E_3E_2)
f_{\#}(\bar{E_1})$. In fact this is a splitting. However, there is
a choice of tightening which first cancels the final $E_1
\bar{E_1}$ and then the subpath $\bar{E_2} E_2$.  Therefore
the splitting $E_3 E_2 \cdot \bar{E_1}$ is not a hard
$1$-splitting.
\end{example}

The following lemma describes the main utility of hard splittings, and the
example above shows that it is not true in general for splittings. 

\begin{lemma} \label{HardSplitLemma}
Suppose that $\sigma_1 \odot \sigma_2$ is a hard splitting, that $\rho_1$ is
a terminal subsegment of $\sigma_1$, and that
$\rho_2$ is an initial subpath of $\sigma_2$.  Then $\rho_1 \odot
\rho_2$ is a hard splitting.

\end{lemma}
\begin{proof}
If there were any cancellation between
images of $\rho_1$ and $\rho_2$ then there would be a possible
tightening between the images of $\sigma_1$ and $\sigma_2$. 
\end{proof}

The following two lemmas will also be crucial for our applications of
hard splittings in Part \ref{Part:BG3}.

\begin{lemma}\label{l:Hsplit}
Every edge-path admits a unique maximal hard splitting into edge-paths.
\end{lemma}
\begin{proof} This follows by an obvious induction on length from the
observation that
if $\rho = \rho_1 \rho_2 \rho_3$,
where the $\rho_i$ are edge-paths,  and if
$\rho = \rho_1 \odot \rho_2 \rho_3$ and $\rho = \rho_1 \rho_2 \odot \rho_3$
then $\rho = \rho_1 \odot \rho_2 \odot \rho_3$. 
\end{proof}

\begin{lemma}
If $\rho = \rho_1 \odot \rho_2$ and $\sigma_1$ and $\sigma_2$
are, respectively, terminal and initial subpaths of $f^k_{\#}(\rho_1)$
and $f^k_{\#}(\rho_2)$ for some $k \ge 0$ then
$\sigma_1 \sigma_2 = \sigma_1 \odot \sigma_2$.
\end{lemma}
\begin{proof}
For all $i \ge 1$, the untightened path $f^i(\sigma_1)$
is a terminal subpath of the untightened path $f^i(f^{k}_\#(\rho_1))$,
while $f^i(\sigma_2)$ is an initial subpath of $f^i(f^{k}_\#(\rho_2))$.

The hardness of the splitting $\rho = \rho_1 \odot \rho_2$ ensures
that no matter how one tightens $f^{k+i}(\rho_1)f^{k+i}(\rho_2)$
there will be no cancellation between $f^{k+i}(\rho_1)$ and $f^{k+i}(\rho_2)$.
In particular, one is free to tighten to obtain $f^i(f^{k}_\#(\rho_1))f^i(f^{k}_\#(\rho_2))$
first, and then
tighten $f^i(\sigma_1)f^i(\sigma_2)$, and 
there can be no cancellation between them.
(It may happen that when one goes to tighten $f^{k+i}(\rho_1)$ completely,
the whole of $f^i(\sigma_1)$ is cancelled, but this does not affect the assertion
of the lemma.)
\end{proof}

The purpose of the remainder of this section is to sharpen results
from the previous section to cover hard splittings \footnote{Bestvina 
{\em et al.} make no explicit mention of the
  distinction between splittings and hard splittings, however
  condition (3) of Proposition 5.4.3 on p.581 (see Lemma
  \ref{SubpathSplitting} below) indicates that they are aware of the
  distinction and that the term `splitting' has the same
  meaning for them as it does here.}.

The following lemma is clear.
\begin{lemma} [cf. Lemma 4.1.1, p.554 \cite{BFH}] \label{HardSplitProperties}
If $\sigma = \sigma_1 \odot \sigma_2$ is a hard splitting, and
  $\sigma_1 = \sigma_1' \odot \sigma_2'$ is a hard splitting then
  $\sigma = \sigma_1' \odot \sigma_2' \odot \sigma_2$ is a hard
  splitting.  The analogous result with the roles of $\sigma_1$ and
  $\sigma_2$ reversed also holds.
\end{lemma}

\begin{remark} The possible existence of an edge-path $\sigma_2$ so that
$f_{\#}(\sigma_2)$ is a single vertex means that
 $\sigma_1 \sigma_2 = \sigma_1 \odot \sigma_2$ and $\sigma_2 \sigma_3 = \sigma_2 \odot \sigma_3$ need {\em not} imply
that  $\sigma_1 \sigma_2 \sigma_3 = \sigma_1 \odot \sigma_2 \odot \sigma_3$.

Indeed if $\sigma_2$ is an edge-path so that 
 $f_{\#}(\sigma_2)$ is a vertex then
$f_{\#}(\sigma_1)$ and $f_{\#}(\sigma_3)$ come together in a tightening of
$f(\sigma_1\sigma_2\sigma_3)$, possibly cancelling.

In contrast, if  $f_{\#}(\sigma_2)$ (and hence each $f_{\#}^k(\sigma_2)$)
contains an edge, then  the hardness of the two splittings ensures that
in any tightening $f_{\#}(\sigma_1\sigma_2\sigma_3)=f_{\#}(\sigma_1)f_{\#}(\sigma_2)
f_{\#}(\sigma_3) $, that is $\sigma_1 \sigma_2 \sigma_3 = \sigma_1 \odot \sigma_2 \odot \sigma_3$.
\end{remark}

The following strengthening of Theorem \ref{MainTrainTrack} ne-(ii) is
a restatement of (a weak form of) \cite[Proposition 5.4.3.(3),
  p.581]{BFH}.

\begin{lemma} \label{SubpathSplitting}
Suppose $f$ is an improved relative train track map and $E$ is a parabolic edge with $f(E) = E u$.  For any initial subpath $w$ of $u$, $E \cdot w$ is
a splitting.
\end{lemma}

\begin{corollary}
Suppose $f$ is an improved relative train track map, $E$ is a
parabolic edge and $f(E) = E u$.  Then $f(E) = E \odot u$.
\end{corollary}

The following lemma is straightforward to prove.
\begin{lemma} \label{BasicHSplit}
Suppose $H_i$ is a parabolic stratum and $\sigma$ is a
path in $G_i$ that intersects $H_i$ nontrivially, and that the endpoints of
$\sigma$ are not contained in the interior of $E_i$.  Then $\sigma$
admits a hard splitting, each of whose pieces is either a basic path
of height $i$ or is contained in $G_{i-1}$. 
\end{lemma}

\begin{lemma} \label{ne(iii)Hard}
If $\sigma$ is a basic path of height $i$ that does not
admit a hard splitting as a concatenation of two basic paths of height 
$i$ or as a concatenation of a basic path of height $i$ with a path of 
weight less than $i$, then either; (i) for some $k$, the path
$f_{\#}^k(\sigma)$ admits a hard splitting into pieces, one of which
is $E_i$ or $\bar{E_i}$; or (ii) $f(E_i) = E_i \odot u_i$, where
$u_i$ is a Nielsen path and, for some $k$, the path $f_{\#}^k(\sigma)$
is an exceptional path of height $i$.
\end{lemma}

\begin{proof}
Follows from the proof of \cite[Lemma 5.5.1, pp.585--590]{BFH}.
\end{proof}

\begin{lemma} [cf. Lemma \ref{ExpSplitting} above] \label{ExpHardSplitting}
Suppose that $f : G \to G$ is a relative train track map, that $H_r$
is an exponentially-growing stratum, that $\sigma$ is a
path or circuit in  $G_r$, and that each $f_{\#}^k(\sigma)$ has the same finite
number of illegal turns in $H_r$.  Then $\sigma$ can be decomposed as
$\sigma = \rho_1 \odot \ldots \odot \rho_k$, where
each $\rho_i$ is either (i) an element of $P_r$;
(ii) an $r$-legal path which starts and ends with edges in $H_r$; or
(iii) of weight at most $r-1$.
\end{lemma}
\begin{proof}
Consider the splitting of $\sigma$ given by Lemma \ref{ExpSplitting}.
The pieces of this splitting are either (i) elements of $P_r$, or (ii)
$r$-legal paths.  By Definition \ref{TrainTrackDef} RTT-(i), any
$r$-legal path admits a hard splitting into $r$-legal paths which
start and end with edges in $H_r$, and paths of weight at most $r-1$.
The turn at the end of a Nielsen path in the splitting of $\sigma$ is
either a mixed turn (with the edge from $H_r$ coming from the Nielsen
path and the other edge being of weight at most $r-1$) or a legal turn
in $H_r$.  In either case, $\sigma$ admits a hard splitting at the
vertex of this turn.
\end{proof}

The next result follows from a consideration of the form of indivisible Nielsen paths, noting Definition \ref{TrainTrackDef} and Lemma \ref{ExpHardSplitting}.

\begin{lemma} \label{NielsenHardSplit}
Any Nielsen path admits a hard splitting into indivisible Nielsen paths.
\end{lemma}

\begin{remark} \label{HardSplitRemark}
If $\rho = \rho_1 \odot \rho_2$ is a hard splitting for the map $f$
then it is a hard splitting for $f_{\#}^k$ for any $k \ge 1$.
\end{remark}

We record a piece of terminology which will be important in Part \ref{Part:BG3}.

\begin{definition} \label{Displayed} \index{path!displayed}
A sub edge-path $\rho$ of a path $\chi$ is {\em displayed} if there is a hard splitting of $\chi$ immediately on either side of $\rho$.
\end{definition}

 \section{A Small Reduction}
\label{ColourCancellation}

In this section we clarify a couple of issues about monochromatic paths,
and state Theorem \ref{MainProp}, which immediately implies the Beaded
Decomposition Theorem.

Our strategy for proving the Beaded Decomposition Theorem is as follows:
given an automorphism $\phi \in$\Aut, we start
with an improved relative train track representative $f :G \to G$ for some iterate $\phi^k$ of $\phi$, as obtained from the conclusion of Theorem \ref{MainTrainTrack}.  We analyse
the evolution of monochromatic paths, and eventually pass to an iterate of $f$
in which we can prove the Beaded Decomposition Theorem.  However, it is crucial
to note that monochromatic paths for $f$ are not necessarily monochromatic paths
for $f^k_{\#}$ when $k > 1$.  See Section \ref{IterateSection} for further discussion about some of these issues.  

These concerns lead to the following definition, where we are
concentrating on a fixed IRTT $f : G \to G$, and so omit mention of $f$ from our notation.

\begin{definition} \label{MonoChiDef} \index{path!monochromatic}
For a positive integer $d$, we define {\em $d$-monochromatic paths}  by  recursion:
edges in $G$ are $d$-monochromatic and if $\rho$ is a
$d$-monochromatic path then every sub edge-path of
$f_\#^{d}(\rho)$ is $d$-monochromatic.
\end{definition}

Note that if $d'$ is a multiple of $d$ then every $d'$-monochromatic path is  $d$-monochromatic
but not {\em vice versa}. Thus if we replace $f$ by an iterate then, for fixed $n$, 
the set of 
$n$-monochromatic paths may get  smaller. The content of the Beaded Decomposition Theorem is that
one need only pass to a bounded iterate in order to ensure that all monochromatic
paths admit a beaded decomposition.  In particular, the Beaded Decomposition
Theorem is an immediate consequence of the following theorem.

\index{bead}\index{Beaded Decomposition Theorem}
\begin{theorem} [Monochromatic paths are beaded] \label{MainProp}
Let $f:G\to G$ be an improved relative train track map.
There exist constants $d$ and $J$, depending only on $f$,
so that every $d$-monochromatic path in $G$ is $(J,f)$-beaded.
\end{theorem}

\begin{definition}[\NF s] \index{future!nibbled}
Let $\rho$ be a (tight) edge-path.  The {\em $0$-step \nf \ of $\rho$} is $\rho$.

For $k \ge 1$, a {\em $k$-step \nf \ of $\rho$} is a sub edge-path of $f_{\#}(\sigma)$, where $\sigma$ is a $(k-1)$-step \nf \ of $\rho$.  A {\em \nf \ of $\rho$} is a $k$-step \nf \ for some $k \ge 0$.

For $k \ge 0$, the {\em $k$-step entire future of $\rho$} is $f_{\#}^k(\rho)$.
\end{definition}

\begin{remark}  Nibbled futures are not assumed to be non-empty.  If
 a path is empty, any statement we claim about the existence of hard splittings
 should be interpreted to hold vacuously.   The $1$-monochromatic paths are precisely
 the \nf s of  single edges. 
\end{remark}
The notion of \nf s is central to Parts \ref{Part:BG2} and \ref{Part:BG3} of this book.  Usually, when proving things about monochromatic paths, we are actually proving things about the \nf s of paths of
bounded length.  In this spirit, rather than just proving Theorem \ref{MainProp},
we prove the following more general theorem about the iterated futures of arbitrary
paths.  We expect this theorem to have applications beyond those
presented  in this work.

\begin{theorem} \label{SuperBDT}
If $f : G \to G$ is an improved relative train track map, then there exists an
integer $d$ with the following property:
for each positive integer $n$, there exists $J>0$ so that for every edge-path $\rho$ with
$|\rho| \le n$ and  every positive integer $k$, 
{\em{every}} $kd$-step \nf \ of $\rho$ is $(J,f)$-beaded.
\end{theorem}

\begin{remark} \label{SuperBDTEnough}
It is clear that Theorem \ref{MainProp} follows immediately from Theorem \ref{SuperBDT}.  Therefore, in order to prove the Beaded Decomposition Theorem,
it suffices to prove Theorem \ref{SuperBDT}.
\end{remark}

\begin{remark} We posted a version of Part \ref{Part:BG2} of this book on the {\em{ArXiv}}
in July 2005. In December 2006, Feighn and Handel posted \cite{FH2}, in which
they develop a powerful refinement of the train track technology.
If one employs their {\em completely split} train track representatives, one
can prove the Beaded Decomposition Theorem with considerably greater ease than
we do here.
One can also streamline significant parts of the proof of Theorem \ref{SuperBDT}.  
However, we feel that the effort that this would save the reader
is offset by the extra machinery that they would be required to accept or absorb. On this
basis, we decided to retain our original proof.
\end{remark}

\section{Nibbled Futures}\label{s:nib}

\begin{notation}
Throughout this section and the rest of Part \ref{Part:BG2}, $f : G \to G$ is an improved 
relative train track map. 

Let $L$ be the maximum of the lengths of the paths $f(E)$ where $E$ ranges over
the edges of $G$.
\end{notation}

Monochromatic paths arise as  {\em nibbled futures} in the sense defined
below. Thus   in order to prove 
Theorem \ref{MainProp} we must understand how {\em nibbled futures} evolve.
The results in this section reduce this challenge to the task of understanding
the nibbled futures of GEPs.

\begin{theorem} [First Decomposition Theorem] \label{DecompTheorem}
For any $n \ge 1$ there exists an integer $V = V(n,f)$ such that if
$\rho$ is an edge-path of length at most $n$ then any \nf \ of $\rho$
admits a hard splitting into edge-paths, each of which is either the
\nf \ of a \gep \ or else has length at most $V$. 
\end{theorem}

The remainder of this section is dedicated to proving Theorem \ref{DecompTheorem}.
We begin by examining the entire future of a path of fixed length (Lemma 
\ref{SplittingLemma}) and then refine the argument to deal with nibbling. In
the proof of the first of these lemmas we require the following observation.
 
\begin{remark} \label{NoOfTurns}
Suppose that $\rho$ is a tight path of weight $r$.  Since $f$ is an improved relative
train track map, the number of $r$-illegal turns in $f_{\#}^l(\rho)$ is a non-increasing
function of $l$, bounded below by $0$.
\end{remark}

\begin{lemma} \label{SplittingLemma}
There is a function $D :  \N \to \N$, depending only on $f$, such that, for any $r \in \{ 1 , \ldots , \omega \}$, if $\rho$ is a path of weight $r$, and $|\rho| \le n$,
then for any $i \ge D(n)$ the edge-path $f_{\#}^{i}(\rho)$ admits a hard splitting into edge-paths, each of
which is either 
\begin{enumerate}
\item   a single edge of weight $r$;
\item an indivisible Nielsen path of weight $r$;
\item a \gep \ of weight $r$; or
\item a path of weight at most $r-1$.
\end{enumerate}
\end{lemma}
\begin{proof}  If $H_r$ is a zero stratum, then $f_{\#}(\rho)$ has
  weight at most $r-1$, and $D(n) = 1$ will suffice for any $n$. 

If $H_r$ is a parabolic stratum, then $\rho$ admits a hard splitting
into pieces which are either basic of height $r$ or of weight at most
$r-1$ (Lemma \ref{BasicHSplit}).  Thus it is sufficient to consider the case where $\rho$ is a
basic path of weight $r$ and $|\rho| \le n$.  By at most $2$ applications
of Lemma \ref{ne(iii)Hard}, we see that there exists a $k$ such
that $f_{\#}^k(\rho)$ admits a hard splitting into pieces which are
either (i) single edges of weight $r$, (ii) exceptional paths of
height $r$, or (iii) of weight at most $r-1$.  By taking the maximum
of such $k$ over all basic paths of height $r$ which are of length at
most $n$, we find an integer $k_0$ so that we have the desired hard splitting 
of $f_{\#}^{k_0}(\rho)$ for all basic paths of height $r$ of length
at most $n$.  Any of the exceptional paths in these splittings which
are not \gep s have bounded length and are either indivisible Nielsen
paths or are decreasing in length.  A crude bound on the length of the exceptional paths which are not \gep s is $L^{k_0}n$ where $L$ is the maximum length of $f(E)$ over all edges $E
\in G$.  Thus, those exceptional paths which are decreasing in length will
become \gep s within less than $L^{k_0}n$ iterations.  Therefore, 
replacing $k_0$ by $k_0 + L^{k_0}n$, we may assume all exceptional
paths in the hard splitting are \gep s. 

Finally, suppose that $H_r$ is an exponential stratum.  As noted in Remark \ref{NoOfTurns}, the number of $r$-illegal turns in $f_{\#}^l(\rho)$ is a non-increasing function
of $l$ bounded below by $0$.  Therefore, there is some $j$ so that the
number of $r$-illegal turns in $f_{\#}^{j'}(\rho)$ is the same for
all $j' \ge j$.  By Lemma \ref{ExpHardSplitting},
$f_{\#}^j(\rho)$ admits a hard splitting into pieces which are either
(i) elements of $P_r$, (ii) single edges in $H_r$, or (iii) paths of weight 
at most $r-1$.  To finish the proof of the lemma it remains to note that 
if $\sigma \in P_r$ then $f_{\#}(\sigma)$ is a Nielsen path by Lemma \ref{PreNielsen}. 

Therefore, the required constant for $H_r$ may be taken to be the maximum of $j+1$ over all the
paths of weight $r$ of length at most $n$. 

To find $D(n)$ we need merely take the maximum of the constants found
above over all of the strata $H_r$ of $G$. 
\end{proof}

In the extension of the above proof to cover nibbled futures, we shall need
the following straightforward adaptation of Lemma \ref{NoTau}.

\begin{lemma} \label{NibNoTau}
Let $\tau$ be a Nielsen path and $\tau_0$ a proper initial (or terminal) sub-path of $\tau$.  No nibbled future of $\tau_0$ contains $\tau$ as a sub-path.
\end{lemma}

\begin{proposition} \label{NibbleSplittingLemma}
There exists a function $D' :  \N \to \N$, depending only on $f$, so that for any $r \in \{ 1 , \ldots , \omega \}$, if $\rho$ is a path of weight $r$ and $|\rho| \le n$, then for any $i \ge D'(n)$ any $i$-step \nf \ of $\rho$ admits a hard splitting into edge-paths, each of which is either
\begin{enumerate}
\item   a single edge of weight $r$;
\item   a \nf \ of a weight $r$ indivisible Nielsen path;
\item\label{nfofgep}   a \nf \ of a weight $r$ \gep ; or
\item   a path of weight at most $r-1$.
\end{enumerate}
Moreover, in Case \eqref{nfofgep}, the \gep \ lies in the $j$-step \nf \ of 
$\rho$ for some $j \le i$.
\end{proposition}

\begin{remark} \label{stable}
Each of the conditions (1) -- (4) stated above is stable in the following sense: once an edge in a $k$-step \nf \ is contained in a path satisfying one of these conditions, then any future of this edge in any further \nf \ will also lie in such a path (possibly the future will go from case (1) to case (4), but otherwise which case it falls into is also stable).  Thus we can split the proof of Proposition \ref{NibbleSplittingLemma} into a number of cases, deal with the cases separately by finding some constant which suffices, and finally take a maximum to find $D'(n)$.  An entirely similar remark applies to a number of subsequent proofs, in particular Theorem \ref{ColourCancelMain}.
\end{remark}

\begin{remark} \label{D'Monotonic}
Since the statement of Proposition \ref{NibbleSplittingLemma} involves all paths
$\rho$ such that $|\rho| \le n$, if the function $D'(n)$ is chosen to be the smallest
function satisfying the conclusion then it is nondecreasing.  We will assume that
the function $D'$ we use is indeed monotonic.
\end{remark}

\begin{proof}[Proof (Proposition \ref{NibbleSplittingLemma}).]
Let $\rho_0 = \rho$ and for $j > 0$ let $\rho_j$ be a sub edge-path of $f_{\#}(\rho_{j-1})$.

If $H_r$ is a zero stratum, then $f_{\#}(\rho)$ has weight at most $r-1$ and it suffices to take $D'(n) = 1$.

Suppose that $H_r$ is an exponential stratum.  By Lemma \ref{SplittingLemma}, the $D(n)$-step entire future of $\rho$ admits a hard splitting of the desired form.  We consider how \nib \ can affect this splitting.  As we move forwards through the \nf \ of $\rho$, cancellation of $H_r$-edges can occur only at $r$-illegal turns and at the ends, where the \nib \ occurs. 

Remark \ref{NoOfTurns} implies that we can trace the $r$-illegal turns forwards through the successive \nf s of $\rho$ (whilst the $r$-illegal continues to exist).  We compare the $r$-illegal turns in $\rho_k$ to those in $f_{\#}^k(\rho)$, the entire future of $\rho$.  We say that the \nib \ {\em first cancels an $r$-illegal turn at time $k$} if the collection of $r$-illegal turns in $\rho_{k-1}$ is the same as the collection in $f_{\#}^{k-1}(\rho)$, but the collection in $\rho_k$ is {\em not} the same as that of $f_{\#}^k(\rho)$.  The first observation we make is that if, at time $k$, the \nib \ has not yet cancelled any $r$-illegal turn then the sequence of $H_r$-edges in $\rho_k$ is a subsequence of the $H_r$-edges in $f_{\#}^k(\rho)$.  Therefore, any splitting of the desired type for $f_{\#}^k(\rho)$ is inherited by $\rho_k$.

Since there is a splitting of the $D(n)$-step entire future of $\rho$ of the desired form, either there is a splitting of $\rho_{D(n)}$, or else $\rho_{D(n)}$ has fewer $r$-illegal turns than $f_{\#}^{D(n)}(\rho)$, and hence than $\rho$.  However, $|\rho_{D(n)}| \le n.L^{D(n)}$.  We apply the above argument to $\rho_{D(n)}$, going forwards a further $D(nL^{D(n)})$ steps into the future.  Since the number of illegal turns in $H_r$ in $\rho$ was at most $n-1$, we will eventually find a splitting of the required form within an amount of time bounded by a function of $n$ (this function depends only on $f$, as required).  Denoting this function by $D_0$, we have that any $D_0(n)$-step \nf \ of any path of exponential weight whose length is at most $n$ admits a hard splitting of the desired form.

Now suppose that $H_r$ is a parabolic stratum.  By Lemma \ref{BasicHSplit}, $\rho$ admits a hard splitting into basic edge-paths.  Therefore we may assume (by reversing the orientation of $\rho$ if necessary) that $\rho = E_r \sigma$ or $\rho = E_r \sigma \overline{E_r}$ where $E_r$ is the unique edge in $H_r$ and $\sigma$ is in $G_{r-1}$.  For the \nf \ of $\rho$ to have weight $r$, the \nib \ must occur only on one side (since the only edges of weight $r$ in any future of $\rho$ occur on the ends).  We assume that all \nib \ occurs from the right.  Once again, the $D(n)$-step entire future of $\rho$ admits a hard splitting of the desired form.  If $\rho = E_r \sigma \overline{E_r}$ then the $D(n)$-step \nf \ of $\rho$ either admits a hard splitting of the required form, or is of the form $E_r \sigma_1$, where $\sigma_1$ is in $G_{r-1}$.  Hence we may assume that $\rho = E_r \sigma$.  Suppose that $f(E_r) = E_r u_r$, and that $u_r$ has weight $s < r$.

Consider first the possibility that $\sigma$ has weight $q > s$ (but less
than $r$ by hypothesis).  We claim that after a bounded amount of time
the \nf \ of $\rho$ admits a splitting into one piece of the form $E_r\sigma'$ where
the weight of $\sigma'$ is strictly less than $q$, and other pieces which are all of the 
form required by the statement of the proposition.  Then, by induction on weight,
we may suppose that we have a splitting into one piece of the form $E_r \sigma''$
where the weight of $\sigma''$ is at most $s$ and all of the other pieces have
the form required by the proposition.  

So, suppose that $\sigma$ has weight $q > s$. There are three cases to consider.
If the weight of $\sigma$ is that of a zero stratum, then it immediately drops in weight
and the claim is proved.

Now suppose that $H_q$ is an exponential stratum.  The future of $E_r$ cannot cancel any edges of weight $q$ or higher in the future of $\sigma$, so the edges of weight $q$ in the \nf \ of $\rho$ are exactly the same as the edges of weight $q$ in the corresponding \nf \ of $\sigma$ (recall we are assuming that \nib \ only occurs from the right).  This $D_0(|\sigma|)$-step \nf \ of $\sigma$ admits a hard splitting into edge-paths which are either\footnote{\gep s have parabolic weight} single edges of weight $q$, the \nf \ of an indivisible Nielsen path of weight $q$, or of weight at most $q-1$.  Let $\sigma_2$ be the subpath of the 
$D_0(|\sigma|)$-step \nf \ of $\rho$ which starts at the right endpoint of $E_r$ up to but not including the first edge of weight $q$.\footnote{In case the \nf\ of $\sigma$ is empty, this is the entire path.}  Then, since mixed turns are legal, the $D_0(n)$-step \nf \ of $\rho$ admits a hard splitting into edge-paths, the leftmost of which is $E_r \sigma_2$.

Suppose now that $H_q$ is a parabolic stratum.  It is easy to see that $\rho$ admits a hard splitting into edge-paths, the leftmost of which is either $E_r \sigma_2$ or $E_r \sigma_2 \overline{E_q}$, where $\sigma_2$ has weight at most $q-1$.  Thus we may suppose that $\rho$ itself has this form.  Again, either the $D(n)$-step \nf \ of $\rho$ admits a hard splitting of the required form, or the $D(n)$-step nibbled future of $\rho$ has the form $E_r \sigma_3$, where $\sigma_3$ has weight at most $q-1$.  The arguments in the previous two paragraphs include the possibility that a \gep \ of weight $r$ occurs as a factor of the hard splitting of the $D(n)$-step \nf \ of $\rho$.  Thus we may assume that in some \nf \ of $\rho$ there will necessarily be a hard splitting on each side of the edge of weight $r$. (Recall by Remark \ref{rem:vacuoussplit} that this includes the case that this edge is an initial or terminal
subsegment.)

As noted above, by induction we have now proved that going forwards into
the \nf \ an amount of time bounded by a function of $n$, we may assume that
$\rho$ has the form $E_r \sigma_4$, where $\sigma_4$ has weight at most $s$
(thus $\sigma_4$ is the path $\sigma''$ from the claim above).  Suppose that
$\sigma_4$ has weight less than $s$.  Then $f_{\#}(E_r\sigma_4) = 
E_r \odot \sigma_5$, where $\sigma_5$ has weight less than $r$.  This is a splitting
of the required form which is inherited by an \nf .   Therefore, we are left with the
case that the weight of $\sigma_4$ is exactly $s$.

We now consider what kind of stratum $H_s$ is.  Suppose that $H_s$ is parabolic.  There are only two ways in which cancellation between weight $s$ edges in the \nf \ of $\rho$ can occur (see Lemma \ref{NoOldCanc}): they might be cancelled by edges whose immediate past is the edge of weight $r$ on the left end of the previous \nf ; alternatively, they can be nibbled from the right.  The $D(n)$-step entire future of $\rho$ admits a hard splitting as $E_r \odot \sigma_6$, where $\sigma_6$ has weight at most $r-1$.  There is no way that \nib \ can affect this splitting.

Finally, suppose that $H_s$ is an exponential stratum.  We follow a similar argument to the case when $H_r$ was an exponential stratum.  Either the $D(n)$-step \nf \ of $\rho$ admits a hard splitting of the desired kind (which means $\rho_{D(n)} = E_r \odot \sigma_7$ where $\sigma_7$ has weight at most $r-1$), or there are fewer $s$-illegal turns in the future of $\sigma_4$ in $\rho_{D(n)}$ than there are $s$-illegal turns in $\sigma_4$.  We then apply the same argument to the \nf \ of $\rho_{D(n)}$ until eventually we achieve a hard splitting of the required form.  

The last sentence in the statement of Proposition \ref{NibbleSplittingLemma} follows
immediately, since in the proof we have only consider paths which arise in the \nf s
of $\rho$.
This completes the proof of Proposition \ref{NibbleSplittingLemma}.
\end{proof}

We are now in a position to prove Theorem \ref{DecompTheorem}.  For this we require the following definition.

\begin{definition}\index{seed}
Suppose that $H_r$ is a stratum, and $E \in H_r$.
An {\em $r$-seed} is a non-empty subpath
$\rho$ of $f(E)$ which is maximal subject to lying in $G_{r-1}$.

If the stratum $H_r$ is not relevant, we just refer to {\em seeds}.
\end{definition}
Note that seeds are edge-paths and that the set of all seeds is finite.  Also,
if $H_r$ is an exponential stratum and $E \in H_r$ then the seeds in $f(E)$
are the sub-paths $b_i$ from Definition \ref{RTT(i)Lemma}.

The following is an immediate consequence of Lemma \ref{ExpHardSplitting} and RTT-(i) of Definition \ref{TrainTrackDef}.

\begin{lemma}
If $E \in H_r$ is an exponential edge and $\rho$ is an $r$-seed in
$f(E)$ then $f(E) = \sigma_1 \odot \rho \odot \sigma_2$ where
$\sigma_1$ and $\sigma_2$ are $r$-legal paths which start and finish
with edges in $H_r$.
\end{lemma}

\medskip

\begin{proof}[Proof (Theorem \ref{DecompTheorem}).]
Suppose that $\rho$ is a path of length $n$ and that $\rho_k$ is a
$k$-step \nf \ of $\rho$.  Denote by $\rho_0 = \rho, \rho_1, \dots ,
\rho_{k-1}$ the intermediate \nf s of $\rho$ used in order to define
$\rho_k$.  

We begin by constructing a 
\index{stack diagram}
van Kampen diagram\footnote{in fact, just a stack of corridors.  Of course, van Kampen diagrams are not required for this proof, but
we find them a convenient way of encoding choices of tightening and nibbling.} $\Delta_k$ which encodes
the $\rho_i$, proceeding by induction
on $k$.  For $k = 1$ the diagram $\Delta_1$ has a single (folded)
corridor with the bottom labelled by $\rho$ and the path $\rho_1$ a
subpath of the top of this corridor.  Suppose that we have associated
a van Kampen diagram $\Delta_{k-1}$ to $\rho_{k-1}$, with a unique
corridor at each time $t = 0 , \ldots , k-2$, such that $\rho_{k-1}$
is a subpath of the top of the latest (folded) corridor.  Then we
attach a new folded corridor to $\Delta_{k-1}$ whose bottom is
labelled by $\rho_{k-1}$.  The path $\rho_k$ is, by definition, a
subpath of the top of this new latest corridor. By convention, we consider $\rho_i$ to occur at time $i$.

Choose an arbitrary edge $\e$ in $\rho_k$ on the (folded) top of the
latest corridor in $\Delta_k$.  We will prove that there is a path
$\sigma$ containing $\e$ in $\rho_k$ so that $\rho_k$ admits a hard
splitting immediately on either side of $\sigma$ and so that $\sigma$
is either suitably {\em short} or a \nf \ of a \gep .  The purpose of
this proof is to find a suitable notion of {\em short}. 

Consider the embedded \index{family forest}
`family forest' $\mathcal F$ for $\Delta_k$, tracing
the histories of edges lying on the folded tops of corridors (see Remark \ref{tree}).  
Let $p$ be the path in $\mathcal F$ which follows the history of
$\e$.  We denote by $p(i)$ the edge  which intersects $p$ and lies on
the bottom of the corridor at time $i$. The edges $p(i)$ form
the {\em past} of $\e$.  We will sometimes denote the edge $\e$ by
$p(k)$.  It will be an analysis of the times at which the weight of
$p(i)$ decreases that forms the core of the proof of the theorem. 

The weights of the edges $p(0), p(1), \ldots , p(k)$ form a
non-increasing sequence.  Suppose this sequence is $W = \{ w_0 ,
\ldots , w_k \}$. A {\em drop} in $W$ is a time $t$ such that $w_{t-1}
> w_t$.  At such times, the edge $p(t)$ is contained in a (folded)
seed in the bottom of a corridor of $\Delta_k$. 

We will show that either successive drops occur rapidly, or else we reach 
a situation wherein each time a drop occurs we lose no essential information by restricting our attention to a small subpath of $\rho_{i}$.

To make this localisation argument precise, we define
 {\em incidents}, which fall into two types.  
 
 An {\em incident of Type A} is a time $t$ which (i) is a drop; and (ii) is such that
there is a hard splitting of $\rho_t$ immediately on either side of
the folded seed containing $p(t)$.  

An {\em incident of Type B} is a
time $t$ such that $p(t-1)$ lies in an indivisible Nielsen path with a
hard splitting of $\rho_{t-1}$ immediately on either side, but $p(t)$
does not;  except that we do not consider this to be an incident if some $\rho_i$, for $i \le t-1$ admits a hard splitting $\rho_i = \sigma_1 \odot \sigma_2 \odot \sigma_3$ with $p(i) \subseteq \sigma_2$ and $\sigma_2$ a \gep .   In case of an incident of Type B, necessarily $p(t)$
lies in the nibbled future of a Nielsen path on one end of $\rho_t$ with a hard splitting of
$\rho_t$ immediately on the other side. 

Define the time $t_1$ to be the last time at which there is an
incident (of Type A or Type B).  If there are no incidents, let $t_1 =
0$.  If this incident is of Type A, the edge $p(t_1)$ lies in a folded 
seed, call it $\pi$, and there is a hard splitting of $\rho_{t_1}$
immediately on either side of $\pi$.  If the incident is of Type B,
the edge $p(t_1)$ lies in the $1$-step \nf \ of a Nielsen path, call
this \nf \ $\pi$ also.  In case $t_1 = 0$, let $\pi = \rho$.  We will
see that there is a bound, $\alpha$ say, on the length of $\pi$ which
depends only on $f$ and $n$, and not on the choice of $\pi$, or the
choice of \nf .  The bound $\alpha$ will be defined solely in terms of Type B incidents.
We postpone the proof of the existence of the bound
$\alpha$ while we examine the consequences of its existence. 

The purpose of isolating the path $\pi$ is that it is a path of
controlled length and the hard splitting \footnote{this splitting is vacuous in case
$t_1 = 0$ and at various other points during this proof which we do not
explicitly mention} of $\rho_{t_1}$
immediately on either side of $\pi$ means that we need only consider the \nf \ or
$\pi$.  Suppose that $\pi$ has weight $r$. 

{\bf Claim 1:}
There exists a constant $\beta = \beta(n,\alpha,f)$ so that one of the following must occur:  
\begin{enumerate}
\item[(i)] for some $t_1 \le i < k$, the edge
$p(i)$ lies in  a \gep \ in $f_{\#}(\rho_{i-1})$
with a hard splitting immediately on either side; 
\item[(ii)] case (i) does not occur; $k - t_1 > \beta$;  and at some time $i \le t_1 + \beta$, the edge $p(i)$ 
lies in an indivisible Nielsen path $\tau$ in $f_{\#}(\rho_{i-1})$
with a hard splitting immediately on either side;  
\item[(iii)] $k - t_1 \le \beta$; or
\item[(iv)] there is a hard splitting of $\rho_k$ immediately on either side of $\e$. 
\end{enumerate}

\medskip

This claim implies the theorem, modulo the bound on $\alpha$, as we shall now explain.
In case (i), for all $j \ge i$, the edge $p(j)$ lies in the \nf \ of a \gep , so in particular this is true for $\e = p(k)$.  If case
(ii) arises then the definition of $t_1$ implies that for $j \ge i$, the edge $p(j)$ always
lies in a path labelled $\tau$ with a hard splitting immediately on
either side, for otherwise there would be a subsequent incident.  Also, the length of this Nielsen path is at most $\alpha L^\beta$.  If case (iii) arises, then the nibbled future of
$\pi$ at time $k$ has length at most $\alpha L^{\beta}$.

To prove the claim, we define two sequences of numbers
$V_{\omega}, V_{\omega - 1}, \ldots , V_1$ and $V'_{\omega},
V'_{\omega - 1}, \ldots , V'_1$, depending on $n$ and $f$, as follows (where $D'(n)$ is the function
  from Proposition \ref{NibbleSplittingLemma}): 
\begin{eqnarray*}
V_{\omega} &  := & D'(\alpha),\\ 
V'_{\omega} & := & V_{\omega} + \alpha L^{V_{\omega}}.
\end{eqnarray*}
For $\omega > i \ge 1$, supposing $V'_{i+1}$ to be defined, 
\[      V_i   :=   V'_{i+1} + D'(\alpha L^{V'_{i+1}}).  \] 
Also, supposing  $V_i$ to be defined, we define
\[      V'_i   :=  V_i +  \alpha L^{V_i} .      \]
The constants $V_i'$ and $V_i$ are defined so that Proposition \ref{NibbleSplittingLemma} may be applied successively to paths
which satisfy Case (4) of the statement of that result.  The key point
is that at time $t_1$ we have control over the length of the path $\pi$,
and so may apply Proposition \ref{NibbleSplittingLemma} to find
a hard splitting of any $i$-step \nf\ of $\rho$ so long as 
$i \ge D'(|\pi|)$.  If we consider the $D'(|\pi|)$-step \nf , then we also
have control of each of the $h$-indivisible paths in this hard splitting,
and thus we may apply Proposition \ref{NibbleSplittingLemma} again.
Note that the paths which satisfy Cases (1)--(3) of Proposition \ref{NibbleSplittingLemma}
also satisfy the requirements of Theorem \ref{DecompTheorem} (given the
as yet unproved bound $\alpha$), so we have to deal with the paths satisfying
Case (4).  We deal with these by successive applications of Proposition
\ref{NibbleSplittingLemma}, considering at each weight the paths satisfying
Case (4).  The constants $V_i$ and $V_i'$ are tuned to allow this induction
on weight to occur.  Note that since $D'$ is a nondecreasing function
(Remark \ref{D'Monotonic}), we have $V_{i+1} < V_{i+1}' < V_i$ for each $i$.

Consider the situation at time $t_1 + V_r$ (recall that $r$ is the
weight of $\pi$).  Possibly $k \le t_1 + V_r$, which is covered by
case (iii) of our claim, so long as $\beta > V_r$.  
Therefore, suppose that $k > t_1 + V_r$. 

According to Proposition \ref{NibbleSplittingLemma}, and the definition of
$t_1$,  at time $t_1 + V_r$ the $V_r$-step \nf \ of $\pi$ which exists
in $\rho_{t_1 + V_r}$ admits a hard splitting into edge-paths, each of
which is either: 
\begin{enumerate}
\item a single edge of weight $r$;
\item   a \nf \ of a weight $r$ indivisible Nielsen path;
\item   a \nf \ of a weight $r$ \gep ; or
\item   a path of weight at most $r-1$.
\end{enumerate}
We need to augment possibility (3) by recalling that Proposition \ref{NibbleSplittingLemma} also 
shows that the \gep \ referred to lies in the $j$-step \nf \ of $\pi$ for some $j \le V_r$.

We analyse what happens when the edge $p(t_1+V_r)$ lies in each of these four types of path.

{\bf Case (1):}
In the first case, by the definition of $t_1$, there will be
a hard splitting of $\rho_k$ immediately on either side of $\e$, since
in this case if there is a drop in $W$ after $t_1 + V_r$ then there is an
incident of Type A, contrary to hypothesis. 

{\bf Case (3):}
If $p(t_1+V_r)$ lies in a
path of the third type then we are in case (i) of our claim, and hence
content.  

The fourth type of path will lead us to an inductive
argument on the weight of the path under consideration.  But first we
consider the \nf s of Nielsen paths. 

{\bf Case (2):}
Suppose that in $\rho_{t_1 + V_r}$ the edge $p(t_1 + V_r)$ lies in the
\nf \ of a Nielsen path of weight $r$, with a hard splitting of
$\rho_{t_1 + V_r}$ immediately on either side.  Suppose that this \nf
\ is $\pi_r$.  If $\pi_r$ is actually a Nielsen path then we lie in
case (ii) of our claim.  Thus suppose that $\pi_r$ is not a Nielsen
path.  It has length at most $\alpha L^{V_r}$, and within  time
$\alpha L^{V_r}$ any \nf \ of $\pi_r$ admits a hard splitting into edge-paths of
types (1), (3) and (4) from the above list.  The required bound on length is 
straightforward, since the length of $\pi$ is at most $\alpha$ and we are
considering a sub-path of a $V_r$-step \nf \ of $\pi$ (recall that $L$ is the maximum
length of paths $L(E)$ for edges $E$ in $G$).

To see that any \nf \ of $\pi_r$ admits a splitting of the required form within time
$\alpha L^{V_r}$, consider the
three types of indivisible Nielsen paths.  If $\tau$ is a Nielsen path
which is a single edge fixed pointwise by $f$, then any \nf \ of
$\tau$ is either a single edge or empty.  

Suppose that $\tau$ is an
indivisible Nielsen path of weight $r$ and $H_r$ is exponential, and
suppose that $\tau'$ is a proper subpath of $\tau$.  Then there is
some iterated image $f_{\#}^l(\tau')$ of $\tau'$ which is $r$-legal. 
 By Proposition
\ref{NibbleSplittingLemma} any $D'(\alpha)$-step \nf \ of $\tau'$ is $r$-legal.
Since $\tau$ has length at most $\alpha$, so does $\tau'$.  Therefore, if
$i \ge D'(\alpha)$ then any $i$-step \nf \ of $\tau'$ admits a hard splitting into paths
of the required form.  Since $V_r > V_\omega = D'(\alpha)$, it is clear that
within time $L^{V_r} > D'(\alpha)$, the \nf s of $\tau'$ admit a hard splitting
of the required form.

Finally suppose that $E_i \tau^k \overline{E_j}$ is an indivisible
Nielsen path of parabolic weight, with $k \ge 0$.  Thus $\tau$ is a Nielsen path of
weight less than $r$, and $E_i$, $E_j$ are edges such that $f(E_i) =
E_i \odot \tau^m$, $f(E_j) = E_j \odot \tau^m$.  A $1$-step \nf \ of
$E_i \tau^k \overline{E_j}$ has one of three forms:  (I) $E_i
\tau^{k_1} \tau'$, where $\tau'$ is a proper sub edge-path of $\tau$;
(II) $\tau' \tau^{k_2} \tau''$ where $\tau'$ and $\tau''$ are proper
sub edge-paths of $\tau$; or (III) $\tau' \tau^{k_3} \overline{E_j}$,
where $\tau'$ is a proper sub edge-path of $\tau$. Note that cases (I) and (III)
are not symmetric because we assume that $k \ge 0$ (and hence $k_1, k_2, k_3 \ge 0$
also).

\smallskip

{\em Case 2(I):}
In this case, $E_i \tau^{k_1} \tau'$ admits a hard splitting into $E_i$ and $\tau^{k_1}
\tau'$, which is of the required sort.  

\smallskip

{\em Case 2(II):}
In this case the path already
had weight less than $r$.  

\smallskip

{\em Case 2(III):}
Suppose we are in case (III), and that $\mu$, the $\alpha L^{V_r}$-step \nf \ of $\tau'\tau^{k_3} \overline{E_j}$  has a copy of $\overline{E_j}$. Lemma \ref{NibNoTau} assures us that no \nf \ of $\tau'$ can contain $\tau$ as a subpath, and therefore there is a splitting of $\mu$ immediately on the right of $\overline{E_j}$, and we are done.  If there is no copy of $\overline{E_j}$ in $\mu$, we are also done, since this \nf \ must have weight less than $r$.

{\bf Case (4):}
Having dealt with cases (1), (2) and (3), we may now suppose that at time $t_1 + V_r + \alpha L^{V_r} = t_1 +
V'_r$ the edge $p(t_1 + V'_r)$ lies in an edge-path of weight at most
$r-1$ with a hard splitting of $\rho_{t_1 + V'_r}$ immediately on
either side.\footnote{Note that again it is possible that $k < t_1 +
  V'_r$, in which case we are in case (iii) of our claim.  We suppose
  therefore that this is not the case.}  Denote this path by $\pi'_r$, chosen to be in the future of $\pi$.
Note that $\pi'_r$ has length at most $\alpha L^{V'_r}$. 

By Proposition \ref{NibbleSplittingLemma} again, either $k < t_1 + V_{r-1}$
or at time $t_1 + V_{r-1}$ the \nf \ of $\pi'_r$ admits a hard
splitting into edge-paths each of which is either: 
\begin{enumerate}
\item   a single edge of weight $r-1$;
\item   a \nf \ of a weight $r-1$ indivisible Nielsen path;
\item   a \nf \ of a weight $r-1$ \gep ; or
\item   a path of weight at most $r-2$.
\end{enumerate}

We continue in this manner.  We may conceivably fall into case (4)
each time until $t_1 + V_1$ when it is not possible to fall into a
path of weight at most $1 - 1$!  Thus at some stage we must fall into
one of the first three cases.  This completes the proof of Claim 1.

\bigskip 

{\bf The existence of $\alpha$}.
We must find a
bound, in terms of $n$ and $f$, on the length of indivisible Nielsen
paths that arise in the \nf \ of $\rho$ with a hard splitting
immediately on either side.\footnote{Recall that the definition of
  Type B incidents excluded the case of Nielsen paths which lie in the
  \nf \ of a \gep \ with a hard splitting immediately on either side.}
To this end, suppose that $\e'$ is an edge which lies in an
indivisible Nielsen path $\tau$ in a $k'$-step \nf \ of $\rho$, and
that there is a hard splitting immediately on either side of $\tau$. 
We again denote the $i$-step nibbled future of $\rho$ by $\rho_i$ for
$0 \le i \le k'$. 

As above, we associate a diagram $\Delta_{k'}$ to
$\rho_{k'}$.\footnote{If we are considering Nielsen paths arising in
  the past of $\e$ above, then we can assume $k' \le k$ and that
  $\Delta_{k'}$ is a subdiagram of $\Delta_k$ in the obvious way.}
Denote by $q$ the path in the family forest of $\Delta_{k'}$ which follows the past of $\e'$. Let $q(i)$ be the edge in $\rho_i$ which intersects $q$.  Let the sequence
of weights of the edges $q(i)$ be $W' = \{ w'_0, \ldots , w'_{k'} \}$. 

Define incidents of Type A and B for $W'$ in exactly the same way as
for $W$, and let $t_2$ be the time of the last incident of Type A for
$W'$.  If there is no incident of Type A for $W'$ let $t_2 = 0$.  Let
$\kappa$ be the folded seed containing $q(t_2)$; in case $t_2 = 0$
let $\kappa = \rho$.  Define $\theta = \mbox{max} \{ n , L \}$ and note
that $| \kappa | \le \theta$.  The path $\tau$ must lie in the \nf \ of
$\kappa$, so it suffices to consider the \nf \ of $\kappa$.  Suppose
that $\kappa$ has weight $r'$. 

We deal with the \nf \ of $\kappa$ in the same way as we dealt with
that of $\pi$.   Let $\kappa_0 = \kappa, \kappa_1, \ldots$ be the \nf s of $\kappa$.

\medskip

{\bf Claim 2:} There exists a constant $\beta' = \beta'(n,f)$ so that one of the following must occur:  
\begin{enumerate}
\item[(i)] for some $t_2 < i < k'$, the edge $q(i)$ lies in a \gep \ in $f_{\#}(\kappa_{i-1})$ that has a hard splitting immediately on either side; 
\item[(ii)] not in case (i), {\em and} at some time $i \le k'$ the edge $q(i)$ lies in an indivisible Nielsen
path $\tau_0$ in $f_{\#}(\kappa_{i-1})$ so that $|\tau_0| \le\theta L^{\beta'}$ and immediately on either side of $\tau_0$ there is a hard splitting, {\em and} there are no incidents of Type B after time $i$;
\item[(iii)] $k' - t_2 \le \beta'$; or 
\item[(iv)] there is a hard splitting of $\kappa_{k'}$ immediately on either side of $\e'$. 
\end{enumerate}

Let us prove that this claim implies the existence of $\alpha$ and hence completes the proof of the theorem.  By definition, $\alpha$ is required to be an upper bound on the length of an arbitrary Nielsen path $\tau$ involved in a Type B incident.  We assume this incident occurs at time $k'$ and use Claim 2 to analyse what happens.  

Case (i) of Claim 2 is irrelevant in this regard.  If case (ii) occurs, the futures of $\tau_0$ are unchanging up to time $k'$, so $\tau = \tau_0$ and we have our required bound.  In case (iii) the length of $\tau$ is at most $\theta L^{\beta'}$, and in case (iv) $\tau$ is a single edge.  It suffices to let $\alpha = \theta L^{\beta'}$.

\medskip

It remains to prove Claim 2.  The proof of Claim 2 follows that of Claim 1 almost verbatim, with $\theta$ in place of $\alpha$ and $\kappa$ in place of $\rho$, etc.,  {\em except}  that the third sentence in Case (2) of the proof becomes invalid because Type B incidents after time $t_2 + V_r$ may occur.  

In this setting, suppose $\pi_r$ (which occurs at time $t_2 + V_r$) is a Nielsen path, but that we are not in case (ii) of Claim 2, and there is a subsequent Type B incident at time $j$, say.  The length of $\pi_r$ is at most $\theta L^{V_r}$. The Nielsen path at time $j-1$ has the same length as the one at time $t_2 + V_r$.  We go forward to time $j$, where the future of $\pi_r$ is no longer a Nielsen path, and continue the proof of Case (2) from the fourth sentence of the proof.

Otherwise, the proof of Claim 2 is the same as that of Claim 1 (the above modification is required at each weight, but at most {\em once} for each weight).  The only way in which the length bounds change is in the replacement of $\theta$ by $\alpha$ (including in the definitions of $V_i$ and $V'_i$). This finally completes the proof of Theorem \ref{DecompTheorem}.
\end{proof}

\section{Passing to an Iterate of $f$} \label{IterateSection}

It is important to be able to replace $f$ by an iterate $f_0 = f_{\#}^k$, for $k \ge 1$.
However, when doing this, it is important to be able to retain control over certain
constants (since which iterate we choose will depend on some of these constants).
In this section we describe what happens to various definitions when we replace 
$f$ by an iterate.  Suppose that $k \ge 1$, and consider the relationship between 
$f$ and $f_0 = f_{\#}^k$.

First, for any integer $j \ge 1$, the set of $kj$-monochromatic paths for $f$ is the same
as the set of $j$-monochromatic paths for $f_0$.  Therefore, once Theorem 
\ref{MainProp} is proved, we will pass to an iterate so that $r$-monochromatic
becomes $1$-monochromatic. However, the story is not quite as simple as that.

It is not hard to see that if $\sigma \odot \nu$ is a 
hard splitting for $f$, then it is also a hard splitting for $f_0$.

When $f$ is replaced by $f_0$, the set of \gep s is unchanged, as are the sets of 
\pep s and indivisible Nielsen paths.  Also, the set of indivisible
Nielsen paths which occur as sub-paths of $f(E)$ for some linear edge $E$ remains unchanged.

With the definition as given, the set of $(J,f_0)$-atoms may be smaller than the set of $(J,f)$-atoms.  This is because an atom is required to be 
$1$-monochromatic.  However, we continue to consider the set of $(J,f)$-atoms
even when we pass to $f_0$, and we also consider paths to be {\em beaded}
if they are $(J,f)$-beaded.

Since we are quantifying over a smaller set of paths the constant $V(n,f_0)$ in Theorem \ref{DecompTheorem} is assumed, without loss of generality, to be $V(n,f)$.  This is an 
important point, because the constant $V$ is used to find the appropriate $J$ when proving 
Theorem \ref{MainProp}.  When passing from $f$ to $f_0$, we need this $J$ to 
remain unchanged, for the appropriate iterate $k$ which we eventually choose depends 
crucially upon $J$ (See Addendum \ref{add:SameJ}).

It is also clear that if $m \le n$ then without loss of generality we may assume that $V(m,f) \le V(n,f)$.  Once again, this is because we are considering a smaller set of paths when defining $V(m,f)$.

We now want to replace $f$ by a fixed iterate in order to control some of the cancellation within monochromatic paths.  The following lemma is particularly useful in the proof of Proposition \ref{gammasingleedge} below, and also for Theorem 
\ref{ColourCancelMain}.  In particular, it will be used to find the value of $d$ in the Beaded Decomposition Theorem.  Lemma \ref{Power1} allows us to tune the improved
relative train track map in order to exclude some troublesome cancellation phenomena
that can otherwise occur in \nf s.

\begin{lemma} \label{Power1}
There exists $k_1 \ge 1$ so that $f_1 = f_{\#}^{k_1}$ satisfies the following. 
Suppose that $E$ is an exponential edge of weight $r$ and that $\sigma$ is an 
indivisible Nielsen path of weight $r$ (if it exists, $\sigma$ is unique up to a change 
of orientation). Then
\begin{enumerate}
\item $|f_1(E)| > |\sigma|$.
\item Moreover, if $\sigma$ is an indivisible Nielsen path of exponential weight 
$r$ and $\sigma_0$ is a proper subedge-path of $\sigma$,  then $(f_1)_{\#}(\sigma_0)$ is $r$-legal.
\item If $\sigma_0$ is a proper initial sub edge-path of $\sigma$ then $(f_1)_{\#}(\sigma_0)$ admits a hard splitting, $f(E) \odot \xi$, where $E$ is the edge on the left end of $\sigma$.
\item Finally, if $\sigma_1$ is a proper terminal sub edge-path of $\sigma$ then $(f_1)_{\#}(\sigma_1) = \xi' \odot f(E')$ where $E'$ is the edge on the right end of $\sigma$.
\end{enumerate}
Now suppose that $\sigma$ is an indivisible Nielsen path of parabolic weight $r$ and that $\sigma$ is a sub edge-path of $f(E_1)$ for some linear edge $E_1$.  The path $\sigma$ is either of the form $E \eta^{m_{\sigma}} \overline{E'}$ or of the form $E \overline{\eta}^{m_{\sigma}} \overline{E'}$, for some linear edges $E$ and $E'$.  Then 
\begin{enumerate}
\item If $\sigma_0$ is a proper initial sub edge-path of $\sigma$ then
\[      (f_1)_{\#}(\sigma_0) = E \odot \eta \odot \cdots \odot \eta \odot \xi''   ,       \]
where there are more than $m_{\sigma}$ copies of $\eta$ visible in this splitting.
\item If $\sigma_1$ is a proper terminal sub edge-path of $\sigma$ then
\[      (f_1)_{\#}(\sigma_1) = \xi' \odot \overline{\eta} \odot \cdots \odot \overline{\eta} \odot \overline{E'},    \]
where there are more than $m_\sigma$ copies of $\overline{\eta}$ visible in this splitting;
\end{enumerate}
\end{lemma}
\begin{proof}
First suppose that  $H_r$ is an exponential stratum, that $\sigma$ is an indivisible Nielsen path of weight $r$, and that $E$ is an edge of weight $r$.  Since $|f_{\#}^j(E)|$ grows exponentially with $j$, and $|f_{\#}^j(\sigma)|$ is constant, there is certainly some $d_0$ so that $|f_{\#}^d(E)| > |\sigma|$ for all $d \ge d_0$.

There is a single $r$-illegal turn in $\sigma$, and if $\sigma_0$ is a proper sub edge-path of $\sigma$.  By Lemma \ref{NoTau}, no future of $\sigma_0$ can contain $\sigma$ as a subpath.  The number of $r$-illegal turns in iterates of $\sigma_0$ must stabilise, so by Lemma \ref{ExpSplitting} there is an iterate of $\sigma_0$ which is $r$-legal.  Since there are only finitely many paths $\sigma_0$, we can choose an iterate of $f$ which works for all such $\sigma_0$.

Suppose now that $\sigma_0$ is a proper initial sub edge-path of $\sigma$, and that $E$ is the edge on the left end of $\sigma$.  It is not hard to see that every (entire) future of $\sigma_0$ has $E$ on its left end.  We have found an iterate of $f$ so that $f_{\#}^{d'}(\sigma_0)$ is $r$-legal.  It now follows immediately that
\[      f_{\#}^{d'+1}(\sigma_0) = f(E) \odot \xi        ,       \]
for some path $\xi$.  The case when $\sigma_1$ is a proper terminal sub edge-path of $\sigma$ is identical.

Now suppose that $H_r$ is a parabolic stratum and that $\sigma$ is an indivisible Nielsen path of weight $r$ of the form in the statement of the lemma.  The claims about sub-paths of $\sigma$ follow from the hard splittings $f(E) = E \odot u_E$ and $f(E') = E' \odot u_{E'}$, and from the fact that $m_{\sigma}$ is bounded because $\sigma$ is a subpath of some $f(E_1)$.

As in Remark \ref{stable}, we can treat each of the cases separately, and finally take a maximum.  
\end{proof}

\section{The Nibbled Futures of \gep s} \label{nfGEPSection}
\index{future!nibbled}

In this section $f$ is an improved relative train track map, although we do not suppose
yet that we have replaced $f$ by an iterate so that Lemma \ref{Power1} holds with
$k_1 = 1$.

The entire future of a \gep \ is a \gep \ but a \nf \ need not be and Theorem \ref{DecompTheorem} tells us that we need to analyse these \nf s.  This analysis will lead us to define \index{proto-\pep}
{\em proto-\pep s}.  In Proposition \ref{gammasingleedge}, we establish a
normal form for proto-\pep s which proves that proto-\pep s are in fact the \pep s which 
appear in the Beaded Decomposition Theorem.

To this end, suppose that
\[ \zeta = E_i \overline\tau^n \overline{E_j}           \]
is a \gep , where $\tau$ is a Nielsen path, $f(E_i) = 
E_i \odot \tau^{m_i}$ and $f(E_j) = E_j \odot \tau^{m_j}$.  As in Definition
\ref{Exceptional}, we consider $E_i \overline\tau^n \overline{E_j}$ to be
unoriented, but here we do not suppose that $j \le i$.  However, we suppose $n > 0$ and thus,
since $E_i \overline\tau^n \overline{E_j}$ is a \gep , $m_j > m_i > 0$.

The analysis of \gep s of the form $E_j {\tau}^n \overline{E_i}$ is entirely similar
to that of \gep s of the form $E_i \overline\tau^n \overline{E_j}$ except that one must reverse all left-right orientations.  Therefore, we ignore this case until Definition \ref{pep} below (and often afterwards also!).

We fix a sequence of \nf s $\zeta = \rho_{-l}, \ldots , \rho_0, \rho_1, \ldots , \rho_k, \ldots$ of $\zeta$, where $\rho_0$ is the first \nf \ which is not the entire future.  Since the entire future of a \gep \ is a \gep , we restrict our attention to the \nf s of $\rho_0$.

There are three cases to consider, depending on the type of sub-path on either end of $\rho_0$.
\begin{enumerate}
\item $\rho_0 = \bar{\tau_0} \bar{\tau}^m \overline{E_j}$;
\item $\rho_0 = \bar{\tau_0} \bar{\tau}^m \bar{\tau_1}$.
\item $\rho_0 = E_i \bar{\tau}^m \bar{\tau_1}$;
\end{enumerate}
where $\tau_0$ is a (possibly empty) initial sub edge-path of $\tau$, and $\tau_1$ is
a (possibly empty) terminal sub edge-path of $\tau$.

In case (1) $\rho_0$ admits a hard splitting
\[      \rho_0 = \overline\tau_0 \odot \overline\tau \odot \cdots \odot \overline\tau \odot \overline E_j.    \]
Since  $\tau_0$ is a sub edge-path of $f(E_i)$, it has length less than
 $L$ and its nibbled futures admit hard splittings as in
Theorem \ref{DecompTheorem} into \nf s of GEPs and paths of length at most
$V(L,f)$. These \gep s 
 will necessarily be of strictly lower weight than $\rho_0$, since $\overline\tau_0$ is. Thus, 
case (1) is easily dealt with  by an induction on weight, supposing that we have a nice 
splitting of the \nf s of lower weight \gep s; this is made precise in Proposition \ref{nfgeptopep}.  Case (2) is entirely similar.

Case (3)  is by far the most troublesome of the three, and it is this case which leads to the definition of {\em proto-\pep s} in Definition \ref{pep} below.  Henceforth assume $\rho_0 = E_i \bar{\tau}^m \bar{\tau_1}$.

Each of the \nf s of $\rho_0$ (up to the moment of death, Subsection \ref{Death}) 
has a \nf \ of 
$\overline\tau_1$ on the right.  If the latter becomes empty at some point, 
the \nf \ of $\rho_0$ at this time has the form $E_i \overline\tau^{n'} \overline\tau_2$, where $\tau_2$ is a proper (but possibly empty) sub edge-path of  $\tau$.  We
 restart our analysis at this moment.  Hence we make the following

\begin{workingassumption} \label{NibAssumption}
We make the following two assumptions on the $k$-step \nf s considered:
\begin{enumerate}
\item $\rho_0 = E_1 \overline{\tau}^m \overline{\tau}_1$;
\item all nibbling of $\rho_k$ occurs on the right; and 
\item the $k$-step \nf \ $\overline\tau_{1,k}$ of $\overline\tau_1$ inherited from $\rho_k$ is non-empty.
\end{enumerate}
\end{workingassumption}

We will deal with the case $m - km_i < 0$ later, in particular with the value of $k$ for 
which $m - (k-1)m_i \ge 0$ but $m - km_i < 0$.  For now suppose that $m - km_i \ge 0$.

In this case, the path $\rho_k$ 
has the form
\[      \rho_k = E_i \overline\tau^{m - km_i} \overline\tau_{1,k} .     \]

There are (possibly empty) Nielsen edge-paths $\iota$ and $\nu$, and an indivisible Nielsen edge-path $\sigma$ so that 
\begin{equation} \label{TauDecomp}
\tau = \iota \odot \sigma \odot \nu \mbox{ and } \tau_1 = \sigma_1 \odot \nu     ,
\end{equation}
where $\sigma_1$ is a proper terminal sub edge-path of $\sigma$.  Now, as in Working Assumption \ref{NibAssumption}, there is no loss of generality in supposing that
\[      \rho_k = E_i \overline\tau^{m-km_i} \bar \nu \bar\sigma_{1,k} , \]
where $\overline\sigma_{1,k}$ is the \nf \ of $\overline\sigma_1$ inherited from $\rho_k$, and that $\overline\sigma_{1,k}$ is non-empty.

Since $|\sigma_1| < L$, by Theorem \ref{DecompTheorem} the path $\sigma_{1,k}$ admits a hard splitting into edge-paths each of which is either the \nf \ of a \gep , or of length at most $V(L,f)$;  we take the (unique) maximal hard splitting of $\sigma_{1,k}$ into edge-paths.  

Let $s = \lfloor m / m_i \rfloor + 1$.  In $\rho_s$ (but not before) there may be some 
interaction between the future of $E_i$ and $\overline\sigma_{1,s}$.   We denote by $\gamma_{\sigma_1}^{k,m}$ the concatenation of those factors in the hard splitting of $\overline\sigma_{1,k}$
which contain edges any part of whose future is eventually cancelled by some edge in the future 
of $E_i$ under any choice of \nf s of $\rho_k$ (not just the $\rho_{k+t}$ chosen earlier) and any choice of tightening.  Below we will analyse more carefully the structure of the paths $\overline\sigma_{1,k}$ and $\gamma_{\sigma_1}^{k,m}$.

We now have $\overline\sigma_{1,k} = \gamma_{\sigma_1}^{k,m} \odot \sigma_{1,k}^\bullet$.  From (\ref{TauDecomp}), we also have
\begin{equation} \label{nfpep}
\rho_k = E_i \overline\tau^{m-km_i}  \overline\nu \gamma_{\sigma_1}^{k,m} \odot \sigma_{1,k}^\bullet .
\end{equation}

\begin{definition}[Proto-\pep s] \label{pep}\index{proto-\pep}
Suppose that $\tau$ is a Nielsen edge-path, $E_i$ a linear edge such that $f(E_i) = E_i 
\odot \tau^{m_i}$ and $\tau_1$ a proper terminal sub edge-path of $\tau$ such that $\tau_1 
= \sigma_1 \odot \nu$ as in (\ref{TauDecomp}).  Let $k,m \ge 0$ be such that $m - km_i \ge 0$ and let $\gamma_{\sigma_1}^{k,m}$ be as in (\ref{nfpep}).
A path $\pi$ is called a {\em proto-\pep} if either $\pi$ of $\overline{\pi}$ is of the form
\[      E_i  \overline\tau^{m-km_i} \overline\nu \gamma_{\sigma_1}^{k,m} .       \]
\end{definition}

\begin{remarks} 
\ \par

\begin{enumerate}
\item The definition of proto-\pep s is intended to capture those paths which remain when a \gep \ is 
partially cancelled, leaving a path which may shrink in size of its own accord.
\item By definition, a proto-\pep \ admits no non-vacuous hard splitting into edge-paths.
\end{enumerate}
\end{remarks}

We now introduce two distinguished kinds of proto-\pep s.

\begin{definition} \label{stabletransient}
Suppose that 
\[      \pi = E_i  \overline\tau^{m-km_i} \overline\nu \gamma_{\sigma_1}^{k,m} ,\] 
is a proto-\pep \ as in Definition \ref{pep}.

The path $\pi$ is a {\em transient} proto-\pep \  if $k=0$.

The path $\pi$ is a {\em stable} proto \pep \ if $\gamma_{\sigma_1}^{k,m}$
is a single edge.
\end{definition}

\begin{lemma} \label{transientispep}
A transient proto-\pep \ is a \pep. 
\end{lemma}
\begin{proof}
With the notation of Definition \ref{pep}, in this case $\gamma^{0,m}_{\sigma_1}$
is visibly a sub-path of $\bar{\tau}$, and the proto-\pep \ is visibly a sub-path of a \gep.
\end{proof}

\begin{lemma} \label{stableprotoispep}
A stable proto-\pep \ is a \pep.
\end{lemma}
\begin{proof}
Since $\bar{\sigma}$ is a Nielsen path, if $\alpha$ is a \nf \ of $\bar{\sigma}$ where
all the nibbling has occurred on the right, then the first edge in $\alpha$ is the
same as the first edge in $\bar{\sigma}$.

On the other hand, $\gamma_{\sigma_1}^{k,m}$ is a \nf \ of $\bar{\sigma}$ where all the
nibbling has occurred on the right.  Therefore, if $\gamma_{\sigma_1}^{k,m}$
is a single edge then it must be a sub-path of $\bar{\sigma}$.  It follows
immediately that any stable proto-\pep \ must be a \pep.
\end{proof}

\begin{remark}
We will prove in Proposition \ref{gammasingleedge} that after
replacing $f$ by a suitable iterate all proto-\pep s are either
transient or stable, and hence are \pep s.
\end{remark}

\subsection{The Death of a proto-\pep } \label{Death}

Suppose that $\pi = E_i \overline\tau^{m-km_i} \overline\nu \gamma_{\sigma_1}^{k,m}$
 is a proto-\pep \ with \nf s satisfying Assumption \ref{NibAssumption}.  Let $q = \lfloor \frac{m-km_i}{m_i}
 \rfloor + 1$, and consider, $\pi_{q-1}$, a $(q-1)$-step \nf \  of $\pi$. As before, we assume 
 that the $(q-1)$-step \nf \ of $\gamma_{\sigma_0}^{k,m}$ inherited from a $\pi_{q-1}$
  is not empty and that the edge labelled $E_i$ on the very left is not nibbled. 

In $\pi_{q}$, the edge $E_i$ has consumed all of the copies of $\overline\tau$ and begins to interact 
with the future of $\overline\nu \gamma_{\sigma_1}^{k,m}$. Also, the future of $\pi$ at time $q$ need not contain a \pep.  Hence we refer to the time $q$ as the {\em death of the \pep}.
Recall that $\tau = \iota \odot \sigma 
\odot \nu$ and that $\gamma_{\sigma_1}^{k,m}$ is a $k$-step \nf \ of $\overline\sigma_1$, where $\sigma_1$ is a proper subpath of $\sigma$.  Let $p = m - (k+q-1)m_i $, so that $0 \le p < m_i$.

The path $\pi_{q-1}$ has the form
\[      \pi_{q-1} = E_i \overline\tau^{p} \overline\nu \gamma_{\sigma_1}^{k+q-1,m} .       \]
Suppose that $\pi_q$ is a $1$-step
\nf \ of $\pi_{q-1}$.  In other words, $\pi_q$ is a subpath of $f_{\#}(\pi_{q-1})$.
Consider what happens when $f(\pi_{q-1})$ is tightened to form $f_{\#}(\pi_{q-1})$ 
(with any choice of tightening).   The $p$ copies of $\overline\tau$ 
(possibly in various stages of tightening) will be consumed by $E_i$, leaving 
$\overline \nu \odot f(\gamma_{\sigma_1}^{k+q-1,m})$ to interact with at least one 
remaining copy of $\tau = \iota \odot \sigma \odot \nu$.  The paths $\nu$ and 
$\overline{\nu}$ will cancel with each other\footnote{The hard splittings imply that 
this cancellation must occur under {\em any} choice of tightening.}.

Lemma \ref{NibNoTau} states that $\gamma_{\sigma_1}^{k,m}$ cannot contain $\sigma$ as a subpath.  Therefore, once $\nu$ and $\overline{\nu}$ have cancelled, not all of $\overline{\sigma}$ will cancel with $f(\gamma_{\sigma_1}^{k+q-1,m})$.  A consequence of this discussion (and the fact that $f(E_i) = E_i \odot \tau^{m_i}$) is the following

\begin{lemma} \label{pepSplitting}
Suppose that $\pi = E_i \overline\tau^{m-km_i} \overline\nu \gamma_{\sigma_0}^{k,m}$ is a proto-\pep , and let $q = \lfloor \frac{m-km_i}{m_i} \rfloor + 1$.  Suppose that $\pi_{q-1}$ is a $(q-1)$-step \nf \ of $\pi$ satisfying Assumption \ref{NibAssumption}.  If $\pi_q$ is an immediate \nf \ of $\pi_{q-1}$ and $\pi_q$ contains $E_i$ then $\pi_q$ admits a hard splitting
\[      \pi_q = E_i \odot \lambda.      \]
\end{lemma}

We now analyse the interaction between $f(\gamma_{\sigma_1}^{k+q-1,m})$ and $\sigma$ more closely.  As usual, there are two cases to consider, depending on whether $\sigma$ has exponential or parabolic weight\footnote{Recall that there are three kinds of indivisible Nielsen paths: constant edges, parabolic weight and exponential weight.  If $\sigma$ has nontrivial proper sub edge-paths, then it is certainly not a single edge, constant or not.}.

In the following proposition, $f_1$ is the iterate of $f$ from Lemma \ref{Power1} and we are using the definitions as explained in Section \ref{IterateSection}.  Also, we assume that proto-\pep s are defined using $f_1$, not $f$.

\smallskip

\begin{proposition} \label{gammasingleedge}
Every proto-\pep \ for $f_1$ is either transient or stable.
In particular, every proto-\pep \  for $f_1$ is a \pep.
\end{proposition}
\begin{proof}
Let $\pi = E_i \overline\tau^{m-km_i} \overline\nu \gamma_{\sigma_1}^{k,m}$ be a 
proto-\pep \ for $f_1$.

Lemma \ref{transientispep} implies that  if $k = 0$ then $\pi$ is a \pep.
Consider Working Assumption \ref{NibAssumption}.
  If Assumption \ref{NibAssumption}.(2) fails to hold at any point, then we can restart our analysis, and 
  in particular we have a transient proto-\pep \  at this moment.
  Thus we may suppose that $\pi$ is an initial sub-path of a $k$-step \nf \ of a \gep , where 
  $k \ge 1$ and we may further suppose that $\pi$ satisfies Assumption \ref{NibAssumption}.(2). We prove that in this case $\pi$ is a stable proto-\pep .

First suppose that $\sigma$ has exponential weight, $r$ say.  If $\overline{\sigma}_0$ is a 
proper initial sub edge-path of $\overline\sigma$ then Lemma \ref{Power1} asserts that 
\[      (f_1)_{\#}(\sigma_0) = f(E) \odot \xi   ,       \]
and $| f(E)| > |\sigma|$.  Note also that $f(E) = E \odot \xi''$ for some path $\xi''$.

Now, at the death of the proto-\pep , the \nf \ of $\gamma_{\sigma_0}^{k,m}$ interacts 
with a copy of $E_i$, and in particular with a copy of $f(\sigma)$ (in some stage
 of tightening).  Now the above hard splitting, and the fact that $\sigma$ is not 
 $r$-legal whilst $f(E)$ is, shows that $\gamma_{\sigma_1}^{k,m}$ must be a 
 single edge (namely $E$).

Suppose now that $\sigma$ has parabolic weight.  Since $\sigma$ has proper 
sub edge-paths, it is not a single edge and so $\sigma$ or $\overline{\sigma}$ 
has the form $E \eta^{m_{\sigma}} \overline{E'}$.  The hard splittings guaranteed
 by Lemma \ref{Power1} now imply that $\gamma_{\sigma_1}^{k,m}$ is a single
  edge in this case also.

Therefore, every proto-\pep \ for $f_1$ is transient or stable, proving the first
assertion of the proposition.  The second assertion follows from the first assertion,
and Lemmas \ref{transientispep} and \ref{stableprotoispep}.
\end{proof}

Finally, we can prove the main result of this section.  In the following, $L_1$ is the maximum length of $f_1(E)$ over all edges $E$ of $G$.

The following statement assumes the conventions of Section \ref{IterateSection}.
\begin{proposition} \label{nfgeptopep}
Under iteration of the map $f_1$ constructed in Lemma \ref{Power1}, any
\nf \ of a \gep \ admits a hard splitting into edge-paths, each of which is either a 
\gep , a \pep , or of length at most $V(2L_1,f)$.
\end{proposition}
\begin{proof}
Suppose that $E_i \overline{\tau}^n \overline{E_j}$ is a \gep \ of weight $r$.  We may suppose
 by induction that any \nf \ of any \gep \ of weight less than $r$ admits a hard splitting of
the required form (the base case $r=1$ is vacuous, since there cannot be a \gep \ of weight 
$1$).

Suppose that $\rho$ is a \nf \ of $E_i \tau^n \overline{E_j}$.  If $\rho$ is the entire future,
it is a \gep \ and there is nothing to prove.  Otherwise, as in the analysis at the beginning of this section, 
we consider the first time when a \nf \ is not the entire future.  Let the \nf \ be $\rho_0$.  
In cases (1) and (2) from that analysis, $\rho_0$ admits a hard splitting into edge-paths, 
each of which is either (i) $\overline{E}_i$; (ii) $\overline{\tau}$; or (iii) a proper 
sub edge-path of $\tau$.  In each of these cases, Theorem \ref{DecompTheorem} asserts
 that there is a hard splitting of $\rho$ into edge-paths, each of which is either of length at
  most $V(L,f)$ or is the \nf \ of a \gep .  Any \nf \ of a \gep \ which occurs in this splitting is 
  necessarily of weight strictly less than $r$, and so admits a hard splitting of the required 
  form by induction.  

Suppose then that $\rho_0$ satisfies Case (3), the third of the cases articulated at 
the beginning of this section.  In this 
case, $\rho_0$ is a transient proto-\pep .  Also, any time that Assumption \ref{NibAssumption}.(2) 
is not satisfied, the \nf \ of $\rho_0$ is a transient proto-\pep .  Thus, we may assume that 
Assumption \ref{NibAssumption} is satisfied. If $m -km_i \ge 0$ then we have
\[      \rho = E_i \overline\tau^{m-km_i} \overline\nu \gamma_{\sigma_1}^{k,m} \odot 
\sigma_{1,k}^\bullet    .       \]
The first path in this splitting is a stable \pep \ by Proposition \ref{gammasingleedge}.  
Once again, Theorem \ref{DecompTheorem} and the inductive hypothesis yield a hard 
splitting of $\sigma_{1,k}^\bullet$ of the required form.

Finally, suppose that Case (3) pertains and 
$m-km_i < 0$.  Let $q = \lfloor \frac{m-km_i}{m_i} \rfloor + 1$ (the significance of this moment 
-- ``the death of the \pep " -- was explained at the beginning of this subsection).  By the 
definition of a \pep \ (Definition \ref{pep}), the $q$-step \nf \ of $\rho_0$ admits a hard 
splitting as
\[      E_i \overline\tau^{m-qm_i} \overline\nu \gamma_{\sigma_1}^{q,m} \odot 
\sigma_{1,q}^\bullet .  \]
By Lemma \ref{pepSplitting}, the immediate future of  $E_i \overline\tau^{m-qm_i}
 \overline\nu   \gamma_{\sigma_1}^{q,m}$ admits a hard splitting as $E_i \odot \xi$.  
 Since $\gamma_{\sigma_1}^{r,m}$ is a single edge, we have a bound of $2L_1$ on the 
 length of $\xi$.  Any \nf \ of $E_i \odot \xi$ now admits a hard splitting into edge-paths,
  each of which is either a \gep , a \pep \ or of length at most $V(2L_1,f)$, by induction on 
  weight and Theorem \ref{DecompTheorem}.
\end{proof}

We highlight one consequence of Proposition \ref{nfgeptopep}:

\begin{corollary} \label{cor:nfpep}
Suppose that $\rho = E_i \overline{\tau}^{m-km_i} 
\overline{\nu} \gamma$ is a 
\pep.  Any immediate \nf \ of $\rho$ (with all nibbling on the right) 
has one of the following two forms:
\begin{enumerate}
\item $\rho' \odot \sigma$, where $\rho'$ is a \pep\ and $\sigma$ admits a hard
splitting into \atom s; or
\item $E_i \odot \sigma$, where $\sigma$ 
admits a hard splitting into \atom s.
\end{enumerate}
In particular, this is true of $f_{\#}(\rho)$.

There are entirely analogous statements in case $\rho$ is a \pep \ where
$\overline{\rho}$ has the above form and all nibbling occurs on the left.
\end{corollary}

 \section{Proof of the Beaded Decomposition Theorem}
\index{Beaded Decomposition Theorem}

In this section, we finally prove Theorem \ref{SuperBDT}.  As noted in Remark 
\ref{SuperBDTEnough}, this immediately implies the Beaded Decomposition 
Theorem.

\begin{proof}[Proof (Theorem \ref{SuperBDT}).]
Take $d = k_1$, the constant from Lemma \ref{Power1}.  Let $L_1$ be  the maximum length of $f_{\#}^{k_1}(E)$ for any edge $E \in G$, let $s = \max\{ 2L_1, n \}$,
and let $J = V(s,f)$, where $V$ is the constant from Theorem \ref{DecompTheorem}.

Suppose that $\rho$ is a path so that $|\rho| \le n$, and let $\rho'$ be a $kd$-step
\nf \ of $\rho$ for some positive integer $k$. Then $\rho'$ is a $k$-step \nf \ of
$\rho$ with respect to $f_1 = f_{\#}^{k_1}$.  By Proposition \ref{gammasingleedge}, every proto-\pep \ for $f_1$ is a \pep .

By Theorem \ref{DecompTheorem}, $\rho'$ admits a hard splitting into edge-paths, each of which is either the \nf \ of a \gep \ or else has length at most $V(n,f)$.  By Proposition \ref{nfgeptopep}, if we replace $f$ by $f_1$ then any \nf \ of a \gep \ admits a hard splitting into edge-paths, each of which is either a \gep , a \pep \ or else has length at most $V(2L_1,f)$.   By Lemma \ref{HardSplitProperties}, the splitting of the \nf \ of a \gep \ is inherited by $\rho$.

We have shown that $\rho$ is $(J,f)$-beaded, as required.
\end{proof}

\begin{proof}[Proof (Addendum \ref{SuperBDT}).]
We have already remarked that, for a fixed $m$, the constant $V(m,f)$ from 
Theorem \ref{DecompTheorem} remains unchanged when $f$ is replaced by an iterate.

As in Section \ref{IterateSection}, we retain the notion of $(J,f)$-beaded with the original
$f$ when passing
to an iterate of $f$

Therefore, when $f$ is replaced by an iterate, Theorem \ref{SuperBDT} remains true
with the same constant $J$.  This immediately implies that the same is true
of the Beaded Decomposition, which is what we were required to prove.
\end{proof}

\section{Refinements of the Beaded Decomposition Theorem} \label{refinements}

The Beaded Decomposition Theorem is the main result of Part \ref{Part:BG2}.  In this section, we 
provide a few further refinements that will be required for future applications.

Throughout this section we suppose that $f$ has been replaced with $f_1$ from 
Lemma \ref{Power1}, whilst maintaining the conventions for definitions from 
Section \ref{IterateSection}.  When we refer to $f$ we mean this iterate $f_1$.  With this in 
mind, a {\em monochromatic path} is a $1$-monochromatic path for $f$.  Similarly, 
armed with Theorem \ref{MainProp}, we refer to $(J,f)$-beads, simply as {\em beads},
and a path which is $(J,f)$-beaded will be referred to simply as {\em beaded}.
The constant $L$ now refers to the maximum length $|f(E)|$ for edges $E \in G$ { \em 
with  the new $f$}.

In the following theorem, the {\em past} of an edge is defined with respect to an 
arbitrary choice of tightening.

\begin{theorem} \label{ColourCancelMain}
There exists a constant $D_1$, depending only on $f$, with the
following properties. Suppose $i \ge D_1$, that $\chi$ is a monochromatic path and
that $\e$ is an edge in $f^{i}_{\#}(\chi)$ of weight $r$ whose past in $\chi$ is 
also of weight $r$.  Then $\e$ is contained in an edge-path $\rho$ so that $f_{\#}^i(\chi)$ has a hard splitting immediately on either side of $\rho$ and $\rho$ is one of the following:
\begin{enumerate}
\item a Nielsen path;
\item a \gep ;
\item  a \pep ; or
\item a single edge.
\end{enumerate}
\end{theorem}

\begin{proof}
Let $\chi$ be a monochromatic path.  For any $k \ge 0$, denote $f_{\#}^k(\chi)$ by $\chi_k$.  In a sense, we prove the theorem `backwards', by fixing an edge $\e_0$ of weight $r$ in $\chi_0 = \chi$ and considering its futures in the paths $\chi_k$, $k \ge 1$.  The purpose of this proof is to find a constant $D_1$ so that if $\e$ is any edge of weight $r$ in $\chi_i$ with past $\e_0$, and if $i \ge D_1$ then we can find a path $\rho$ around $\e$ satisfying one of the conditions of the statement of the theorem.

Fix $\e_0 \in \chi_0$.  By Theorem \ref{MainProp}, there is an edge-path $\pi$ containing $\e_0$ so that $\chi$ admits a hard splitting immediately on either side of $\pi$ and $\pi$ either (I) is a \gep ; (II) has length at most $J$; or (III) is a \pep .  In the light of Remark \ref{stable}, it suffices to establish the existence of a suitable $D_1$ in each case.  To consider the futures of $\e_0$ in the futures $f_{\#}^k(\chi)$ of $\chi$, it suffices to consider the futures of $\e_0$ within the (entire) futures of $\pi$.  Therefore, for $k \ge 0$, let $\pi_k = f_{\#}^k(\pi)$.  Suppose that we have chosen, for each $k$, an edge $\e_k$ in $\pi_k$ such that: (i) $\e_k$ lies in the future of $\e_0$; (ii) $\e_k$ has the same weight as $\e_0$; and (iii) $\e_k$ is in the future of $\e_{k-1}$ for all $k \ge 1$.

{\bf Case (I):} $\pi$ is a \gep .  In this case, the path $\pi_k$ is a \gep \ for all $k$, any future of $\e_0$ lies in $\pi_k$, and there is a hard splitting of $\chi_k$ immediately on either side of $\pi_k$.  Therefore, the conclusion of the theorem holds in this case with $D_1 = 1$.

\medskip

{\bf Case (II):} $|\pi| \le J$.  Denote the weight of $\pi$ by $s$.  Necessarily $s \ge r$.  By Lemma \ref{SplittingLemma} the path $\pi_{D(J)}$ admits a hard splitting into edge-paths, each of which is either
\begin{enumerate}
\item   a single edge of weight $s$; \label{DJ1}
\item   an indivisible Nielsen path of weight $s$; \label{DJ2}
\item   a \gep \ of weight $s$; or \label{DJ3}
\item   a path of weight at most $s-1$. \label{DJ4}
\end{enumerate}
We consider which of these types of edge-paths our chosen edge $\e_{D(J)}$ lies in.  In case (\ref{DJ1})  there is a hard splitting of $\pi_{D(J)}$ immediately on either side of the edge $\e_{D(J)}$, so for all $i \ge D(J)$ there is a hard splitting of $\pi_i$ immediately on either side of $\e_i$, since $\e_i$ and $\e_{D(J)}$ both have the same weight as $\e_0$.  For cases (\ref{DJ2}) and (\ref{DJ3}), $\e_{D(J)}$ lies in an indivisible Nielsen path or \gep \ with a hard splitting of $\pi_{D(J)}$ immediately on either side, so for all $i \ge D(J)$ any future of $\e_0$ in $\pi_i$, and in particular $\e_i$, lies in an indivisible Nielsen path of \gep \ immediately on either side of which there is a hard splitting of $\pi_i$.

Finally, suppose we are in case \eqref{DJ4} and not in any of cases \eqref{DJ1}--\eqref{DJ3}.  Then $\e_{D(J)}$ lies in an edge-path $\tilde\rho$ with a hard splitting of $\pi_{D(J)}$ immediately on either side, and that $\tilde\rho$ is not a single edge, an indivisible Nielsen path, or a \gep \footnote{In this case necessarily $s \le r-1$}.  We need only consider the future of $\tilde\rho$.  For $k \ge 0$, let $\rho_{D(J) + k} = f_{\#}^k(\tilde\rho)$ be the future of $\tilde\rho$ in $\pi_{D(J)+k}$.  Now, $|\tilde\rho| \le J L^{D(J)}$ so by Lemma \ref{SplittingLemma} the edge-path $\rho_{D(J) + D(J L^{D(J)})}$ admits a hard splitting into edges paths, each of which is either
\begin{enumerate}
\item   a single edge of weight $s-1$;
\item an indivisible Nielsen path of weight $s-1$;
\item   a \gep \ of weight $s-1$; or
\item   a path of weight at most $s-2$.
\end{enumerate}

We proceed in this manner.  If we ever fall into one of the first three cases, we are done.  Otherwise, after $s-r+1$ iterations of this argument, the fourth case describes a path of weight strictly less than $r$.  Since the weight of each $\e_i$ is $r$, it cannot lie in such a path, and one of the first three cases must hold.  Thus we have found the required bound $D_1$ in the case that $|\pi| \le J$.

\medskip

{\bf Case (III):} $\pi$ is a \pep . 

Let $\pi = E_i \overline\tau^{m-km_i} \overline\nu  \gamma_{\sigma_1}^{k,m}$ as in Definition \ref{pep}.  We consider where in the path $\pi$ the edge $\e_0$ lies. First of all, suppose that $\e_0$ is the unique copy of $E_i$.  Since $\e_0$ is parabolic, it has a unique weight $s$ future at each moment in time.  Let $q = \lfloor \frac{m-km_i}{m_i} \rfloor +1$, the moment of death. For $1 \le p \le q-1$, the edge $\e_p$ is the leftmost edge in a \pep \ and there is a hard splitting of $\pi_p$ immediately on either side of this \pep .  For $p \ge q$, Lemma \ref{pepSplitting} ensures that there is a hard splitting of $\pi_p$ immediately on either side of $\e_p$. Therefore in this case the conclusion of the theorem holds with $D_1 = 1$.

Now suppose that the edge $\e_0$ lies in one of the copies of $\overline\tau$ in $\pi$, or in the visible copy of $\overline\nu$.  Then any future of $\e_0$ lies in a copy of $\tau$ or $\nu$ respectively, which lies in a \pep \ with a hard splitting immediately on either side, until this copy of $\overline\tau$ or $\overline\nu$ is consumed by $E_i$.  Again, the conclusion of the theorem holds with $D_1 = 1$.

Finally, suppose that $\e_0$ lies in $\gamma_{\sigma_1}^{k,m}$.  For ease of notation, for the remainder of the proof $\gamma$ will denote $\gamma_{\sigma_1}^{k,m}$.  By Proposition \ref{gammasingleedge} $\gamma$ is a single edge.  Until the $q$-step \nf \ of $\pi$, any future of $\gamma$ of the same weight is either $\gamma$ or will have a splitting of $\pi$ immediately on either side. 

Since $\sigma$ is an indivisible Nielsen path, and $\gamma$ is a single edge, $\gamma$ is the leftmost edge of $\overline{\sigma}$.  Therefore $[\sigma \gamma]$ is a proper sub edge-path of $\sigma$. 

Suppose that $\sigma$ has exponential weight (this weight is $r$).  By Lemma \ref{Power1} and the above remark, $f_{\#}(\sigma \gamma)$ is $r$-legal.  Therefore, any future of $\gamma$ which has weight $r$ will have, at time $q$ and every time afterwards, a hard splitting immediately on either side.

Suppose now that $\sigma$ has parabolic weight $r$.  Since $[\sigma E]$ is a proper sub edge-path of $\sigma$, and since there is a single edge of weight $r$ in $f(E)$ and this is cancelled, it is impossible for $\gamma$ to have a future of weight $r$ after time $q$.
\end{proof}

Recall that the number of strata for the map $f : G \to G$ is $\omega$.  Recall also the definition of {\em displayed} from Definition \ref{Displayed}

\begin{lemma}
Let $\chi$ be a monochromatic path.  Then the number of displayed \pep s in $\chi$ of length more than $J$ is less than $2 \omega$.
\end{lemma}
\begin{proof}
Suppose that $\chi$ is a monochromatic path, and that $\rho$ is a subpath of $\chi$, with a hard splitting immediately on either side, such that $\rho$ is a \pep , and $|\rho| > J$.  Then, tracing through the past of $\chi$, the past of $\rho$ must have come into existence because of nibbling on one end of the past of $\chi$.  Suppose this nibbling was from the left.  Then all edges to the left of $\rho$ in $\chi$ have weight strictly less than that of $\rho$, since it must have come from a proper subpath of an indivisible Nielsen path in the \nf \ of the \gep \ which became $\rho$.  Also, any \pep \ to the left of $\rho$ must have arisen due to nibbling from the left. Therefore, there are at most $\omega$ \pep s of length more than $J$ which came about due to nibbling from the left.  The same is true for \pep s which arose through nibbling from the right.
\end{proof}

\begin{lemma} \label{Displayededge}
Let $D_1$ be the constant from Theorem \ref{ColourCancelMain}, and let $f_2 = (f_1)_{\#}^{D_1}$.  If $\rho$ is an \atom , then either $(f_2)_{\#}^{\omega}(\rho)$ is a beaded path all of whose beads are Nielsen paths and \gep s, or else there is some displayed edge $\e \subseteq (f_2)_{\#}^{\omega}(\rho)$ so that all edges in $(f_2)_{\#}^{\omega}(\rho)$ whose weight is greater than that of $\e$ lie in Nielsen paths and \gep s.
\end{lemma}
\begin{proof}
Suppose that $\rho$ is an \atom \ of weight $r$. If $H_r$ is a zero stratum and $(f_2)_{\#}(\rho)$ has weight $s$ then $H_s$ is not a zero stratum.  Thus, by going forwards one step in time if necessary, we suppose that $H_r$ is not a zero stratum, so $(f_2)_{\#}(\rho)$ has weight $r$.

By Theorem \ref{ColourCancelMain}, all edge of weight $r$ in $(f_2)_{\#}(\rho)$ are either displayed or lie in Nielsen paths or \gep s (since we are considering the entire future of an \atom , \pep s do not arise here).  If all edge of weight $r$ in $(f_2)_{\#}(\rho)$ lie in Nielsen paths or \gep s then we consider the \atom s in $(f_2)_{\#}(\rho)$ of weight less than $r$ (this hard splitting exists since $\rho$ and hence $(f_2)_{\#}(\rho)$ are monochromatic paths).  We now consider the immediate future of these \atom s in $(f_2)_{\#}^2(\rho)$, etc.  It is now clear that the statement of the lemma is true.
\end{proof}

Finally, we record an immediate consequence of the Beaded Decomposition
Theorem and Proposition \ref{nfgeptopep}:

\index{future!nibbled}
\begin{theorem} \label{t:nfbead}
Suppose that $\sigma$ is a beaded path.  Any \nf \ of $\sigma$ is also beaded.
\end{theorem}

\renewcommand{\thesection}{\ref{Part:BG3}.\arabic{section}}

\setcounter{section}{0}

\part{The General Case} \label{Part:BG3}

In Part \ref{Part:BG3}, we bring together the techniques developed in
Parts \ref{Part:BG1} and \ref{Part:BG2} to prove the main result of this book.

\renewcommand{\thethmspec}{{\bf{Main Theorem. \kern-.3em}}}

\smallskip

\noindent
\begin{thethmspec} 
{\em If $F$ is a finitely generated
free group and $\phi$ is
an automorphism of $F$ then $F \rtimes_{\phi} \mathbb Z$ satisfies a
quadratic isoperimetric inequality.}
\end{thethmspec}

\smallskip

In Part \ref{Part:BG1} we proved the Main Theorem  in the case of
{\em positive} automorphisms. That proof proceeded via an analysis
of van Kampen diagrams in the universal cover of the mapping
torus $R\times[0,1]/\langle(x,0)\sim (f(x),1)\rangle$, 
where $R$ is a 1-vertex graph with fundamental
group $F$ and  $f$ is the obvious
homotopy equivalence with $f_*=\phi$.
 
Such  $f$ are the prototypes for the 
{\em improved relative train track maps} of
Bestvina, Feighn and Handel
\cite{BFH}. In Part \ref{Part:BG2} we refined the train track technology in pursuit
of topological representatives of arbitrary automorphisms that share with the 
prototypes $f$ features that proved
crucial in Part \ref{Part:BG1}. We identified {\em{beads}}
as the basic units of an edge-path that  play the
role in the general setting that single edges (letters)
played in the case of positive automorphisms. The claim of beads to this 
 role was underscored by the
Beaded Decomposition Theorem.

With these technical innovations in hand, we now set about the task of
adapting the arguments of Part \ref{Part:BG1} to the general case, following
the proof  from Part \ref{Part:BG1} as closely as possible and providing the
(often fierce) technical details needed to translate each step into
the more general context provided by the topological representatives constructed in
Part \ref{Part:BG2}.  We shall not repeat the proofs of technical lemmas 
from Part \ref{Part:BG1} when the adaptation is obvious. Nor shall
we repeat our account of the intuition underlying our overall strategy
of proof and intermediate strategies at  key stages. 

Unfortunately, the adaptation to the general case is not entirely smooth. Thus
at times we are obliged to break from the narrative
that parallels Part \ref{Part:BG1} in order to deal with phenomena that do not
arise in the case of positive automorphisms ---
Section \ref{TrappedHNP}, for example. But we as far as
possible we have organised matters so that, having taken
account of the new phenomena, we can return to the main
narrative with the new phenomena controlled and
packaged into concise terminology. Thus, with
considerable technical exertions in our wake, we are
able to arrange  
matters so that  the
final stages of the proof of our Main Theorem
consist only of references to the corresponding
sections of Part \ref{Part:BG1} with a brief explanation of what changes, if any,
must be made in the general setting.

We have already noted that,
from the analysis of improved relative train tracks in Part \ref{Part:BG2}, it
emerged that beads are the correct analogue for the
role played by `letters' in the positive case. An important
manifestation of this is that the Main
Theorem can be reduced to a statement concerning the 
existence of a linear bound (in terms of $|\partial\Delta|$) on the
number of beads along the bottom of any corridor in a van Kampen diagram
$\Delta$ in the universal cover of the mapping tori that we consider.
In contrast to
the positive case, however, the existence of such a bound does not immediately
imply  the Main Theorem,  because there 
is no global bound
on the length of a bead.  

Nevertheless,
proving a bound on the number of beads is
by far the bulk of our work, occupying Sections 
\ref{FastBeads}--\ref{BonusSection}, which closely 
follow Sections \ref{NonConstantSection}--\ref{BonusScheme} (with
different numbering and modified structure).  In Section 
\ref{LongGepsandPepsSection} we explain how the bound on the
number of beads, together with the ideas from the {\em{Bonus Scheme}}
in Section \ref{BonusSection}, finally gives the Main Theorem.

In Section 
\ref{BracketingSection} we explain how to deduce estimates on the
geometry of van Kampen diagrams for {\em all}
mapping tori of free group automorphisms
from the specially-crafted  ones that we work
with during our main proof. The key estimate -- the linear
bound on the length
of $t$-corridors -- when reformulated algebraically,   yields the
{\em Bracketing Theorem} stated in the introduction. 

In Section \ref{BrinkSection} we explain how
our proof of the Main Theorem  allows one to reprove
the main result of 
\cite{BrinkDyn}.

We suggest that readers approach Part \ref{Part:BG3} as follows.
First, they must be familiar with the structure of the argument
in Part \ref{Part:BG1} and the vocabulary of beads in Part \ref{Part:BG2}.
This will enable them to skim smoothly through Sections
\ref{DiagSection}--\ref{Section:Iterate} of the current paper.
Next, they can gain an accurate overview of the proof of the Man Theorem 
 reading the introduction to each of Sections 
\ref{DiagSection}--\ref{LongGepsandPepsSection} together
with the titles of their subsections (and the introductions to
subsections when they exist).
There is then no alternative but to delve into the details
of the proof.

Section \ref{BracketingSection} can be read independently.
The argument in Section
\ref{BrinkSection} is easy to understand in outline, but the
proof appeals to detailed results from Sections \ref{FastBeads},
\ref{TeamsSection} and \ref{BonusSection}.

\section{The Structure of Diagrams} \label{DiagSection}

Associated to any finite group-presentation $\G=\AR$ one has the
standard
combinatorial 2-complex $K(\Cal A:\Cal R)$
with fundamental group $\G$ and directed edges labelled by the
$a\in\A$.  
There is a 1-1 correspondence between words in the letters $\A^{\pm 1}$ and
combinatorial loops in the 1-skeleton of  $K(\Cal A:\Cal R)$. Words such that
$w=1$ in $\G$ correspond to loops that are null-homotopic.
Van Kampen's Lemma explains the connection\footnote{For a complete account of
the equivalences in this subsection, see \cite{steer}.} between free equalities
demonstrating the membership $w\in\langle\!\langle \Cal R\rangle\!\rangle$
and combinatorial null-homotopies for the corresponding loops. 

Such a null-homotopy is given by a van Kampen diagram over  $\AR$, which is
a 1-connected, combinatorial planar 2-complex $\Delta$ in $\mathbb R^2$ with a 
basepoint; each oriented edge is labelled by a generator $a_i^{\pm 1}$ with $a_i\in\A$
and the boundary label on each face is some $r_j^{\pm 1}$ with $r_j\in\Cal R$
(read from a suitable basepoint).  There is a  unique label-preserving map
from the 1-skeleton of $\Delta$ to the 1-skeleton of the standard
2-complex $K(\Cal A: \Cal R)$, and this extends to a combinatorial map
$\Delta\to K(\Cal A: \Cal R)$.  

Van Kampen's Lemma implies that 
the number of faces in a least-area van Kampen diagram 
with boundary label $w$ is the least number $N$ of factors among free
equalities $w=\prod_{j=1}^Nu_jr_ju_j^{-1}$.  Thus the {\em Dehn
  function} of $\AR$ can be defined to be the minimal function
$\delta(n)$ such that every null-homotopic edge-loop of length at most $n$
in $K(\Cal A: \Cal R)$ is the  restriction to $\partial\Delta$ of a
combinatorial map $\Delta\to K(\Cal A: \Cal R)$ where $\Delta$ is a
$1$-connected, planar combinatorial 2-complex.  When described in this
manner, it is natural to call  the \index{Dehn function}Dehn function  the \index{isoperimetric 
inequality} {\em
combinatorial isoperimetric function} of $K(\Cal A: \Cal R)$; the
combinatorial isoperimetric function of an arbitrary compact
combinatorial 2-complex is defined in the same way.

There is a standard  diagrammatic argument for showing that the Dehn
functions of quasi-isometric groups are $\simeq$ equivalent --- see
\cite{Alonso}. In that argument, it is unimportant that the
complexes considered have only one vertex. Thus if $K$ is any compact
combinatorial 2-complex with  fundamental group $\G$, then the
combinatorial  isoperimetric function of $K$ is $\simeq$ equivalent to
the Dehn function of $\G$. We shall exploit the freedom stemming from
this equivalence.  Specifically, we shall prove the Main Theorem 
by establishing a quadratic upper bound on the combinatorial
isoperimetric function of a carefully-crafted 2-complex $M$ with
fundamental group $F\rtimes_{\phi^r}\mathbb Z$, where $r>0$. In other
words, we identify a constant $C>0$ such that  every null-homotopic
combinatorial loop of length at most $n$ in $M^{(1)}$ is the boundary of a
combinatorial map to $M$ from a 1-connected planar 2-complex  with  at
most $Cn^2$ 2-cells.  In fact, we prove something more refined than this (see 
Section \ref{Linear} below).

\begin{remark}\label{f.index}
Note that we are free to pass from $F\rtimes_\phi\mathbb Z$ to the
finite-index subgroup  $F \rtimes_{\phi^r}\mathbb Z$ because the
$\simeq$ class of the Dehn function of a group is an  invariant  of
commensurability.
\end{remark}

Henceforth we shall use the term \index{van Kampen diagram}
 {\em van Kampen diagram} to refer to the
domain of a combinatorial map to $M$ from a 1-connected planar
2-complex, with
oriented edges {\em labelled} by letters representing the oriented edges of
the target. (Note that this agrees with the standard terminology
in the special case $M=K(\Cal A:
\Cal R)$.) Such a diagram is said to be {\em least-area} if it has the
least number of 2-cells among all diagrams with the same boundary
label.

\subsection{The Mapping Torus}

Let $G$ be a compact graph and let $f:G \to G$ be a continuous
map that sends each
edge $e_i$ of $G$ to an immersed edge-path $u_i=\e_1\dots\e_m$ in $G$.  We
attach to each vertex $v\in G$ a new edge $t_v$ joining $v$ to
$f(v)$. We then attach one 2-cell to this augmented graph for each
edge $e_i$; the 2-cell is attached  along the edge path
$t_v^{-1}e_it_{v'}u_i^{-1}$, where $v$ and $v'$ are the initial
and terminal vertices of $e_i$ and where the inverse is taken in the path groupoid
(i.e. $u_i^{-1}$ is $u_i$ traversed backwards). The resulting 2-complex
is the \index{mapping torus}
{\em mapping torus} of $f$, which we shall denote $M(f)$.

In this part of the book we are primarily concerned with van Kampen diagrams over
$M(f)$, where $f$ is a homotopy equivalence representing a given
free-group automorphism $\phi$.  In this case $\pi_1(M(f)) \cong
\pi_1(G) \rtimes_{\phi} \Z$.  The $1$-cells in such a diagram
$\Delta_0$ are either labelled by some $t_u$ or by an edge $e \in G$.
We will refer to all of the edges $t_u$ as {\em $t$-edges} and, when
it does not cause confusion, denote them simply by $t$.  For the other
edges in $\Delta_0$, it is necessary to distinguish between the edge
and its {\em label} in $G$.

\begin{notation}[Labels $\widecheck\rho$]
If an edge $\e$ in a van Kampen diagram  over
$M(f)$ is labelled by an edge in $G$, then we write  $\widecheck{\e}$
to denote that label. More generally, if an edge-path $\rho$ in such a
diagram contains no $t$-edges,  we write $\widecheck{\rho}$
to denote the path in $G$ that labels $\rho$.
\end{notation}

\subsection{Time, folded $t$-corridors, singularities and bounded
  cancellation} \label{Folding}

Assume we are in the setting of the previous paragraph.
A {\em $t$-corridor} (more simply, \index{corridor}
{\em corridor}) is then defined exactly
as in Section \ref{sec:Time,corridors}, and we have the corresponding
notion of \index{time}{\em{time}} (which may be thought of as a map to $\mathbb R$
that is constant on non-$t$ edges, integer-valued on vertices, and
sends the endpoints of each $t$-edge to integers that differ by $1$).
As in  Subsections \ref{ss:Cond} and \ref{ss:Fold},
 we see that each least-area
diagram is the union of its corridors, and we may assume that the tops
of all corridors are {\em folded}. (In Subsection \ref{ss:foldBead}
we shall specify how this folding is to be done, but for the
results in this subsection it is not necessary to prescribe it.)

We write $\bot(S)$ and $\top(S)$ to denote
the {\em top} and {\em bottom} of a (folded) corridor,
respectively.  \index{singularity}{\em Singularities} are defined exactly as in
Part \ref{Part:BG1}.

We restrict our attention to least-area disc diagrams. The argument
used to prove Lemma \ref{l:onepoint} applies {\em verbatim}  in the
present setting  to prove:

\begin{lemma}If $S$ and $S'$ are distinct corridors in a least-area
diagram, then $\bot(S)\cap\bot(S')$ consists of at most one point.
\end{lemma}

Let $L$ be the maximum length of $f(E)$ for $E$ an edge in $G$.
As in Proposition \ref{SingularityProp} we have

\begin{proposition}[Bounded singularities] \label{SingularityProp:BG3}

\ 
 
\begin{enumerate} 
\item[1.] If the tops of two corridors in a  least-area 
diagram meet, then their intersection is a singularity. 
\item[2.] 
There exists a constant $B$ depending only 
on $\phi$  such that less than $B$ 2-cells 
hit each singularity in any  least-area diagram over $M(f)$. 
\item[3.]  
If $\Delta$ is a least-area diagram over $M(f)$, 
then there are less than $2|\partial \Delta|$ non-degenerate singularities 
in $\Delta$, and each has length at most $LB$.  
\end{enumerate} 
\end{proposition} 
\begin{proof} Except for one minor difficulty,  the proof from
Part \ref{Part:BG1} translates directly to the current setting. The minor
difficulty is that in the current context the map
$f$ is a homotopy equivalence rather than a group automorphism,  
 and $f^{-1}$ is not defined as a topological map. 
Thus, given a path $\rho$,
we need a canonical path $\sigma$ in $G$ such that $f_\#(\sigma)=\rho$,
where $ f_\#$ is tightening rel endpoints.  

Consider
$\wt{M(f)}$, the universal cover of $M(f)$.  Its 1-skeleton consists of
a collection of trees (copies of the universal cover of  $G$)   joined by
$t$-edges.  Consider a lift to $\wt{M(f)}$ of the unique edge-path
$\tau_0\rho \tau_1^{-1}$ such that the $\tau_i$ are $t$-edges.
Both endpoints of this lift lie in one
of the
 trees
$T\cong \tilde G$; define
$\wt{\sigma}$ to be the unique injective path 
which joins them in $T$, and define
$\sigma$ to be the image of $\wt \sigma$ in $M(f)$. 
\end{proof}

As in Lemma \ref{BCL}, the above result yields as a
special case (cf. \cite{Cooper} and \cite[Lemma 2.3.1,
  pp.527--528]{BFH}):

\index{Bounded Cancellation Lemma}
\begin{lemma}[Bounded Cancellation Lemma] \label{BCL:BG3} There is a
  constant $B$, depending only on $f$, so that if 
$I$ is an interval consisting of $|I|$ edges 
on the bottom of a (folded) corridor $S$ in a least-area diagram over
$M(f)$, and every edge of $I$ dies in $S$, then $|I| < B$. 
\end{lemma}

\subsection{Past, Future and Colour in Diagrams}\label{s:colour}

These concepts, for edges and 2-cells in van Kampen diagrams $\Delta$,
are defined exactly as in Section \ref{time}.
The {\em immediate past} (or 
\index{ancestor}{\em ancestor}) of an edge at the top of a
corridor in any diagram
 is the unique edge at the bottom of the corridor 
that lies in
the same 2-cell; the {\em entire past} of an edge is defined by taking the
transitive closure of the relation ``is the immediate past of''. The
past of a 2-cell is defined similarly. The \index{future!of an edge}{\em future} of an edge $e_0$ is
the set of edges that have $e_0$ in their past. The future of 2-cells
is defined similarly. The evolution of edges is
described by a graph $\Cal F$  whose vertices are the 1-cells $e$ of $\Delta$,
which has an edge connecting each $e$ to its  immediate ancestor.
Note that \index{family forest}$\mathcal F$ is a forest.
Its connected components  define  \index{colour!of an edge}{\em colours} in $\Delta$;
each edge not labelled $t$ is
assigned a unique colour, as is each 2-cell. Note that colours are in bijection with
  a subset of the edges of the boundary of the
diagram. 
The union of the 2-cells in a corridor $S$ that have colour
$\mu$ will be denoted $\mu(S)$.

As in Part \ref{Part:BG1},  simple separation arguments yield the following
observations.

\begin{lemma} Each $\mu(S)$ is
 connected and intersects each of $\top(S)$ and $\bot(S)$ in an interval.
\end{lemma}

\begin{lemma} [cf. Lemma \ref{perfect}]
Let $\e_1, \e_2$ and $\e_3$ be three (not necessarily adjacent) edges
that appear in order of increasing subscript as one reads from left to
right along the bottom of a corridor.  If the future of $\e_2$
contains an edge of $\partial \Delta$ or of a singularity, then no 
edge in the future of $\e_1$ can cancel with any edge in the future of
$\e_3$. 
\end{lemma}

Again following Part \ref{Part:BG1}, given a diagram $\Delta$ we define $\Cal Z$ to
be the set of pairs $(\mu,\mu')$ such that the coloured regions
$\mu(S)$ and $\mu'(S)$ are adjacent in some corridor $S$. The proof of
Lemma \ref{NoOfAdjacencies} establishes: 

\begin{lemma} \label{NoOfAdj}
$$|\Cal Z| \le 2\,|\partial\Delta| - 3.$$
\end{lemma}

\section{Adapting Diagrams to the Beaded Decomposition}

We refer the reader to Part \ref{Part:BG2} for the definitions and results
which we require here about improved relative train track maps,
nibbled futures, monochromatic paths,
hard splittings and the language of {\em{beads}}
 --- including
$(J,f)$-atoms, {\gep}s and {\pep}s and what it means for a path
to be $(J,f)$-beaded.
We shall proceed under the assumption that the reader is
familiar with each of these terms, and work axiomatically with
the following outputs from Part \ref{Part:BG2}.

\begin{theorem}[Beaded Decomposition Theorem, Part \ref{Part:BG2}] \label{BDT}
For every $\phi\in{\rm{Out}}(F_r)$,
there exist positive integers $k$, $r$ and $J$ such that $\phi^k$ has an improved
relative train-track representative $f_0:G\to G$ 
with the property that every $(f_0)_{\#}^r$-monochromatic
path in $G$ is $(J,f_0)$-beaded.
\end{theorem}

We remind the reader that
beads are either monochromatic paths (in case they are \atom s) or else \gep s
or \pep s (which may be monochromatic, but do not have to be). 
Thus, by the above theorem and Proposition \ref{nfgeptopep}, any nibbled
future of a $(J,f_0)$-bead is $(J,f_0)$-beaded. Any hard splitting of an edge-path is 
inherited by its (nibbled) futures, by definition. And if one refines a hard splitting
by decomposing the factors in a hard splitting, the result is again
a hard splitting (Lemma \ref{HardSplitProperties}). Thus we have:

\begin{coroll} [Theorem \ref{t:nfbead}] 
Let $f = (f_0)_{\#}^r$ be as in the Beaded Decomposition Theorem above. If an edge-path $\sigma$ in $G$
is $(J,f_0)$-beaded, then any $f$-nibbled future of $\sigma$ is $(J,f_0)$-beaded.
In particular, $f_{\#}(\sigma)$ is also $(J,f_0)$-beaded.
\end{coroll}

\begin{remark} An important point to recall from Part \ref{Part:BG2} is that the
decomposition of an edge-path into $(J,f_0)$-beads is canonical.
\end{remark}

The value of the constant $J$ in the Beaded Decomposition Theorem will be
of no importance in what follows, so we drop it from the terminology.
Similarly, we will fix the map $f_0$.  Once we have passed to the power $f = (f_0)_{\#}^r$, the above results remain true when $f$ is replaced by an iterate.  Therefore, we
refer simply to  \index{bead}``beads" and ``beaded paths".
\medskip

\subsection{Refolding corridors according to the Beaded Decomposition}\label{ss:foldBead}

\ 
\medskip

{\em{Henceforth\footnote{There are exceptions
to this in Theorem \ref{BoundS:BG3}, Section \ref{BracketingSection} and Section 
\ref{BrinkSection}}, we consider only diagrams over the \index{mapping torus}
mapping torus of $M(f)$, where
$f$ is an iterate of $(f_0)_{\#}^r$ as in the Beaded Decomposition Theorem.} In Section \ref{Section:Iterate}, we will fix the map $f$ once and for all.}

\medskip

We return to the matter of how best to fold the tops of corridors in
least area diagrams over $M(f)$. Given an arbitrary least-area
diagram, we refold the tops
of corridors in order of increasing time. The process begins with
edges at the minimal time on the boundary of the diagram, where
there is no folding to be done provided the boundary label is
reduced.

Focussing on a
particular corridor $S$, our folding up to  ${\rm{time}}(S)$
defines the histories of all edges up to this time and hence
assigns colours to the edges on $\bot(S)$, decomposing
it as a concatenation of monochromatic paths, one for
each of the colours $\mu(S)$. Theorem \ref{BDT} decomposes
each of these labels as a hard splitting of beads $\sigma_i$.
The hardness of the splitting means that after tightening the
$f(\sigma_i)$, their concatenation will be a tightening of 
$f_{\#}(\check\mu(S))$. We insist that the first step in 
the tightening of the naive top of $S$, is that determined by
the tightening of labels just described: i.e.~we first tighten
beads {\em{within}} colours, each according to a left-to-right
convention (which labels inherit
from the orientation of the corridors within
the diagram). Then, as a second step, we tighten (again with
a left-to-right convention) the concatenation of the tightened
images of the colours.  A diagram which is folded according to
these conventions will be called {\em well-folded}.
\index{well-folded diagram}

The key point of this convention is that the hard splitting 
of the label on each colour is carried into the future --- of course the futures of
the original beads may split into a concatenation of several
beads, and some beads at the ends of each colour may be cancelled by interaction 
with neighbouring colours, but {\em{each bead (more precisely\footnote{we shall
generally drop this cumbersome distinction in the sequel},
bead-labelled arc) in the beaded decomposition
of each coloured interval on $\top(S)$ is contained into the future of a 
unique bead-labelled arc of the same colour on $\bot(S)$.}} Thus 
$\top(S)$ is a concatenation of beads, each with a definite colour,
where neighbouring beads are separated by a hard splitting if
they are of the same colour but perhaps not if they are of a 
different colour.
 (It also becomes sensible to discuss
the future of a bead in a [well-folded] diagram.)

\medskip
{\em{
We henceforth suppose (usually without comment) that our diagram 
has been refolded according to this convention.}}
\smallskip
 
\begin{definition}\label{BeadLength} [cf.~Definition \ref{beadNorm}]
The {\em bead length} of $[S]_\beta$, of a corridor $S$ in a well-folded
diagram is the number of beads along $\bot(S)$. \index{corridor!bead-length of}
\end{definition}

\begin{remark} 
It is important to note that the decomposition of $\bot(S)$
and $\top(S)$ into coloured intervals is
{\em not} a hard splitting in general.  Indeed it
is the analysis of 
the cancellation between these intervals as one
flows $S$ forwards in time  that forms the meat of
this part of the book.
\end{remark}

\subsection{Abstract Futures of Beads}\label{stackDiags}
Given an edge-path $\rho$ in $G$, expressed as a concatenation
of monochromatic edge-paths $\rho=\rho_1\dots\rho_m$, consider
the van Kampen diagram $\Delta(l,\rho)$ with boundary
label equal to $t^{-l}\rho t^l \overline{f_{\#}^l(\rho)}$; 
this is a simple stack of corridors. The above convention dictates
how we should fold the corridors of $\Delta$ and
determines the future at each time up to $l$ for each  bead in the
beaded decompositions of the $\rho_i$. 

\index{stack diagram}
\index{future!of a bead}
We define the (full) {\em{abstract future of a bead in $\rho$}}   to be
(the label on) its future in $\Delta(l,\rho)$.

\section{Linear Bounds on the Length of Corridors} \label{Linear}

In any least-area diagram, 
each corridor has at least two edges on the boundary, namely its $t$-edges.
The {\em length} of a  \index{corridor!length of}
corridor $S$ is defined to be the number of 2-cells that it 
contains. The area of a least-area diagram is the sum of the lengths of
its corridors, and therefore our Main Theorem 
 is an immediate consequence of:

\begin{theorem} \label{BoundS:BG3}
Let $\phi$ be an automorphism of a finitely generated free group
and let $f$ be a topological representative for a positive power of $\phi$.  
There is a 
constant $K$, depending only on $f$, so that each corridor in a least-area 
diagram $\Delta$ over $M(f)$ has length at most $K\n$.
\end{theorem}

Note that the Main Theorem actually depends only on establishing 
Theorem \ref{BoundS:BG3}
 for a single  topological representative $f^k$ of a suitable
 power of our given free group automorphism
$\phi$; in the next section we shall articulate what that suitable
power is.
The bulk of this part of the book will then be devoted to 
 proving the existence of the constant $K$ for this
particular $f^k$. (In Section \ref{BracketingSection} we shall 
deduce Theorem \ref{BoundS:BG3} from this special case.)

Having restricted attention to a particular $f^k$, we may
further restrict our attention to diagrams that are well-folded
in the sense of Subsection \ref{ss:foldBead}, since
refolding the corridors of an arbitrary a diagram does
not change the configuration of corridors or their length.
In a well-folded diagram,
the top of each corridor $S$ is a concatenation of beads, and
the vast majority of our work (up to and including Section \ref{BonusSection}) goes into proving  the following result.

\begin{theorem}\label{noName} If $f$ and $k$ are as  above, then
there is a constant $K_1$ such that  all
corridors $S$ in  well-folded, least-area diagrams $\Delta$ over  $M(f_{\#}^k)$,
have bead length $[S]_\beta\le K_1\ |\partial\Delta|$. 
\end{theorem}

The linear bound on the length of $S$ that we require for
Theorem \ref{BoundS:BG3} does not follow directly from this estimate because
there is no uniform bound on the length of certain beads, namely
\gep s and \pep s. However, we shall see in Section \ref{LongGepsandPepsSection} 
that the
ideas developed in Part \ref{Part:BG1} to implement the Bonus Scheme adapt
 to the current setting to provide the following estimate:

\begin{proposition}\label{longGEPS} There are constants $J$ and $K_2$,
depending only on $f$, such that the beads $\beta$ on $\bot(S)$ of length
greater than $J$ satisfy
$$
\sum_{\beta} |\beta| \le K_2\, |\partial\Delta|.
$$
\end{proposition} 
The constant $J$ in the above statement is the one from Theorem \ref{BDT}.

\section{Replacing $f$ by a Suitable Iterate} \label{Section:Iterate}

In order to establish the  bound on the 
length of corridors required to prove Theorem \ref{BoundS:BG3}, 
we must analyse how 
corridors grow as they flow into the future and 
assess what cancellation can take place to inhibit this 
growth.  This is much more difficult than in Part \ref{Part:BG1} because
now we must cope with the cancellation that takes place
within colours. But in common with our approach in Part \ref{Part:BG1},
we can appeal to Remark \ref{f.index}
 repeatedly in order to replace our
topological representative $f$ by some iterate of $f$ that
affords a more stable situation in which 
cancellation phenomena are more amenable to analysis.

In the present setting, we have to be a little careful about specifying
what we mean by ``an iterate", because we wish to consider only
topological  representatives  whose restriction to
each edge is an immersion, and this property is not inherited by
powers of the map. To avoid this problem, we deem the phrase\footnote{and obvious
variations on it}
{\em replacing $f$ by an iterate}, to mean that for fixed
$k\in \mathbb N$, we pass from consideration of $f : G \to G$
to consideration of  the map $f_{\#}^k : G \to G$ that sends
each edge $E$ in $G$ to  the tight edge-path $f_{\#}^k(E)$ that is 
homotopic rel endpoints to $f^k(E)$.

When we replace $f$ by $f_{\#}^k$, we leave behind the   
mapping torus $M(f)$ and consider instead $M(f_{\#}^k)$, which although 
homotopic to a $k$-sheeted covering of  $M(f)$ is distinct
from it. 

A corridor in a van Kampen diagram over $M(f_{\#}^k)$ can be divided into
a stack of $k$ corridors in order to yield a van Kampen diagram over 
$M(f)$. This observation will play little role in our arguments, but it highlights
one reason for hoping to simplify diagrams by passing to an iterate of $f$: the van Kampen
diagrams over $M(f_{\#}^k)$ are a proper subset (after subdivision\footnote{the obvious subdivision of a diagram $\Delta$ is called the {\em{$k$-refinement}}}
of $\Delta$) of 
the diagrams over $M(f)$; in the diagrams of this subset, corridors  
flow unhindered for at least $k$ steps in time.

\subsection{Finding the desired iterate} \label{Iterate}

We have already passed to a large iterate in order to obtain
the Beaded Decomposition Theorem. In the present subsection we
pass to further iterates in order to control the behaviour of
the images of beads.

Before settling on a specific $f$ for the remainder of the
paper, we must remove an irritating ambiguity concerning the ordering
of strata in the filtration associated to the train track structure.
This is required in order to render the choices in Section 
\ref{s:PrefFuture} coherent.

\begin{definition} \label{interchangeable}\index{strata!interchangeable}
Suppose that $f : G \to G$ is an improved relative train track map, and
that $H_i, H_j$ are strata for $f$.  We say that $H_i$ and $H_j$ are
{\em interchangeable} if one can reorder the strata, so that one
still has an improved relative train track structure, but
the order of $H_i$ and $H_j$ is reversed.
\end{definition}

If $H_i$ and $H_j$ are interchangeable, and $i > j$, then no iterate
of any edge in $H_i$ crosses an edge in $H_j$ (and neither do the
iterates of any edges occurring in the iterated images of edges in $H_i$).

\begin{convention} \label{convention:interchangeable}
We suppose that for any improved relative train track map that
we consider, if $H_i$ and $H_j$ are interchangeable strata so
that $H_i$ is an exponential stratum and $H_j$ is a parabolic
stratum then $i > j$.  

We further assume that if $H_i
= \{ E_i \}$ and $H_j = \{ E_j \}$ are interchangeable
 parabolic strata  and $n\mapsto |f^n(E_i)|$
grows exponentially while $n\mapsto |f^n(E_j)|$
grows polynomially, then $i > j$.  And if both these functions
grow polynomially, then the degree of polynomial growth of
the former is at least as great as the latter. 
\end{convention}

In the following lemma, $\omega$ is the number of strata
in the train track structure for $f$. Also recall that an
edge $\e$ in a path $\sigma$ is said to be \index{edge!displayed}
 {\em displayed} if there is a hard
splitting $\sigma = \sigma_1\odot \e \odot \sigma_2$.  The definition of
a \index{path!displayed}
{\em displayed sub edge-path} is entirely analogous, and will be used later.

\begin{lemma} \label{fofatom:BG3}
One can replace $f$ by an iterate to ensure  that
 if $\rho$ is any atom then
either the beads of  $f_{\#}^\omega (\rho)$ are
Nielsen paths and \gep s only, or else there is a displayed edge
$\epsilon $ in $ f_{\#}^\omega (\rho)$ so that 
\begin{enumerate}
\item $\e$ is 
 of highest weight amongst all displayed
edges in all $f_{\#}^k(\rho)$, for $k \ge 1$, and 
\item the growth of $n\mapsto |f^n_{\#}(\epsilon)|$ is at least 
as large as that of
any displayed edge in any $f_{\#}^k(\rho)$.
\end{enumerate}
\end{lemma}

\begin{proof} Lemma \ref{Displayededge} contains
all but statement (2), whose validity is assured by
 Convention 
\ref{convention:interchangeable}.
\end{proof}

Our next two results capture the \index{end stability}
 {\em{end
stability}} that Proposition \ref{power} provided in the
case of positive automorphisms. This is the first stage in our
analysis at which we encounter an awkward point that does not
arise in Part \ref{Part:BG1}, namely   there may exist beads (more
specifically atoms) $\rho$ such that $f_{\#}(\rho)$ is a single vertex.

\begin{definition} \label{vanishing} A \index{atom!vanishing}
{\em vanishing bead (atom)} \index{bead!vanishing}
$\rho$ is one with $f_{\#}(\rho)$  a single vertex. 
\end{definition}

\begin{lemma} \label{EndStab1}
There exists a constant $k_0$, depending only on $f$ so that the
map $f_0 = f_{\#}^{k_0}$ satisfies the following properties.  Let
$\rho$ be a non-vanishing bead,
let $i \in \{ 1 ,\ldots , \omega \}$, and let  $\sigma_i$ be the 
leftmost bead in $(f_0)_{\#}(\rho)$ of weight at least $i$.
\begin{enumerate}
\item If $\sigma_i$ is not a \gep\ or a \pep\ then the leftmost bead
of weight at least $i$ in $(f_0)_{\#}^j(\rho)$ is the same for all $j \ge 1$.
Furthermore, in this case $\sigma_i$ is a single (displayed) edge or
a Nielsen bead.
\item If $\sigma_i$ is a \gep\ or a \pep\ then the leftmost bead of weight
at least $i$ in $(f_0)_{\#}^j(\rho)$ is contained in the (abstract)
 future of $\sigma_i$
for all $j \ge 1$.
\end{enumerate}
\end{lemma}
\begin{proof}
If $\sigma$ is a bead then
all iterated images of $\sigma$ are beaded paths, 
and a simple finiteness
argument shows that there is a
bound on the number of beads which are not
\gep s or \pep s.
\end{proof}

An entirely similar argument applies to rightmost beads, of course.
In order to deal with the different types of beads, we also
need the following variant.

\begin{lemma} \label{EndStab2}
There exists a constant $k_1$, depending only on $f$, so that
the map $f_1 = f_{\#}^{k_1}$ satisfies the following properties.   Let
$\rho$ be a non-vanishing bead 
and let $\sigma$ be the leftmost bead in 
$(f_1)^j_{\#}(\rho)$ which is not a Nielsen bead.
\begin{enumerate}
\item  If $\sigma$ is not a \gep \ or a \pep\ then for all $j \ge 1$
the leftmost bead in $(f_1)^j_{\#}(\rho)$ which is not a Nielsen bead
is $\sigma$.  Furthermore, in this case $\sigma$ is a (displayed)
edge.
\item  If $\sigma$ is a \gep\ or a \pep\ then for all $j \ge 1$
the leftmost bead in $(f_1)_{\#}(\rho)$ which is not a Nielsen bead
is in the future of $\sigma$.
\end{enumerate}
\end{lemma}

We are finally in a position to articulate all of the properties
that we want to arrange for $f$ by replacing it with an iterate.

\begin{proposition} \label{FinalPower} There is a constant $D_2$
that depends only on $f$, so that if we replace $f$ by 
 $f_{\#}^{D_2}$ then,
\begin{enumerate}
\item\label{ExpLong}  the conclusion of Lemma \ref{Power1} holds with $k_1 = 1$: in particular, if $\e$ is an exponential edge of weight $i$, then $f(\e)$
is longer than the unique indivisible Nielsen path of weight $i$
(if it exists);
\item the conclusion of Theorem \ref{ColourCancelMain} holds with $D_1 = 1$;
\item the conclusion of Lemma \ref{fofatom:BG3} holds;
\item the conclusions of Lemmas \ref{EndStab1} and \ref{EndStab2} hold; and
\item\label{DisplayExp} if $\rho$ is a bead  then
$f_{\#}(\rho)$ contains at least three displayed copies of any
exponential edge that is displayed 
in any  $f_{\#}^j(\rho),\ j\ge 1$. Moreover, the leftmost (and rightmost)
such displayed edge $\e$ is contained in a displayed path
of the form $f(\e)$.
\end{enumerate}
\end{proposition}

\medskip
\noindent{\bf{Power Decree:}}{\em{
For the remainder of the paper, we will assume
that $f : G \to G$ is an improved
relative train track map that satisfies the properties 
in Proposition \ref{FinalPower}. We shall also
operate under Convention \ref{convention:interchangeable}.
}}

{\em Let $L$ be the maximal length of $f(E)$, for
edges $E \in G$.}

\section{Preferred Futures of Beads} \label{s:PrefFuture}

The reader who is comparing our progress to Part \ref{Part:BG1}
will find that we are now in the position that we were at the
start of Section \ref{PrefFutSec}. Thus we now want to define
the preferred future of a bead $\rho$ (in three senses\footnote{in $f_{\#}(\rho)$,
in a diagram, and in a concatenation of beaded paths})
 and then begin a study of  fast beads.

Unfortunately, the definition of the preferred future of a bead 
in a diagram is much more 
cumbersome than the analogue in Part \ref{Part:BG1}.

\subsection{Abstract Preferred Futures and Growth}
 
First we note that if beads (or more generally edge paths
in $G$) are ever going to vanish in the sense
of Definition \ref{vanishing}, then they do so immediately.

\begin{lemma}
If $\sigma$ is an edge path in $G$ and
$f_{\#}^k(\sigma)$ is a vertex for some $k\ge 1$,
then $f_{\#}(\sigma)$ is already a vertex.
\end{lemma}
\begin{proof}
For all vertices $v \in G$, $f(v)$ is a fixed point of $f$.  Therefore,
the endpoints of $f_{\#}^j(\sigma)$ are the same for all $j \ge 1$.
If $f^k_{\#}(\sigma)$ is a point, then the endpoints of $f_{\#}^k(\sigma)$
are equal, hence the tight path $f_{\#}(\sigma)$ is  a loop. Since $f$ is a 
homotopy equivalence, this loop must be trivial.
\end{proof}

\begin{definition} [Abstract preferred futures]\label{PrefFuture}
The (immediate) {\em preferred future} of a \index{future!preferred}
non-vanishing bead $\sigma$ is a particular
 bead in the beaded decomposition of $f_{\#}(\sigma)$, as
defined below. The {\em $k$-step preferred future} is then
defined by an obvious recursion.
\begin{enumerate}
\item If $\sigma$ is a \gep\  then $f_{\#}(\sigma)$ is
also a \gep, and we define the preferred future of $\sigma$
to be $f_{\#}(\sigma)$.
\item If $\sigma$ is a \pep\ then either $\sigma$ or $\overline{\sigma}$
has the form $\sigma = E \overline{\tau}^k\nu\gamma$.
If it is $\sigma$, then
by Corollary \ref{cor:nfpep}, $f_{\#}(\sigma)$ is
either of the form $\sigma' \odot \xi$, where $\sigma'$
is a \pep\ (which has the same weight as $\sigma$), 
or else of the form $E \odot \xi$, where $E$
has the same weight as $\sigma$ and is the unique highest
weight edge in $f_{\#}(\sigma)$.  In the first case, the preferred
future of $\sigma$ is $\sigma'$.  In the second case, the 
preferred future of $\sigma$ is $E$. The preferred future of a \pep \ $\sigma$
where $\overline{\sigma}$ has the above form is
defined in an entirely analogous way.
\item If $\sigma$ is a Nielsen path then 
the preferred future of
$\sigma$ is $f_{\#}(\sigma) = \sigma$. 
\item Finally, we consider  a non-vanishing
\atom\ $\sigma$.

\noindent
(a) If the beaded decomposition of
 $f_{\#}(\sigma)$  consists entirely of Nielsen paths
and \gep s, 
 then we fix a highest weight \gep\ to
be the preferred future of $\sigma$;  
otherwise, we fix a highest weight
Nielsen path.

\noindent
(b) If not, then let $\e$ be the edge described in Lemma \ref{fofatom:BG3},
 fix a displayed occurrence of $\e$ in $f_{\#}(\sigma)$ (in case $\e$
 is exponential, choose a displayed occurrence that is neither
leftmost nor rightmost\footnote{this exists by Proposition \ref{FinalPower}}) 
and define this to be the preferred future of
$\e$.
\end{enumerate}
\end{definition}

\begin{remark}
Suppose that $\e$ is an edge in $G$, considered as a bead, and suppose
that $\e$ is not contained in a zero-stratum.  Then $\e$ has a preferred
future, which is an edge contained in the same stratum as $\e$.  We always
assume that the preferred future of $\e$ is a (fixed) occurrence of $\e$ in 
$f_{\#}(\e)$ which satisfies the requirements of the above definition.  This
situation is very close in spirit to the definition of preferred future in Part \ref{Part:BG1}.
\end{remark}

We now divide the beads into classes according to the growth
of the paths $f_{\#}^k(\sigma)$, $k = 1,2,\ldots$.  Specifically,
we define left-fast and left-slow beads in accordance with
Subsection \ref{ss:speed}.

\begin{definition} [Left-fast beads] \label{LeftFast}\index{bead!fast, slow}
\gep s and Nielsen paths are {\em left-slow}.

Suppose that $\alpha$ is an \atom\ or a \pep.  Then $\alpha$ 
is {\em left-fast} if the distance between the left end of 
$f_{\#}^k(\alpha)$ and the left end of the 
preferred future of $\alpha$ in 
$f_{\#}^k(\alpha)$ grows at least quadratically with $k$, and {\em left-slow}
otherwise.
\end{definition}
Note that if a \pep\ $\sigma$ is left-fast then it is $\overline{\sigma}$
which it is of the form $E \overline{\tau}^k \nu\gamma$.

\begin{remark}
We only care that fast growth be super-linear, but it happens that
this is the same as being at least quadratic  (cf.~\cite{BG-Growth}).
\end{remark}

The concepts of {\em right-fast} and {\em right-slow} beads are
entirely analogous.

\subsection{Preferred future in diagrams} \label{PrefFutureDiag}

\index{future!preferred}
In this subsection we define the notion of `preferred futures' within
van Kampen diagrams.  We also
define `biting' and `consumption', which are the analogues in this
paper of `consumption' from Section \ref{PrefFutSec}.

The folding convention of Subsection \ref{Folding} expresses
$\bot(S)$ as the concatenation of coloured paths $\mu(S)$, each labelled
by a monochromatic path in $G$. The Beaded Decomposition Theorem
gives us a hard splitting into beads
\[      \widecheck{\mu(S)} = \widecheck{\beta}_1 \odot 
\widecheck{\beta}_2 
\odot \cdots \odot \widecheck{\beta}_{m_\mu}    ,      \]
and it is convenient to refer to the sub-paths 
$\beta_i \subseteq \bot(S)$ carrying the labels $\widecheck\beta_i$
as beads, as we did in Subsection \ref{Folding}.

If $\mu_1, \ldots , \mu_k$ are the colours appearing in $S$, in order,
then the label on $\top(S)$ is obtained by tightening
\[ f_{\#}(\widecheck{\mu_1(S)}) \cdots f_{\#}(\widecheck{\mu_k(S)}).    \]
The path $f_{\#}(\widecheck{\mu_1(S)}) \cdots 
f_{\#}(\widecheck{\mu_k(S)})$ \index{future!semi-naive}
is called the {\em semi-naive future} of $S$.

We have adopted a left-to-right convention to remove any
ambiguity in how one tightens the semi-naive
future to obtain the label of $\top(S)$.

We previously defined the (immediate) future of a 
bead $\beta \subset \bot(S)$ to consist of those edges of $\top(S)$ 
whose immediate past lies  in $\beta$. Since it is integral to what
we shall do now, we re-emphasize:

\begin{lemma}
The immediate future of a bead $\beta\subset\bot(S)$ is a (possibly
empty) interval equipped with a hard-splitting into beads.
\end{lemma}

If $\rho$ is the immediate
future of $\beta$, then $\rho$ is also an interval in the semi-naive
future of $S$, and hence its label
$\widecheck{\rho}$ is a specific sub-path of 
$f_{\#}(\widecheck{\beta})$.   [Note that one has more than the
path $\widecheck{\rho}$ here, one also has its position within
$f_{\#}(\widecheck{\beta})$; thus, for example, 
we would distinguish between the two visible copies
of $\widecheck{\rho}$ in  $f_{\#}(\widecheck{\beta})=
\widecheck{\rho}\sigma\widecheck{\rho}$.]

\begin{definition}[Preferred and tenuous futures in $\Delta$] \label{PreferredFutureDef} 
Consider a bead $\beta\subset\mu(S)\subset\bot(S)$
 in $\Delta$ \index{future!preferred} \index{future!tenuous}
whose immediate future $\rho\subset\top(S)$ determines the
subpath $\widecheck{\rho_0}$ of $\widecheck{\beta}$ in $G$.

If the (abstract) preferred future $\widecheck\beta_+$ of  $\widecheck\beta$, as
defined in Definition \ref{PrefFuture}, is entirely contained in  
$\widecheck{\rho_0}$, then the corresponding sub-path
$\beta_+$ of $\rho$ 
is the {\em preferred future} of $\beta$. 

If $\widecheck{\rho_0}$ does not contain $\widecheck\beta_+$,
 then $\beta$ does not have a preferred future. In this situation
we say that the future of $\beta$ is {\em{tenuous}}.
\end{definition}

\begin{remark} \label{rem:Pref}
Note that, if it exists, the preferred future of a bead $\beta\subset\mu(S)$ is 
a bead in the beaded decomposition of both $\rho$ and the $\mu$-coloured
interval of $\top(S)$.

Also, if a bead happens to be a single edge  $\e$ whose label is not contained in 
a zero stratum, the preferred future is a single (displayed) edge, with 
the same label as $\e$.
\end{remark}

\begin{definition} [Biting and consumption] If the future of a bead
$\beta\subset\bot(S)$ is tenuous,  \index{biting}
we say that $\beta$ is {\em bitten in $S$}. If, in the notation of
(\ref{PreferredFutureDef}),  {\em no edge} of the preferred future of 
$\widecheck{\beta}$ appears in $\widecheck{\rho}$, then we say that
$\beta$ is {\em consumed} in $S$.
\index{bead!bitten} \index{bead!consumed}
\end{definition}

\begin{remark}
The above definition says in particular that any bead whose label
is a vanishing atom is consumed.
\end{remark}

Let $\beta'\subset\bot(S)$ be a bead whose label is non-vanishing.
If $\beta'$ is bitten in $S$, 
there is a specific edge $\e$ in the semi-naive future of $S$ 
that, during the tightening process, is the first to  cancel with
an edge $\e'$ in the interval labelled by the preferred future of 
$\widecheck{\beta'}$.  The edge $\e$ is in the immediate future
of a bead $\beta$, necessarily of a different colour than $\beta'$.

\begin{definition} \label{d:bite}
In the above situation,
we say that $\beta$ {\em bites $\beta'$ from the left} if $\beta$
lies to the left of $\beta'$ in $S$, and that $\beta$ {\em bites
$\beta'$ from the right} if $\beta$ lies to the right of $\beta'$ in $S$.
We say that the edges $\e$ and $\e'$ discussed above {\em exhibit} the biting.
\end{definition}

The above concepts
of biting and consumption replace the single, simpler, notion
of consumption from Section \ref{PrefFutSec}: there, since
the preferred future was a single edge,  if it was bitten 
it was consumed.   In Part \ref{Part:BG1}, a frequently used concept was
for an edge to be `eventually consumed'.  In this part of the book, we need 
the following replacement:

\begin{definition} \label{d:eventbite}
Suppose that $\rho_1 \subset \mu_1(S)$ and $\rho_2 \subset
\mu_2(S)$ are beads in $\bot(S)$.  We say that $\rho_1$ is
{\em eventually bitten} by $\rho_2$ if there is a corridor $S'$
which contains a preferred future $\beta_1$ of $\rho_1$ and a bead $\beta_2$
in the future of $\rho_2$ so that $\beta_2$ bites $\beta_1$ in $S'$.
\end{definition}

With these definitions in hand, we have the following, which is an appropriate
replacement for \ref{C_0}

\begin{lemma} [cf.  Lemma \ref{C_0}]  \label{C0Lemma}
There exists a constant $C_0$ with the following property: if $\rho$
is a bead such that $f_{\#}(\rho)$ contains a left-fast displayed
edge $E$ and if $UV\rho$ is a (tight) path with $V\rho = V \odot \rho$ 
 and $|V | \geq C_0$ then for all $j \geq 1$ the preferred
future of $E$ is not bitten when $f^j(UV\rho)$ is tightened.
Moreover, $|f_{\#}^j(UV\rho)| \to \infty$ as $j \to \infty$.
\end{lemma}
\begin{proof}
We first prove the result in the special case that $V\rho$ is a \nf \ of
a left-fast edge $E_1$, where $\rho$ is the preferred future of $E_1$.
In other words, we will prove the existence of a constant $C_0'$ so that
if $|V| \ge C_0'$ then the statement of the lemma holds for the particular
path $U V\rho$.
(We will later reduce to this special case.)

Note that $V$ and $V\rho$ are monochromatic paths, and thus admit a
beaded decomposition.  Suppose first that $V$ does not contain any
beads of length greater than $J$.  In this case, the proof is entirely
parallel to that of Lemma \ref{C_0}, where we count using the number 
of non-vanishing beads rather than the number of edges.

In case $V$ contains long \gep s or long \pep s, we note that the cancellation
by $U$ on the left, and possibly by one of the edges in the \gep \ or \pep \ on
the right can only decrease the length of a \gep \ or \pep \ by at most $2B$
at each iteration.  Thus it is straightforward to include long \gep s and \pep s
into the above calculation.  We now turn to the general case.

Suppose that $V$ is an arbitrary path so that $V\rho = V \odot \rho$.
Then $V$ can shrink of its own accord (it needn't be beaded), and can be cancelled
by the future of $U$.  However,  there is certainly a constant $C_0$ so that
if $|V| \ge C_0$ then by the time this shrinking of $V$ combined with cancelling
by the future of $U$ can have reduced $V$ to the empty path, the future of the edge $E$ has
at least $C_0'$ edges to the left of its preferred future.  
We are then in the
special case that we dealt with first.
\end{proof}

The following two lemmas are proved in an entirely similar manner to Lemma \ref{NoOldCanc}. Recall that displayed edges are particular types of beads,
and the (abstract) preferred futures of beads were defined in
Definition \ref{PrefFuture}.  Recall from Remark \ref{rem:Pref} that the 
preferred future of a displayed edge whose label is not contained in a 
zero stratum is a single displayed edge.

\begin{lemma} \label{OldLemma} Let $\chi_1\sigma\chi_2$
be a tight path in $G$.
Suppose that $\chi_1$ and $\chi_2$ are monochromatic 
and that, for $i = 1,2$, 
the edge $E_i$ is displayed in $\chi_i$ and that $E_i$ is not
in a zero stratum.  Suppose that $\sigma$ is a concatenation of beaded paths.
Then the preferred futures of $E_1$ and $E_2$ cannot cancel
each other in any tightening of $f_{\#}(\chi_1)f_{\#}(\sigma)f_{\#}(\chi_2)$.

Suppose that $S$ is a corridor in a well-folded diagram, and that $\mu_1(S)$
and $\mu_2(S)$ are non-empty paths in $\bot(S)$, where $\mu_1$ and $\mu_2$ are
colours.  Suppose further that for $i = 1,2$ there is a  displayed edge $\e_i$
such that $\widecheck{\e_i}$ is not contained in a zero stratum.  Then the edges
in the semi-naive future of $S$ corresponding to the preferred futures of
$\e_1$ and $\e_2$ do not cancel each other when folding the semi-naive future
of $\bot(S)$ to form $\top(S)$.
\end{lemma}

\begin{lemma} \label{ParabolicOldEdges}
Let $S$ be a corridor and suppose that $\e_1$ and $\e_2$ are edges in
$\bot(S)$ whose labels lie in parabolic strata. In the naive future of
each $\e_i$ (that is, before even the beads have been tightened),
there is a unique edge $\e_i'$ with the same label as
$\e_i$. At no stage during the tightening of $\top(S)$
can $\e_1'$ cancel with $\e_2'$.
\end{lemma}

\begin{corollary} \label{ParaCancelledbyHigher}
A displayed edge in any coloured interval $\mu(S)$ which is 
 labelled by a parabolic edge $\widecheck E_i \in H_i$ can
only be consumed by an edge whose label is in 
$\overline{G \smallsetminus G_i}$.
\end{corollary}

\subsection{Abstract paths, futures and biting}
\label{AbstractFolding}

In many of the arguments in later sections, we wish to work with
concatenations of beaded paths in $G$
 rather than sides of  corridors in diagrams.  This is done as in 
Subsection \ref{stackDiags} by associating to such a path
$\rho = \rho_1\dots\rho_m$, with the $\rho_i$ beaded, 
the van Kampen diagram $\Delta(l,\rho)$ with boundary
label $t^{-l}\rho t^l \overline{f_{\#}^l(\rho)}$. But we modify
the usual definition of colour by defining the colours on the bottom
of the first (earliest) corridor not to be single edges but
rather to be intervals labelled $\rho_i$.
We then use the 
definitions of the previous subsection (biting, preferred future {\em{etc.}})
to define the associated
concepts {\em{for beads in $\rho$}}.

We emphasize, $\rho$ itself need not
beaded; only the $\rho_i$ are. We also emphasize that edges do not
have preferred futures, only beads do.

However, some beads are single, displayed edges, and when considered
as beads they do have a preferred future.

\section{Counting Fast Beads} \label{FastBeads}

This section is the analogue of Section \ref{NonConstantSection}; it is here
that the proof of the Main Theorem  begins in earnest.

Let $\Delta$ be a minimal area van Kampen diagram, folded
according to the convention of Section \ref{Folding},
and fix a corridor $S_0$ in $\Delta$. As explained in 
Section \ref{Linear}, the core of our task is to bound
the number of beads in the decomposition of $\bot(S_0)$.
In order to do so, we must undertake a detailed study of
the preferred futures of these beads.

First we dispense 
with the case that $\widecheck{\beta}$ is a
vanishing \atom.  

\begin{lemma} \label{VanishLemma}
Suppose that $\mathcal{S}$ is the collection of beads in $S_0$
which are not vanishing \atom s.  If $\sum_{\beta \in \mathcal{S}}
|\beta| = D$ then $|S_0| \le B(D+1)$.
\end{lemma}
\begin{proof} This follows immediately from the Bounded Cancellation
Lemma.
\end{proof}
  
Narrowing our focus in the light of this lemma, we define:

\begin{definition} [Bead norm]\label{beadNorm} \index{bead norm}
Given a concatenation
 $\rho=\rho_1\dots\rho_m$ 
 of  beaded paths, we
define the {\em{bead norm}} of $\rho$, denoted $\| \rho \|_\beta$, 
to be the number of non-vanishing beads in the concatenation.
(This is poor notation, since the norm depends
on the decomposition into the $\rho_i$ and not just the edge-path $\rho$. But
in the contexts we shall use it, specifically $\bot(S_0)$, it will
always be clear which decomposition we are considering.)
\end{definition}

\begin{remark}
All beads have length at least $1$.  Thus bead norm is dominated by
length.  In particular, estimates concerning
Bounded Singularities and Bounded Cancellation
remain true when distance is replaced by bead norm; cf.~Lemma 
\ref{EasyBits}.
\end{remark}

\begin{remark} \label{BeadsDontVanish} An important advantage
of bead norm over edge-length is that when one takes the
repeated images $f_{\#}^k(\chi)$ of a monochromatic path,
its length can decrease,
due to cancellation within beads, whereas bead norm cannot.
\end{remark}

In Definition \ref{BeadLength} we defined the bead length $[S]_\beta$ of a
corridor $S$ in a well-folded diagram.  It is convenient for our future
arguments to concentrate on non-vanishing \atom s, and hence on
bead norm rather than bead length.  However,
an immediate consequence of the Bounded Cancellation Lemma
is the following bi-Lipschitz estimate:

\begin{lemma} \label{BLengthNorm}
Suppose $S$ is a corridor in a well-folded corridor.  Then
\[	\| S \|_{\beta} \le [S]_\beta \le B \| S \|_{\beta}	.	\]
\end{lemma}

\subsection{The first decomposition of $S_0$} [cf. 
Subsection \ref{decomp}]
\index{corridor!first decomposition of}

Let $\beta$ be a bead in $S_0$ that is not a vanishing \atom. 
As we follow the preferred future of $\beta$ forwards in time, one of the
following events must occur:
\begin{enumerate}
\item[1.] The last preferred future of $\beta$ intersects the boundary
of $\Delta$ nontrivially.
\item[2.] The last preferred future of $\beta$ intersects a singularity
nontrivially.
\item[3.] The last preferred future of $\beta$ is bitten in a corridor $S$.
\end{enumerate}
We remark that, unlike in Part \ref{Part:BG1}, these events are not
mutually exclusive; this is because a bead can consist of
more than one edge.

We shall bound the bead norm of $S_0$ by finding a bound on the 
number of non-vanishing beads  in each of the
three cases.

We divide Case (3) into two sub-cases:
\begin{enumerate}
\item[3a.] The preferred future of $\beta$ is bitten by a bead
that is not in the future of $S_0$.
\item[3b.] The preferred future of $\beta$ is bitten by a bead
that is in the future of $S_0$.
\end{enumerate}

\subsection{Bounding the easy bits} [cf. Subsection \ref{EasyBounding}]

Label the non-vanishing beads which fall into the above classes
$S_0(1)$, $S_0(2)$, $S_0(3a)$ and $S_0(3b)$, respectively.  We shall
see, just as in Part \ref{Part:BG1}, that $S_0(3b)$ is by far the most 
troublesome of these sets.

The following lemma is proved in an entirely similar way to
Lemmas \ref{bound1and2} and \ref{bound3a}, using the Bounded Cancellation
Lemma and simple counting arguments.

\begin{lemma} \label{EasyBits}
\ 

\begin{enumerate}
\item $\| S_0(1) \|_\beta \le \n$.
\item $\| S_0(2) \|_\beta \le 2B\n$.
\item $\| S_0(3a) \|_\beta \le B \n$.
\end{enumerate}
\end{lemma}

We have thus reduced our task of bounding $\| S_0 \|_\beta$
to bounding the numbers of beads in $S_0(3b)$, i.e. to
understanding cancellation {\em within} the future of $S_0$.
The bound on the number of beads in $S_0(3b)$ is proved
in an analogous way to Part \ref{Part:BG1}, and takes up a large
part of the remainder of this part of the book (through Section 
\ref{BonusSection}).

\subsection{The chromatic decomposition} \label{Chromatic}
[cf. Subsection \ref{chromatic}]
\index{corridor!chromatic decomposition of}

Fix a colour $\mu$ and consider the interval $\mu(S_0)$ in
$\bot(S_0)$ consisting of beads coloured $\mu$.

We shall subdivide $\mu(S_0)$ into five (disjoint but possibly
empty) subintervals according to
the fates of the preferred futures of the beads.

Let $l_\mu(S_0)$ be the rightmost bead $\beta$ in $\mu(S_0)$
such that $f_{\#}(\widecheck{\beta})$ contains a left-fast displayed edge
$\epsilon$
so that the preferred future of $\epsilon$ is eventually
bitten from the left from within the future of $S_0$.
Let $A_1(S_0,\mu)$ be the set of beads in $\mu(S_0)$ from the left
end up to and including $l_\mu(S_0)$.

Let $A_2(S_0,\mu)$ consist of those beads which are not in 
$A_1(S_0,\mu)$ but whose preferred futures are bitten from 
the left from within the future of $S_0$.

Let $A_3(S_0,\mu)$ denote those beads which do not lie in
$A_1(S_0,\mu)$ or $A_2(S_0,\mu)$ and which fall into
the set $S_0(1) \cup S_0(2) \cup S_0(3a)$.

All of the beads which are not in $A_1(S_0,\mu)$, $A_2(S_0,\mu)$
or $A_3(S_0,\mu)$ must have their preferred future bitten from
the right from within the future of $S_0$.

Analogous to the definition of $l_\mu(S_0)$, we define
a bead $r_\mu(S_0)$: the bead $r_\mu(S_0)$ is the
leftmost bead $\beta'$ so that $f_{\#}(\widecheck{\beta'})$ contains
a right-fast displayed edge whose preferred future is eventually bitten
from the right from within the future of $S_0$.

Let $A_4(S_0,\mu)$ denote those beads which are not in
$A_1(S_0,\mu)$, $A_2(S_0,\mu)$ or $A_3(S_0,\mu)$ and which
lie strictly to the left of $r_\mu(S_0)$. 

Finally, let $A_5(S_0,\mu)$ denote those edges not in 
$A_1(S_0,\mu)$, $A_2(S_0,\mu)$, $A_3(S_0,\mu)$ or
$A_4(S_0,\mu)$ which lie to the right of $r_\mu(S_0)$
(include $r_\mu(S_0)$ in $A_5(S_0,\mu)$ if it has not 
already been included in one of the earlier sets).

Now Lemma \ref{EasyBits} immediately implies

\begin{lemma}
\[      \sum_{\mu} \| A_3(S_0,\mu) \|_\beta \le (3B+1) \n       .       \]
\end{lemma}

We also have

\begin{lemma} \label{A1Short}
Let $C_0$ be the constant from Lemma \ref{C0Lemma} above.  Then
\begin{enumerate}
\item $\| A_1(S_0,\mu) \|_\beta, \| A_5(S_0,\mu) \| \le C_0$; and
\item $| A_1(S_0,\mu) \smallsetminus l_\mu(S_0)|,
| A_5(S_0,\mu) \smallsetminus r_{\mu}(S_0)|  \le C_0$.
\end{enumerate}
\end{lemma}

\begin{proof}
We prove the bounds only for $A_1(S_0,\mu)$, the proofs for
$A_5(S_0,\mu)$ being entirely similar.

The entire future of beads in 
$A_1(S_0,\mu)$ other than $l_\mu(S_0)$ must be eventually
consumed from the left from within the future of $S_0$;
cf. Lemma \ref{perfect}.

If $\| A_1(S_0,\mu) \|_{\beta}$ or $| A_1(S_0,\mu) \smallsetminus 
l_\mu(S_0)|$ were greater than $C_0$ then we would conclude 
from Lemma \ref{C0Lemma} that
no left-fast bead in the immediate future of $l_\mu(S_0)$
could be bitten at any stage from the left from within the future of $S_0$,
contrary to the definition of $l_\mu(S_0)$.
\end{proof}

As we continue to follow the proof from Part \ref{Part:BG1}, our next goal is to reduce
the task of bounding the bead norm of $S_0$ to that of bounding
the number of Nielsen beads contained in $A_2(S_0,\mu)$
and $A_4(S_0,\mu)$.  We focus exclusively on $A_4(S_0,\mu)$, the arguments
for $A_2(S_0,\mu)$ being entirely similar.

In outline, our argument proceeds  in analogy with the subsections
beginning with  Subsection \ref{ss:furtherdecomp}, commencing with the decomposition
of $A_4(S_0,\mu)$ into subintervals $C_{(\mu,\mu')}$. But we quickly
encounter  a new phenomenon that requires an additional section
of argument --
HNP cancellation -- which does not arise in the case of positive
automorphisms.

\subsection{The decomposition of $A_4(S_0,\mu)$ into the $C_{(\mu,\mu')}$.}

All beads in $A_4(S_0,\mu)$ are eventually bitten from the right
from within the future of $S_0$.  For a colour $\mu' \neq \mu$, 
define a subset $C_{(\mu,\mu')}$ of $A_4(S_0,\mu)$ as follows:
given a bead $\sigma \in A_4(S_0,\mu)$, there is a bead $\sigma'$
in $S_0$ so that $\sigma$ is eventually bitten by $\sigma'$.
If $\sigma'$ is coloured $\mu'$ then $\sigma \in C_{(\mu,\mu')}$. 

The sets $C_{(\mu,\mu')}$ form intervals in $S_0$.

\section{HNP-Cancellation and Reapers} 
\label{TrappedHNP}

The results of the previous section
 reduce the task of bounding $\| S_0 \|_\beta$ to
that of establishing a bound on the sum
of the bead norms of the monochromatic
intervals $C_{(\mu,\mu')}$.  In 
Part \ref{Part:BG1}, the corresponding intervals (also labelled $C_{(\mu,\mu')}$)
contained no exponential edges.  In the current context, however, 
there
may be exponential edges {\em{trapped}} in 
Nielsen paths, which may themselves be contained in beads of any type. 
This raises the concern that our attempts to control the length of the
$C_{(\mu,\mu')}$ in the manner of Part \ref{Part:BG1} will be undermined by
the {\em{release}} of these trapped edges when the Nielsen path
is bitten, leading to  rapid growth in subsequent
nibbled futures of the Nielsen path. Our purpose in this section
is to develop  tools to control  this situation, specifically Lemmas
\ref{DealWithExp} and \ref{WholeNielsen}.

We must also deal with a second threat 
that arises from the phenomenon described in Example \ref{Ex:HNP};
we call this {\em Half Nielsen Path (HNP-) cancellation}.

Recall that a {\em{\pep}} is an edge path $\rho$ in $G$;
it is associated to a \gep\  and
 either $\rho$ or $\bar{\rho}$ is of the
form $E \bar{\tau}^k \bar{\nu} \gamma$
where $E$ is an edge with
$f_\#(E) = E\odot\tau^m$, where $\tau$ and $\nu$
are  Nielsen paths, and $\bar{\gamma}\nu$ is a terminal segment of $\tau$
(and $m, k > 0$).
These are the prototypes of the following types of paths.

\begin{definition}
Suppose that $E$ is a linear edge with $f_{\#}(E) = E \odot \tau^m$, where $\tau$ is a 
Nielsen path and $m > 0$.  Suppose further that $\nu$ is a Nielsen path and $\gamma$
an edge-path so that $\bar{\gamma}\nu$ is a terminal segment of $\tau$.  \index{PEP}

A {\em \rpep} is a path $\rho$ so that either $\rho$ or $\bar{\rho}$ has the form
$E \bar{\tau}^k \bar{\nu} \gamma$ where $k > 0$.
\end{definition}

\begin{remark} Every
\pep\ is  a \rpep, but an arbitrary
 \rpep\  has no \gep\  associated to it.
\end{remark}

It is important to note that in the following definition the \rpep \ 
being discussed is {\em{not}} assumed
to be a bead in the decomposition of $\bot(S)$. (Beads along
$\bot(S)$ are monochromatic whereas we want to discuss HNP
cancellation, as in Definition \ref{HNPBite}, in the context
of adjacent colours interacting.)

\begin{definition}[HNP cancellation] \label{HNPDef} \index{HNP cancellation}
Let $S$ be a corridor in
a well-folded diagram, let
$\e$ and $\e'$ be edges in the 
naive (unfolded) future of $\bot(S)$ that cancel in the passage to $\top(S)$
and assume that  $\e$ is to the
left of $\e'$.

Suppose further that the past of $\e$ is $e$ with
label $\check e = E$ a linear edge 
and that $\e'$ is in the future of an edge $e_\gamma$
whose label is an edge $\gamma$.

We call the cancellation of $\e$ and $\e'$ {\em left HNP-cancellation}
and write $\e \lhnp \e'$
if the interval from $e$ to $e_\gamma$ in $\bot(S)$ (inclusive) is
labelled by a \rpep\
of the form $E \bar{\tau}^k \bar{\nu}\text{\em\o}\gamma$, where
$\tau$ is a Nielsen path so that $\tau = \xi \nu$, where $\xi$ and $\nu$
are Nielsen paths, and 
$\overline{\text{\em\o}\gamma}$ is a terminal sub edge-path of $\xi$.

{\em Right HNP-cancellation} is defined by reversing the roles of $\e$ and 
$\e'$
and insisting upon a \rpep\ in $\bot(S)$ of the form
$\bar{\gamma}\bar{\text{\em\o}} \nu \tau^k \bar{E}$.  It is denoted $\e \rhnp \e'$.

When we are unconcerned about the distinction between left and right, we refer
simply to {\em HNP-cancellation}.

We extend this definition to concatenations of beaded paths in $G$ by using the
obvious stack-of-corridors diagram as in Subsection \ref{stackDiags}.
\end{definition}

\begin{remark}
HNP-cancellation occurs at the `moment of
death' of the \rpep; see Section \ref{nfGEPSection} for an explanation of the
significance of this moment  and an analysis of it (in
the language of \pep s).
\end{remark}

\begin{lemma} \label{l:oEmpty}
Suppose that $E \overline{\tau}^k \overline{\nu} \text{\em\o}\gamma$ is a \rpep \
which exhibits an HNP-cancellation, as in Definition \ref{HNPDef}.
Then $\text{\em\o}$ is empty, so $\gamma$ is the first edge of $\overline{\xi}$.
\end{lemma}
\begin{proof}
The assumption that HNP-cancellation occurs means that we can
restrict our attention to cancellation when tightening
\[	f(E\overline{\tau}^k\overline{\nu}\text{\em\o}\gamma).	\]
This can be written as
\[	E \tau^m f(\overline{\tau}^k \overline{\nu}) f(\text{\em\o}\gamma).	\]
The path $\overline{\tau}^k\overline{\nu}\text{\em\o}\gamma$ admits a
hard splitting $\overline{\tau} \odot \cdots \odot  \overline{\tau}
\odot \overline{\nu} \odot \text{\em\o}\gamma$.
Therefore, under any choice of tightening, the $m$ copies of 
$\tau$ cancel with the $k$ copies of $f(\overline{\tau})$ (partially
tightened), then with $f(\overline{\nu})$; they then begin to
interact with $f(\text{\em\o}\gamma)$.  Just as in the proof of 
Proposition \ref{gammasingleedge}, under the assumptions of 
Lemma \ref{Power1}, there is only a single edge in $\text{\em\o}\gamma$
whose future can interact with $f(E)$ when tightening. 
\end{proof}

We now present the deferred example that explains the
need to consider HNP-cancellation. This will also lead
us to a further definition --- {\em{HNP biting}} --- that
encodes a genuinely troublesome situation where HNP
cancellation must be accounted\footnote{We usually
account for it by excluding it from our definitions.  When it
cannot be excluded, we often sidestep it, using the notions of `robust future'
and `robust past' given in Definitions \ref{Robust} and \ref{RobustPast}
below.}  for. Fortunately, many
other instances of HNP-cancellation are swept-up 
by our general cancellation and finiteness arguments, allowing
us to avoid a detailed analysis of the possible outcomes.

The problem at the heart of the following example did
not arise in Part \ref{Part:BG1} because the natural
realisation of a positive automorphism does not map
any linear edge   across other linear edges.

\begin{example} \label{Ex:HNP}
Suppose that $u$ is a Nielsen path, and that $E_1$ and 
$E_2$ are edges so that $f(E_i) = E_i u^k$ for $i = 1,2$
and some integer $k > 0$.  For any integer $j$, the path
$\tau_j = E_1 u^j \overline{E_2}$ is an indivisible Nielsen path.

Suppose that $E_3$ is an edge so that $f(E_3) = E_3 \tau_j^l$,
for some integers $j$ and $l$ (with $l > 0$).  For ease of 
notation, we will assume that $l =1$.

Consider the path $\rho = E_3 \overline{\tau}_j^r E_2$,
for some $r > 0$.  Then $\rho$ is a \rpep.

In the iterated images $f_{\#}(\rho)$, the visible copy of $E_2$
has a unique future labelled $E_2$, which we will call the `preferred
future' of $E_2$ for the purposes of this example.
After $r+1$ iterations of $\rho$ under $f_{\#}$ (and any choice
of tightening at each stage), the future of $E_3$ cancels the preferred future
of the visible copy of $E_2$.  If we encode the evolution of $\rho$
in a \index{stack diagram}stack diagram as in Subsection \ref{stackDiags} then the cancellation
of $E_2$ is HNP-cancellation.
\end{example}

In the following discussion, we assume that the reader is
familiar with  Part \ref{Part:BG1}, in particular the vocabulary of
teams and reapers.

The phenomenon described in
the above example causes problems when 
the sub-path $\rho_1 = \overline{\tau}_j^r E_2$ of $\rho$ is
monochromatic and $E_2$ is displayed in $\rho_1$.  In
this situation, it shows that the most obvious
adaptation of Lemma \ref{C1Lemma} would be false.
It is for this reason that we must exclude
HNP-biting in Definition \ref{DefEi}.

Similarly, because  Example \ref{Ex:HNP}
renders a naive version of the results of Section \ref{ConstantSection} false,
HNP-biting must be excluded from
the Two Colour Lemma and the associated results
in Section \ref{PincerSection}.

A situation in which we cannot exclude HNP-biting by decree
arises in the analysis of teams and in particular
the definition of a {\em{reaper}} (Subsection \ref{ss:reaper}).
Suppose that  $\rho$ labels some interval in the bottom
of a corridor, with many copies of $\overline{u}$ to its immediate right.
In this case, the edge $\e_2$ labelled
$E_2$ will consume copies of $\overline{u}$
in the first $r$ units of time, but its future
will then be cancelled (assuming no
other cancellation occurs from either side, and that there are no
singularities, etc.).  Since $\e_2$ was acting as the reaper
of a team, we must find a continuing manifestation of it
at subsequent times, for otherwise we
will lose control over the length of teams
($r$ being arbitrary) and the structure of our main argument
will fail. This problem is solved by introducing the
{\em{robust future}} of $\e_2$
(Definition \ref{Robust}), which in this case is an edge labelled
$E_1$ that `replaces' the preferred future of $\e_2$ when it
is cancelled.

\begin{definition} \label{HNPBite}\index{HNP biting}
Suppose that $\chi_1$ and $\chi_2$ are beaded paths in $G$
and $\chi_1\chi_2$ is tight.
Suppose that there is a bead $\rho_1 \subset \chi_1$
and a bead $\rho_2 \subset \chi_2$ so that
\begin{enumerate}
\item\label{rho1} either $\rho_1$ is a displayed edge $\gamma$ in $\chi_1$ which is linear or else
$\rho_1$ is a displayed \pep \ in $\chi_1$ of the form
$E \bar{\tau}^k \bar{\nu} \gamma$, where $\gamma$ is a linear edge;
\item when tightening $f_{\#}(\chi_1)f_{\#}(\chi_2)$ to form
$f_{\#}(\chi_1\chi_2)$, $\rho_1$ bites $\rho_2$ and the edge $\e'$
in the exhibiting pair $( \e',\e)$ (see Definition \ref{d:bite}) is
in the future of $\gamma$;
\item moreover\footnote{The \rpep\ implicit in the symbol $\rhnp$ is not
the \pep\ in \eqref{rho1}.}, $\e' \rhnp \e$.
\end{enumerate}
Under these circumstances we say that $\rho_2$ is {\em left-HNP-bitten} by $\rho_1$
and we write $\rho_1 \lbite \rho_2$.
There is an entirely analogous definition of {\em right-HNP-biting} $\rho_1 \rbite \rho_2$,
and when we are unconcerned about the direction we will refer
simply\footnote{We swap orientation in Definition \ref{d:evHNPb}
so as to emphasize  this point immediately.}
 to {\em HNP-biting}.

We make the analogous definition for HNP-biting within diagrams.
\end{definition}

\begin{definition}\label{d:evHNPb}
Suppose that $\chi_1$ and $\chi_2$ are beaded paths
and that $\rho_1$ is a bead in $\chi_1$.  We say that $\rho_1$
is {\em eventually HNP-bitten} by $\chi_2$ if $\rho_1$ is eventually
bitten by $\chi_2$ (Definition \ref{d:eventbite}) and this biting is HNP-biting.

We make the analogous definition within diagrams.
\end{definition}

\begin{definition}
Suppose that $E$ and $E'$ are edges in $G$.  We say that \index{edge!indistinguishable}
$E$ and $E'$ are {\em indistinguishable} if there is a Nielsen path
$\tau$ and an integer $s > 0$ so that $f(E) = E \tau^s$ and
$f(E') = E' \tau^s$.
\end{definition}

The edges $E_1$ and $E_2$ in Example \ref{Ex:HNP} are
indistinguishable.

\subsection{Parabolic HNP-cancellation and robust futures}

The following is a simple (but key) observation, and has an obvious
application to HNP-cancellation of edges of parabolic weight.

\begin{lemma} \label{HNPreplace}
Suppose that $\tau$, $\nu$, $\nu'$ and $\sigma$ are Nielsen paths, 
with $\sigma$ irreducible and
$\tau = \nu' \overline{\sigma} \nu$.  Suppose further that $\gamma$ is the initial edge
of $\sigma$, and that $f(\gamma) = \gamma \odot \xi^l$ for some Nielsen
path $\xi$.  Then $\sigma$ has the form $\gamma \xi^r \overline{\gamma'}$
where $r$ is some integer and $\gamma'$ is an edge so that
$\gamma$ and $\gamma'$ are indistinguishable.

Moreover, suppose that $E$ is an edge so that $f(E) = E \odot \tau^m$,
and let $\rho = E \overline{\tau}^i \overline{\nu} \gamma$ be a \rpep \ with
$0 \le i < m$. Then $f_{\#}(\rho)$ has the form
$E \odot \tau^{m-i-1} \nu' \gamma' \overline\xi^j$ where
$\gamma$ and $\gamma'$ are indistinguishable.
\end{lemma}
\begin{proof}
The first assertion is an immediate consequence of the structure of indivisible Nielsen
paths of parabolic weight, and the second is then obvious (a detailed analysis
of the Nielsen paths of parabolic weight is undertaken in Section \ref{TrainTracks}).
\end{proof}

\begin{definition}\label{ParabPrefFut} In general, non-displayed edges 
$\e$ in diagrams do not 
have preferred futures. But if $\check\e$ has parabolic weight, there is a
unique edge of the same weight in $f_\#(\check\e)$, and it is natural
to define the (immediate) preferred future of $\e$ to be the corresponding
edge in the immediate future of $\e$. (If $\e$ happens to be displayed,
this agrees with our earlier definition.)
\end{definition}

In Section \ref{PincerSection}, when proving the Pincer Lemma,
we will have to exclude HNP-biting.  This will also be the case
in the applications of the Pincer Lemma in Sections \ref{TeamsSection}
and \ref{BonusSection}.  Thus, in following the future of  a linear
edge $\gamma$ when HNP-cancellation occurs, we would like to ignore the 
preferred future (which disappears), and rather follow the future of 
the indistinguishable edge $\gamma'$ from Lemma \ref{HNPreplace}
above.  Thus we make the following

\begin{definition} [Robust Futures for Parabolic
Edges] \label{Robust} \index{edge!parabolic} \index{future!robust}
Suppose that $\e$ is a (not necessarily displayed) edge in a colour
$\mu(S)$, and that $\widecheck{\e}$ is contained in a parabolic
stratum.  If the preferred future of $\e$ is cancelled from the left 
[resp.~right] by
HNP-cancellation in $\top(S)$, then
Lemma \ref{HNPreplace} provides an edge $\gamma'$
that is indistinguishable from $\widecheck\e$ and survives in the
tightened path $f_{\#}(E \bar{\tau}^k\bar{\nu}\text{\em\o} \gamma)$
[resp.~its reverse] considered in Definition \ref{HNPDef}.

We define the {\em robust
future} of an edge $\e \subseteq \bot(S)$ as follows.  If the preferred
future of $\e$ survives in $\top(S)$, then the robust future of $\e$
is just the preferred future of $\e$.  If the preferred future is
cancelled by HNP-cancellation, then the robust future of $\e$ is the
above edge labelled $\gamma'$, provided this survives in
$\top(S)$.  Otherwise there is no robust future.
\end{definition}

\begin{definition} [Robust Pasts for Linear Edges]\label{RobustPast}\index{edge!linear!robust past of}
Let $\e'$ be an edge of $\top(S)$ and suppose
that both it and its immediate past are labelled by linear edges. If
$\e'$ is not the robust future of any edge then the robust past of
$\e'$ is the past of $\e'$.  But if $\e'$ {\em{is}} the (immediate) robust
future of $\e$ then the robust past of $\e'$ is $\e$.
\end{definition}

Just as for preferred futures, the notions of robust future and robust 
past can be extended arbitrarily many steps forwards or backwards
in time by iterating the definition.

\subsection{A setting where we require cancellation lemmas}
\label{ss:setting}

Consider the following situation.  Let $\chi_1\sigma\chi_2$
be a tight path in $G$ with $\chi_1$ and $\chi_2$
monochromatic and  $\sigma$  a path with a preferred
decomposition into monochromatic paths (each of which comes
equipped with a beaded decomposition).  
We will analyse the possible
interaction between $\chi_1$ and $\chi_2$ in iterates of
$\chi_1\sigma\chi_2$ under $f$ (where the tightening follows
the convention of Subsection \ref{AbstractFolding}).

As ever, the following lemma remains valid with left/right orientation
reversed.

 \begin{lemma} \label{OneHNP}
Suppose that $\chi_1$, $\chi_2$ and $\sigma$ are as above, and 
suppose that each non-vanishing bead in $\chi_2$ is 
eventually bitten by a bead from $\chi_1$ in some
iterated image  $f_{\#}^k(\chi_1\sigma\chi_2)$
of $\chi_1\sigma\chi_2$.

Suppose further that $\rho$ is a bead in $\chi_2$ so that
$f_{\#}(\rho)$ has parabolic weight,
and that $\rho$ is eventually left-HNP-bitten 
by a bead from $\chi_1$ in the evolution of $\chi_1\sigma\chi_2$. Then $\rho$
is the rightmost non-vanishing bead in $\chi_2$.
\end{lemma}
\begin{proof}
Pass to the iterate $f_{\#}^{k-1}(\chi_1\sigma \chi_2)$ so that
the preferred future of $\rho$ lies in a \rpep\  $\pi$, which exhibits the
(eventual) HNP-biting of $\rho$ in the tightening to form
$f_{\#}^{k}(\chi_1\sigma \chi_2)$.  Let $\rho_1$ be the preferred
future of $\rho$ in $f_{\#}^{k-1}(\chi_1\sigma \chi_2)$.  Since
$f_{\#}(\rho)$ has parabolic weight, $\rho_1$ has parabolic weight,
and is either a displayed edge or a displayed \pep\ or \gep.
 We must prove
that no bead to the right of $\rho_1$ is eventually bitten by
the future of $\chi_1$.

By Definition \ref{HNPBite} and Lemma \ref{l:oEmpty} the
\rpep \ $\pi$ has the form $\gamma \bar{\tau}^k\bar{\nu} \e$, where
\begin{enumerate}
\item $\gamma$ is an edge so that $f(\gamma) = \gamma \odot \tau^m$;
\item $\gamma$ is either a displayed edge in the future of $\chi_1$ in 
$f^{k-1}_{\#}(\chi_1\sigma\chi_2)$ or else if the rightmost edge in a 
displayed \pep; and
\item $\e$ is contained in $\rho_1$.
\end{enumerate}
Let $\alpha$ be the displayed edge or \pep\ containing $\gamma$.

Let $\rho_1'$ be the terminal part of $\rho_1$ from $\e$ to its right end,
and let $\chi_2'$ be the terminal part of the future of $\chi_2$ in
$f_{\#}^{k-1}(\chi_1\sigma\chi_2)$, from $\e$ to its right end.

Since $\rho_1$ is displayed, we have $\chi_2' = 
\rho_1' \odot \beta$
for some path $\beta$.  

By Lemma \ref{HNPreplace}, when tightening to form
$f_{\#}^{k}(\chi_1\sigma \chi_2)$, the edge $\e$ is replaced
by an indistinguishable edge $\e'$ which comes from the
future of $\alpha$.  Suppose that $\delta$ is that part
of $f_{\#}(\alpha \rho_1')$ from $\e'$ to the right end.
Since $\alpha$ is a (linear) edge or a \pep,
the edge $\e'$ survives in all iterates of $\alpha$ (under
any choices of cancellation).  Similarly, since $\e$ and
$\e'$ are indistinguishable, $\e'$ survives in all iterates
of $\delta$ (under any choices of tightening).  This implies that
we have a hard splitting 
$f_{\#}(\alpha \chi_2') =  f_{\#}(\alpha \rho_1') \odot f_{\#}(\beta)$,
and the fact that $\alpha$ is displayed implies that
no bead in $\beta$ can be eventually bitten by the future
of $\chi_1$, as required.
\end{proof}

In applications of Lemma \ref{OneHNP} (and of Lemmas \ref{DealWithExp} and \ref{WholeNielsen} below), we usually take
$\chi_1 = \widecheck{\mu_1(S)}$ and $\chi_2 = \widecheck{\mu_2(S)}$,
where $\mu_1$ and $\mu_2$ are colours and $S$ is some corridor,
and we will choose $\sigma$ to be the label of that part of $\bot(S)$
which lies strictly between $\mu_1(S)$ and $\mu_2(S)$.\footnote{
However, it will also be convenient sometimes to take $\chi_1$ to
be a subinterval of $\widecheck{\mu_1(S)}$ consisting of an interval of beads.}
Since the folding conventions of Subsections \ref{Folding} and 
\ref{AbstractFolding} are compatible, and because of the hardness of our splittings,
the interaction between $\mu_1$
and $\mu_2$ in the future of $S$ can be analysed by studying the
interaction between the futures of $\chi_1$ and $\chi_2$ in iterated images
of $\chi_1\sigma \chi_2$ under $f$.

\subsection{Reapers} \label{ss:reaper}

In Part \ref{Part:BG1} proving the existence of reapers was straightforward
(see Section \ref{teamSec}).
In the current context, however, we have to work harder to
prove that a suitable incarnation of
a reaper exists, because of the phenomena discussed
in the preceding subsection. At the heart of our difficulties
is the fact that Nielsen atoms need not be single edges.

\begin{definition} \index{Nielsen path!beaded}
A {\em beaded Nielsen path} in a corridor $S$ is a subinterval
$\sigma \subset \bot(S)$ so that $\widecheck{\sigma}$ is a beaded
path all of whose beads are Nielsen paths.
\end{definition}

Note that in the above definition we do not assume that $\sigma$
is a single colour, or even that each bead in $\widecheck{\sigma}$
is contained in a single colour.  
Examples of beaded Nielsen paths include
that part of a \gep\ between the extremal edges,
and the sub-paths $\overline{\tau}^i$ of a \rpep\ 
$E\overline{\tau}^k\bar{\nu}\text{\em\o}\gamma$.

Although the beads in a beaded Nielsen path might not be displayed in
a path $\widecheck{\mu(S)}$, it is still possible to define the future
of a bead in a beaded Nielsen path, and the notions of preferred
future and biting still make sense.  We will use this observation in 
the sequel.

The following notion is parallel to that of Definition \ref{def:swol},
which was pivotal in the bonus scheme (cf. Section
\ref{BonusSection} below). Here, it plays a more central role.

\begin{definition}[Swollen present and swollen future]
\label{Swollen}\index{future!swollen}\index{swollen present}
Suppose $S$ is a corridor and that $I \subseteq \bot(S)$ is a
beaded Nielsen path in $S$. 
The {\em swollen present} of $I$ is the\footnote{Uniqueness is immediate from
the observation that if a terminal sub-path $\sigma$ of a Nielsen path $\tau$ is itself Nielsen
then   $\sigma$ is a concatenation of beads in $\tau$.}
maximal subinterval $I' \subseteq \bot(S)$ such that (i) $I \subseteq I'$;
(ii) $I'$ is a beaded Nielsen path in $S$; and (iii) the beads of $I$ 
are beads of $I'$.

The {\em left swollen present} of $I$ is that part of the swollen present from the
left end up to the right end of $I$, whilst the {\em right-swollen present} goes
from the left end of $I$ to the right end of the swollen present.

If the actual future of $I$ is a beaded Nielsen path the
(immediate) {\em swollen future} $sw_1(I)$ of $I$ is the swollen present of the
(actual) future of $I$. With a similar 
qualification, the {\em swollen future} $sw_k(I)$ at $\time(S) + k$ is
defined to be $sw_1(sw_{k-1}(I))$.

With the same qualifications, the left and right swollen futures are defined in
the obvious ways.
\end{definition}

The first qualification in the above definition is required because it is possible
that the immediate future of a beaded Nielsen path is not a beaded Nielsen
path.  Thus we must be careful only to apply
this concept in cases where we know the swollen future to exist.

\begin{definition}[Reapers] \label{ReaperDef}\index{reaper}
Suppose that $S$ is a corridor and $I \subset \bot(S)$ is
a beaded Nielsen path in $S$ with nonempty swollen future $sw_1(I)$.
Suppose that $\alpha$ is an edge in $\bot(S)$
immediately adjacent to $I$ on the left.  We say that $\alpha$ is
a {\em left-reaper for $I$} if (i) $\widecheck{\alpha}$ is a linear edge;
(ii) $\widecheck\alpha$ bites some of the future of $\widecheck{I}$ in 
$f_{\#}(\widecheck{\alpha I})$; and
(iii) the robust future of $\alpha$ is immediately adjacent to $sw_1(I)$ in $\top(S)$.

There is an entirely analogous definition of {\em right-reapers}.
As usual, when we are unconcerned about the direction we will refer
to {\em reapers}.
\end{definition}

\begin{definition}[Left-edible] \label{d:left-edible} \index{path!edible}
Let $S$ be a corridor in a well-folded diagram, and $I \subset \bot(S)$
a beaded Nielsen path.  
We say that $I$ is {\em left-edible}
if each bead in $I$ is eventually bitten by a bead coloured $\mu$ in the
future of $S$, where $\mu(S)$ lies to the left of $I$.

{\em Right-edible} paths are defined with a reversal of the left-right orientation.
\end{definition}

In the remainder of this section we work towards proving
Propositions \ref{OneSideStable} and \ref{Reaper}.

\begin{proposition} \label{OneSideStable}
Let $S$ be a corridor in a well-folded diagram and $I \subset \bot(S)$ a
left-edible path so that $|I| \ge B+J$.  Then the immediate future of
$I$ in $\top(S)$ is left-edible.
\end{proposition}

The following lemma is straightforward, and allows us to focus our attention
on the time when cancellation between colours begins.

\begin{lemma} \label{l:gotofuture}
Let $S$ be a corridor in a well-folded diagram and let $I \subset \bot(S)$
be a left-edible colour, all of whose beads are eventually bitten by
beads coloured $\mu$.  Let $S^I$ be the corridor in the future of $S$ so that
the first  biting of a bead in the left swollen future of $I$ by something
coloured $\mu$ occurs in $S^I$.  Then the left swollen future of $I$ in $\bot(S^I)$
is left-edible.
\end{lemma}

 In the following statement $B$ is the Bounded Cancellation Constant
from Proposition \ref{BCL:BG3} and $J$ is the constant from the Beaded
Decomposition Theorem \ref{BDT}.  The corridor $S^I$ is as in Lemma 
\ref{l:gotofuture} above, and $I^\lambda$ is the left swollen future of $I$ in
$S^I$.

\begin{proposition} \label{Reaper}
Suppose that $S$ is a corridor in a well-folded diagram and $I \subset \bot(S)$
is a left-edible path, all of whose beads are eventually bitten by beads coloured
$\mu$.  Suppose also that $|I| \ge B+J$.  Then
\begin{enumerate}
\item the immediate future of $I^\lambda$ in $\top(S^I)$ has an associated
left reaper $\alpha$, which is coloured $\mu$; and
\item for each bead in the immediate future of $I^\lambda$, when it is eventually
bitten the biting is by the robust future of $\alpha$.
\end{enumerate}
\end{proposition}

\subsection{Two Cancellation Lemmas}

The following lemma is useful in the proof of Lemma
\ref{C1Lemma:BG3} below. We record it now because a variation
on it (Lemma \ref{WholeNielsen}) is needed in the proof of 
Proposition \ref{Reaper}.

We revert to the setting described in Subsection \ref{ss:setting}.

\begin{lemma} \label{DealWithExp}
Assume that in the iterates of $\chi_1\sigma\chi_2$
(i.e.~forward-images  under $f_{\#}$)
each bead in $\chi_2$ is eventually bitten by a bead in $\chi_1$.
Suppose that $\chi_2$ has weight $i$, where $H_i$ is an
exponential stratum, and that all beads of weight $i$ in $\chi_2$
are Nielsen beads.
Let $\rho$ be a bead of weight $i$ in $\chi_2$.
\begin{enumerate}
\item  If $\rho$ is not bitten in
$f_{\#}(\chi_1\sigma\chi_2)$ but is eventually bitten
in the image $f_{\#}^k(\chi_1\sigma\chi_2)$  
then $\rho$ is entirely consumed in $f_{\#}^k(\chi_1\sigma\chi_2)$.
\item If $\rho$ is bitten but not entirely consumed in
$f_{\#}(\chi_1\sigma\chi_2)$ then
$\rho$ is the rightmost bead in $\chi_2$.
\end{enumerate}
\end{lemma}
\begin{proof}
There is at most one indivisible Nielsen path of weight $i$
and the lemma is vacuous unless there
is exactly one.

Let $\beta$ be a bead in $\chi_2$ of weight $i$, and suppose that
an edge $\eta$ in the future of $\chi_1$ is the edge which cancels the
rightmost edge in the preferred future of $\beta$ to exhibit the
biting of $\beta$ by $\chi_1$.
Since $\beta$ is an indivisible Nielsen path, it has edges of weight $i$
on both ends, as does its preferred future, 
and so $\eta$ has weight $i$.  Suppose that the 
past of $\eta$ in $\chi_1\sigma\chi_2$ has weight $i$.  Then by 
Theorem \ref{ColourCancelMain} and Assumption \ref{FinalPower}, $\eta$
is either a displayed edge in the future of $\chi_1$, or else is
contained in a Nielsen bead.  
Suppose first that $\eta$ is contained in a Nielsen bead $\tau$.
Since $\eta$ is to cancel with an edge in $\beta$, the path $\tau$
must have weight $i$.  Hence $\tau = \bar\beta$, and $\beta$ is 
entirely consumed when it is bitten.

Suppose then that $\eta$ is displayed in the future of $\chi_1$.  By
Assumption \ref{FinalPower}.(\ref{DisplayExp}) we may assume
that the edge $\eta$ is contained in a displayed path of the form
$f(\eta)$.  Since $f(\eta)$ is $i$-legal, and $\beta$ is not, it is not
possible for the illegal turn in $\beta$ (of weight $i$) to be cancelled
by any iterates of $\eta$.  However, $|f(\eta)| > |\beta|$, by Assumption
\ref{FinalPower}\eqref{ExpLong}, so it is not possible for the displayed
copy of $f(\eta)$ to be cancelled by the future of $\beta$.  Therefore,
in this case $\beta$ must be the rightmost bead in $\chi_2$.

Furthermore, suppose that $\beta$ and $\eta$ are as above, and
the past of $\eta$ in $\chi_1\sigma\chi_2$ has weight $i$, 
and suppose moreover that 
$\beta$ is not bitten in $f_{\#}(\chi_1\sigma\chi_2)$.  Then $\beta$ is 
bitten by $\eta$ in some  $f_{\#}^k(\chi_1\sigma\chi_2)$, and $k \ge 2$.  
Thus we may assume that the immediate past of $\eta$ is also displayed
and is $\eta$.  By applying Lemma 
\ref{EndStab1} and noting that the rightmost edge of $\beta$ must be
$\bar\eta$, we see that the sub-path between the immediate
past of $\beta$ and the immediate 
past of $\eta$ has the form $\cdots \bar\eta \omega \eta \cdots$
for some path $\omega$.  The path $\omega$ must start and finish at he
same vertex, and in
order for the written copy of $\bar\eta$ to cancel with the written
copy of $\eta$ it must be that $f_{\#}(\omega)$ is a point.  However,
$\omega$ is not a point, because otherwise the past of $\beta$ and
the past of $\eta$ would already cancel.  This contradicts the fact
that $f$ is a homotopy equivalence.  The same argument
shows that if $\eta$ is contained in a Nielsen bead and $\beta$
is not bitten in $f_{\#}(\chi_1\sigma\chi_2)$
then $\beta$ cannot be bitten by $\eta$.

Therefore, if $\beta$ is bitten by an edge $\eta$ whose past in 
$\chi_1\sigma\chi_2$ has weight $i$ then $\beta$ is close to
the left end of $\chi_2$, and is either entirely consumed when
bitten or is the rightmost bead in $\chi_2$.

We may now assume that the bead $\rho$ is cancelled by an edge
$\eta$ whose past in $\chi_2$ has weight greater than $i$.  The above
arguments show that we may assume that the immediate past of
$\eta$ also has weight greater than $i$, and by Lemma \ref{EndStab1}
we may assume that this past is contained in a displayed edge,
a \gep, or a \pep.  It is easy to see that the immediate past of 
$\eta$ cannot have exponential weight and cannot be a \gep.  
Thus we may assume that the immediate past of $\eta$ is either
the edge on the left end of a \pep\ of the form $\gamma\nu\tau^k\overline{E}$,
 (and that the edge $\gamma$
is parabolic) or else is displayed and parabolic.

Lemma \ref{EndStab1} and the above arguments imply that this
immediate past of $\eta$ must be a linear edge, and the above
arguments now imply that if $\rho$ is bitten in a corridor it must be
entirely consumed.
 \end{proof}

 The following variant of Lemma \ref{DealWithExp} is the one
we need in the  proof of Proposition \ref{Reaper}.
 We continue to study $\chi_1\sigma\chi_2$ as in 
Subsection \ref{ss:setting}.

\begin{lemma} \label{WholeNielsen}
Suppose that $\chi_2$ is a beaded
Nielsen path and each of its beads is eventually bitten
by a bead in $\chi_1$ in some iterated image of $\chi_1
\sigma\chi_2$ under $f$.

Let $\rho$ be a bead in $\chi_2$ which is not bitten in $f_{\#}(\chi_1\sigma\chi_2)$.
If $\rho$ is bitten but not consumed in some iterated image of $\chi_1\sigma\chi_2$
then $\rho$ is the rightmost bead in $\chi_2$.
\end{lemma}
\begin{proof}
We follow the proof of Lemma \ref{DealWithExp} above, with
the added wrinkle that there may be parabolic weight Nielsen
paths to consider in $\chi_2$.  In this case there needn't be 
a unique Nielsen path of weight $i$.

Suppose that $\rho$ is as in the statement of the Lemma.  If
$\rho$ has exponential weight, then the arguments of the proof
of Lemma \ref{DealWithExp} give the required properties.  
If $\rho$ has parabolic weight,  Lemma 
\ref{ParabolicOldEdges} implies that when $\rho$ is bitten by
an edge $\eta$ in the future of $\chi_1$,  the immediate
past of $\eta$ has weight greater than that of $\rho$.  Also,
this immediate past must be parabolic. Arguing as in the proof of 
Lemma \ref{DealWithExp}, one sees  that either $\rho$ is entirely 
consumed when bitten, or else $\rho$ is the rightmost bead
in $\chi_2$.
\end{proof}

\begin{corollary} \label{cor:reaper}
Suppose that $I$ is a beaded Nielsen path in $\bot(S)$ for some corridor
$S$ of a well-folded diagram, and suppose that all beads of $I$ are eventually 
bitten from the left  by beads in a single colour $\mu$.  Then, with the possible
exception of $B$ beads 
on the left end and one bead on the right (the final one bitten),
whenever $\mu$ bites a Nielsen bead in the future of $I$, 
it consumes it entirely.
\end{corollary}

\bigskip

\noindent{\bf{Proof of the Proposition \ref{OneSideStable}}}
\smallskip

\begin{proof}
If the immediate future of $I$ in $\top(S)$ were not left-edible, then Corollary
\ref{cor:reaper} would ensure that no bead in $I$ which is not bitten in $S$ is ever
bitten by $\mu$.  However, the assumption on the length of $I$ (and the Bounded 
Cancellation Lemma) ensure that there {\em are} beads in $I$ not bitten in $S$.
The fact that $I$ {\em is} left-edible therefore ensures that the future of $I$ in
$\top(S)$ is also left-edible.
\end{proof}

\bigskip

\noindent{\bf{Proof of the Proposition \ref{Reaper}}}
\smallskip

\begin{proof}
Let $S'$ be the corridor containing
the immediate past of $I^\lambda$.
Lemma
\ref{WholeNielsen} implies that in $\top(S')$ there is an edge $\rho$ in
$\mu$ which cancels a whole Nielsen path in the future of $I$.

Since $|I| \ge B +J$, there is a bead in $I$ not bitten in $\top(S)$.  The proof of Lemma \ref{WholeNielsen} now implies that there is a reaper as in the statement of the proposition.
\end{proof}

\section{Non-fast and Unbounded Beads} \label{C1Section}

With the technical exertions of the previous section behind us,
we are now able to return to the  main argument,
picking up the flow of Part \ref{Part:BG1} at Subsection \ref{NonConstantSubsect}. 
Thus our
next purpose is to reduce the task of bounding
the bead norm of the intervals $C_{(\mu,\mu')}$ to that of bounding
the lengths of certain long blocks of Nielsen atoms.  These blocks
are the analogue of the intervals $C_{(\mu,\mu')}(2)$ from
Part \ref{Part:BG1}, and will be the building blocks of the {\em{teams}} 
introduced in Section
\ref{TeamsSection} (in analogy with Section \ref{teamSec}).

\begin{definition} \label{Slowpep} \index{PEP!slow}
Suppose that $\rho = \gamma \nu \tau^k \overline{E_i}$ is a
\rpep \ (with $k \ge 0$).  We say that $\rho$ is {\em left-slow} if $\gamma$ is empty
 or a concatenation of left-slow beads.

There is an entirely analogous definition of {\em right-slow} \rpep s of the form
$\rho = E_i \overline{\tau}^k \overline{\nu}  \overline{\gamma}$.
\end{definition}

Often, we will just speak of {\em slow} \rpep s, since a single \rpep \ 
can only be left-slow or right-slow, but not both.

\begin{definition} \index{bead!unbounded}
Suppose that the bead $\rho$ is such that $f_{\#}(\rho)$ is not a Nielsen
bead.  Then the function $n \mapsto |f_{\#}^n(\rho)|$ grows at least linearly.
In this case, we call $\rho$ an {\em unbounded bead}.
\end{definition}

\begin{definition} \label{tame} \index{tameness}
A beaded path is called {\em right-tame} if all of its beads are
\gep s, slow \pep s, Nielsen paths and \atom s which do not
have a right-fast displayed edge in their immediate future.
\end{definition}

The next lemma follows immediately from the definition.
\begin{lemma}
$A_4(S_0,\mu)$ is a right-tame path.
\end{lemma}

\begin{lemma}
Suppose that $\alpha$ is a non-vanishing atom which is not right-fast.  Then either
all of the beads in $f_\#(\alpha)$ are Nielsen paths and \gep s, or else
the preferred future of $\alpha$ is parabolic.
\end{lemma}
\begin{proof}
The only modification to Lemma \ref{fofatom:BG3} is the 
exclusion of exponential edges in the second case, which is 
valid because such an edge would obviously contradict the fact 
that $\alpha$ is not right-fast.
\end{proof}

\begin{definition} \label{untrapped}\index{path!weight of!untrapped}
Suppose that $\sigma$ is a right-tame path.  The {\em untrapped
weight} of $\sigma$ is the largest $j$ so that $f_{\#}(\sigma)$ contains a
bead of weight $j$ which is not Nielsen.
\end{definition}

\begin{definition} \label{DefEi}
Suppose that, for some pair $(\mu,\mu') \in \mathcal Z$, the
untrapped weight of $C_{(\mu,\mu')}$ is $j$.
For each $1 \le i \le j$, define $\rho_i$
to be the leftmost bead in $C_{(\mu,\mu')}$ so that $f_{\#}(\rho_i)$
has an unbounded bead of weight at least $i$ that is not
HNP-bitten in the future of $S_0$.\footnote{Note that it is
possible that $\rho_i = \rho_{i+1}$ for some $i$.}

Let $\mathcal E_i$ denote those beads in $C_{(\mu,\mu')}$
from the right end up to and including $\rho_i$, and let
$\mathcal D_i = \mathcal E_i \smallsetminus \mathcal E_{i+1}$.
\end{definition}

The following is the analogue of Lemma \ref{C1Lemma}
\begin{lemma} \label{C1Lemma:BG3}
For all $1 \le i \le \omega$ there is a constant $C_1(i)$ so that
for each of the paths $C_{(\mu,\mu')}$ and decomposition into intervals
$\mathcal D_i$ as above, we have
\[	\| \mathcal D_i \|_\beta \le C_1(i) .	\]
\end{lemma}
\begin{proof} As far as possible, we try to
follow the proof of
Lemma \ref{C1Lemma}.  However, due to the phenomena
described in Section \ref{TrappedHNP}, the proof here is
somewhat more complicated.

We go forward to the time, $t$ say, which is one step before the moment
when $\mu'$ first starts to bite the preferred futures.  By virtue
of Remark \ref{BeadsDontVanish}, and the definition of 
$\mathcal D_i$, there are at least as many beads in the future
of $\mathcal D_i$ at time $t$ as there are in $S_0$.  Therefore,
it is sufficient to bound the number of beads in the future of 
$\mathcal D_i$ at time $t$; to ease the
notation, we write $\mathcal D_i$ for this future,
i.e.~pretend that $t=\time(S_0)$.

It is possible that there exist beads $\rho \in \mathcal D_i$ so that
$f_{\#}(\rho)$ has weight greater than $i$.  In such a case,
all of the beads in $f_{\#}(\rho)$ of weight greater than $i$ are Nielsen beads. 

Consider the highest weight $k$ for which there is a 
bead $\rho$ in $\mathcal D_i$ with $f_{\#}(\rho)$ of weight $k$,
and suppose that $k > i$.  Suppose first that $\rho$ has exponential
weight.  Then by Lemma \ref{DealWithExp} either $\mathcal D_i$ has
bead norm at most $B$ (and length at most $\ell=JB(B+1)$), or else
$\rho$ is entirely consumed when it is bitten.  In the first case
$\rho$ is the leftmost bead in $\mathcal D_i$, and also in $C_{(\mu,\mu')}$.
A similar argument applies when $\rho$ has parabolic weight.

Thus, excluding cases where $|\mathcal D_i|<\ell$, we may
 treat  the Nielsen beads of weight higher than $i$ as indivisible
units, which are entirely consumed when bitten.  We are therefore
in the situation of the proof of Lemma \ref{C1Lemma}, where the
unbounded beads in $\mathcal C_i$ grow apart at a linear rate,
and so must be cancelled quickly.  Otherwise, the proof is entirely
parallel to the one from Part \ref{Part:BG1}.
\end{proof}

We are trying to reduce the task of bounding the bead norm to
that of bounding the size of intervals consisting entirely
of Nielsen beads, which are each consumed by a reaper.  In order
to make this reduction, we still have some HNP-biting to deal with.
In order to deal with this, we need an analogue of Lemma \ref{G34pics}.

Recall that $L$ is the maximal length of $f(E)$ where $E$ is an edge
in $G$.

\begin{proposition} [cf. Lemma \ref{G34pics}] \label{C4Lemma}
There is a constant $C_4$ depending only on $f$ which satisfies the
following properties.
If $I$ is an interval on $\top(S)$ labelled by a beaded path all of whose beads
are Nielsen \atom s, then the path labelling the past of $I$ in $\bot(S)$
is of the form $u\alpha v$ where $\alpha$ is a beaded path all of whose beads
are Nielsen \atom s and $|u|$ and $|v|$ are less than $C_4$. 

If the past of $I$ begins (respectively ends) with a point fixed by
$f$, then $u$ (respectively $v$) is empty.

In particular, $| I | \le |\alpha| + 2LC_4$.
\end{proposition}
\begin{proof} 
The interval $I \subset \top(S)$ is a beaded path, all of whose beads are  
Nielsen paths
of length at most $J$.  Therefore, along $I$ there are points where
$I$ admits a hard splitting and these points occur with a frequency of
at least one every $J$ edges.  Since these points are vertices, the
set of labels of points at which the splitting occurs is finite.
Consider the path from $\top(S)$ to $\bot(S)$ starting from one of
these vertices.  The label of this path is $w\bar{t_i}$ where
$w$ is a (possibly empty) path in $G$ of length at most $L$, and
$t_i$ is one of the edges from the mapping torus $M(f)$.  (We are about
to use a finiteness argument and it will be important that the
repetition we infer includes the labels of the points on
$\bot(S)$.  Thus it is important which of the $t$-edges this path
includes.)

Since the data we record --- the label of the vertex on
$\top(S)$, the path $w\bar{t_i}$ and the label of the end of this path
on $\bot(S)$) --- range over a finite set, there is a constant $C'$ such
that in the interval within $C'$ vertices of the left end of
$I$ there will be repetition of these data.  Since the vertices
occur at least every $J$ edges, this repetition occurs within $C'J$ of
the left end of $I$.

Once we have found this repetition, we have an interval $\lambda
\subset \bot(S)$, an interval $\eta \subset \top(S)$ and a path $w_0$
of length at most $L$ such that $f_{\#}(\lambda) = w_0 \eta
\bar{w_0}$. Therefore, the free homotopy class of $f_{\#}(\lambda)$ is
the same as that of $\eta = f_{\#}(\eta)$, since $\eta$ is a
beaded path all of whose beads are Nielsen paths. Since $f$ is a homotopy equivalence,
the free homotopy class of $\lambda$ must be the same as that of
$\eta$.

Suppose that $\eta = p_1 \dots p_m$ where each $p_i$ is an indivisible
Nielsen path.  Now, $\lambda$ is tight, so $\lambda = \sigma
p_ip_{i+1} \dots p_m p_1 \dots p_{i-1} \bar{\sigma}$, for some path
$\sigma$. Thus, if `$\sim$' denotes free homotopy,
\[ f(\lambda) \sim f_{\#}(\sigma) p_i \dots p_{i-1}
f_{\#}(\bar{\sigma}), \] 
which tightens to
\[ w_0 p_1 \dots p_m \bar{w_0}.  \]
By the Bounded Cancellation Lemma, tightening the path $f(\lambda)$ as
written above reduces the length of $f_{\#}(\sigma)$ by less than $B$,
and the result has length at most $2L + |\eta|$.  This implies that
$|f_{\#}(\sigma)| < L+B$.  Therefore, $\| \sigma
\|$ is bounded, and by a small increase we may also assume that
$i=1$.  By considering only one vertex out of every $B(L+B)$,
we can find such a path $\eta$ where there is some $p_j$ in the
middle of $\lambda$ such that the path from the copy of $p_j \subset \top(S)$ to the copy of
$p_j \subset \bot(S)$ is a single edge labelled $t$, for some $j$.

We have argued that, for some path $\eta$ of bounded length which lies
on the left end of $I$, the past of $\eta$ is of the form $u \eta u'$
where $|u|$ and $|u'|$ are bounded, and the paths from the splitting
points in $\eta \subset I$ to $\bot(S)$ consist of single edges
labelled $t$.

Consider the analogous situation on the right end of $I$.  We can
find a path $\eta' \subset I$ lies at the right end of $I$
such that the past of $\eta'$ is of the form $v' \eta' v$
where $|v|$ and $|v'|$ are bounded and the paths from the
vertices of $\eta' \subset I$ to $\bot(S)$ consist of single
edges labelled $t$.  

Consider the paths along $\bot(S)$ and $\top(S)$ from the left
end of $\eta$ to the right end of $\eta'$.  We have a path $\rho
\subset \bot(S)$ with fixed points of $f$ on either end which maps to
a Nielsen path $f_{\#}(\rho) \subset I \subset \top(S)$.  The same
argument as in the proof of Lemma \ref{PreNielsen} then shows that
$\rho = f_{\#}(\rho)$.  Hence $\rho$ is a beaded path, all of whose
beads are Nielsen
paths, and the paths $u$ and $v$ on either side of $\rho$ are
of bounded length as required.  This proves the first assertion in the
statement of the lemma.

The second assertion follows similarly, and the final assertion
follows immediately from the first. 
\end{proof}

Consider a pair $(\mu,\mu') \in \mathcal Z$, and recall the
definition of the subintervals $\mathcal E_i$ from Definition \ref{DefEi}.

\begin{proposition} \label{GetC2}
There is a constant $C_5$, depending only on $f$ so that the
following holds.  For each $(\mu,\mu') \in \mathcal Z$, the
interval $C_{(\mu,\mu')} \smallsetminus \mathcal{E}_1$ 
in $A_4(S_0,\mu)$
has the form $u N v$ where $u$ and $v$ are such that
$\| u \|_{\beta}, \| v \|_\beta \le C_5$ and $N$ is a beaded path all
of whose beads are Nielsen beads.
\end{proposition}
\begin{proof}
By Lemma \ref{OneHNP}, for each adjacency of colours $(\mu,\mu')$
there can only be one bead in $\mu(S)$ which is eventually HNP-bitten
by $\mu'$.

The result now follows from Proposition \ref{C4Lemma} and the 
definition of $\mathcal E_1$.
\end{proof}

\begin{definition}
For $(\mu,\mu') \in \mathcal Z$, define $C_{(\mu,\mu')}(2):=
N$, the beaded Nielsen path from Proposition \ref{GetC2}.
\end{definition}

The sum of our arguments to this point has reduced the
task of bounding the sum of the bead norms of the 
intervals $\mu(S_0)$ in $S_0$ to that of bounding the
sum of the lengths of the intervals $C_{(\mu,\mu')}(2)$ for
pairs $(\mu,\mu') \in \mathcal Z$.

We summarise
the results from this section as follows.

\begin{proposition} \label{PointofC1}
There is a constant $C_1$, depending only on $f$, so that
\[      \| C_{(\mu,\mu')}\|_\beta \le \|C_{(\mu,\mu')}(2)\|_{\beta} + C_1.      \]
\end{proposition}

\begin{remark} \label{BeadtoLength}
Since the intervals $C_{(\mu,\mu')}(2)$ consist entirely of Nielsen
beads, we have the following obvious relationship between
length and bead norm:
\[ |C_{(\mu,\mu')}(2)| \le \|C_{(\mu,\mu')}(2) \|_\beta
\le J|C_{(\mu,\mu')}(2)|.       \]
Therefore, in order to finish the bound on bead norm, it is sufficient
to bound the total lengths of the intervals $C_{(\mu,\mu')}(2)$.
\end{remark}

It is important for the remainder of the paper that the path
$C_{(\mu,\mu')}(2)$ is a beaded path that consists entirely of
Nielsen atoms.  This is a stronger statement than just asserting
it is a Nielsen path, since we require a decomposition into beads
of uniformly bounded size, each of which is a Nielsen path.
This makes the path $C_{(\mu,\mu')}(2)$ very similar to the long
blocks of constant letters which played such a prominent role in
Part \ref{Part:BG1}

At this point the reader may benefit from consulting 
Section \ref{A4sec}, which outlines the strategy for the remainder 
of the proof of the Main Theorem (the strategy from the positive
case still holds here).  For the remainder of this part of the book, we will
mostly continue without reminding the reader of this strategy.

\section{The Pleasingly Rapid Disappearance of Colours} 
\label{PincerSection}

We are now at the point in our arguments where we need to formulate
and prove the Pincer Lemma, as in Section \ref{ConstantSection}.
In Part \ref{Part:BG1} the Pincer Lemma was proved by counting colours
which {\em essentially vanished}, which is to say they came to consist
entirely of constant letters.  For positive automorphisms, this is a
well-defined event and can only occur once for each colour.  For
general automorphisms, the analogues of constant letters are
indivisible Nielsen paths.  However, since Nielsen paths can contain
non-constant edges, indivisible Nielsen paths are not
indivisible in an absolute sense (the terminology refers to the fact that an indivisible
Nielsen path cannot be split into two Nielsen paths).  Thus, it is
possible that a colour can be labelled by a Nielsen path at some time
$t$ but not at some later time $t+k$. There are two ways to circumvent
this problem.  The first is to concentrate on the times when a colour
decreases in weight, whilst the second is to focus on the times when a
colour becomes Nielsen and seek compensation when a colour 
subsequently ceases to
be Nielsen. We mostly pursue the second idea but there are aspects of the first also.

The version of the Pincer Lemma which we need in this part of the book is
Theorem \ref{PincerLemma:BG3}.

The ideas in the proof of the Pincer Lemma here are very similar to
those in Part \ref{Part:BG1} but the execution is somewhat different.

\begin{definition} \label{StablyNielsen}
Suppose that $I$ is a non-empty beaded Nielsen path
 and that $U$ and $V$ are
beaded 
paths.  We say that $I$ is \index{Nielsen path!stability of}{\em stably Nielsen} in the
path $UIV$ if the future\footnote{as defined in (\ref{stackDiags})}
of $I$ in $f_{\#}(UIV)$ is also a non-empty Nielsen beaded path.

Suppose that $\mu_1, \mu_2$ and $\mu_3$ are colours in a 
well-folded diagram and that the intervals
$\mu_1(S), \mu_2(S)$ and  $\mu_3(S)$ are non-empty and
adjacent in $\bot(S)$. If
$\widecheck{\mu_2(S)}$ is a non-empty Nielsen path, then we say 
that $\mu_2(S)$ is {\em stably Nielsen} if, in the above sense,
 $\widecheck{\mu_2(S)}$ is stably
Nielsen in $\widecheck{\mu_1(S)}
\widecheck{\mu_2(S)} \widecheck{\mu_3(S)})$.
\end{definition}

\begin{lemma}[Relative Buffer Lemma] \label{RelBuffer}\index{Buffer Lemma!relative}
Let $i \in \{ 1 , \dots , \omega-1 \}$ and let $I \subset \bot(S)$ be an
edge-path labelled by edges in $G_i$. Suppose that the colours
$\mu_1(S)$ and $\mu_2(S)$ lie either side of $I$, adjacent to it.
Provided that the whole of $I$ does not die in $S$, no edge in the future 
of $\mu_1(S)$  with label in
$G \smallsetminus G_i$ will ever cancel
with an edge in the future of $\mu_2(S)$  with label in
$G \smallsetminus G_i$.
\end{lemma}
\begin{proof}
Given Lemmas \ref{EndStab1} and \ref{EndStab2},  the proof of Lemma \ref{BufferLemma} applies modulo changes of terminology.
\end{proof}

We now need the following `two-sided' version of Proposition \ref{OneSideStable}.

\begin{lemma} \label{StableNielsenLemma}
Let $\mu_1$, $\mu_2$, $\mu_3$ and $S$ be as in Definition
\ref{StablyNielsen}, and suppose that $\mu_2(S)$ is stably Nielsen.  Then
for all corridors $S'$ in the future of $S$, if $\mu_1(S')$ and
$\mu_3(S')$ are nonempty then $\mu_2(S')$ is a (possibly empty)
Nielsen path.
\end{lemma}
\begin{proof} 
Whilst $\mu_1(S')$ and $\mu_3(S')$ are non-empty, any bead in $\mu_2$
which is bitten must be bitten by a bead coloured either $\mu_1$ or $\mu_3$.  Let
$I_1$ be the set of (Nielsen) beads in $\mu_2(S)$ which are eventually bitten
by a bead coloured $\mu_1$ (and are bitten whilst $\mu_1(S')$ and $\mu_3(S')$ are non-empty).
Define $I_2$ to be those beads in $\mu_2(S)$ which are bitten by a bead coloured $\mu_3$
(with the same proviso).

Suppose that $I_1$ and $I_2$ are non-empty.  They form intervals, and $I_1$
is to the left of $I_2$.

Proposition \ref{Reaper}, and the fact that $\mu_2(S)$ is stably Nielsen, implies
that unless $I_1$ is immediately consumed there is a left reaper coloured
$\mu_1$ associated to $I_1$, and similarly there is a right reaper coloured
$\mu_3$ associated to $I_2$.  The properties of reapers in Definition \ref{ReaperDef}
imply the result.  

In case one or both of $I_1$ and $I_2$ are empty (or immediately consumed), 
there is at most one reaper to consider, but the result follows in the same way.
\end{proof}

\begin{lemma}[Buffer Lemma] \label{BufferLemma:BG3} \index{Buffer Lemma}
Suppose, for some corridor $S$ in a well-folded diagram, that
$I \subset \bot(S)$ is a beaded Nielsen path and
that $\mu_1(S)$ and $\mu_2(S)$ lie either side of $I$,
 immediately adjacent to it.
 Suppose further that $\check{I}$ is stably Nielsen in $\check{\mu_1(S)}\check{I}
\check{\mu_2(S)}$.  Provided that the whole of $I$ does not die in $S$, no
bead in $\mu_1(S)$ can be eventually bitten by a bead coloured $\mu_2$ (and vice versa), 
unless it is (eventually) HNP-bitten.
\end{lemma}
\begin{proof}
Given Lemmas \ref{EndStab1}, \ref{EndStab2}
and \ref{StableNielsenLemma}, and the exclusion of HNP-biting, the proof 
of Lemma \ref{BufferLemma} applies.
\end{proof}

The proof of the following lemma follows that
of Lemma \ref{BufferLemma}.

\begin{lemma}[Weighted Buffer Lemma] \label{WtBuffer} \index{Buffer Lemma!weighted}
Suppose, for some corridor $S$ in a well-folded diagram, 
that $I \subset \bot(S)$ is a beaded path consisting
of Nielsen beads and beads of weight at most $i$, and 
that $\mu_1(S)$ and $\mu_2(S)$ lie either side of $I$,
 immediately adjacent to it.
Suppose further that the only  beads of $f_{\#}(\check{\mu_1(S)}\check{I}\check{\mu_2(S)})$ 
that are in the future of  $I$ and have weight greater than $i$ are
Nielsen beads.  

Then, provided that the whole of $I$ does not die in $S$, no bead
in $\mu_1(S)$ can be eventually bitten by a bead coloured $\mu_2$ (and vice versa), 
unless it is (eventually) HNP-bitten.
\end{lemma}

\subsection{The Two Colour Lemma}

Example \ref{Ex:HNP} can be used to construct examples where the
above two results are false if HNP-biting is not excluded. The same
is true of the results in this section. This accounts for the caution
that the reader will note in Sections
\ref{TeamsSection}, \ref{BonusSection} and 
\ref{LongGepsandPepsSection}, where we are careful to ensure
that the Pincer Lemma is applied only to  
pincers that  involve no HNP-biting.

\begin{definition} [Stable $f$-neutering]  \index{neutering!stable}
Suppose that $U$ and $V$ are beaded paths, that for some $k$ the futures of $V$ in
$f_{\#}^k(UV)$ and $f_{\#}^{k+1}(UV)$ are Nielsen, but that the future of $V$ in
$f_{\#}^{k-1}(UV)$ contains a non-Nielsen bead.  

Denote the futures of $U$ and $V$ in $f_{\#}^{k-1}(UV)$ by $U^{k-1}$ and $V^{k-1}$,
respectively.  Let $\beta$ be the  rightmost non-Nielsen bead in $f_{\#}(V^{k-1})$.
If the biting of $\beta$ in the tightening of $f_{\#}(U^{k-1})f_{\#}(V^{k-1})$ to form $f_{\#}^k(UV)$
is not HNP-biting then we say that $U$ {\em stably left $f$-neuters} $V$
in $k$ steps.

The definition of {\em stable right $f$-neutering} is identical with the roles of $U$ and $V$
reversed, and when we are unconcerned about the direction we will refer simply
to {\em stable $f$-neutering}.
\end{definition}
In the light of Proposition \ref{OneSideStable}, once stably $f$-neutered, the subsequent futures
of $V$ remain beaded Nielsen paths.

\begin{proposition} [Two Colour Lemma, cf. Proposition \ref{TwoColourLemma}]
\label{TwoColour} \index{Two Colour Lemma}
There exists a constant $T_0$, depending only on $f$, so that if $U$
and $V$ are beaded
paths and $U$ stably $f$-neuters $V$ then it does so in at most $T_0$ steps.
\end{proposition}

\begin{proof}
Denote the future of $U$ in $f_{\#}^i(UV)$ by $U^i$ and the future
of $V$ by $V^i$.

As in the proof of Proposition \ref{TwoColourLemma}, we will decompose each
of the paths $V^i$ into an {\em{unbounded part}} and a 
{\em{bounded part}}.  The bounded
part will be an interval on the right end of $V^i$ whose immediate (abstract)
future is a beaded Nielsen path.  The unbounded interval lies on the left end
of $V^i$, and we will bound its length.  

This would
be a straightforward adaptation of the proof from Part \ref{Part:BG1} if Proposition
\ref{PointofC1} provided a bound of the length of that part of $C_{(\mu,\mu')}$
not contained in $C_{(\mu,\mu')}(2)$.  However, the bound in Proposition
\ref{PointofC1} is just a bound on bead norm.
Thus, we need to deal with the possibility of long \gep s and \pep s.

The following enumerated claims will together yield an upper bound
on the length of the unbounded part of $V^i$, which in the course of 
the proof will be decomposed into $V^i_{\text{fast}}$ and $V^i_{\text{nc}}$ 

Three of the claims concern the existence of a constant $k_j$ that
depends only on $f$; we use the abbreviation $\exists k_j=k_j(f)$.

\smallskip

{\bf Claim 1:}  $\exists k_1=k_1(f)$ such that any \gep \ in $V^i$
has length less than $k_1$.
\smallskip

This follows in a straightforward way from the Buffer Lemma \ref{BufferLemma:BG3}
and the fact that the obvious preferred future of the rightmost edge in any \gep \ in
$V^i$ must eventually cancel with an edge from the future of $U^i$.

\smallskip

Next we consider long \pep s in $V^i$.  Suppose that $\rho$ is a \pep \ in
$V^i$.  Then the label on ${\rho}$ or ${\bar{\rho}}$ 
has the form $E\bar{\tau}^k\bar{\nu}\gamma$,
where $\tau$ is Nielsen path, $f(E) = E \odot \tau^m$ and $\bar{\gamma}\nu$ is
a terminal segment of $\tau$.  We consider a number of different cases.
First we dismiss a case that 
follows immediately from Lemma \ref{C0Lemma} and from the fact that
exponential edges are left-fast:

\smallskip

{\bf Claim 2:}  If $\check{{\rho}} = E\bar{\tau}^k\bar{\nu}\gamma$ and $\gamma$ is
an exponential edge then the right end of $\rho$ lies within $C_0$
of the left end of $V^i$.

\smallskip

Next we consider $V^i_{\text{fast}}$, which is defined 
to consist of those beads from the left end of $V^i$
up to and including the rightmost bead in $V^i$ whose immediate (abstract) future
contains a left-fast bead.  
\smallskip

{\bf Claim 3:} $\exists k_2=k_2(f)$ such that 
$|V^i_{\text{fast}}|\le k_2$.
\smallskip

This follows immediately from Lemma \ref{C0Lemma} unless the rightmost
bead in $V^i_{\text{fast}}$ is a \pep .  (Note that this rightmost bead is not a \gep ,
since a \gep \ does not have a left-fast bead in its immediate abstract future.)

Suppose, then, that the rightmost bead in $V^i_{\text{fast}}$ is a \pep , say $\rho$.
If $\check{{\rho}} = E\bar{\tau}^k\bar{\nu}\gamma$, then 
we are done by Claim 2. So
suppose that $\check{{\rho}} = \bar{\gamma} \nu \tau^k \bar{E}$.  Let
$\e$ be the edge in $\rho$ whose label is $\bar{E}$.  The preferred
future of $\e$ is to be cancelled by an edge in the future of $U^i$.  
By an obvious finiteness argument (as in the proof of Proposition \ref{TwoColourLemma}), 
there is a constant $p$ so that the path $V^p$ contains
no left-fast beads.  This gives a bound on the amount of time before the future
of $\rho$ is bitten, and hence a bound on the amount that the future of $\rho$ can
shrink before then.  Suppose that $V^j$ is the first future of $V^i$ in which the future 
of $\rho$ has been bitten.  Because the preferred future of $\e$ is to be cancelled, 
Proposition \ref{Reaper} and the Buffer Lemma \ref{BufferLemma:BG3} imply that
the length of the future in $V^j$ of $\rho$ is bounded above by a 
constant depending only on $f$.

The required bound on $|V^i_{\text{fast}}|$
is now at hand:  Lemma \ref{C0Lemma} bounds the length of 
$V^i_{\text{fast}} \smallsetminus \rho$, and the combination of the bound on $j$
and the bound on the length of the future of $\rho$ in $V^j$ gives a bound on
the length of $\rho$.  This completes the proof of Claim 3.  We remark that the
above argument also gives a bound on the amount of time it takes for $V^1_{\text{fast}}$
to be entirely consumed.

\smallskip

We now define a set $V^i_{\text{nc}}$ as follows: Let $\rho_{\text{nc}}$ be the rightmost
bead in $V^i$ whose immediate abstract future is not Nielsen.  We define
$V^i_{\text{nc}}$ as follows:
\begin{enumerate}
\item if $\rho_{\text{nc}} \in V^i_{\text{fast}}$ then $V^i_{\text{nc}} = \emptyset$;
\item if $\rho_{\text{nc}}$ is not a \pep , then $V^i_{\text{nc}}$ consists of those
beads from (but not including) the rightmost bead in $V^i_{\text{fast}}$ up to
and including $\rho_{\text{nc}}$;
\item if $\rho_{\text{nc}}$ is a \pep \ with label of the form $\bar{\gamma} \nu \tau^k \bar{E}$
or $\rho_{\text{nc}}$ is a \pep \ with label of the form $E\bar{\tau}^k\bar{\nu}\gamma$
and $\gamma$ is not a Nielsen path, then $V^i_{\text{nc}}$ consists of those
beads in $V^i$ from (but not including) the rightmost bead in $V^i_{\text{fast}}$
up to and including $\rho_{\text{nc}}$;
\item \label{LooseEndRight} finally, if $\rho_{\text{nc}}$ is a \pep \ with label of the form
$E \bar{\tau}^k\bar{nu} \gamma$
and $\gamma$ is either empty or a Nielsen path, then $V^i_{\text{nc}}$ consists of
that interval from (but not including) the rightmost bead in $V^i_{\text{fast}}$ up
to and including the leftmost edge in $\rho_{\text{nc}}$ (the label of this
leftmost edge is $E$).
\end{enumerate}
Note that in Case \ref{LooseEndRight} the bead $\rho_{\text{nc}}$ is certainly
not contained in $V^i_{\text{fast}}$.

\smallskip

{\bf Claim 4:} $\exists k_3=k_3(f)$  such that  $|V^i_{\text{nc}}|\le k_3$.
\smallskip

The proof of Claim 3 above established an upper bound on the time before
all of $V^i_{\text{fast}}$ is entirely consumed, and hence also on the time
before the future of $V^i_{\text{nc}}$ begins to be consumed.
We now follow the proof of Lemma \ref{C1Lemma:BG3}, which
establishes an upper bound 
 on the time that can elapse before the final non-constant bead in $V^i$ is bitten.
We will be done if we can bound this time from below by a positive
constant times $|V^i_{\text{nc}}|$. 

In the current setting, we have non-constant beads in $V^i_{\text{nc}}$ 
that may not be growing apart like those in the proof of Lemma \ref{C1Lemma:BG3}.\footnote{This is because we are now measuring length rather than bead-norm.}
But  there {\em{is}} a lower
bound on the rate at which the surviving futures of these
beads can come together.
Hence the length of $V^i_{\text{nc}}$ provides a lower bound on 
the amount of time that must elapse before $V^j$ becomes  stably Nielsen,
since the future of $V^i_{\text{nc}}$ must be entirely consumed before
this time. (Note that in Case \ref{LooseEndRight},
the preferred future of the edge $\bar{E}$ in $\rho_{\text{nc}}$ must be 
eventually consumed by the future of $U^i$.)  This proves Claim 4.

\smallskip

The {\em unbounded part} of $V^i$ is the union of $V^i_{\text{fast}}$ and
$V^i_{\text{nc}}$, whilst the {\em bounded part} is the remainder of $V^i$.
The sum of the previous four claims bound the length of
the unbounded part of $V^i$ by a constant that depends only on $f$.

There is a similar bound on the number of edges in $U^i$ that have
an edge in their future that cancels with an edge in the future of $V^i$.
(Here we need the hypothesis that the path $V^k$ becoming stably Nielsen
does not arise from HNP-biting.)

At this stage, we can follow the proof of Proposition \ref{TwoColourLemma} directly.
After an amount of time bounded by a constant
that depends only on $f$, either the future of $V$ 
becomes stably Nielsen or empty, or else there is a repetition 
of the following data: (i) the unbounded part of
$V^i$ plus the leftmost $B+J$ edges of the bounded part; (ii) a
terminal segment of $U^i$ containing all of the edges that can ever
interact with the future of $V$.  Once we have such a repetition, if
the future of $V$ has not become stably Nielsen or vanished then it never will,
contrary to hypothesis.
\end{proof}

We need a weighted version of neutering and the two-colour lemma.

\begin{definition} [$(f,i)$-neutering]\index{neutering!weighted}
Fix $i \in \{ 1, \ldots , \omega \}$ and let $U$ and $V$ be beaded paths.
Suppose  that for some $k$ the future of $V$ in $f_{\#}^k(UV)$ has weight
less than $i$, but that the future of $V$ in $f_{\#}^{k-1}(UV)$ has weight at 
least $i$.

Denote the futures of $U$ and $V$ in $f_{\#}^{k-1}(UV)$ by $U_{k-1}$ and
$V_{k-1}$, respectively.  Let $\beta$ be the rightmost bead in $f_{\#}(V_{k-1})$
of weight at least $i$.  If the biting of $\beta$ in the tightening of
$f_{\#}(U_{k-1})f_{\#}(V_{k-1})$ to form $f_{\#}^k(UV)$ is not HNP-biting
then we say that $U$ {\em $(f,i)$-neuters} $V$ in at most $k$ steps.
\end{definition}

\begin{proposition}[Weighted Two Colour Lemma] \label{WeightTwo}\index{Two Colour Lemma!weighted}
There exists a constant $T_0'$, depending only on $f$, so that for any
$i \in \{ 1 , \ldots , \omega \}$, if $U$ and $V$ are beaded
paths and $U$ $(f,i)$-neuters $V$ then it does so in at most $T_0'$ steps.
\end{proposition}
\begin{proof}
We decompose the futures of $U$ and $V$ in $f_{\#}^k(UV)$
as in Lemma \ref{TwoColour}.

The proof is similar to that of Lemma \ref{TwoColour}, 
except that when we appeal to the proof of
Proposition \ref{C1Lemma:BG3} we assume that we have a path
$\mathcal E_j$ with $j \ge i$.
Otherwise, the proof of Lemma \ref{TwoColour} above and that
of Proposition \ref{TwoColourLemma} can now be followed {\em mutatis
mutandis}.
\end{proof}

By replacing $T_0$ by $T_0'$ if necessary, we may assume that $T_0 \ge
T_0'$. We henceforth make this assumption. 

\subsection{The disappearance of colours: Pincers and implosions}

\begin{definition} \label{PincerDef:BG3}
Consider a pair of non-constant edges $\e_1$ and $\e_2$ which cancel
in a corridor $S_t$ of $\Delta$, and suppose that, for $i = 1,2$,  the
immediate past of $\e_i$ lies in a bead of some $\mu_i(S_t)$ that is 
either a unbounded atom, a
\gep \ or a \pep .  Suppose further that the cancellation of
$\e_1$ and $\e_2$ is not HNP-cancellation, and that
$\mu_1 \neq \mu_2$.  Consider the paths $p_1,
p_2$ in $\F \subset \Delta$ tracing the histories of $\e_1$ and
$\e_2$.  Suppose that at time $\tau_0$ the paths $p_1$ and $p_2$ lie in a
common corridor $S_b$.  Under these circumstances, we define the
\index{pincer} {\em
  pincer} $\Pin = \Pin(p_1,p_2,\tau_0)$ to be the sub-diagram of
$\Delta$ enclosed by the chains of $2$-cells along $p_1$ and $p_2$,
and the chain of $2$-cells connecting them in $S_b$.

We define $S_{\Pin}$ to be the earliest corridor of the pincer in
which $\mu_1(S_{\Pin})$ and $\mu_2(S_{\Pin})$ are adjacent.  Define
$\wt{\chi}(\Pin)$ to be the set of colours $\mu \not\in \{ \mu_1 ,
\mu_2 \}$ such that there is a $2$-cell in $\Pin$ coloured
$\mu$. Finally, define\index{pincer!life of}
\[      \life(\Pin) = \time(S_{\Pin}) - \time(S_b) . \]
\end{definition}
See Section \ref{ConstantSection} for illustrative pictures.

\begin{proposition}[Unnested Pincer Lemma, cf. Proposition \ref{prePincerLemma}] \label{prePincer}
\index{Pincer Lemma}

\ 

\noindent There exists a constant $T_1$, depending only on $f$, 
such that for any pincer $\Pin$
\[      \life(\Pin) \le T_1(1 + |\wt{\chi}(\Pin)|)      .       \]
\end{proposition}

In the proof of Proposition \ref{prePincerLemma} (Regular Implosions) 
the strategy was to identify a constant $T_1$ such
that over each period of time of length $T_1$ within a pincer, at
least one colour became constant.  There are a number of impediments
to implementing this strategy in the current situation.   The first is
that Nielsen paths can consist of edges which are not constant edges,
so if a colour {\em becomes Nielsen} then it may cease to be
Nielsen at some stage in the future.  In order to overcome this
impediment, we make the following

\begin{definition}
Suppose that for some colour $\mu$ and some corridor $S$, the path
$\check{\mu(S)}$ is stably Nielsen, and let $\nu_1$ and $\nu_2$ be the colours
immediately on either side of $\mu$ in $S$.  If there is some
corridor $S'$ in the future of $S$ in which $\check{\mu(S')}$ is not
Nielsen and $S'$ is the earliest such corridor, then we say that \index{colours!resuscitated}
$\mu$ is {\em resuscitated} in $S'$.  By Lemma
\ref{StableNielsenLemma}, at least one of $\nu_1$ and $\nu_2$ is not adjacent
to $\mu$ in $S'$,  so either $\nu_1(S')$ or $\nu_2(S')$ is empty.
If  $\nu_i(S')$ is empty, we say that $\nu_i$  {\em sacrifices
itself} for $\mu$.  
\end{definition}

\begin{remark}
A colour can sacrifice itself for at most one colour.

A colour may become stably Nielsen and be resuscitated a number of
times, but a different colour must sacrifice itself for each resuscitation.

The concept of `becoming stably Nielsen' is analogous to that of a colour
`essentially vanishing' in Section \ref{ConstantSection}.  However, the concept
of `resuscitation' does not have an analogue in Part \ref{Part:BG1}.
\end{remark}

Fix a pincer $\Pin$ and assume that $\life(\Pin) > 1$.  The strategy
to prove Proposition \ref{prePincer} is to identify a constant $T_1$
so that  during the life of $\Pin$, in each
$T_1/2$ steps of time there is a colour that becomes
stably Nielsen (perhaps vanishing) 
In order to obtain the bound in the statement of
Proposition \ref{prePincer},
we then count the colours which
become stably Nielsen or vanish, and the colours which sacrifice themselves 
for those that are
resuscitated.  A colour can therefore be counted twice
-- once for disappearing (or for the last time it becomes stably
Nielsen), and once as a sacrifice -- but no colour is counted more than
twice.  Thus Proposition \ref{prePincer} is an immediate consequence of
the following result whose proof will occupy the remainder of this
subsection.

\begin{proposition} \label{pSmall}
There is a constant $T_1$, depending only on $f$, so that for any
pincer $\Pin$ in a minimal area van Kampen diagram over $M(f)$, in any
interval of time of length $T_1 / 2$, at least one colour in
$\tilde{\chi}(\Pin)$ becomes stably Nielsen or vanishes.
\end{proposition}

\begin{definition} [$p$-implosive arrays] \label{implosive}
Let $p$ be a positive integer and $S$ a corridor. \index{implosive array}  A {\em
  $p$-implosive array} of colours in  $S$ is an ordered tuple $A(S) = [
  \nu_0(S), \ldots ,   \nu_r(S) ]$, with $r > 1$, such that
\begin{enumerate}
\item each pair of colours $\{ \nu_j, \nu_{j+1} \}$ is separated in
  $S$ only by a stably Nielsen (or empty) path;
\item in each of the corridors $S = S^1, S^2, \ldots , S^p$ in the
  future of $S$, no $\nu_j(S^i)$ is empty or a stably Nielsen path,  $j = 1
  , \ldots , r-1$;
\item in $S^p$, {\em either} an edge coloured $\nu_0$ from a
   unbounded \atom , a \gep \ or a
   \pep \ cancels with an edge coloured $\nu_r$ from a
   unbounded \atom , a \gep \ or a
  \pep \ (and hence the colours $\nu_j$ with $j = 1, \dots , r-1$ are
  consumed entirely), or {\em else} each of the colours $\nu_j$ ($j =
  1, \dots , r-1$) become stably Nielsen or vanish, while $\nu_0$ and $\nu_r$
  are not Nielsen in $f_{\#}(\widecheck{\nu_0(S^p)} \cdots
  \widecheck{\nu_r(S^p)})$ (although they may nevertheless become stably Nielsen
  or even disappear in $S^p$ because of colours external to the
  array).
\end{enumerate}
\end{definition}
Arrays satisfying the first of the conditions in (3) are said to be of
{\em Type I}, and those satisfying the second condition are said to be of
{\em Type II}.  (These types are not mutually exclusive).

The {\em residual block} of an array of Type II is the stably Nielsen
path which lies between $\nu_0(S^p)$ and $\nu_r(S^p)$ (if either
$\nu_0(S^p)$ begins or $\nu_r(S^p)$ ends with an interval of
Nielsen \atom s include these  in the residual block). Note that the
residual block may be empty.  The {\em
  enduring block} of the array is the set of stably Nielsen paths in
$\bot(S)$ that have a future in the residual block.

Note that there may exist some \index{colours!unnamed}
{\em unnamed colours} between $\nu_j(S)$ and
$\nu_{j+1}(S)$; if they exist, these form a stably Nielsen path.

\begin{remark}
Let $[ \nu_0(S), \dots , \nu_r(S)]$ be a $p$-implosive array.
\begin{enumerate}
\item Any $q$-implosive sub-array of $[ \nu_0(S), \dots , \nu_r(S)]$
  has $q = p$.
\item If an edge of $\nu_i$ cancels with an edge of $\nu_j$ and $j-i >
  1$, then this cancellation can only take place in $S^p$.  If the
  edges cancelling come from displayed unbounded \atom s, \gep s or
  \pep s, then the sub-array $[ \nu_i(S), \dots , \nu_j(S)]$ is
  $p$-implosive of Type I.
\item If $u, v$ and $w$ are beaded
  edge-paths such that $u$, $v$ and $f_{\#}(uwv)$ are Nielsen paths
  then $w$ is a Nielsen path.  It follows that the residual block of
  any array of Type II contains edges from at most two of the colours
  $\nu_j$, and if there are two colours then they 
are consecutive,  $\nu_j,\,\nu_{j+1}$.
\item Likewise, the enduring block of an implosive array of
  Type II is an interval involving at most two of the $\nu_j$ and if
  there are two such colours they must be consecutive.
\end{enumerate}
\end{remark}

\begin{lemma} \label{ExistsImpArray}
Let $\Pin$ be a pincer.  The ordered list of colours along each
corridor before $\time(S_{\Pin})$ in a pincer $\Pin$ must contain
a $p$-implosive array for some $p$.
\end{lemma}

\begin{proof}
The definition of $p$-implosive array is designed so that when a
colour becomes stably Nielsen (or disappears) in a pincer there is a
$p$-implosive array.  See the proof of Lemma \ref{haveImp} for
more details.
\end{proof}

\begin{definition} \label{HNP-ImpArray}\index{implosive array!HNP}
Suppose that $A(S) = [\nu_0(S), \ldots , \nu_r(S)]$ is a $p$-implosive
array. We say that $A(S)$ is an {\em HNP-implosive array} if either
\begin{enumerate}
\item $A(S)$ is of Type I and in $S^p$ the cancellation between
  $\nu_0$ and $\nu_r$ is HNP-biting, or
\item $A(S)$ is of Type II and in $S^p$, for some $0 < i < r$,
  $\nu_0$ and $\nu_i$ are involved in HNP-biting or
  for some $0 < j < r$, $\nu_j$ and $\nu_r$ are involved in
  HNP-biting.
\end{enumerate}
\end{definition}

In order to follow the
arguments from Part \ref{Part:BG1}, we need to sharpen
Lemma \ref{ExistsImpArray}:
HNP-cancellation can beget $p$-implosive arrays with $p$
arbitrarily large, and therefore we must
argue for the frequent occurrence of 
$p$-implosive arrays that are not HNP-implosive. A first step in this
direction is given by the following

\begin{lemma} \label{NoHNP}
Let $\Pin$ be a pincer, and let $\mu_1$ and $\mu_2$ be the
colours associated to the bounding-paths $p_1$ and $p_2$ 
of $\Pin$.  Then there is no HNP-biting between beads in 
$\mu_1$ and $\mu_2$ within $\Pin$.
\end{lemma}

\begin{proof} Follows from Lemmas \ref{OneHNP} and \ref{DealWithExp}.
\end{proof}

When we are unconcerned about $p$ in a $p$-implosive array, we refer
merely to an {\em implosive array}.  The first restriction to note
concerning implosive arrays is this:

\begin{lemma} \label{NotManyColours}
If $[\nu_0(S), \dots , \nu_r(S)]$ is implosive of Type I, then $r \le
B$.  If it is implosive of Type II, then $r \le 2B$.
\end{lemma}

\begin{proof}
In Type I arrays, the interval $\nu_1(S^p)\cdots \nu_{r-1}(S^p)
\subset \bot(S^p)$ is to die in $S^p$, so the bound is an immediate
consequence of the Bounded Cancellation Lemma.  For Type II arrays,
one applies the same argument to the intervals joining $\nu_0(S^p)$
and $\nu_r(S^p)$ to the residual block.
\end{proof}

\begin{proof}[Proof of Proposition \ref{pSmall}]
We give a suitable formulation of `short' so that in
any corridor $S$ within $\Pin$, $S$ contains a short
$p$-implosive array.  Proposition \ref{pSmall} then follows from an
obvious finiteness argument.

Let $A(S) = [\nu_0(S) , \ldots , \nu_r(S)]$ be the implosive array
guaranteed to exist by Lemma \ref{ExistsImpArray}, and suppose that $p
\ge 2T_0$ (if not then a colour becomes stably Nielsen or vanishes within $2T_0$
of $\time (S)$). 

We can decompose each of the colours $\nu_j(S)$ in 
analogy with Part \ref{Part:BG1}, using the decomposition in
Section \ref{Chromatic} above.

We fix a constant $\Lambda_1$ so that if $\| A(S) \| >
\Lambda_1$ then one of the following must occur in $S^{T_0}$:
\begin{enumerate}
\item there is a block of displayed Nielsen \atom s in some
  $\nu_j(S^{T_0})$ of length at least $J+4B$,
\item there is a displayed \gep \ in some $\nu_j(S^{T_0})$ of length
  at least $J + 4B + 2$,
\item there is a displayed \pep \ in some $\nu_j(S^{T_0})$ of length
  at least $J + 4B + L + 1$, or
\item there is an interval of unnamed colours in $A(S)$ (which form a
  stably Nielsen block) of length at least $J + 4B$ between
  $\nu_0(S^{T_0})$ and $\nu_r(S^{T_0})$.
\end{enumerate}

In the remainder of the proof, we shall use the term {\em{block}}
to refer generically to the identified interval in whichever
of the above cases we find ourselves. Increasing $\Lambda_1$
if necessary, we may assume 
that the past of the block in $S$ satisfies the
relevant condition from (1) -- (4)  with the bound increased
by $2BT_0$.

For such a block $I$ in $S^{T_0}$, consider the first edge on either side
of this block which is not contained in a Nielsen path.  These edges
may be on one end of a  \gep \ or a \pep \ (including the \gep
\ or \pep\ from condition (2) or (3)), or may be contained
in unbounded \atom s.  Call these edges $\e_1$ and $\e_2$.

The Buffer Lemma \ref{BufferLemma:BG3} implies that either (i) one of
$\e_1$ and $\e_2$ must be `stabbed in the back' -- we do not exclude
the possibility that this stabbing happens by HNP-biting, or
(ii) there is HNP-cancellation across the above block.

We first dispose of case (ii).  Suppose, for ease of notation,
that the edge $\e_1$ HNP-bites the edge $\e_2$ across the above
block $I$.  Let $\e_1$ have weight $k$.  Then all edges in $I$ and
$\e_2$ must have weight less than $k$.  Let $\e_2'$ be the first edge
to the right of $I$ that has weight at least $k$.  Then the Relative
Buffer Lemma \ref{RelBuffer} implies that either $\e_1$ or $\e_2'$
must be stabbed\footnote{Note that if there is no such edge $\e_2'$
in $\Pin$ then $\e_1$ must be stabbed in the back, by Lemmas
  \ref{RelBuffer} and \ref{NoHNP}.}
in the back (again, this could be by
HNP-biting). 

We have argued that some edge must be stabbed in the back.
Suppose that this stabbing is of an edge $\e$ in $S^{T_0}$ and that $\e$ has
weight $k_1$.  Consider first the possibility that $\e$ is stabbed in
the back via HNP-biting.  Then this occurs by an edge $\e'$ of
weight at least $k_1 + 1$.  Now, either this stabbing in the back
occurs within $T_0$ of $S^{T_0}$, or by the Weighted Two Colour Lemma
(\ref{WeightTwo}) there is another block as in (1) -- (4) above.  This
block has higher weight than the previous block, and as above leads to
another stabbing in the back.  If this stabbing is HNP-biting,
pass to a yet higher weight stabbing, and so on.

Eventually (after less than $\omega$ iterations of this argument), we 
get an edge $\e$ stabbed in the back with the stabbing not 
HNP-biting.  Suppose that $\e$ has weight $k_2$.  Suppose for
ease of notation that $\e$ is to the left of the long block, and
suppose that $\e$ is coloured $\nu_i$.  Because of the block of Nielsen
\atom s to the non-stabbing side of $\e$, the Two Colour Lemma
(Proposition \ref{TwoColour}) implies that if the edge $\e'$ which
stabs $\e$ in the back is coloured by $\nu_j$ then $i - j > 1$;
we then write $\nu_j \searrow \nu_i$.

Passing to an innermost pair $\nu_{l_1} \searrow \nu_{l_2}$ between
$\nu_i$ and $\nu_j$ we can see that there are no blocks in $S^{T_0}$
satisfying any of (1) -- (4) above, for otherwise there would be a
further stabbing, leading to a related pair of colours between our
innermost pair, contradicting the innermost nature of this pair.

Once there are no such blocks, we have a bound on the length of the
$p$-implosive array implicit in the
relation $\nu_{l_1} \searrow \nu_{l_2}$. An obvious
finiteness argument now finishes the proof.
\end{proof}

We have already seen how Proposition \ref{pSmall} implies 
Proposition \ref{prePincer}.  Just as in Section \ref{ConstantSection},
we must now deal with the 
possibility of `nested pincers'.  

\subsection{Super-buffers}
\index{super-buffers}

\begin{definition} We consider sequences of $5$-tuples of tight edge-paths in $G$.
\[ U_k := \Big( u_{k,1}, u_{k,2}, u_{k,3}, u_{k,4}, u_{k,5} \Big),
\ \ k= 1,2,...  \]
with $|u_{k,1}|$ and $|u_{k,2}|$ at most $C_0 + C_1 + 2B(B+1) + 1$, 
while $|u_{k,2}|$ and
$|u_{k,4}|$ are at most $C_0+C_1+J$ and $|u_{k,3}| \le 4B(B+1)+1$.\footnote{The purpose of these constants is just as in Definition \ref{T1'Lemma}, with appropriate changes due to Lemmas \ref{A1Short}
and \ref{VanishLemma} and Proposition \ref{PointofC1}.}   
We fix an integer $T_1'$ sufficiently large to ensure that 
for any sequence of length $T_1'$ there will be a repetition,
i.e.~some
$t_1 < t_2 \le T_1'$ with
\[  \Big( u_{t_1,1}, u_{t_1,2}, u_{t_1,3}, u_{t_1,4}, u_{t_1,5} \Big) =
\Big( u_{t_2,1}, u_{t_2,2}, u_{t_2,3}, u_{t_2,4}, u_{t_2,5} \Big) . \]
We also choose
$T_1' \ge T_1$.
\end{definition}

With appropriate changes of terminology and the results
of the previous subsection in hand, the
proof of 
Proposition \ref{NoDoubleNeuter} yields:

\begin{lemma} \label{SuperBuffer}
Let $V = V_1 V_2 V_3$ be a tight concatenation of three beaded paths in $G$.
If the future of $V_2$ is not stably Nielsen in $f_{\#}^{T_1'}(V)$ then
the future of $V_2$ is not stably Nielsen in $f_{\#}^k(V)$ for any
$k \ge 0$.
\end{lemma}

\subsection{Nesting and the Pincer Lemma}

Let $\lambda_0 = J+2B(T_0+1)+1$, which is the obvious analogue
of the constant of the same name in Section \ref{ConstantSection}.  As in
Remark \ref{decreell}, it is convenient to assume that
$LC_4 < \lambda_0$, and we increase $\lambda_0$ 
to make this so. (This makes certain statements in Section \ref{TeamsSection}
easier, but has no serious affect.)

\begin{definition} \label{NestDef}
Consider one pincer $\Pi_1$ contained in another $\Pi_0$.
Suppose that in the corridor $S \subseteq \Pi_0$ at the top of
$\Pi_1$ (where its boundary paths $p_1(\Pi_1)$ and $p_2(\Pi_1)$
come together) the future in $\top(S)$ of at least one of the edges
containing $p_1(\Pi_1) \cap \top(S)$ or $p_2(\Pi_1)\cap \top(S)$
is not contained in any stably Nielsen path and this future
\footnote{We allow this future to be empty, in which case ``contained
in" means that the immediate past of the long stably Nielsen path
is not separated from $\Pi_1$ by any edge that has a future in
$\top(S)$.} lies in a beaded path consisting of Nielsen beads and
beads of weight strictly less than the weight of the edges containing
$p_1(\Pi_1) \cap \top(S)$ and $p_2(\Pi_1)\cap \top(S)$, and that this
beaded path has at least $\lambda_0$ non-vanishing beads.
Then we say that $\Pi_1$ is \index{pincer!nested} {\em nested} in
$\Pi_0$. 
\end{definition}

\begin{remark} \label{remark:nesteddef}
Besides the obvious translations, the above differs from
Definition \ref{NestingDef} in that the path at the top
of the pincer may now consist of Nielsen beads and lower 
weight beads, whereas in Part \ref{Part:BG1} it consisted entirely
of constant letters.  This more general setting does not make
any of the proofs in this section harder (because of the
Weighted Two Colour Lemma), but is needed because
of the more complicated definition of the `cascade of pincers'
below (Definition \ref{CascadeDef}).
\end{remark}

\begin{definition}\index{pincer!colours $\chi(\Pi)$}
For a pincer $\Pi_0$, let $\{ \Pi_i \}_{i \in I}$ be the set of all
pincers nested in $\Pi_0$.  Then define
\[      \chi(\Pi_0) = \tilde{\chi}(\Pi_0) \smallsetminus \bigcup_{i \in I}
\tilde{\chi}(\Pi_i)     .       \]
\end{definition}

The corridor $S_t$ was defined in Definition \ref{PincerDef:BG3}.

\begin{lemma} [cf. Lemma \ref{nestLife}]
If the pincer $\Pi_1$ is nested in $\Pi_0$ then
$\time(S_t(\Pi_1)) < \time(S_{\Pi_0})$.
\end{lemma}
\begin{proof}
The existence of the beaded path at the top of the pincer
$\Pi_1$ makes this an immediate consequence of the 
Weighted Buffer Lemma \ref{WtBuffer}.
\end{proof}

Define $T_1 = T_1' +2T_0$.  The following theorem is the main result
of this section, and is the strict analogue of Theorem \ref{PincerLemma}.
The proof in the current context follows the proof from Part \ref{Part:BG1} {\em mutatis mutandis}.

\begin{theorem} [Pincer Lemma] \label{PincerLemma:BG3}
For any pincer $\Pi$ \index{Pincer Lemma} \index{pincer!life of}
\[ \life(\Pi) \le T_1 (1 + |\chi(\Pi)|) .       \]
\end{theorem}

\section{Teams} \label{TeamsSection}
\index{teams}

By virtue of Lemma \ref{PointofC1}, Remark \ref{BeadtoLength} 
and the results of Section \ref{FastBeads},
we have reduced the task of bounding the bead norm 
of $S_0$ to that of 
bounding the lengths of certain blocks $C_{(\mu,\mu')}(2)$ which
consist of Nielsen beads coloured $\mu$ all of which are to
be eventually bitten by beads coloured $\mu'$ in the future of $S_0$.  
By Proposition \ref{Reaper}, if
such a block has length at least $B+J$, then there is an 
associated reaper, which consumes Nielsen beads in
$C_{(\mu,\mu')}(2)$ at a constant rate (and entirely consumes
any bead it bites, up to the final bead).  Note that to each
pair $(\mu,\mu')$ there is at most one associated reaper.

This puts us in the situation where we can develop the technology
of {\em teams} as in Section \ref{teamSec}. However, there are
a number of key differences to Part \ref{Part:BG1}: we
already had to work hard in 
Section \ref{TrappedHNP} to establish  the existence of 
a {\em reaper} for
$C_{(\mu,\mu')}(2)$, and now we have to work harder 
to identify the times
$\hat{t}_1(\mu,\mu')$ and $t_1(\mathcal{T})$ attached to a pair
$(\mu,\mu') \in \mathcal Z$ and a team $\mathcal T$, using
the {\em robust} past of the reaper instead of the actual past;
this is required  in order that the 
 Pincer Lemma apply to teams of genesis
(G3).  It is worth remarking that once we have identified the pincer
$\Pi_{\mathcal{T}}$ associated to a team $\mathcal T$ of genesis (G3),
we revert to an analysis of actual pasts (as in the definition of pincer).

Note that the colour of the edges in the robust future of an edge may 
not always be the same, contrary to the actual future. In fact, whenever
the robust past is not the actual past, the colour changes.
 This explains a 
slight difference between Definition \ref{preTeam} below
and Definition \ref{preteam}.

Consider an interval $C_{(\mu,\mu')}(2)$ so that
$|C_{(\mu,\mu')}(2)| > B+J$, and let $\epsilon^\mu$
be the reaper associated to $C_{(\mu,\mu')}(2)$ in 
Proposition \ref{Reaper} above.  Let $t_0$ be the 
time at which $\epsilon^\mu$ first bites a Nielsen bead
in $C_{(\mu,\mu')}(2)$, and let $\beta_\mu$ be the
rightmost bead in the future of $C_{(\mu,\mu')}(2)$ at 
this time.  Note that $\beta_\mu$ is a Nielsen bead.
Let $\epsilon_\mu$ be the rightmost edge in $\beta_\mu$.

\begin{remark}
Since $|C_{(\mu,\mu')}(2)| > B+J$, and each bead of 
$C_{(\mu,\mu')}(2)$ is to be bitten by $\mu'$, the colour
 of $\epsilon^\mu$ is $\mu'$.
\end{remark}

\begin{lemma}
Suppose that the immediate past of $\epsilon_\mu$
exists (i.e.~that $\epsilon_\mu$ does not lie on
$\partial \Delta$).  Then the immediate past of $\epsilon_\mu$ 
lies in some bead $\sigma$, and $\sigma$ contains
the immediate past of each edge in $\beta_\mu$.
\end{lemma}

The above lemma, applied at each stage in the past,
implies that we can follow the past of the edge $\epsilon_\mu$
and deduce consequences about the past of all edges
in $\beta_\mu$.

\index{teams!times $t_i(\T)$}
We now define a time $\hat{t}_1(\mu,\mu')$ as follows:
We go back to the last point in time when
(i) the past
of $\epsilon_\mu$ and the robust past of $\epsilon^\mu$
lay in a common corridor; and
(ii) $\epsilon_\mu$ is contained in a beaded Nielsen path
whose swollen present is immediately adjacent to the
robust past of $\epsilon^\mu$.

{\em{We denote this corridor $S_{\uparrow}$.}}

\begin{definition} \label{preTeam}\index{reaper}\index{pre-team}
The robust past of $\epsilon^\mu$ at time $\hat{t}_1(\mu,\mu')$
is called the {\em reaper}, and is denoted $\hat\rho(\mu,\mu')$.
The interval $\hat{\mathfrak T}(\mu,\mu')$ is the maximal
beaded Nielsen path in $\bot(S_{\uparrow})$ all of whose
beads are eventually bitten by $\hat{\rho}(\mu,\mu')$.
The {\em pre-team} $\hat{\mathcal T}(\mu,\mu')$ is defined to be the 
set of pairs $(\mu_1,\mu_2) \in \mathcal Z$ so that (i)
the robust past of $\epsilon^\mu$ is coloured $\mu_2$ at some time
between $\hat{t}_1(\mu,\mu')$ and $t_0$; and (ii) 
$\hat{\mathfrak T}(\mu,\mu')$ contains some edges coloured
$\mu_1$.
The number of beads in $\hat{\mathfrak T}(\mu,\mu')$ is denoted
$\| \hat{\mathcal T} \|$.
\end{definition}

As in Section \ref{teamSec}, we will define {\em teams} to be pre-teams
satisfying a certain maximality condition (see Definition \ref{TeamDef}
below).

\begin{remark}
Just as in Remark \ref{rem:virt}, if $\hat{t}_1(\mu,\mu') < \time(S_0)$ then near
the right-hand end of $\hat{\mathfrak{T}}(\mu,\mu')$ one may have an
interval of colours $\nu$ for which $\nu(S_0)$ is empty.
\end{remark}

\begin{lemma} [cf. Lemma \ref{t1high}] \label{TeamBound}
If $\hat{t}_1(\mu,\mu') \ge \time(S_0)$ then
\[      \sum_{(\mu_1,\mu_2) \in \hat{\mathcal{T}}(\mu,\mu')} |C_{(\mu,\mu')}(2)|
\le \| \hat{\mathcal{T}}(\mu,\mu') \| + B(B+1). \]
\end{lemma}
\begin{proof}
The extra $B(B+1)$ is to account for the beads consumed
before the reaper comes into play. Otherwise
the proof is just as in Part \ref{Part:BG1}.
\end{proof}

\subsection{The Genesis of pre-teams} [cf. Subsection \ref{genesis}]
\index{teams!genesis of}

We consider the various events that may occur at $\hat{t}_1(\mu,\mu')$ which prevent
us pushing the pre-team back one step in time.  Recall that $S_{\uparrow}$ is the corridor
at time $\hat{t}_1(\mu,\mu')$ which contains $\hat{\mathcal T}(\mu,\mu')$.  Suppose that
$\mu_2$ is the colour of $\hat{\rho}(\mu,\mu')$.

There are four types of events:
\begin{enumerate}
\item[(G1)]  The immediate past of $C_{(\mu,\mu_2)}(S_\uparrow)$ is separated from the 
robust past
of $\hat{\rho}(\mu,\mu')$ by an intrusion of $\partial \Delta$.
\item[(G2)]  We are not in Case (G1), but the immediate past of $C_{(\mu,\mu_2)}(S_\uparrow)$
is separated from the robust past of $\hat{\rho}(\mu,\mu')$ because of a singularity.
\item[(G3)]  The immediate past of $C_{(\mu,\mu_2)}(S_\uparrow)$ is still in the same corridor
as the robust past of $\hat{\rho}(\mu,\mu')$, but the swollen present of the immediate
past of $C_{(\mu,\mu_2)}(S_\uparrow)$ is not immediately adjacent to the robust past
of $\hat{\rho}(\mu,\mu')$.
\item[(G4)]  We are not in any of the above cases, but the immediate past of the rightmost
edge in $C_{(\mu,\mu_2)}(S_\uparrow)$ is not contained in a beaded Nielsen path.
\end{enumerate}

We now make the definition of a team.

\begin{definition}[cf. Definition \ref{shortDef}]\label{TeamDef}\index{teams}
All pre-teams $\hat{\mathcal T}(\mu,\mu')$ with $\hat{t}_1(\mu,\mu') \ge \time(S_0)$ are
defined to be teams, but the qualification criteria for pre-teams with $\hat{t}_1(\mu,\mu') < \time(S_0)$ are more selective.

If the genesis of $\hat{\mathcal T}(\mu,\mu')$ is of type (G1) or (G2), then the rightmost
component of the pre-team may form a pre-team at times before $\hat{t}_1(\mu,\mu')$.  In
particular, it may happen that $(\mu_1,\mu_2) \in \hat{\mathcal T}(\mu,\mu')$ but
$\hat{t}_1(\mu,\mu') > \hat{t}_1(\mu_1,\mu_2)$ and hence $(\mu,\mu') \not\in
\hat{\mathcal T}(\mu_1,\mu_2)$.  To avoid double-counting in our estimates on $\| \mathcal T \|$
we disqualify the (intuitively smaller) pre-team $\hat{\mathcal T}(\mu_1,\mu_2)$ in these
settings.

If the genesis of $\hat{\mathcal T}(\mu,\mu')$  is of type (G4), then again it may happen that
what remains to the right of $\hat{\mathcal T}(\mu,\mu')$ at some time before $\hat{t}_1(\mu,\mu')$
is a pre-team.  In this case, we disqualify the (intuitively larger) pre-team $\hat{\mathcal T}(\mu,\mu')$.

The pre-teams that remain after these disqualifications are now defined to be {\em teams}.
\end{definition}

A typical team will be denoted $\mathcal T$ and all hats will be dropped from the notation
for their associated objects (just as in Section \ref{teamSec}).

A team is said to be {\em short} if $\| \mathcal T \| \le \lambda_0$ or
$\sum_{(\mu_1,\mu_2) \in \mathcal T} |C_{(\mu_1,\mu_2)}(2)| \le \lambda_0$.
Let $\Sigma$ denote the set of short teams.

\begin{lemma} [cf. Lemma \ref{G4lemma}]
Teams of genesis (G4) are short.
\end{lemma}

We wish our ultimate definition of a team to be such that every pair $(\mu,\mu')$
with $C_{(\mu,\mu')}(2)$ non-empty is assigned to a team.  The above definition
fails to achieve this because of two phenomena:  first, a pre-team $\mathcal T(\mu,\mu')$
with genesis of type (G4) may have been disqualified, leaving $(\mu,\mu')$
teamless; second, in our initial discussion of pre-teams we excluded pairs $(\mu,\mu')$ with $| C_{(\mu,\mu')}(2)| \le B + J$.  The following
definitions remove these difficulties.

\begin{definition} [Virtual team members] \label{Virtual:BG3}
If a pre-team $\pTmm$ of type (G4) is disqualified under the terms of Definition \ref{TeamDef} 
and the smaller team necessitating disqualification is $\pT(\mu_1,\mu_2)$,  
then we define $(\mu,\mu')\vin\pT(\mu_1,\mu_2)$ and $\pTmm\subset_v\pT(\mu_1,\mu_2)$. 
We extend the relation $\subset_v$ to be transitive and extend $\vin$ correspondingly. 
If $(\mu,\mu')\vin\T$ then $(\mu_2,\mu')$ is said to be a
\index{teams!virtual members}{\em virtual member} of  
the team $\T$. 
\end{definition}

\begin{definition} 
If  $(\mu,\mu')$ is such that  $1\le |C_{(\mu,\mu')}(2)|\leq B+J$ and 
$(\mu,\mu')$ is neither a member nor a virtual member of any previously 
defined team, then we define $\T_{(\mu,\mu')}:=\{(\mu,\mu')\}$ to be a
(short) team with $\|\T_{(\mu,\mu')}\|=|C_{(\mu,\mu')}(2)|$.  
\end{definition}
 
\begin{lemma} [cf. Lemma \ref{allIn}]\label{allIn:BG3} 
Every   $(\mu,\mu')\in\vecZ$ with $C_{(\mu,\mu')}(2)$ non-empty is a member 
or a virtual member of exactly one team, and there are less than $2\n$ teams. 
\end{lemma}  
\begin{proof} 
The first assertion is an immediate consequence of the preceding 
three definitions, and the second  follows 
from the fact that $|\vecZ| < 2\n$. 
\end{proof} 
\medskip 

\subsection{Pincers associated to teams of genesis (G3)} [cf. Subsection \ref{ss:PinG3}]

In this subsection we describe a pincer $\Pi_\T$ canonically associated to each team
of genesis (G3), as in Subsection \ref{ss:PinG3}.  The only real difference between
the definitions here and those in Part \ref{Part:BG1} is the use of robust past and 
beaded Nielsen paths.  Sadly, this variation leads to complications
in the cascade of pincers; see
Definition \ref{CascadeDef} and Remark \ref{remark:nesteddef}.

\begin{definition} [cf.~Definition \ref{narrowPast}]\index{past!narrow}
The {\em narrow past} of a team $\T$ at time $t$
consists of those beaded Nielsen paths whose beads are displayed in their
colour and whose future is contained in $\mathfrak T$.
The narrow past may have several components at each time, the set of which
are ordered left to right according to the ordering in $\mathfrak T$ of their
futures.  We call these components {\em sections}.
\end{definition}

{\em For the remainder of this subsection we consider only long teams of
genesis (G3).}

The following lemma follows from the definition of teams
of genesis (G3) in a straightforward manner.

\begin{lemma} \label{BeadsForPincer}
Let $\T$ be a team of genesis (G3).
There exist beads $y(\T)$ and $y_1(\T)$ of different colours,
both lying strictly between
the immediate past of the swollen present of $\T$
and the robust past of $\hat\rho(\mu,\mu')$, so that
$y(\T)$ is bitten by $y_1(\T)$ and this is not HNP-biting.
\end{lemma}

\begin{definition} [The Pincer $\tilde{\Pi}_\T$] \label{prePincerDef}
Choose a leftmost pair of beads $y(\T), y_1(\T)$ satisfying Lemma
\ref{BeadsForPincer}, and let $x(\T)$ be the leftmost \index{pincer!of a team}
edge in $y(\T)$.  Let $x_1(\T)$ be the edge in $y_1(\T)$ which
is the past of the edge which cancels with the leftmost edge 
in the immediate future of $x(\T)$.

Define $\tilde{p}_l(\T)$ to be the path in the family forest $\mathcal F$
that traces the history of $x(\T)$ to $\partial\Delta$, and let $\tilde{p}_r(\T)$
be the path that traces the history of $x_1(\T)$.

Define $\tilde{t}_2(\T)$ to be the earliest time at which the paths
$\tilde{p}_l(\T)$ and $\tilde{p}_r(\T)$ lie in the same corridor.
\end{definition}

\begin{remark}
Since the pair $y(\T), y_1(\T)$ in Definition \ref{prePincerDef}
are the leftmost pair satisfying Lemma 
\ref{BeadsForPincer}, any non-vanishing beads which lie
between $\mathfrak T$ and this pair are involved
in HNP-biting and are of lower weight than $y_1(\T)$, by the
Weighted Buffer Lemma \ref{WtBuffer}.
\end{remark}

\begin{lemma} \label{Lemma:GetPincer}
The segments of the paths $\tilde{p}_l(\T)$ and $\tilde{p}_r(\T)$, 
together with the  path joining them along the bottom of the corridor 
at time  $\tilde{t}_2(\T)$ form a pincer.
\end{lemma}
\begin{proof}
Note that when choosing the beads $y(\T)$
and $y_1(\T)$ we excluded HNP-cancellation.  That the
paths in the statement of the lemma form a pincer then follows
immediately from the definition of pincers.
\end{proof}
We denote the pincer described in Lemma 
\ref{Lemma:GetPincer} above by $\tilde{\Pi}_\T$.

\subsection{The cascade of pincers} \index{pincer!cascade}
The Pincer Lemma argues for the regular disappearance of colours
within a pincer during those times when more than two colours continue
to survive along its corridors.  However, when there are only two colours,
the situation is more complicated.

Recall that the constant $T_0$ is as in Proposition \ref{TwoColour}, subject
to the requirement that $T_0 \ge T_0'$ as in the assumption immediately
after Proposition \ref{WeightTwo}.  The pincer $S_{\Pin}$ associated to a pincer
$\Pin$ is defined in Definition \ref{PincerDef:BG3}.

\begin{lemma} \label{CascadeLemma}
One of the following must occur: 
\begin{enumerate}
\item \label{GetPincer} $\time(S_{\tilde{\Pi}_\T}) > t_1(\T) - T_0$;
\item the path $\tilde{p}_l(\T)$ and the entire narrow past
of $\T$ are not in the same corridor at time $t_1(\T) - T_0$; or
\item \label{GetCascade} at time $t_1(\T) - T_0$ the path 
$\tilde{p}_l(\T)$ and the
narrow past of $\T$ are separated by a path which does not split
as a beaded path whose beads are either
 Nielsen paths or of weight less than $\tilde{p}_l(\T)$.
\end{enumerate}
\end{lemma}
\begin{proof}
If not, the Weighted Two Colour Lemma (Lemma \ref{WeightTwo}) would give a contradiction, since there {\em is} to be interaction between
the beads $y(\T)$ and $y_1(\T)$ at time $t_1(\T)$, and this interaction
is not HNP-biting.
\end{proof}

We now consider each of the three cases in turn, seeking
a definition of times $t_2(\T)$ and $t_3(\T)$ and (possibly)
a pincer $\Pi_\T$.  The following definition is entirely
analogous to Definition \ref{PincerDef}, with the appropriate
translations.

\begin{definition} [cf. Definition \ref{PincerDef}] \label{CascadeDef} 
\ 
\begin{enumerate}\index{teams!times $t_i(\T)$}
\item \label{Case1} Suppose some section of the narrow past of $\T$ 
is not in the same corridor as $\tilde{p}_l(\T)$ at time $t_1(\T)
- T_0$:  In this case\footnote{this includes the possibility
that $\tilde{p}_l(\T)$ does not exist at time $t_1(\T) - T_0$}
we define $t_2(\T) = t_3(\T)$ to be the earliest time
at which the entire narrow past of $\T$ lies in the same
corridor as $\tilde{p}_l(\T)$ and has length at least $\lambda_0$.
\item \label{Case2} Suppose that Case (\ref{Case1}) does not
occur and $\time(S_{\tilde{\Pi}_\T}) > t_1(\T) - T_0$.  We define
$\Pi_\T = \tilde{\Pi}_\T$ and $t_3(\T) = \time(S_{\Pi_\T})$.  If the
narrow past of $\T$ at time $t_1(\T) - T_0$ has length less than
$\lambda_0$, we define $t_2(\T) = t_3(\T)$, and otherwise
$t_2(\T) = \tilde{t}_2(\T)$.
\item \label{Case3} Suppose that neither Case (\ref{Case1}) or
Case (\ref{Case2}) occurs:  In this case,
Lemma \ref{CascadeLemma}(\ref{GetCascade}) pertains.  We pass to the
latest time at which there is a path between $\tilde{p}_l(\T)$ and
the narrow past of $\T$ which has an edge of at least the same
weight as $\tilde{p}_l(\T)$ at this time and is not contained
in a Nielsen path.  Choose a pair of beads $y'(\T)$, $y_1'(\T)$
as in Lemma \ref{BeadsForPincer}, as
well as edges $x'(\T)$, $x_1'(\T)$.  Let $\tilde{p}'_l(\T)$ 
be the path 
tracing the history of $x'(\T)$.  Let
$\tilde{p}'_r(\T)$ trace the history of the edge $x_1'(\T)$ that
cancels $x'(\T)$.  Let $\tilde{t}_2'(\T)$ be the earliest time at which the
paths $\tilde{p}'_l(\T)$ and $\tilde{p}'_r(\T)$ lie in the same
corridor and consider the pincer formed by these paths after time 
$\tilde{t}_2'(\T)$ and the path joining them along the bottom
of the corridor at time $\tilde{t}_2'(\T)$.

We now repeat our previous analysis with the primed objects 
$\tilde{p}_l'(\T), \tilde{t}_2'(\T)$, {\em etc.} in place of
$\tilde{p}_l(\T), \tilde{t}_2(\T)$, {\em etc.}, checking whether
we now fall into Case (\ref{Case1}) or (\ref{Case2}); if we do not
then we pass to $\tilde{p}''_l(\T)$, {\em etc.}. We iterate this
analysis until we  fall into Case (\ref{Case1}) or 
(\ref{Case2}), at which point we acquire the desired definitions
of $\Pi_\T, t_2(\T)$ and $t_3(\T)$.
\end{enumerate}
Define $p_l(\T)$ (resp. $p_r(\T)$) to be the left (resp. right)
boundary path of the pincer $\Pi_\T$ extended backwards in
time through $\mathcal F$ to $\partial \Delta$.  Define 
$p_l^+(\T)$ to be the sequence of edges (one at each time) lying on
the leftmost of the primed $\tilde{p}_l(\T)$ from the top of
$\pi_{\T}$ to time $t_1(\T)$.
\end{definition}

\begin{definition} [cf. Definition \ref{def:chiP}]
Let $\T$ be a long team of genesis (G3).  
We define $\chi_P(\T)$ to
be the set of colours containing the paths $\tilde{p}_l(\T),
\tilde{p}'_l(\T), \tilde{p}''_l(\T), \ldots$ that arise in 
Case (\ref{Case3}) of Definition 
\ref{CascadeDef} but do not become $p_l(\T)$.
\end{definition}

\begin{lemma} [cf. Lemma \ref{disj1}]

\ 

\begin{enumerate}
\item If $\T$ is a long team of genesis (G3),
\[ t_1(\T) - t_3(\T) \le T_0(|\chi_P(\T)| + 1 ).  \]
\item If $\T_1$ and $\T_2$ are distinct teams then
$\chi_P(\T_1) \cap \chi_P(\T_2) = \emptyset$.
\end{enumerate}
\end{lemma}

\subsection{The length of teams}

This subsection follows Subsection \ref{TeamLemmas}.  We
consider the lengths of arbitrary teams.

\begin{definition} [cf. Definition \ref{down1}]
\label{Down1Def} \index{teams!$\down_i(\T)$}
Let $\T$ be a team.  Define $\mbox{\rm{down}}_1(\T) \subset 
\partial \Delta$ to consist of those edges $e$ that are labelled
by some $t_i$ and satisfy one of the following conditions:
\begin{enumerate} 
\item[1.] $e$ is at the left end of a corridor containing a section of the narrow 
past of $\T$ that is not leftmost at that time; 
\item[2.]  $e$ is at the right end of a corridor containing a section of the narrow 
past of $\T$ that is not rightmost at that time;  
\item[3.]  $e$ is at the right end of a corridor which 
contains the rightmost section of the narrow past of $\T$ at that time but which does 
not intersect $p_l(\T)$.\\ 
\end{enumerate} 
\end{definition}

\begin{definition} [cf. Definition \ref{def:partialT}]\index{teams!$\partial^\T$}
Define $\partial^\T \subset \partial \Delta$ to be the intersection
of the narrow past of $\T$ with $\partial \Delta$.
\end{definition}

\begin{lemma} [cf. Lemma \ref{disj2}]

\ 

\begin{enumerate}
\item For distinct teams $\T_1$ and $\T_2$, the sets
$\partial^{\T_1}$ and $\partial^{\T_2}$ are disjoint.
\item For distinct teams $\T_1$ and $\T_2$, the 
sets $\text{\rm{down}}_1(\T_1)$ and $\text{\rm{down}}_1(\T_2)$
are disjoint.
\end{enumerate}
\end{lemma}

\begin{definition} [cf. Definition \ref{QT}]
Suppose that $\T$ is a team of genesis (G3).  We define \index{teams!$\QT$} $\QT$ be the
set of edges $\e$ with the following properties:
$p_l(\T)$ passes through $\e$  before time $t_3(\T)$,  the corridor $S$
with $\e\in\bot(S)$ contains the entire narrow past of $\T$, and
this narrow past has length at least $\lambda_0$.
\end{definition}

The following lemma reduces the task of bounding the total length
of teams to that of bounding the size of the sets $\QT$.  Its proof
follows that of Lemma \ref{TeamAgeLemma}.

\begin{lemma} [cf. Lemma \ref{TeamAgeLemma}]

\

\begin{enumerate}
\item If the genesis of $\T$ is of type (G1) or (G2), then
\[ \| \T \| \le 2LC_4 |\text{\rm{down}}_1(\T)| + |\partial^\T|.      \]
\item If the genesis of $\T$ is of type (G3), then
\[ \| \T \| \le 2C_4 |\text{\rm{down}}_1(\T)| + |\partial^\T|
+ 2LC_4 |\QT| + 2LC_4T_0 (|\chi_P(\T)| + 1) + \lambda_0.        \]
\end{enumerate}
\end{lemma}

\subsection{Bounding the size of $\QT$}

Let $\mathcal{G}_3$ be the set of long teams of genesis
(G3) for which $\QT$ is nonempty. Our goal for the remainder
of this section is to find a bound for $\sum_{\T \in \mathcal{G}_3}
|\QT|$.

\begin{lemma} [cf. Lemma \ref{t1-t2Lemma}]
For all $\T \in \mathcal{G}_3$
\[ t_3(\T) - t_2(\T) = \life(\Pi_\T) \le T_1(|\chi(\Pi_\T)| + 1).       \]
\end{lemma}

\begin{lemma} [cf. Lemma \ref{nesters}]
If $\T_1, \T_2 \in \mathcal{G}_3$ are distinct teams then
$\chi(\Pi_{\T_1}) \cap \chi(\Pi_{\T_2}) = \emptyset$.
\end{lemma}
\begin{proof}
The pincers $\Pi_{\T_i}$ are disjoint or else one is 
contained in the other.  In the latter case, say $\Pi_{\T_1}
\subset \Pi_{\T_2}$, the definition of nesting (Definition
\ref{NestDef}), and of the pincer associated to a team
(Definition \ref{CascadeDef}) ensure that $\Pi_{\T_1}$
is actually nested in $\Pi_{\T_2}$ (cf. Remark 
\ref{remark:nesteddef}).
\end{proof}

\begin{corollary} [cf. Corollary \ref{t1-t2Corr}]
\[\sum_{\T \in \mathcal{G}_3} t_3(\T) - t_2(\T)
\le 3T_1 |\partial \Delta|.     \]
\end{corollary}

We have now reduced our task for this section to bounding
the number of edges in the $\QT$ which occur before
$t_2(\T)$; this is the cardinality of the following set.

\begin{definition} [cf. Definition \ref{down2}]
\label{Down2Def}\index{teams!$\down_i(\T)$}
For a team $\T \in \mathcal{G}_3$ we define $\text{\rm{down}}_2(\T)$
to be the set of edges in $\partial \Delta$ that lie at the right-hand
end of a corridor containing an edge in $\QT$ before time $t_2(\T)$.
\end{definition}

Just as in Part \ref{Part:BG1}, it is not necessarily the case that 
the sets $\text{\rm{down}}_2(\T)$ are disjoint for distinct
teams, and we must deal with the possibility of `double-counting'.

The left-to-right ordering defined on paths in $\mathcal F$
in Section \ref{teamSec} is defined in the current context exactly
as in Part \ref{Part:BG1}.

\smallskip

\noindent{\bf Notation:} Let $\Gthree'$ be the set of teams $\T
\in \Gthree$ with $\down_2(\T) \neq \emptyset$.

\smallskip

\begin{lemma} [cf. Lemma \ref{separate}]
Consider $\T \in \Gthree'$.  If a path $p$ in $\mathcal F$ is to
the left of $p_l(\T)$ and a path $q$ is to the right of $p_r(\T)$,
then there is no corridor connecting $p$ to $q$ at any time
$t < t_2(\T)$.
\end{lemma}

\begin{definition} [cf.~Definition \ref{depthDef}]\index{teams!depth of}
$\T_1 \in \Gthree'$ is said to be {\em below} $\T_2 \in \Gthree'$
if $p_l(\T_1)$ and $p_r(\T_1)$ both lie between $p_l(\T_2)$ and
$p_r(\T_2)$ in the left-to-right ordering.

$\T_1$ is {\em to the left} of $\T_2$ if both $p_l(\T_1)$ and 
$p_r(\T_2)$ lie to the right of $p_r(\T_1)$.

We say that $\T$ is at {\em depth} $0$ if there are no teams above
it.  Then, inductively, we say that a team $\T$ is at depth $d+1$ if $d$ 
is the maximum depth of those teams above $\T$.

A {\em final depth} team is one with no teams below it.
\end{definition}

Note that there is a complete left-to-right ordering of those teams
in $\Gthree'$ at any given depth.

\begin{lemma} [cf. Lemma \ref{trapS_0}]
If there is a team from $\Gthree'$ below a team $\T \in \Gthree'$,
then $t_1(\T) \ge \time(S_0) \ge t_2(\T)$.
\end{lemma}
\begin{proof}
The proof from Part \ref{Part:BG1} works almost verbatim.  In particular,
the same proof shows that $\time(S_0) \ge t_2(\T)$.

To see that $t_1(\T) \ge \time(S_0)$, suppose that $\T'$ is a team
below $\T$.  Associated to
the team $\T'$ we have the beaded Nielsen path $\mathfrak T'$,
which is to be consumed by some reaper.  
The definitions of nesting and of the pincer $\Pi_{\T'}$ 
ensure that this
consumption of $\mathfrak T'$ must occur before time $t_1(\T)$.
On the other hand, $\mathfrak T$ has 
a non-empty future or past in $S_0$.
\end{proof}

With the preceding results in hand, a direct translation
of the proof of Lemma \ref{chiT} finishes the
work of this section:

\begin{lemma} [cf. Lemma \ref{chiT}]\index{teams!colours $\chi_*(\T)$}
There exist sets of colours $\chi_c(\T)$ and $\chi_\delta(\T)$
associated to each team $\T \in \Gthree'$ such that the sets
associated to distinct teams are disjoint and the following 
inequalities hold.

For each fixed team $\T_0 \in \Gthree'$ (of depth $d$ say), the
teams of depth $d+1$ that lie below $\T_0$ may be described
as follows:
\begin{enumerate}
\item[$\bullet$] There is at most one distinguished team
$\T_1$, and
\[ \| \T_1 \| \le 2B \Big( T_1(1+ |\chi(\Pi_{\T_0})|) +
T_0(|\chi_P(\T_0)| + 1) \Big).  \]
\item[$\bullet$] There are some number of final-depth teams.
\item[$\bullet$] For each of the remaining teams $\T$ we
have
\[ |\down_2(\T_0) \cap \down_2(\T)| \le T_1 \Big(
1+|\chi_c(\T)| \Big) + T_0 \Big( |\chi_\delta(\T)|+2 \Big). \]
\end{enumerate}
\end{lemma}

\begin{corollary} [cf.~Corollary \ref{Firstdown2sum}]
Summing over the set of teams $\T \in \Gthree'$ that are
not distinguished, we get
\[ \sum_\T \Big| \down_2(\T) \Big| \le
2 \Big| \bigcup_\T \down_2(\T)\Big| + \sum_\T T_1\Big(
1 + |\chi_c(\T)| \Big) + \sum_\T T_0 \Big( |\chi_\delta (\T)|
+ 2 \Big).      \]
\end{corollary}

Summing over the same set of teams again, we finally obtain:

\begin{corollary}
\[      \sum_\T |\down_2(\T)| \le |\partial \Delta| (2 + 3T_1 + 5T_0).  \]
\end{corollary}

\section{The Bonus Scheme} \label{BonusSection}\index{bonus scheme}

This section closely follows Section \ref{BonusScheme}.  
We have at last reached a stage where the proofs from
Part \ref{Part:BG1} can be translated without significant
modification.

In the previous section we defined teams and obtained a 
global bound on $\sum \| \T \|$.  If $C_{(\mu,\mu')}(2)$ is
non-empty then $(\mu,\mu')$ is a member or virtual member
of a unique team.  If the team is such that $t_1(\T) \ge \time(S_0)$,
then no member of the team is virtual and we have the inequality
\[ \| \T \| \ge \sum_{(\mu_1,\mu_2) \in \T} |C_{(\mu_1,\mu_2)}| - B(B+1),
\]
established in Lemma \ref{TeamBound}.  
This inequality might fail in case $t_1(\T) < \time(S_0)$.  The {\em
bonus scheme} assigns additional edges to teams in order to 
compensate for this failure.

By definition, at time $t_1(\T)$ the reaper $\rho = \rho_\T$ lies
immediately to the right of $\mathfrak T$.  The beads of $\mathfrak T$
not consumed from the right by $\rho$ by $\time(S_0)$ have a 
preferred future in $S_0$.  This preferred future, if contained in a single
colour, lies in $C_{(\mu_1,\mu_2)}(2)$ for some member 
$(\mu_1,\mu_2) \in \T$.  It could also intersect more than one colour
\footnote{Since Nielsen beads have bounded length, and 
there is a bound on the number of adjacencies of colours, there are relatively few such beads.}.
However, not all beads in the
 $C_{(\mu_1,\mu_2)}(2)$ need arise in this way:  some may not
 have a Nielsen bead as an ancestor at time $t_1(\T)$.  And if 
 $(\mu_1,\mu_2)$ is only a virtual member of $\T$, then no bead of
 $C_{(\mu_1,\mu_2)}(2)$ lies in the future of $\mathfrak T$.  The
 {\em bonus} beads in $C_{(\mu_1,\mu_2)}(2)$ are a certain subset of 
 those that do not have a Nielsen bead as an ancestor at time $t_1(\T)$.
 They are defined as follows.

\begin{definition}
Let  $\T$ be a team with  $t_1(\T) < \time(S_0)$ and consider a time
$t$ with $t_1(\T) < t < \time(S_0)$.

The \index{future!swollen}
{\em swollen future} of $\T$ at time $t$ is defined as in
Definition \ref{Swollen} with respect to the interval $\mathfrak{T}$,
which lies at time $t_1(\T)$.

Let $\epsilon$ be a non-Nielsen bead that lies immediately to the
left of the swollen future of $\T$, but whose immediate
ancestor is not a right
linear edge in this position.  If the path from $\epsilon$ to the reaper
$\rho_\T$ of $\T$ is a \gep, then we say that $\epsilon$ is a
\index{rascal}{\em rascal}.  Otherwise, if $\epsilon$ provides
more Nielsen beads than the reaper consumes, then $\epsilon$ 
is a {\em terror}.\index{terror}

In both cases, the {\em bonus provided by $\epsilon$} is the 
set of beads in the swollen future of $\T$ in $S_0$ that have
$\epsilon$ as their most recent ancestor which is not a 
Nielsen bead,
and which are eventually consumed by $\rho_\T$.

The set $\bonus(\T)$ is the union of the bonuses provided
to $\T$ by all rascals and terrors.
\end{definition}

\begin{lemma} [cf. Lemma \ref{C1toTeamLength}]
For any team $\T$,
\[      \sum_{(\mu_1,\mu_2) \in \T\ or\ (\mu_1,\mu_2) \vin \T}
| C_{(\mu_1,\mu_2)}(2)| \le \| \T \| + |\bonus(\T)| +B +J.      \]
\end{lemma}
 Note that the \gep \ which contains
a rascal in the above definition is not displayed.
We now proceed to bound the total bonus provided 
to teams by all
rascals and terrors.  Terrors
are straightforward to deal with.

\begin{lemma} [cf. Lemma \ref{lem:TerrorCont}]
The sum of the lengths of the bonuses provided to all
teams by terrors is less than $2L\n$.
\end{lemma}
\begin{proof}
Let $\epsilon$ by a terror, associated to a team $\T$.
Since the region from $\epsilon$ to the reaper of $\T$ is not
a \gep, $\epsilon$ must be right-fast.  Therefore, it will be 
separated from the team to which it is associated after one unit
of time.  Hence the bonus that $\epsilon$ provides is at most $L$.

That there can be at most one terror per 
adjacency of colours follows in a straightforward manner from
Lemma \ref{EndStab2} and the definition of terror.

Thus the total contribution of all terrors is less than $2L \n$.
\end{proof}

  In
parallel with Definition \ref{RascalPin}, we make the following

\begin{definition}
Fix a team $\T$ with $t_1(\T) < \time(S_0)$ and consider the
interval of time $[\tau_0(\epsilon),\tau_1(\epsilon)]$, where 
$\tau_0(\epsilon)$ is the
time at which a rascal $\epsilon$ appears at the left end of the 
swollen future of $\T$, and $\tau_1(\epsilon)$ is the time at 
which the robust future of $\epsilon$ is no longer to the immediate
left of the future of the swollen future of $\T$.

In the case where the robust future $\hat{\epsilon}$ of $\epsilon$ at 
time $\tau_1(\epsilon)$ is cancelled from the left by an edge
$e'$, we define $\tau_2(\epsilon)$ to be the earliest time when the
pasts of $\hat{\epsilon}$ and $e'$ are in the same corridor.  The path
in $\mathcal F$ that traces the past of $\hat{\epsilon}$ is denoted
$p_\epsilon$ and the past following the ancestors of $e'$ from 
$\tau_2(\epsilon)$ to $\tau_1(\epsilon)$ is denoted $p'_\epsilon$.  The 
pincer\footnote{we include the degenerate case here where the
``pincer" has no colours other than those of $\epsilon$ and $e'$.} 
formed by $p_\epsilon$, $p'_\epsilon$ and the corridor joining
them at time $\tau_2(\epsilon)$ is denoted $\Pi_\epsilon$.
\end{definition}

The only essential difference between the above definition and 
Definition \ref{RascalPin} is the use of the robust future of 
$\epsilon$ rather than the pp-future. 

With this definition in hand, the remaining results from
Section \ref{BonusScheme} may be translated directly,
yielding in particular:

\begin{proposition} [cf. Lemma \ref{BonusBound}]
Summing over all teams that are not short, we have
\[ \sum_{\T}|\bonus(\T)| \le \Big( (B+3)(3T_1 + 2T_0)L +
6BT_1 + 4BT_0 + 2\lambda_0 + 2B + 5L + 1 \Big) \n .\]
\end{proposition}

\section{From Bead Norm to Length} \label{LongGepsandPepsSection}

The output of the results up to now is a
bound for the bead norm of our corridor $S_0$.  In order to
complete the proof of Theorem \ref{BoundS:BG3} in the case of the
specified IRTT $f$ (which implies our Main
Theorem) we need to turn this into a bound
on the length of $S_0$.  For this we need to bound the total
length of the \gep s and \pep s in $S_0$ which have length more
than $J$ (or indeed any other fixed length).  In this section we 
explain how the techniques of the bonus scheme can be used
to establish such a bound. 

If a bead $\rho$ in $\mu(S_0)$  has length greater 
than $J$, it is either a \gep\ or a \pep.  
If it is a \pep\ then we may trace its past: at each time,
this past is either of length at most $J$ or 
else is a \pep\ or a \gep.  Whilst
this past remains a \pep, the number of Nielsen paths will decrease with
each 
backwards step in time, so at some point in the past of $\rho$, 
it must become
a \gep.  

Suppose now that $\rho$ is a \gep.  The past of a \gep\ is either a 
\gep\ or else has length at most $J$.  Thus, the length of the
\gep\ decreases as we go into the past until eventually it is
of length at most $J$.

There is a strong analogy between teams of genesis (G4) and long 
\gep s and \pep s.  On
one end of a long bead is a linear edge which consumes the
Nielsen beads in the middle.  This linear edge can be considered
as a reaper.  On the other end of a \gep\ is a linear edge which can
be considered as a rascal.   The moment when the past of a \pep\ becomes a
\gep\ is analogous to $\tau_1(\epsilon)$ from the bonus scheme,
and so a \pep\ in $S_0$ can be thought of as a team with a rascal 
$\epsilon$ with $\tau_1(\epsilon) \le \time(S_0)$.  Similarly, a long 
\gep\ in $S_0$ can be thought of as a team with a rascal $\epsilon$
so that $\tau_1(\epsilon) > \time(S_0)$.  

We can define the bonus associated to such a rascal exactly
as we did in the previous section.  Since we
are in the setting of  genesis  type (G4), all of the Nielsen beads in a long
\gep\ or \pep\ are in the bonus.  Thus it is enough to bound the 
total of the bonuses associated to long \gep s and \pep s.

The only thing we need to be able to follow the bonus
scheme directly is a bound on the number of long \gep s
and \pep s in $S_0$.

\begin{lemma}
The number of beads of length greater than $J$ in
$S_0$ is less then $4 \n$.
\end{lemma}
\begin{proof}
Let $\rho$ be a bead in $S_0$ of length greater than $J$, and
assign a time $\tau_1(\rho)$ to $\rho$ as described above.  If
$\rho$ is a \gep\ then $\tau_1(\rho) > \time(S_0)$, whilst if 
$\rho$ is a \pep\ then $\tau_1(\rho) \le \time(S_0)$.

Let $\rho'$ be the past or future of $\rho$ at time $\tau_1(\rho) -1$.
Consider the `event' at time $\tau_1(\rho)$ which stops the robust
future of $\rho'$ being a \gep.  

This `event' is either an intrusion of the boundary, a singularity, or else
there is an associated pincer caused by a cancellation from another
colour.  There are less than $\n$ events of each of the first two
types.

The Buffer Lemma ensures that there is at most one event of the
third type for each adjacency of colours.  An application of 
Lemma \ref{NoOfAdj} completes the proof.
\end{proof}

A bound on the total length of long beads in $S_0$ now follows
exactly as in the bonus scheme from Section \ref{BonusSection}
(the detailed arguments being in Section \ref{BonusScheme}).

\subsection{The end of the main road}

In Section \ref{Linear} we discussed how our Main Theorem 
follows from Theorem \ref{noName} and Proposition \ref{longGEPS}.
The bound that we just established
on the total length of long beads in $S_0$ proves
Proposition \ref{longGEPS}. The output of our estimates in 
the previous sections bounded the bead norm of $S_0$ by a linear
function of $|\partial\Delta|$, and Theorem \ref{noName} follows
from this because 
\[	[S]_{\beta} \le B \| S \|_\beta	,	\]
(see Lemma \ref{BLengthNorm}).

Thus the proof of the Main
Theorem  is finally at an end, and the reader can
join us in wondering why a statement as simple and engaging
as this theorem should require such a complicated proof.

\section{Corridor Length Functions and Bracketing} \label{BracketingSection}

In this section we prove Theorem \ref{BoundS:BG3} in full generality and
deduce the Bracketing Theorem from it. Our proof of Theorem \ref{BoundS:BG3}
proceeds via a discussion of {\em{corridor length functions}} for
more general semidirect products and mapping tori. 
Such functions should be regarded as measuring the complexity of
van Kampen diagrams in the spirit of isoperimetric and isodiametric
functions. We prove the following results (see Subsection \ref{ss:defL}
for precise definitions of the terms involved).

\begin{proposition}\label{CLchange} Let $G_1$ and $G_2$ be 
compact combinatorial complexes with
fundamental group $\Pi$,
and for $i=1,2$ let $f_i:G_i^{(1)}\to G_i^{(1)}$ be
an edge-path map of 1-skeleta inducing $\phi\in\rm{Out}(\Pi)$. Then the
$t$-corridor length function for the mapping torus $M(f_1)$
is $\simeq$ equivalent to that of $M(f_2)$.
\end{proposition}

\begin{proposition}\label{CLfi} If $\Pi$ is finitely generated and 
$\G=\Pi\rtimes_\phi \Z$ is finitely presented, then for every positive
integer $p$, the corridor length function of $\Pi$ is $\simeq$ equivalent
to that of
$\G_p=\Pi\rtimes_{\phi^p} \Z$
\end{proposition}

In the previous section we completed the proof of 
Theorem \ref{BoundS:BG3} in the case of one particular IRTT representative $f$
of a certain power of an arbitrary free-group automorphism $\phi$.  
The above results complete the proof in the general case.
Before turning to the proof of these results, we explain
how the Bracketing Theorem stated in the introduction 
is obtained by applying Theorem \ref{BoundS:BG3} to the most
naive topological representation of a free group automorphism $\phi$.

\subsection{The Bracketing Theorem}

The terms in   the following theorem were defined in the
introduction.

\medskip \index{Bracketing Theorem}
\noindent{\bf Theorem \ref{Bracket}.}
{\em
There exists a constant $K = K(\phi,\mathcal B)$
such that any word $w \equiv e_1 \dots e_n$ that represents
the identity in $F\rtimes_\phi\Z$ 
admits a $t$-complete bracketing $\beta_1,\dots , \beta_m$
such that the content $c_i$ of each $\beta_i$
satisfies $d_F(1,c_i)\le Kn$.
}

\smallskip

\begin{proof}
We work with the mapping torus $M$
of the obvious realisation of $\phi$ on
the graph with one vertex whose edges are indexed by $\mathcal B$.
Given a word $w$, we consider
a minimal-area van Kampen diagram over $M$ with boundary
label $w$. 
 We insert a bracket 
$w_1(w_2)w_3$ if and only if there is a $t$-corridor whose
ends are labelled by the initial and terminal letters of $w_2$.
(One must allow $t$-corridors of zero length in this description;
one would exclude them by making the easy reduction to words that
have no proper sub-words that are null-homotopic.)

These brackets are pairwise compatible because distinct $t$-corridors
cannot cross. And
because every $t$-edge in the boundary of a van Kampen
diagram is the end of a (perhaps zero-length) corridor, the
bracketing is complete. The content of the bracket is the freely
reduced form of the label  along the top or bottom of the corridor
(according to the orientation of the sentinels). In the former case, the
length of the corridor bounds the length of this label, and in the
latter case one has to multiply the length by at most $L=\max\{|\phi(b)|:b\in
\mathcal B\}$.
\end{proof}

\subsection{Corridor length functions} \label{ss:defL}
If $\Pi$ is a group with finite generating
set $\A$ and $\phi\in\rm{Aut}(\Pi)$ is such that $\G=\Pi\rtimes_\phi\Z$
is finitely presented, then $\G$ has a finite presentation of the
form 
\[	\big\langle \A,\,t \mid \mathcal R,\, t^{-1}at=\hat\phi(a)\ (a\in\A)\big\rangle ,	\]
where $t$ is the generator of the visible $\Z$, the relations $\mathcal R$
involve only the letters $\A$, and $\hat\phi(a)\in F(\A)$ is equal to
$\phi(a)$ in $\Pi$. 

\index{corridor length function}
We are concerned with the geometry of $t$-corridors
in van Kampen diagrams over such presentations. Thus we associate to the
presentation the {\em{$t$-corridor length function}} $\Lambda :\mathbb N\to
\mathbb N$, which is defined as follows. For each $w\in F(\A\cup\{t\})$
with $w=1$ in $\G$, we choose a van Kampen diagram for $w$ in which the
length of the longest $t$-corridor is as small as possible, and we
define $\lambda_t(w)$ to be this length. We then define
$$
\Lambda(n) := \max \{\lambda_t(w) \mid w=_\G 1,\ |w|\le n\}.
$$

More generally, since we have a well-defined notion of van Kampen
diagram and $t$-corridor in the setting of mapping tori of
{\em{edge-path maps}}\footnote{an {\em{edge-path map}} is a cellular
map that sends edges to edge-paths} of
combinatorial complexes, we can define the {\em{$t$-corridor length function}}
for such a complex.

\subsection{Invariance under change of topological representative}
\index{topological representative}

The scheme of the following proof follows the standard method of
 showing that
features of the geometry of \index{van Kampen diagram}
van Kampen diagrams are preserved under
quasi-isometry. However, one has to be careful to deal only with fibre-preserving
maps in order to retain control over the $t$-corridor structure.

\smallskip

\begin{proof}[Proof of Proposition \ref{CLchange}.]
We have a cocompact action of $\G =\Pi\rtimes_\phi\Z$
on the universal cover $X_i=\tilde M(f_i)$ for $i=1,2$, where the
action of $\Pi$ leaves invariant the connected components $C_{i,m}$
of the preimage
of $G_i\subset M(f_i)$ and the generator $t$
of $\Z$ acts so that $t^r.C_{i,m} = C_{i,m+r}$.

The cocompactness of the actions means that there exist constants $\delta_1, \delta_2$
so that every vertex in $C_{i,m}$ is within a distance $\delta_i$ of 
any $\Pi$-orbit of vertices in $C_{i,m}$, where distance is measured
in the combinatorial metric on the 1-skeleton (unit edge lengths).

We define $\G$-equivariant quasi-isometries between the 1-skeleta
of the $X_i$  as follows. First we
pick  base vertices $x_i\in C_{i,0}$ and define
$g_1:\gamma.x_1\mapsto \gamma.x_2$ and $g_2:\gamma.x_2\mapsto \gamma.x_1$.
Then, for each  vertex  $v\in C_{i,m}\smallsetminus \Gamma.x_i$ we choose
a closest element  $v'\in \Gamma.x_i\cap C_{i,m}$
and define $g_i(v) := g(v')$. Next, we extend to the edges in 
$C_{i,m}$ by sending each to a shortest edge path connecting the
images of its vertices. Finally, we extend  $g_i$ to $t$-edges
in $X_i$ so that it sends each such homeomorphically onto
the $t$-edge joining the images of its endpoints.

With the maps $g_1, g_2$ in hand, we can now
push van Kampen diagrams back and forth between $X_1$ and
$X_2$ as in the standard proof of the qi-invariance of Dehn functions
(cf.~\cite{BH}, page 143).
Thus, 
given a loop $\ell$ in the 1-skeleton of $X_1$, labelled $u_1t^{\e_1}u_2
\dots u_lt^{\e_l}$
we consider the loop $ g_1\circ \ell$ in $X_2^{(1)}$
and fill it with a
van Kampen diagram $\Delta$ so as minimize the length
of the longest $t$-corridor. We will be done if we can bound
$\lambda_t(\ell)$ by a linear function of this length.

Viewing $\Delta$ as a
map from a cellulated 2-disc to $X_2$, we compose it with $g_2$
to obtain a map to $X_1$. This new map is obtained from $\Delta$
by simply changing the labels on the edges: the $t$-edges are unchanged
while the edges labelled by 1-cells in $G_2$ are now labelled by 
edge-paths in the 1-skeleton of $G_1$ whose length is
bounded by the constants of the quasi-isometry $g_2$; the boundary
label of the diagram will be $\ell'=v_1t^{\e_1}v_2
\dots v_lt^{\e_l}$, where the $v_j$ are edge-paths of uniformly bounded
length and each $v_j$ is contained in the same component $C_{1,m_j}$
as $u_j$. (This is the point at which we use the fact that we chose
our quasi-isometries to respect fibres.) The faces of this diagram
can be filled with van Kampen diagrams in $X_1$; in the case of 2-cells
with no $t$-labels, we use only lifts of 2-cells from $G_1$; in the
case of 2-cells labelled $t^{-1}\rho t \sigma$ we divide them into 
(short) $t$-corridors in the obvious manner. 
The result\footnote{A familiar problem in this type of argument
arises from degeneracies that threaten the planarity of the
diagram; such problems are removed by surgery \cite{LS}. In
the current setting these surgeries take place only in the regions
between the $t$-corridors and therefore do not affect our discussion.} is a van
Kampen diagram for $\ell'$ in $X_1$ whose $t$-corridors are in bijection
with those of $\Delta$ and whose length is bounded by $k$ times the
length of those in $\Delta$, where $k$ is a constant that depends only
on our quasi-isometries.

To complete the desired diagram filling our original loop $\ell$, 
we need an annular diagram between $\ell$ and $\ell'$ that does not
 disrupt the structure of $t$-corridors in $\Delta'$. To this end,
 we join the vertices of $u_j$ to those of $v_j$ by paths in $C_{i,m_j}$
 of minimal length and fill the resulting loop with a diagram mapping to
 $C_{i,m_j}$; this gives a diagram
 $\Delta''$ with holes corresponding to the
 occurrences of $t^{\pm 1}$ in $\ell$.
 Next, if the
 arc joining the termini of $u_j$ and $v_j$ is labelled $\rho_i$,
 then we  insert a $t$-corridor into the
 hole associated to $\dots u_jt u_{j+1}\dots$, where the bottom
 of the $t$-corridor is labelled $\rho_j$. (If $t$ is replaced by $t^{-1}$,
 the bottom of the corridor is the arc $\sigma_{j+1}$ joining the
   initial vertex of $u_{j+1}$
 to that of $v_{j+1}$.) To complete the construction of $\Delta$, one
 uses 2-cells in $C_{i,m_{j+1}}$ to fill the loop formed by the top of
 the $t$-corridor and $\sigma_{j+1}$.
 \end{proof}

\begin{corollary} If $\Pi$ is finitely generated and $\G=\Pi\rtimes_\phi\Z$
is finitely presented then,
up to $\simeq$ equivalence, the $t$-corridor length function of 
$\Pi\rtimes_\phi\Z$ depends only on the semidirect product (i.e.~
although it
depends on the form of the finite presentation, it does
 not depend on the
choice of $\A$ and $\hat\phi$).
\end{corollary}

\subsection{Passing to Powers} 

The purpose of this subsection is to prove Proposition \ref{CLfi}. 

Let $\At$ be as above. Identifying
 $\G_p=\Pi\rtimes_{\phi^p}\Z$ with the subgroup 
$\Pi\rtimes p\Z$ of $\G$, we take generators $\A\cup\{\tau\}$
where $\tau=t^p$ in $\G$.
To each word $w\in\FAt$ that
equals $1\in\G$ we associate a word $w_p$ in the free group
on $\A\cup\{\tau\}$ according to the following
scheme. First
 we draw a path on the integer lattice in $\mathbb R^2$ that begins at the
origin and proceeds up one space as we read $t$, down one as we read $t^{-1}$
and moves one space to the right as we read a letter from $\Apm$.
We shall modify $w$ by replacing certain open 
segments of this path that lie in the vertical
 intervals $[mp,(m+1)p]$; these segments
are of two types, called {\em{bumps}} and {\em{steps}}. 

If both endpoints of the subpath are at height $mp$
and none of its edge are at height $(m+1)p$, then the segment
is called an {\em{up-bump}}. If the initial
endpoint is at height $mp$, the terminus at
height $(m+1)p$ and all other vertices are at 
heights in $(mp, (m+1)p)$,
then the segment is called an {\em{up-step}}. A {\em{down-bump}}
and {\em{down-step}} are defined similarly.

When we have replaced all steps and bumps from the path defined by $w$,
  the horizontal segments of the
resulting path will all run
at heights divisible by $p$.

To this end, we write $w=u_1v_1u_2v_2\dots$ where
$u_1$ is the first non-trivial prefix of $w$
whose exponent sum in $t$ is $0\mod p$ and $v_1$ is the (possibly empty)
subword  before the next $t^{\pm 1}$, then $u_2$
 is the first non-trivial prefix of $w$
whose exponent sum in $t$ is $0\mod p$, and so on.
Each $u_i$ labels either a bump or a step.

If $u_i$ labels a bump
then  we replace it by the reduced word $
U_i\in \FA$ that is equal in $\G$ to $u_i$.
 If $u_i=t^\e u_i',\, \e=\pm 1,$ is a step, then
we replace it by the unique reduced word $t^{\e p}U_i$ with
$U_i\in\FA$ and $t^\e U_i=u_i$ in $\G$.

Let $\tilde w_p\in\FAt$ be the word obtained from $w$ by the above process and
let $w_p\in\FAtp$ be the word obtained from $\tilde w_p$ by (starting from
the left) replacing sub-words labelled $t^{\pm p}$ by $\tau^{\pm p}$
and then freely reducing.

As usual, in the following lemma $L=\max\{|\phi(a)| : a\in\mathcal A\}$.

\begin{lemma}\label{lengthwp} $w=\tilde w_p=w_p$ in $\G$
and $ |w_p| \le |\tilde w_p| \le L^{p-1} |w|$.
\end{lemma}

\begin{proof} The bound on $|\tilde w_p|$ comes from the following observation.
For  a bump labelled $u_i$,
one can pass from $u_i$ to $U_i$ by deleting all letters $t^{\pm 1}$
from $u_i$ and replacing each occurrence of $a\in \A$ in $u_i$, 
say $u_i=\alpha a\beta$, by the freely
reduced word in $\FA$ representing $\phi^r(a)$, where $-r$ is the exponent
sum of $t$ in $\alpha$. Similarly, if a step is labelled $u_i=t^\e u_i'$,
then $U_i$ is obtained by deleting all $t$ from $u_i'$ and replacing
 each occurrence of $a\in \A$ in $u_i$,
say $u_i'=\alpha a\beta$, by the freely
reduced word in $\FA$ representing $\phi^{\e(p-r)}(a)$, 
where $\e r$ is the exponent
sum of $t$ in $\alpha$.
\end{proof}

The replacement scheme described in the preceding proof corresponds
to the construction of a singular-disc diagram  $A(w)$ exhibiting
the equality $w = \tilde w_p$ in $\G$. 
 Specifically, for each bump or
step, one draws the vertical line joining each vertex to the height where
it will be pushed, one labels it by the appropriate power of $t$, and then
one
fills-in the resulting line of rectangles  with 2-cells whose boundary
labels have the form $t^{-1}at\phi^{-1}(a)$.
(Starting from this
specific planar embedding one will in general
have to flip some of the components of the interior in order to get
an embedded diagram $A(w)$ with boundary cycle $\tilde w_p w_p^{-1}$.)

\begin{lemma}\label{corrinAw} $A(w)$ is a union of $t$-corridors;
each has at most one of its ends on the boundary arc labelled $\tilde w_p$,
and the length of a $t$-corridor in $A(w)$ is at most
  $L^{p-1}\max |u_i|$, where the $u_i$ are the sub-words
of $w$ labelling  bumps and steps.
\end{lemma}

\begin{proof} The diagram $A(w)$ consists of a string of disc diagrams, one
for each bump or step. A  $t$-corridor in a disc corresponding to
a bump labelled $u_i$ has both of its ends on the arc labelled $u_i$,
while a $t$-corridor in a disc corresponding to a step
labelled $tu_i'$ may have one end on the corresponding arc labelled $t^p$
in $\tilde w_p$ and one on the arc labelled $u_i'$ or (if the change in height
along $u_i'$ is not monotone) both ends on the arc labelled $u_i'$. In all
cases, the label on the bottom side of the corridor is a concatenation of
 less than $|u_i|$ words of the form $\phi^r(a)$ with $a\in\A$ and $|r|\le
p-1$.
\end{proof}

\noindent{\bf{Proof of Proposition \ref{CLfi}.}}
As we discussed immediately before subsection \ref{Iterate},
 the set of diagrams for $\G_p$ is,
after $p$-refinement, a subset of the diagrams over $\G$,
and hence the corridor length function of the latter
$\preceq$-dominates that of the former. (There are some
constants to take account of here, such as a factor
of $p$ in length
coming from the $p$-refinement, and an $L^{p-1}$ needed
to estimate the area of a $t$-corridor in terms of the corresponding
$\tau$-corridor, but these are trivial matters.)
Thus the true content of the proposition is that  
the corridor length function of $\G$ is
$\preceq$-bounded above by that of the $\G_p$.

For each freely-reduced
 word $W\in \FAtp$ that is null-homotopic in $\G_p$
we fix a van Kampen diagram $\Delta(W)$ whose $\tau$-corridors have length
at most $\Lambda(|W|)$. Then, 
for each freely-reduced $w\in \FAt$ that is 
null-homotopic in $\G$ we define a van Kampen diagram $\Delta_p(w)$ 
as follows. First, we replace
 $\Delta(w_p)$ by its $p$-refinement (which
has boundary label $\tilde w_p$). We then attach
to this the singular-disc diagram $A(w)$ along the portion of
its boundary labelled $\tilde w_p$.

We claim that the length of each
$t$-corridor in $\Delta_p(w)$ is at most
$$
L^{p-1}\, (2+ \Lambda(L^{p-1}|w|)).
$$

It follows from Lemma \ref{corrinAw} that 
each of the 
$t$-corridors in $\Delta_p(w)$ is either contained in 
the annular diagram $A(w)$, or else is a layer in the
$p$-refinement of a $\tau$-corridor from $\Delta(w_p)$,
possibly augmented on each end by a $t$-corridor in $A(w)$.
(The fact that there are no $t$-corridors in $A(w)$ with
both ends on the boundary arc labelled $\tilde w_p$ is
crucial here.)

The length of a $t$-corridor in $A(w)$ is at most
$L^{p-1}|v|$. The length of a $\tau$-corridor from
$\Delta(w_p)$
is at most $\Lambda(|w_p|) \le \Lambda(L^{p-1}|w|)$, and the
length of each layer in its refinement is therefore at most
$L^{p-1}\, \Lambda(L^{p-1}|w|)$.
\hfill$\square$
  
\section{On a Result of Brinkmann} \label{BrinkSection}

The following theorem is the main result in
\cite{BrinkDyn}. It plays a vital role in the first proof
that the conjugacy problem is solvable for free-by-cyclic groups
\cite{BMMV}  (our Corollary \ref{Conjug}).

\begin{theorem} \label{Brink} \cite[Theorem 0.1]{BrinkDyn}
Let $\phi : F \to F$ be an automorphism of a finitely generated
free group.  Then there exists a constant $K \ge 1$ such that for
any pair of exponents $N, i$ satisfying $0 \le i \le N$, the following
two statements hold:
\begin{enumerate}
\item\label{Brink1} If $w$ is a cyclic word in $F$, then
\[ \| \phi^i(w) \| \le K \Big( \| w \| + \| \phi^N(w) \| \Big) ,        \]
where $\| w \|$ is the length of the cyclic reduction of $w$ with
respect to some word metric on $F$.
\item \label{Brink2} If $w$ is a word in $F$, then
\[ | \phi^i(w)| \le K \Big( |w| + |\phi^N(w)| \Big),    \]
where $|w|$ is the word length of $w$.
\end{enumerate}
\end{theorem}

The purpose of this section is to explain how to extract
Theorem \ref{Brink} from our proof of the Main Theorem.
We regard words and cyclic words in $F_n$
 as, respectively, based and unbased loops
in the graph $R$ with one vertex and $n$ edges; the assertions
of Theorem \ref{Brink} are then statements about how the 
lengths of the tightened
images of such loops grow when one applies the obvious
topological realisation $\overline\phi$ of $\phi$.
As in the previous subsection, these assertions will follow
if we can establish the corresponding bounds with $\overline\phi:R\to R$
replaced by a topological (IRTT) representative $f:G\to G$ of
a power of $\phi$ satisfying Assumption
\ref{FinalPower}.

\begin{remark} The proof given below shows that the constant $K$
of Theorem \ref{BoundS:BG3} suffices for Theorem \ref{Brink}.
Brinkmann \cite{BrinkDyn}
states that (his constant) $K$ can be computed
effectively, but we do not see how to prove this. Indeed,
given his approach (and ours), this assertion would
seem to require an effective construction of an improved
relative train track representative for $\phi$, and a
proof that such a construction exists does not
seem to be available at the moment.
\end{remark}

 The following lemma allows a proof of the assertions in
\eqref{Brink1} and \eqref{Brink2} to be undertaken simultaneously.

\begin{lemma}
If $\sigma$ is a nontrivial loop in $G$, then
for some $j\ge 1$, the loop $f_{\#}^j(\sigma)$ admits a 
splitting at a vertex.
\end{lemma}
\begin{proof}
According to \cite[Lemma 4.1.2, p.554]{BFH}, $\sigma$ admits a splitting
$\sigma = \sigma_1$, where $\sigma_1$ is a path, but we argue further
to arrange for this splitting to be at a vertex.

We divide the argument into a number of cases, 
depending on the largest $i$ so that the stratum $H_i$
contains an edge of
$\sigma_1$.  If this $H_i$
is a zero stratum,  $f_{\#}(\sigma_1) \subset G_{i-1}$
and an obvious induction applies.     If $H_i$
parabolic, then we apply
\cite[Lemma 4.1.4]{BFH} to the circuit $\sigma$ to obtain a splitting into paths,
at least one of which is a basic path, and so has a vertex at one end. 
If $H_i$ is an exponential stratum, then there is a positive integer
$K$ so that the number of $i$-illegal turns in $f^k_{\#}(\sigma_1)$ is the 
same for all $k \ge K$.  In this case, since all Nielsen paths of exponential
weight are edge-paths and all periodic paths are Nielsen, \cite[Lemma 4.2.6]{BFH} 
implies that $f^K_{\#}(\sigma_1)$ admits a splitting into sub-paths which are either 
$r$-legal or pre-Nielsen paths. 
If all sub-paths of $f^K_{\#}(\sigma_1)$ are pre-Nielsen paths, then 
$f^{K+1}_{\#}(\sigma_1)$ is a Nielsen path, and we ensured in 
Section \ref{TrainTracks} that all Nielsen paths are edge-paths.

Suppose, then, that $f^K_{\#}(\sigma_1)$ contains an $r$-legal path $\rho$
of weight $r$ in its splitting.  Then an iterate $f^i_{\#}(\rho)$ of $\rho$ contains a 
displayed edge $\e$ of weight $r$, and the path $f^{K+i}_{\#}(\sigma_1)$
splits immediately on either side of $\e$.  
Since $\sigma$ has weight $i$, the splitting of $f^{K+i}_{\#}(\sigma_1)$ 
induces a splitting of $f^{K+i}_{\#}(\sigma)$ at a vertex, as required.
\end{proof}

In order to prove the statements \eqref{Brink1} and \eqref{Brink2}, 
we analyze the van Kampen diagram $\Delta$ over the
mapping torus of $f:G\to G$ that has boundary label
$t^{-k}\sigma t^k f_{\#}^k(\sigma)^{-1}$. This is a
simple stack of corridors as consider in Subsection \ref{stackDiags}.

In the restricted setting of stack diagrams, 
many of the difficulties that had to be overcome in the proof of 
Main Theorem  do not arise
(there are no singularities, for example), but
there remain difficulties that one does not encounter
in the context of positive automorphisms.

The number of edges in  $\partial\Delta$ not
labelled
$t$ is the quantity that determines the upper bound we seek,
 $n := |\sigma| + |f^N(\sigma)|$).
We must bound the length of each corridor in $\Delta$
linearly in terms of $n$.  Theorem \ref{BoundS:BG3} provides a
bound in terms of  $\n$, so we must argue is that in the context
of stack diagrams, one can dispose of the contribution of the
$t$-edges to this bound.
In order to do so, we make an exhaustive list of those
places in the  proof of Theorem
\ref{BoundS:BG3} where $t$-edges were accounted for, 
and we explain why, in each case, 
they are not required in the setting of simple stack diagrams.

(1) The $t$-edges contributed to the bound on the size of 
$S_0(2)$ and $S_0(3a)$ in Section \ref{FastBeads}, but
these sets do not arise in stack diagrams.

(2) The $t$-edges were required in determining the
sets $\down_1(\T)$ used to bound the
lengths of teams (see Definition \ref{Down1Def}).
But $\down_1(\T)$ was
used only to bound the lengths of those
teams whose narrow past had several components at 
some time in the past, and this cannot happen in a stack
diagram.

(3) The $t$-edges entered the definition of $\down_2(\T)$,
which was
used to bound the
number of edges in $\QT$ before time $t_2(\T)$ (see Definition 
\ref{Down2Def}). But
there are no such edges in a stack of corridors, so
we do not have to worry about double-counting, and an
improved bound
on the lengths of teams can be derived directly from the
Pincer Lemma, noting that there are less than
$2|\partial \Delta|$ adjacencies of colours.

(4) In the bonus scheme, the set $\partial^e$ is used to 
bound the size of the interval of time $[\tau_0(e),\tau_2(e)]$,
but in a stack of corridors it is clear that $\tau_0(e)
=\tau_2(e)$, so the edges $\partial^e$ are not required.

(5) Likewise, when bounding the size of the bonuses provided by
rascals, we do not need to use
the edges $\down_2(e)$ if our diagram is simply a stack of corridors

(6) A final use of $t$-edges is hidden in our references to Part \ref{Part:BG1}
in the implementation of the Bonus scheme, specifically
the bound on the sum of the lengths of blocks satisfying condition
(iv) of the `tautologous tetrad'. This is unnecessary in stack
diagrams because there are no singularities and
no edges that are cancelled by edges from outside the future of
$S_0$, so the paths $\pi_l$ and $\pi_r$ travel forwards
in time  until they hit the boundary and
$\sum |\text{bdy}(\mathfrak{B})|<n$ bounds the size of the sum of
all such blocks.
$\square$

\backmatter

\printindex
 
 \end{document}